\documentclass{amsart}
\usepackage{graphicx,psfrag}

\numberwithin{equation}{section}

\begin{document}


\newtheorem{example}{Example}[section]
\newtheorem{note}[example]{Note}
\newtheorem{theorem}[example]{Theorem}
\newtheorem{corollary}[example]{Corollary}
\newtheorem{definition}[example]{Definition}
\newtheorem{proposition}[example]{Proposition}
\newtheorem{algorithm}[example]{Algorithm}
\newtheorem{lemma}[example]{Lemma}
\newtheorem{problem}[example]{Problem}
\newtheorem{conjecture}[example]{Conjecture}


\newcommand{\infinity}{\infty}
\renewcommand{\mod}{\text{mod}\,}

 
\font\twelvesym=msbm10 at 12pt
\font\tensym=msbm10
\font\sevensym=msbm7
\font\fivesym=msbm5
\newfam\ssymfam
\textfont\ssymfam=\tensym
\scriptfont\ssymfam=\sevensym
\scriptscriptfont\ssymfam=\fivesym
\def\ssym{\fam\ssymfam\tensym}


\newcommand{\Z}{{\ssym Z}}
\newcommand{\N}{{\ssym N}}

\newcommand{\T}{{\mathcal T}}
\newcommand{\U}{{\mathcal U}}
\newcommand{\tT}{\tilde{\mathcal T}}
\newcommand{\tU}{\widetilde{\mathcal U}}
\newcommand{\Y}{{\mathcal Y}}
\newcommand{\B}{{\mathcal B}}
\newcommand{\D}{{\mathcal D}}
\newcommand{\M}{{\mathcal M}}
\renewcommand{\P}{{\mathcal P}}
\newcommand{\R}{{\mathcal R}}

\newcommand{\Proof}{\medskip\noindent {\it Proof: }}
\newcommand{\cqfd}{\hfill $\Box$ \medskip}
\newcommand{\run}{{\mathcal X}}
\newcommand{\wt}{{\rm wt\,}}
\newcommand{\owt}{{wt}}
\newcommand{\mwt}{\tilde{wt}}
\newcommand{\ochi}{{\chi}}
\newcommand{\mchi}{\tilde{\chi}}
\newcommand{\pchi}{\ddot\chi}
\newcommand{\tkappa}{\tilde{\kappa}}
\newcommand{\bari}{\overline\imath}
\newcommand{\barj}{\overline\jmath}
\newcommand{\barh}{\overline h}
\newcommand{\boldm}{\boldsymbol{m}}
\newcommand{\oboldm}{\overline{\boldm}}
\newcommand{\boldn}{\boldsymbol{n}}
\newcommand{\boldu}{\boldsymbol{u}}
\newcommand{\oboldu}{\overline{\boldu}}
\newcommand{\bolde}{\boldsymbol{e}}
\newcommand{\boldB}{\boldsymbol{B}}
\newcommand{\boldC}{\boldsymbol{C}}
\newcommand{\boldCC}{\overline{\boldsymbol{C}}}
\newcommand{\boldQ}{\boldsymbol{Q}}
\newcommand{\boldQQ}{\overline{\boldsymbol{Q}}}
\newcommand{\boldN}{\boldsymbol{N}}
\newcommand{\boldmu}{\boldsymbol{\mu}}
\newcommand{\boldnu}{\boldsymbol{\nu}}
\newcommand{\boldDelta}{\boldsymbol{\Delta}}
\newcommand{\wombat}{\rule[-6pt]{0pt}{46pt}}
\newcommand{\qbinom}[2]{{\genfrac{[}{]}{0pt}{}{#1}{#2}}_q}


\newcommand{\muIndUp}[1]{\boldmu^{(#1)\phantom{*}}}
\newcommand{\muIndDn}[1]{\boldmu^{(#1)*}}
\newcommand{\nuIndUp}[1]{\boldnu^{(#1)\phantom{*}}}
\newcommand{\nuIndDn}[1]{\boldnu^{(#1)*}}
\newcommand{\muIndUpA}[1]{\boldmu^{(#1)\phantom{*}}_\uparrow}
\newcommand{\muIndDnA}[1]{\boldmu^{(#1)*}_\uparrow}
\newcommand{\nuIndUpA}[1]{\boldnu^{(#1)\phantom{*}}_\uparrow}
\newcommand{\nuIndDnA}[1]{\boldnu^{(#1)*}_\uparrow}


\renewcommand{\subjclassname}{%
          \textup{2000} Mathematics Subject Classification} 


\DeclareRobustCommand{\SkipTocEntry}[5]{} 

\makeatletter
        \newcommand\@dotsep{4.5}
        \def\@tocline#1#2#3#4#5#6#7{\relax
          \ifnum #1>\c@tocdepth 
          \else
            \par \addpenalty\@secpenalty\addvspace{#2}%
            \begingroup \hyphenpenalty\@M
            \@ifempty{#4}{%
              \@tempdima\csname r@tocindent\number#1\endcsname\relax
            }{%
              \@tempdima#4\relax
            }%
            \parindent\z@ \leftskip#3\relax \advance\leftskip\@tempdima\relax
            \rightskip\@pnumwidth plus1em \parfillskip-\@pnumwidth
            #5\leavevmode\hskip-\@tempdima #6\relax
            \leaders\hbox{$\m@th
              \mkern \@dotsep mu\hbox{.}\mkern \@dotsep mu$}\hfill
            \hbox to\@pnumwidth{\@tocpagenum{#7}}\par
            \nobreak
            \endgroup
          \fi}
        \makeatother 



\hyphenation{boson-ic 
             ferm-ion-ic 
	     para-ferm-ion-ic
             two-dim-ension-al
	     two-dim-ension-al}

\title[Fermionic expressions for Virasoro characters]{Fermionic
            expressions for minimal model\\ Virasoro characters}
\thanks{Research supported by the Australian Research Council (ARC)}
\author[Trevor Welsh]{Trevor~A.~Welsh}
\address{Department of Mathematics and Statistics,
             The University of Melbourne,
             Park\-ville, Victoria 3010, Australia.}
\email{trevor@ms.unimelb.edu.au}

\begin{abstract}

Fermionic expressions for all minimal model Virasoro
characters $\chi^{p, p'}_{r, s}$ are stated and proved. 
Each such expression is a sum of terms of
{\em fundamental fermionic form} type.
In most cases, all these terms are written down using certain trees which
are constructed for $s$ and $r$ from the Takahashi lengths and truncated
Takahashi lengths associated with the continued fraction of $p'/p$.
In the remaining cases, in addition to such terms, the fermionic expression
for $\chi^{p, p'}_{r, s}$ contains a different character
$\chi^{\hat p, \hat p'}_{\hat r,\hat s}$, and is thus recursive in nature.

Bosonic-fermionic $q$-series identities for all characters
$\chi^{p, p'}_{r, s}$ result from equating these fermionic
expressions with known bosonic expressions.
In the cases for which $p=2r$, $p=3r$, $p'=2s$ or $p'=3s$,
Rogers-Ramanujan type identities result from equating these
fermionic expressions with known product expressions
for $\chi^{p, p'}_{r, s}$.

The fermionic expressions are proved by first obtaining fermionic
expressions for the generating functions $\chi^{p, p'}_{a, b, c}(L)$
of length $L$ Forrester-Baxter paths, using various combinatorial
transforms. In the $L\to\infty$ limit, the fermionic expressions for
$\chi^{p, p'}_{r, s}$ emerge after mapping between the trees that
are constructed for $b$ and $r$ from the Takahashi and truncated
Takahashi lengths respectively.

\end{abstract}

\maketitle

\newpage

\tableofcontents

\newpage

\setcounter{secnumdepth}{10}

\setcounter{section}{0}

\section{Prologue}\label{PrologueSec}

\subsection{Introduction}\label{IntroSec}

The rich mathematical structure of the two-dimensional conformal
field theories of Belavin, Polyakov, and Zamolodchikov \cite{bpz,dms-book} 
is afforded by the presence of infinite dimensional Lie algebras
such as the Virasoro algebra in their symmetries.
Indeed, the presence of these algebras ensures the solvability
of the theories.
In the conformal field theories known as the minimal models,
which are denoted $\M^{p, p'}$ where $1<p<p'$ with $p$ and $p'$ coprime,
the spectra are expressible in terms of
the Virasoro characters $\hat\chi^{p,p'}_{r,s}$ where
$1\le r<p$ and $1\le s<p'$. The irreducible highest weight module
corresponding to the character $\hat\chi^{p,p'}_{r,s}$ has central
charge $c^{p,p'}$ and conformal dimension $\Delta^{p, p'}_{r, s}$
given by:
\begin{equation}
c^{p,p'} = 1 - \frac{6(p' - p)^2}{p p'},
\qquad
\Delta^{p, p'}_{r, s} =
\frac{ (p' r - p s)^2 - (p' - p)^2 }{4 p p'}.
\end{equation}
In \cite{feigin-fuchs,rocha-caridi,felder}, it was shown that
$\hat\chi^{p,p'}_{r,s}=q^{\Delta^{p, p'}_{r, s}}\chi^{p,p'}_{r,s}$,
where the (normalised) character $\chi^{p,p'}_{r,s}$ is given by:
\begin{equation}\label{RochaEq}
\chi^{p, p'}_{r, s}=
{\frac{1}{(q)_\infty}}\sum_{\lambda=-\infty}^\infty
(q^{\lambda^2pp'+\lambda(p'r-ps)}-q^{(\lambda p+r)(\lambda p'+s)}),
\end{equation}
and as usual, $(q)_{\infty}=\prod_{i=1}^\infty (1-q^i)$.
Note that $\chi^{p,p'}_{r,s}=\chi^{p,p'}_{p-r,p'-s}$.

An expression such as (\ref{RochaEq}) is known as
bosonic because it arises naturally via the
construction of a Fock space using bosonic generators. 
Submodules are factored out from the Fock space using an
inclusion-exclusion procedure, thereby leading to an expression
involving the difference between two constant-sign expressions.

However, there exist other expressions for $\chi^{p, p'}_{r, s}$
that provide an intrinsic physical interpretation of the states of the
module in terms of quasiparticles. 
These expressions are known as fermionic expressions
because the quasiparticles therein are forbidden to occupy identical
states. Two of the simplest such expressions arise when $p'=5$ and $p=2$:
\begin{align}
\label{RRinfferm}
\chi^{2,5}_{1,2}=
\sum_{n=0}^\infty \frac{q^{n^2}}{(q)_n},
\qquad
\chi^{2,5}_{1,1}=
\sum_{n=0}^\infty \frac{q^{n(n+1)}}{(q)_n}.
\end{align}
Here, $(q)_0=1$ and $(q)_{n}=\prod_{i=1}^n (1-q^i)$ for $n>0$.
Equating expressions of this fermionic type with the corresponding
instances of (\ref{RochaEq}) yields what are known as bosonic-fermionic
$q$-series identities.
In this paper, we give fermionic expressions for all $\chi^{p,p'}_{r,s}$.
In doing so, we obtain a bosonic-fermionic identity for each character
$\chi^{p,p'}_{r,s}$.

{}From both a mathematical and physical point of view, further interest
attaches to the fermionic expressions for $\chi^{p, p'}_{r, s}$
because in certain cases, another expression is available for
these characters.
This expression is of a product form and is obtained from (\ref{RochaEq})
by means of Jacobi's triple product identity
\cite[eq.\ (II.28)]{gasper-rahman}
or Watson's quintuple product identity \cite[ex.\ 5.6]{gasper-rahman}.
The former applies in the cases where $p=2r$ or $p'=2s$, giving:
\begin{equation} \label{ProdEq1}
\chi^{2r, p'}_{r, s}=
\hskip-1mm
\prod_{\begin{subarray}{c}
  n=1\\n\not\equiv0,\pm rs\,(\mod rp')\end{subarray}}^{\infty}
\hskip-2mm
\frac{1}{1-q^n},
\qquad
\chi^{p, 2s}_{r, s}=
\hskip-1mm
\prod_{\begin{subarray}{c}
  n=1\\n\not\equiv0,\pm rs\,(\mod sp)\phantom{'}\end{subarray}}^{\infty}
\hskip-2mm
\frac{1}{1-q^n}.
\end{equation}
The latter applies in the cases where $p=3r$ or $p'=3s$, giving:
\begin{equation}
\label{ProdEq3}
\chi^{3r, p'}_{r, s}=
\hskip-1mm
\prod_{\begin{subarray}{c}
  n=1\\n\not\equiv0,\pm rs\,(\mod 2rp')\\n\not\equiv\pm2r(p'-s)\,(\mod 4rp')\\
       \end{subarray}}^{\infty}
\hskip-3mm\frac{1}{1-q^n},
\qquad
\hskip-1mm
\chi^{p, 3s}_{r, s}=
\prod_{\begin{subarray}{c}
  n=1\\n\not\equiv0,\pm rs\,(\mod 2sp)\\n\not\equiv\pm2s(p-r)\,(\mod 4sp)\\
       \end{subarray}}^{\infty}
\hskip-3mm\frac{1}{1-q^n}.
\end{equation}
It may be shown (see \cite{christe,bytsko-fring}) that apart from these
and those resulting from identifying
$\chi^{3r,p'}_{2r,p'-s}=\chi^{3r,p'}_{r,s}$ and
$\chi^{p,3s}_{p-r,2s}=\chi^{p,3s}_{r,s}$ in (\ref{ProdEq3}),
there exist no other expressions for $\chi^{p,p'}_{r,s}$
as products of terms $(1-q^n)^{-1}$.

In the case where $r=1$ and $p'=5$, the first expression in (\ref{ProdEq1})
yields via (\ref{RRinfferm}), the celebrated Rogers-Ramanujan identities
\cite{rogers,rogers-ramanujan}:
\begin{align}
\label{RR1}
\sum_{n=0}^\infty \frac{q^{n^2}}{(q)_n}
&=\prod_{j=1}^{\infty}
\frac{1}{(1-q^{5j-4})(1-q^{5j-1})},\\
\label{RR2}
\sum_{n=0}^\infty \frac{q^{n(n+1)}}{(q)_n}
&=\prod_{j=1}^{\infty}
\frac{1}{(1-q^{5j-3})(1-q^{5j-2})}.
\end{align}

Via (\ref{ProdEq1}) or (\ref{ProdEq3}),
the fermionic expressions for $\chi^{p,p'}_{r,s}$ given in this
paper thus lead to generalisations of the Rogers-Ramanujan identities
in the cases where $p=2r$, $p=3r$, $p'=2s$ or $p'=3s$.

In fact, prior to the advent of conformal field theory, the search for
generalisations of the Rogers-Ramanujan identities led to many expressions
that are now recognised as fermionic expressions for certain
$\chi^{p,p'}_{r,s}$. Slater's compendium \cite{slater} contains,
amongst other things, fermionic expressions for all characters when
$(p,p')=(3,4)$, when $(p,p')=(3,5)$ and, of course, when $(p,p')=(2,5)$.
The multisum identities of Andrews and Gordon \cite{andrews74,gordon} deal
with all cases when $(p,p')=(2,2k+1)$ for $k>2$.
In \cite{andrews84}, Andrews showed how the Bailey chain may be used
to yield further identities, and in particular, produced
fermionic expressions for the characters
$\chi^{3,3k+2}_{1,k}$, $\chi^{3,3k+2}_{1,k+1}$,
$\chi^{3,3k+2}_{1,2k+1}$ and $\chi^{3,3k+2}_{1,2k+2}$ for $k\ge1$.
The Bailey chain together with the Bailey lattice was further exploited
in \cite{foda-quano} to yield more infinite sequences of characters
(see \cite[p1652]{foda-quano} for these sequences).

Another route to the Virasoro characters is provided by two-dimensional
statistical models or one-dimensional quantum spin
chains that exhibit critical behaviour
\cite{baxter-book,takahashi-suzuki}.
In this paper it suffices for us
to concentrate on the restricted solid-on-solid (RSOS) statistical
models of Forrester-Baxter \cite{forrester-baxter,abf} in regime III,
because all Virasoro characters $\chi^{p,p'}_{r,s}$ arise in the
calculation of the one-point functions
of such models in the thermodynamic limit.

The RSOS models of \cite{forrester-baxter} are parametrised by
integers $p'$, $p$, $a$, $b$ and $c$, where
$1\le p<p'$ with $p$ and $p'$ coprime, $1\le a,b,c<p'$ and $c=b\pm1$.
(The special cases where $p=1$ or $p=p'-1$ were first dealt
with in \cite{abf}: the $p=1$ case may be associated with regime II
of these models.)
Here, the use of the corner transfer matrix method \cite{baxter-book}
naturally leads to expressions for the one-point functions
in terms of generating functions of certain lattice paths.
Such a path of length $L$ is a sequence $h_0, h_1, h_2, \ldots, h_L,$
of integers such that:
\begin{enumerate}
\item $1 \le h_i \le p'-1$  for $0 \le i \le L$,
\item $h_{i+1} = h_i \pm 1$ for $0 \le i <   L$,
\item $h_0 = a, h_L = b.$
\end{enumerate}
The set of all such paths is denoted $\P^{p, p'}_{a, b, c}(L)$.

\begin{figure}[ht]
\includegraphics[scale=1.00]{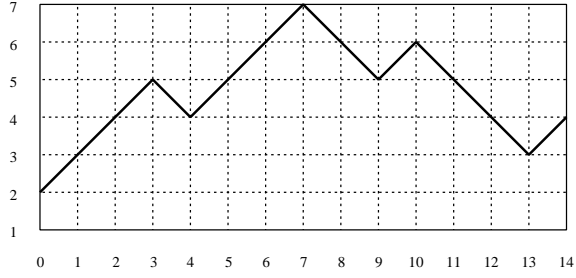}
\caption{Typical path.}
\label{TypicalBasicFig}
\medskip
\end{figure}

Each element of $\P^{p, p'}_{a, b, c}(L)$ is readily depicted on
a two-dimensional $L \times (p'-2)$ grid by connecting the points
$(i,h_i)$ and $(i+1,h_{i+1})$ for $0\le i<L$. A typical element
of $\P^{3,8}_{2, 4, 3}(14)$ is shown in Fig.~\ref{TypicalBasicFig}.

Each path $h\in \P^{p, p'}_{a, b, c}(L)$ is assigned a weight%
\footnote{This weighting function is described in Section \ref{WeightSec}.
It makes use of the values of $p$ and $c$.
The weighting function is not required for this introduction.}
$wt(h)\in\Z_{\ge0}$, whereupon the generating function
$\chi^{p,p'}_{a,b,c}(L)$ is defined by:
\begin{equation}\label{GenFunDefEq}
\chi^{p,p'}_{a,b,c}(L)
=\sum_{h\in{\P}^{p,p'}_{a,b,c}(L)} q^{wt(h)}.
\end{equation}

By setting up recurrence relations for $\chi^{p,p'}_{a,b,c}(L)$
(see Appendix \ref{BosonicSec}), it may be verified that:
\begin{equation}\label{FinRochaEq}
\begin{split}
\chi^{p,p'}_{a,b,c}(L)&=
\sum_{\lambda=-{\infinity}}^{\infinity}
q^{\lambda^{2} p p'+ \lambda (p'r-pa)}
\qbinom{L}{\frac{L+a-b}{2}-p'\lambda}\\[0.5mm]
&\hskip20mm
-\sum_{\lambda=-\infinity}^\infinity
q^{(\lambda p+r)(\lambda p'+a)}
\qbinom{L}{\frac{L+a-b}{2}-p'\lambda-a},
\end{split}
\end{equation}
where
\begin{equation}\label{groundstatelabel}
r=\left\lfloor\frac{pc}{p'}\right\rfloor+\frac{b-c+1}2,
\end{equation}
\noindent and, as usual, the Gaussian polynomial
$\qbinom{A}{B}$ is defined to be:
\begin{equation}\label{gaussian}
\qbinom{A}{B}=
\left\{
  \begin{array}{cl}
     \displaystyle \frac{(q)_A}{(q)_{A-B}(q)_B} &
	  \quad\mbox{if } 0\le B\le A;\\[3mm]
     0 &
	  \quad\mbox{otherwise}.
  \end{array} \right.
\end{equation}

So as to be able to take the $L\to\infty$ limit in (\ref{FinRochaEq}),
we assume that $|q|<1$.
Then, if $r=0$ or $r=p$,
\begin{align}\label{ChiLim0Eq}
\lim_{L\to\infty}\chi^{p, p'}_{a,b,c}(L)&=0.\\
\intertext{Otherwise, if $0<r<p$,}
\label{ChiLimEq}
\lim_{L\to\infty}\chi^{p, p'}_{a,b,c}(L)&=\chi^{p, p'}_{r, a}.
\end{align}
In view of this, we refer to $\chi^{p, p'}_{a,b,c}(L)$ as a finitized
character.

Since the expression (\ref{FinRochaEq}) is of a similar form
to (\ref{RochaEq}) and the former yields the latter in the $L\to\infty$
limit, we refer to (\ref{FinRochaEq}) as a bosonic expression for
$\chi^{p, p'}_{a,b,c}(L)$.
In this paper we first derive fermionic expressions
for $\chi^{p,p'}_{a,b,c}(L)$. The fermionic expressions
for $\chi^{p, p'}_{r, a}$ then arise in the $L\to\infty$ limit.
For example, for $L$ even, we have (c.f.~\cite[p43]{macmahon}):
\begin{align}
\label{RRfinferm}
\chi^{2,5}_{2,2,3}(L)=
\sum_{n=0}^\infty q^{n^2} \qbinom{L-n}{n},
\qquad
\chi^{2,5}_{1,3,2}(L)=
\sum_{n=0}^\infty q^{n(n+1)} \qbinom{L-1-n}{n}.
\end{align}
Using (\ref{ChiLimEq}) and (\ref{groundstatelabel}),
in the $L\to\infty$ limit, these yield the fermionic expressions
for the characters $\chi^{2,5}_{1,2}$ and $\chi^{2,5}_{1,1}$ given
in (\ref{RRinfferm}).

Note that en route to the bosonic-fermionic $q$-series identities
for the characters $\chi^{p,p'}_{r,s}$,
we will have proved bosonic-fermionic polynomial identities
for the finitized characters $\chi^{p, p'}_{a,b,c}(L)$.

{}From a physical point of view, a systematic study of fermionic expressions
for $\chi^{p,p'}_{r,s}$ was first undertaken in \cite{stony-brook1}. An
excellent overview of this and related work is given in \cite{stony-brook2}.
In \cite{stony-brook1}, fermionic expressions were conjectured
for all the unitary characters $\chi^{p,p+1}_{r,s}$,
as well as for characters $\chi^{p,kp+1}_{1,k}$ where $k\ge2$,
and characters $\chi^{p,p+2}_{(p-1)/2,(p+1)/2}$ for $p$ odd.
In fact, the expressions for $\chi^{p,p+1}_{r,s}$ and
$\chi^{p,p+1}_{p-r,p+1-s}$ in \cite{stony-brook1} are different, but since
these two characters are equal, two different fermionic expressions
for each character $\chi^{p,p+1}_{r,s}$ are thus provided.
A further two expressions for each character $\chi^{p,p+1}_{r,s}$
were conjectured in \cite{melzer}, and proofs of these and the
expressions of \cite{stony-brook1} were given in the cases
of $p=3$ and $p=4$.
A proof of the expressions of \cite{stony-brook1} in the special
case where $s=1$ was given in \cite{berkovich} using the
technique of \lq\lq telescopic expansion\rq\rq.
A similar technique was used in \cite{schilling} to give a proof of
all four expressions for each $\chi^{p,p+1}_{r,s}$.
In \cite{warnaar1,warnaar2}, an analysis of the lattice paths
described above, provided a combinatorial proof of two of the
fermionic expressions for $\chi^{p,p+1}_{r,s}$.
In \cite{foda-welsh}, a combinatorial proof along the lines of that
used in the current paper, was given for all four fermionic expressions
for $\chi^{p,p+1}_{r,s}$.
A proof for the $\chi^{p,p+1}_{1,1}$ case using Young tableaux and
rigged-configurations was presented in \cite{dasmahapatra-foda}.

Fermionic expressions for $\chi^{p,p'}_{r,s}$ for all $p$ and $p'$,
and certain $r$ and $s$, were stated in \cite{berkovich-mccoy}.
To specify the restrictions on $r$ and $s$, let
$\T$ be the set of Takahashi lengths and $\tT$ the set of truncated
Takahashi lengths associated with the continued fraction of $p'/p$
(these values are defined in Section \ref{TakSec}).
Then either both $p'-s\in\T$ and $r\in\tT$,
or $s\in\T$ and $r\in\tT$ and either $s\ge s_0$ or $r\ge r_0$,
where $r_0$ and $s_0$ are the smallest positive integers such that
$|p'r_0-ps_0|=1$.
{}From the point of view of path generating functions, the reason
for this strange restriction is expounded in \cite{foda-welsh-kyoto}.

The expressions of \cite{berkovich-mccoy} were proved and vastly
extended in \cite{bms} to give
fermionic expressions for all $\chi^{p, p'}_{r, s}$ where $s\in\T$
and $r$ is unrestricted except for $1\le r<p$.
In the most general cases, these expressions are positive sums
of {\em fundamental fermionic forms} whereas in all previous expressions
a single fundamental fermionic form was involved.
The expressions of \cite{bms} are obtained through first using
telescopic expansion to derive fermionic
expressions for the finitizations $\chi^{p,p'}_{a,b,c}(L)$:
these finitized expressions are linear combinations, not necessarily
positive sums, of fundamental fermionic forms. 
An interesting feature of the resulting fermionic expressions is that
they make use of a modified definition of the Gaussian polynomial.
As explained in \cite{foda-welsh-kyoto}, the modified definition of the
Gaussian polynomial serves to account for terms that correspond
to a character $\chi^{\hat p,\hat p'}_{\hat r,\hat s}$ where
$\hat p$ and $\hat p'$ are smaller than $p$ and $p'$ respectively.
In this article, we provide and prove fermionic expressions for
all $\chi^{p, p'}_{r, s}$ and $\chi^{p,p'}_{a,b,c}(L)$.
Here, the modified Gaussian polynomial does not correctly account
for the extra terms that sometimes arise.
Thus we revert to the original definition of the Gaussian polynomial,
and deal with the extra terms, which again correspond to
a character $\chi^{\hat p,\hat p'}_{\hat r,\hat s}$, by writing
a recursive expression for $\chi^{p, p'}_{r, s}$:
a sum of a number of fundamental fermionic forms plus
$\chi^{\hat p,\hat p'}_{\hat r,\hat s}$.
Since $\hat p'<p'$, this recursive expression terminates,
and in fact, the degree of recursion is small compared to $p'$.

As described later, obtaining the fundamental fermionic forms that occur in
the fermionic expression for $\chi^{p, p'}_{r, s}$ requires the construction
of two trees: the first depends on $s$ and the Takahashi lengths for $p'/p$,
and the second depends on $r$ and the truncated Takahashi lengths for $p'/p$.
The first of these tree constructions is similar to the recursion
underlying the \lq\lq branched chain of vectors\rq\rq\ of \cite[p340]{bms}.
Obtaining the fundamental fermionic forms that occur in
the fermionic expression for $\chi^{p,p'}_{a,b,c}(L)$
also requires the construction of two trees --- one for $a$ and
one for $b$, with each constructed using the Takahashi lengths for $p'/p$.
In fact, in our proof, these expressions are derived first, and
the fermionic expressions for $\chi^{p, p'}_{r, s}$ emerge
after taking the ${L\to\infty}$ limit.


\subsection{Structure of this paper}
In the remainder of this prologue, we describe our fermionic expressions
for $\chi^{p, p'}_{r, s}$ and $\chi^{p,p'}_{a,b,c}(L)$ in detail.
Central to these expressions is the continued fraction of $p'/p$,
for which the notation is fixed in Section \ref{ContFSec}.
The fermionic expressions for the Virasoro characters
$\chi^{p, p'}_{r, s}$ are stated in Section \ref{CharSec} as a
positive sum over {\em fundamental fermionic forms} $F(\boldu^L,\boldu^R)$
where $\boldu^L$ and $\boldu^R$ are certain vectors.
The $F(\boldu^L,\boldu^R)$ are manifestly positive definite.
The construction of the vectors $\boldu^L$ and $\boldu^R$, and
the notation required to define $F(\boldu^L,\boldu^R)$ is
described and discussed in the subsequent Sections \ref{TakSec}
through \ref{MNsysSec}.

The fermionic expressions for $\chi^{p, p'}_{r, s}$ and
$\chi^{p,p'}_{a,b,c}(L)$ have a recursive component in that they
(sometimes) involve a term
$\chi^{\hat p, \hat p'}_{\hat r, \hat s}$ or
$\chi^{\hat p,\hat p'}_{\hat a,\hat b,\hat c}(L)$ respectively,
where $\hat p'<p'$ and $\hat p<p$.
These terms are specified in
Section \ref{ExtraSec},
where a limit to the degree of recursion is also given.

The fermionic expression for a particular $\chi^{p,p'}_{a,b,c}(L)$ is
given in Section \ref{FinFermSec} or \ref{FermLikeSec} depending on the
value of $b$.
In either case, $\chi^{p,p'}_{a,b,c}(L)$ is expressed in terms
of fundamental fermionic forms $F(\boldu^L,\boldu^R,L)$ which are
also manifestly positive definite.
In the case in Section \ref{FinFermSec},
the expression for $\chi^{p,p'}_{a,b,c}(L)$ is a positive sum of such terms.
However, as described in Section \ref{FermLikeSec}, in the other case
we use linear combinations of such terms with both positive and negative
coefficients. 


Complete examples of the constructions that are described in this
section are provided in Appendix \ref{ExampleSec}.
Maple programs that construct (and evaluate) the fermionic expressions
for any case are available from the author.

The proof of these expressions occupies Sections \ref{CombinSec}
through \ref{FermCharSec}.
In Section \ref{DissSec}, we discuss various interesting aspects of
these expressions and their proof.

Appendices \ref{AppBSec}, \ref{DireApp} and \ref{TakApp} give some technical
details that are required in the proof.
Appendix \ref{BosonicSec} provides a proof of the bosonic expression
(\ref{FinRochaEq}).

\subsection{Continued fractions}\label{ContFSec}

The continued fraction of $p'/p$ lies at the heart of the constructions
of this paper. Let $p$ and $p'$ be coprime integers with $1\le p<p'$.
If
\begin{equation*}
\frac{p'}{p}=
{c_0+\frac{\displaystyle\strut 1}{\displaystyle c_1 +
\frac{\displaystyle\strut 1}{
\genfrac{}{}{0pt}{}{\lower-5pt\hbox{$\vdots$}}{\displaystyle c_{n-2}+
\frac{\displaystyle\strut 1}
{\displaystyle\strut c_{n-1} +
\frac{\displaystyle\strut 1}{\displaystyle c_n}}}}}}
\end{equation*}
with $c_i\ge1$ for $0\le i<n$, and $c_n\ge2$,
then $[c_0,c_1,c_2,\ldots,c_n]$ is said to be the
{\it continued fraction} for $p'/p$.

We refer to $n$ as the {\em height} of $p'/p$.
We set $t=c_0+c_1+\cdots+c_n-2$ and refer to it as the {\em rank} of $p'/p$.

We also define:\footnote{The $t_{n+1}$ defined here
differs from that defined in \cite{berkovich-mccoy,bms}, and that
defined in \cite{flpw}.}
\begin{equation}\label{ZoneEq}
t_k=-1+\sum_{i=0}^{k-1} c_i,
\end{equation}
for $0\le k\le n+1$.
Then $t_{n+1}=t+1$ and $t_n\le t-1$.
If the non-negative integer $j\le t_{n+1}$ satisfies
$t_{k}<j\le t_{k+1}$, we say that $j$ is {\em in zone $k$}.
We then define $\zeta(j)=k$.
Note that there are $n+1$ zones and that for 
$0\le k\le n$, zone $k$ contains $c_k$ integers.

\subsection{Fermionic character expressions}\label{CharSec}

Here, we present fermionic expressions for all Virasoro characters
$\chi^{p,p'}_{r,s}$ where $1<p<p'$ with $p$ and $p'$ coprime,
$1\le r<p$ and $1\le s<p'$.

Using notation that will be defined subsequently,
when $p'>3$ the expression for $\chi^{p,p'}_{r,s}$ takes the form:
\begin{equation}\label{FermCEq}
\chi^{p,p'}_{r,s}=
  \sum_{\begin{subarray}{c} \boldu^L\in\U(s)\\
                            \boldu^R\in\tU(r)
        \end{subarray}}
\hskip-2mm
F(\boldu^L,\boldu^R)
\:+\:
  \left\{ \begin{array}{cl}
    \displaystyle
    \chi^{\hat p,\hat p'}_{\hat r,\hat s}
   &\mbox{if } \eta(s)=\tilde\eta(r) \mbox{ and }\hat s\ne0\ne\hat r;\\[1.5mm]
   0 &\mbox{otherwise},
         \end{array} \right.
\end{equation}
where $\U(s)$ and $\tU(r)$ are
sets of $t$-dimensional vectors that are described below.
For $t$-dimensional vectors $\boldu^L$ and $\boldu^R$, the
{\em fundamental fermionic form} $F(\boldu^L,\boldu^R)$ is defined by:
\begin{equation}\label{F2Eq}
\begin{split}
&\hskip-4mmF(\boldu^L,\boldu^R)
=\sum_{\begin{subarray}{c}
           n_1,n_2,\ldots,n_{t_1}\in\Z_{\ge0}\\
           {\boldm}\equiv\boldQQ(\boldu^L+\boldu^R)
        \end{subarray}}
  \hskip-5mm
  q^{\tilde{\boldn}^T\boldB\tilde{\boldn}
   +\frac{1}{4}{\boldm}^T{\boldCC}{\boldm}
   -\frac{1}{2}(\oboldu^L_\flat+\oboldu^R_\sharp)\cdot\boldm
   +\frac{1}{4}\gamma(\run^L,\run^R)}
\\[1mm]
&\hskip24mm\times\:
  \frac{1}{(q)_{m_{t_1+1}}}
  \prod_{j=1}^{t_1} \frac{1}{(q)_{n_j}}
  \prod_{j=t_1+2}^{t-1}
  \left[
 {m_j-\frac{1}{2}
     ({\boldCC}^*{\boldm}\!-\!\oboldu^L\!-\!\oboldu^R)_j\atop m_j}
  \right]_q\!\!,\hskip-1mm
\end{split}
\end{equation}
where the sum is over non-negative integers $n_1,n_2,\ldots,n_{t_1}$,
and $(t-t_1-1)$-dimen\-sional
vectors%
\footnote{All vectors in this paper will be column vectors: for
typographical convenience, their components will be expressed in
row vector form.}
$\boldm=(m_{t_1+1},m_{t_1+2},\ldots,m_{t-1})$
each of whose integer components is congruent, modulo 2,
to the corresponding component of the vector $\boldQQ(\boldu^L+\boldu^R)$
which we define below.
The $t_1$-dimensional vector $\tilde{\boldn}$ is defined by:
\begin{equation}\label{ntildeEq}
\textstyle
\tilde{\boldn}=
(n_1-\frac{1}{2}u_1,n_2-\frac{1}{2}u_2,\ldots,n_{t_1-1}-\frac{1}{2}u_{t_1-1},
 n_{t_1}-\frac{1}{2}u_{t_1}+\frac{1}{2}m_{t_1+1}),
\end{equation}
having set $(u_1,u_2,\ldots,u_{t})=\boldu^L+\boldu^R$.
We also set $\run^L=\run(s,\boldu^L)$
and $\run^R=\tilde{\run}(r,\boldu^R)$, where this
and all the other notation occurring in the expressions (\ref{FermCEq})
and (\ref{F2Eq}) will be defined in the subsequent sections.
In the particular cases for which $p=2$, we find that $t=t_1+1$ and yet
$m_t$ is undefined. We obtain the correct expressions by setting $m_t=0$.

The exceptional case of $(p,p')=(2,3)$ is dealt with by:
\begin{equation}\label{FermExceptEq}
\chi^{2,3}_{1,1}=\chi^{2,3}_{1,2}=1,
\end{equation}
which follows from (\ref{RochaEq}) and Jacobi's triple product identity,
or simply from the first case of (\ref{ProdEq1}).
The expression (\ref{FermExceptEq}) is sometimes required when
iterating (\ref{FermCEq}).

We note that in the cases where $p'=p+1$, certain choices mean that up
to four fermionic expressions for $\chi^{p,p+1}_{r,s}$ arise from
(\ref{FermCEq}). These are the expressions of \cite{melzer}.
One choice is as explained in footnote \ref{ChFN} in Section \ref{FormVSec},
and an analogous choice is available in Section \ref{TruncSec}.
In cases other than $p'=p+1$, only one fermionic expression
for $\chi^{p,p'}_{r,s}$ arises from (\ref{FermCEq}).

We also note that when $r=1$, the expression obtained from (\ref{FermCEq})
does not always appear to coincide with a known expression for the same
character \cite{berkovich-mccoy,bms,flpw}.
However, the discrepancy appears only in the $(t_1+1)$th component
of $\boldu^R$.
It is easily seen (see Lemma \ref{ShiftSigmaLem}) that this discrepancy
has no effect on the summands of (\ref{F2Eq}).
In these cases, the known expressions may be reproduced by setting
$\tilde\sigma_1=t_1$ in the description of Section \ref{TruncSec}
instead of the stated $\tilde\sigma_1=t_1+1$.


\subsection{The Takahashi and string lengths}\label{TakSec}

For $p$ and $p'$ coprime with $1\le p<p'$,
we now use the notation of Section \ref{ContFSec}
to define the corresponding
set $\{\kappa_i\}_{i=0}^t$ of {\em Takahashi lengths},
the multiset $\{\tilde\kappa_i\}_{i=0}^t$ of
{\em truncated Takahashi lengths},
and the multiset $\{l_i\}_{i=0}^t$
of {\em string lengths}.
These definitions are based on those of \cite{takahashi-suzuki,bms}.
First define $y_k$ and $z_k$ for
$-1\le k\le n+1$ by:
\begin{align*}
y_{-1}&=0;&z_{-1}&=1;\\
y_0&=1;&z_0&=0;\\
y_k&=c_{k-1}y_{k-1}+y_{k-2};&z_k&=c_{k-1}z_{k-1}+z_{k-2},
\quad (1\le k\le n+1).
\end{align*}
In particular, $y_{n+1}=p'$ and $z_{n+1}=p$.
Now for $t_k<j\le t_{k+1}$ and $0\le k\le n$, set:
\begin{align*}
\kappa_j&=y_{k-1}+(j-t_k)y_k;\\
\tilde\kappa_j&=z_{k-1}+(j-t_k)z_k;\\
l_j&=y_{k-1}+(j-t_k-1)y_k.
\end{align*}
Note that $\kappa_{j}=l_{j+1}$ unless $j=t_k$ for some $k$,
in which case $\kappa_{t_k}=y_k$ and $l_{t_k+1}=y_{k-1}$.
It may be shown (see Lemma \ref{CFLem}) that if $j$ is in zone $k$ then
$\kappa_j/\tkappa_j$ has continued fraction $[c_0,c_1,\ldots,c_{k-1},j-t_k]$.
In particular, $y_k/z_k$ has continued fraction
$[c_0,c_1,\ldots,c_{k-1}]$ for $0\le k\le n+1$.
We define $\T=\{\kappa_i\}_{i=0}^{t-1}$
and $\T'=\{p'-\kappa_i\}_{i=0}^{t-1}$
(we do not include $\kappa_t$ in the former, nor $p'-\kappa_t$ in the latter).
Apart from the cases in which $p=1$ or $p=p'-1$,
we find that $\T\cap\T'=\emptyset$.
We also define $\tT=\{\tkappa_i\}_{i=t_1+1}^{t-1}$
and $\tT'=\{p-\tkappa_i\}_{i=t_1+1}^{t-1}$.
(We see that these latter two sets are the sets
$\T$ and $\T'$ that would be obtained
for the continued fraction $[c_1,c_2,\ldots,c_{n}]$.)

For example, in the case $p'=223$, $p=69$, which yields the continued
fraction $[3,4,3,5]$, so that $n=3$,
$(t_1,t_2,t_3,t_4)=(2,6,9,14)$ and $t=13$,
we obtain:
\begin{align*}
\{y_k\}_{k=-1}^4&=\{0,1,3,13,42,223\},\\
\{z_k\}_{k=-1}^4&=\{1,0,1,4,13,69\},\\
\{\kappa_j\}_{j=0}^{13}&=
\{1,2,3,4,7,10,13,16,29,42,55,97,139,181\},\\
\{\tilde\kappa_j\}_{j=0}^{12}&=
\{1,1,1,1,2,3,4,5,9,13,17,30,43,56\},\\
\{l_j\}_{j=1}^{13}&=
\{1,2,1,4,7,10,3,16,29,13,55,97,139\}.
\end{align*}
Then $\T=\{1,2,3,4,7,10,13,16,29,42,55,97,139\}$
and $\tT=\{1,2,3,4,5,9,13,17,$ $30,43\}$.

\subsection{Takahashi trees}\label{FormSec}

This and the following two sections describe the sets
$\U(s)$ and $\tU(r)$ that occur in (\ref{FermCEq}).

Given $a$ with $1\le a<p'$, we use the Takahashi lengths
to produce a {\em Takahashi tree} for $a$.
This is a binary tree of positive integers
for which each node except the root node,
is labelled $a_{i_1i_2\cdots i_k}$
for some $k\ge1$, with $i_j\in\{0,1\}$ for $1\le j\le k$.%
\footnote{This notation will be abused by letting
$a_{i_1i_2\cdots i_k}$ denote both the particular node of the tree, and
the value at that node: the interpretation will be clear from the context.}
The root node of this tree is unlabelled.
Each node is either a {\em branch-node}, a {\em through-node}
or a {\em leaf-node}.
These have 2, 1, or 0 children respectively.
Each child of each branch-node is either a branch-node or a through-node,
but the child of each through-node is always a leaf-node.
Naturally, each child of the node labelled $a_{i_1i_2\cdots i_k}$
is labelled $a_{i_1i_2\cdots i_k0}$ or $a_{i_1i_2\cdots i_k1}$.
In fact, our construction has
$a_{i_1i_2\cdots i_k0} < a$ and $a_{i_1i_2\cdots i_k1} > a$ for
non-leaf-nodes and $a_{i_1i_2\cdots i_k}=a$ for leaf-nodes.

\begin{figure}[ht]
\includegraphics[scale=1.00]{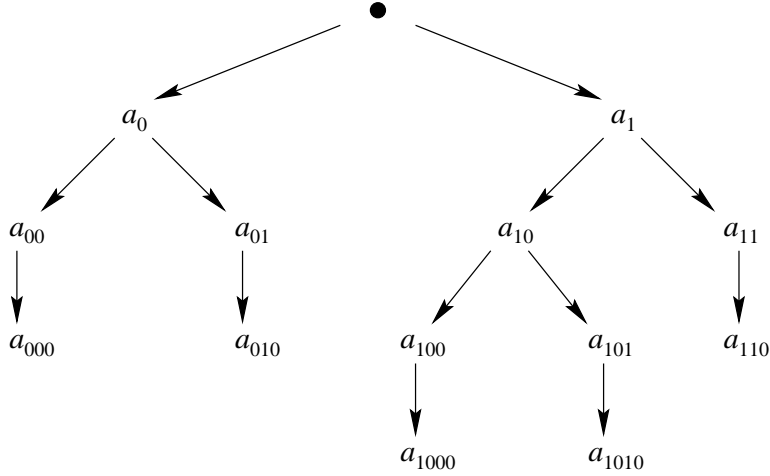}
\caption{Typical Takahashi tree.}
\label{TakahashiTreeFig}
\medskip
\end{figure}

We obtain the Takahashi tree for $a$ as follows.
In the case that $a\in\T\cup\T'$, the Takahashi tree comprises
a single leaf-node in addition to the root node.
This leaf-node is $a_0=a$, and the root node is designated a
through-node.
Otherwise the root node is a branch-node: $a_0$ is set to be
the largest element of $\T\cup\T'$ smaller than $a$,
and $a_1$ the smallest element of $\T\cup\T'$ larger than $a$.
We now generate the tree recursively.
If $a_{i_1i_2\cdots i_k}\ne a$,
and $\vert a_{i_1i_2\cdots i_k}-a\vert\in\T$ then
we designate $a_{i_1i_2\cdots i_k}$ a through-node and define 
$a_{i_1i_2\cdots i_k0}=a$ to be a leaf-node.
Otherwise, when $\vert a_{i_1i_2\cdots i_k}-a\vert\notin\T$,
we make $a_{i_1i_2\cdots i_k}$ a branch-node.
We take $\kappa_x\in\T$ to be the largest element of
$\T$ smaller than $\vert a_{i_1i_2\cdots i_k}-a\vert$.
Then, for $i_{k+1}\in\{0,1\}$, set
$a_{i_1i_2\cdots i_ki_{k+1}}
 =a_{i_1i_2\cdots i_k}+(-1)^{i_k}\kappa_{x+\vert i_k-i_{k+1}\vert}$.
This ensures that $a_{i_1i_2\cdots i_k0}< a$ and $a_{i_1i_2\cdots i_k1}> a$.
It is easy to see that the tree is finite.
In fact, each leaf-node occurs no deeper than $n+1$ levels
below the root node. Thus, there are at most $2^n$ leaf-nodes.
A typical Takahashi tree is shown in Fig.~\ref{TakahashiTreeFig}.

In fact, precisely this Takahashi tree arises in the case for which
$p'=223$, $p=69$ and $a=66$: the nodes take on the values
$a_0=55$, $a_1=84$, $a_{00}=65$, $a_{01}=68$, $a_{10}=55$,
$a_{11}=68$, $a_{100}=65$, $a_{101}=68$ and
$a_{000}=a_{010}=a_{1000}=a_{1010}=a_{110}=66$.

Note that the above construction ensures that
$\vert a_{i_1i_2\cdots i_{k-1}}-a_{i_1i_2\cdots i_k}\vert\in\T$ for each
node $a_{i_1i_2\cdots i_k}$ that is present in the Takahashi tree for $a$.
In addition, it is readily seen that if $a_{i_1i_2\cdots i_k}$ is
a branch-node, then $a_{i_1i_2\cdots i_k1}-a_{i_1i_2\cdots i_k0}\in\T$.

For each leaf-node $a_{i_1i_2\cdots i_k0}$, we also set
$a_{i_1i_2\cdots i_k1}=a_{i_1i_2\cdots i_k0}=a$ for later convenience.

\subsection{Takahashi tree vectors}\label{FormVSec}

Using the Takahashi tree for $a$, we now define a set $\U(a)$
of $t$-dimensional vectors.
The cardinality of this set is equal to the number of
leaf-nodes of the Takahashi tree.
First, for $0\le j\le t+1$,
define the vector $\bolde_j=(\delta_{1j},\delta_{2j},\ldots,\delta_{tj})$
where, as usual, the Kronecker delta is defined by $\delta_{ij}=1$
if $i=j$ and $\delta_{ij}=0$ otherwise.
Then, for $0\le i< j\le t+1$, define:
\begin{equation}\label{uijEq}
{\boldu}_{i,j}
=\bolde_{i}-\bolde_{j}-\sum_{k:i\le t_k<j} \bolde_{t_k}.
\end{equation}

Given a leaf-node $a_{i_1i_2\cdots i_d}$, we obtain the
corresponding vector $\boldu\in\U(a)$ as follows.
We first define a {\em run} $\run=\{\tau_j,\sigma_j,\Delta_j\}_{j=1}^d$.
Set $\bari=1-i$ for $i\in\{0,1\}$.
If $a_{\bari_1}\in\T$, set $\Delta_1=-1$ and
define $\sigma_1$ such that $\kappa_{\sigma_1}=a_{\bari_1}$,
otherwise if $a_{\bari_1}\in\T'$, set $\Delta_1=1$ and
define $\sigma_1$ such that $\kappa_{\sigma_1}=p'-a_{\bari_1}$.%
\footnote{\label{ChFN}In the cases where $p=1$ or $p=p'-1$, we have
$\T\cup\T'=\{2,3,\cdots,p'-2\}$. If $1<a_{\bari_1}<p'-1$, we may
proceed by considering either $a_{\bari_1}\in\T$ or
$a_{\bari_1}\in\T'$.
The two alternatives lead to different equally valid expressions.}
For $2\le j\le d$, set $\Delta_j=-(-1)^{i_{j-1}}$ and
define $\sigma_j$ such that
$\kappa_{\sigma_j}=\vert a_{i_1i_2\cdots i_{j-1}}
                        -a_{i_1i_2\cdots i_{j-1}\bari_j}\vert$.
Define $\tau_1=t+1$.
For $2\le j\le d$, define $\tau_j$ such that
$\kappa_{\tau_j}=a_{i_1i_2\cdots i_{j-2}1}-a_{i_1i_2\cdots i_{j-2}0}$.
Finally define:
\begin{equation}\label{uEq}
\boldu(\run)=
\sum_{m=1}^d \boldu_{\sigma_m,\tau_m}
+
\left\{\begin{array}{l@{\quad\quad}l}
 0  & \mathrm{if\ } \Delta_1=-1;\\
\bolde_t & \mathrm{if\ } \Delta_1=1.
 \end{array}\right.
\end{equation}
The set $\U(a)$ comprises all vectors $\boldu(\run)$ obtained
in this way from the leaf-nodes of the Takahashi tree for $a$.

To illustrate this construction, again consider the case
$p'=223$, $p=69$ and $a=66$. In the case of the leaf-node
$a_{1010}$, we have $\vert a_{101}-a_{1011}\vert=2=\kappa_1$,
$\vert a_{10}-a_{100}\vert=10=\kappa_5$,
$\vert a_{1}-a_{11}\vert=16=\kappa_7$
and $a_0=55=\kappa_{10}$,
and hence $\sigma_4=1$, $\sigma_3=5$, $\sigma_2=7$, and $\sigma_1=10$.
Also, we have $a_{101}-a_{100}=3=\kappa_2$,
$a_{11}-a_{10}=13=\kappa_6$
and $a_{1}-a_{0}=29=\kappa_8$,
and thus $\tau_4=2$, $\tau_3=6$, $\tau_2=8$ and $\tau_1=14$.
Since $a_0\in\T$, we have $\Delta_1=-1$ leading
to the run $\run_1=\run=\{\{14,8,6,2\},\{10,7,5,1\},\{-1,1,-1,1\}\}$,
from which using (\ref{uEq}), we obtain the vector
$\boldu(\run_1)=(1,-1,0,0,1,-1,1,-1,0,1,0,0,0)$.
In the case of the leaf-node $a_{000}$, we obtain the run
$\run_2=\run=\{\{14,8,2\},\{12,6,0\},\{1,-1,1\}\}$, 
which yields the vector $\boldu(\run_2)=(0,-1,0,0,0,0,0,-1,0,0,0,1,1)$.
After performing similar calculations for the leaf-nodes
$a_{010}$, $a_{1000}$ and $a_{110}$, we obtain the set:
\begin{align*}
\U(66)=\{ & (1,-1,0,0,1,-1,1,-1,0,1,0,0,0), \\
          & (0,-1,0,0,0,0,0,-1,0,0,0,1,1), \\
          & (1,-1,0,0,1,-1,0,-1,0,0,0,1,1), \\
          & (0,-1,0,0,0,0,1,-1,0,1,0,0,0), \\
          & (1,-1,0,0,0,-1,0,0,0,1,0,0,0) \}.
\end{align*}

Given $a$, the run $\run=\{\tau_j,\sigma_j,\Delta_j\}_{j=1}^d$
may be recovered from the vector $\boldu\in\U(a)$,
and thus we define $\run(a,\boldu)=\run$.
We also define $\Delta(\boldu)=\Delta_d$.
For example, in the case of the vector
$\boldu=(1,-1,0,0,1,-1,1,-1,0,1,0,0,0)$ that arises in the above example,
we have $\run(66,\boldu)=\run_1$ and $\Delta(\boldu)=1$.

\subsection{Truncated Takahashi tree}\label{TruncSec}

For $1\le r<p$, we now mirror the above constructions of the
Takahashi tree, the set $\U(a)$ and the vectors
$\run(a,\boldu)$ using the truncated Takahashi lengths in
place of the Takahashi lengths.

The truncated Takahashi tree for $r$ is obtained by using
the prescription given in Section \ref{FormSec} for $a=r$,
after replacing each occurrence of $\T$ and $\T'$ with $\tT$ and $\tT'$
respectively, and using truncated Takahashi lengths
$\tilde\kappa_j$ in place of Takahashi lengths $\kappa_j$.

\begin{figure}[ht]
\includegraphics[scale=0.90]{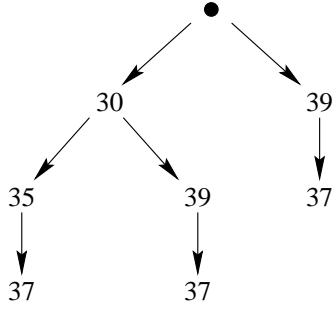}
\caption{Truncated Takahashi tree for $r=37$, $p'=223$, $p=69$.}
\label{TruncTreeFig}
\medskip
\end{figure}

For instance, in the case $p'=223$ and $p=69$,
where $\tT=\{1,2,3,4,5,9,13,17,$ $30,43\}$
and $\tT'=\{68,67,66,65,64,60,56,52,39,26\}$,
the truncated Takahashi tree for $r=37$ is given in Fig.~\ref{TruncTreeFig}.

We now use the truncated Takahashi tree for $r$ to produce a
set $\tU(r)$ of $t$-dimensional vectors, the cardinality of this set
being equal to the number of leaf-nodes of the tree.

For each leaf-node $r_{i_1i_2\cdots i_d}$ of the truncated
Takahashi tree for $r$, we obtain the corresponding vector
$\boldu\in\tU(r)$ in a similar
way to the prescription given in Section \ref{FormVSec} as follows.
We first define the run
$\tilde\run=\{\tilde\tau_j,\tilde\sigma_j,\tilde\Delta_j\}_{j=1}^d$.
If $r_{\bari_1}\in\tT$, set $\tilde\Delta_1=-1$ and
define $\tilde\sigma_1$ such that $\tkappa_{\tilde\sigma_1}=r_{\bari_1}$,
otherwise if $r_{\bari_1}\in\tT'$,
set $\tilde\Delta_1=1$ and define $\tilde\sigma_1$ such that
$\tkappa_{\tilde\sigma_1}=p-r_{\bari_1}$.
For $2\le j\le d$, set $\tilde\Delta_j=-(-1)^{i_{j-1}}$ and
define $\tilde\sigma_j$ such that
$\tkappa_{\tilde\sigma_j}=\vert r_{i_1i_2\cdots i_{j-1}}
                        -r_{i_1i_2\cdots i_{j-1}\bari_j}\vert$.
Define $\tilde\tau_1=t+1$.
For $2\le j\le d$, define $\tilde\tau_j$ such that
$\tkappa_{\tilde\tau_j}=r_{i_1i_2\cdots i_{j-2}1}-r_{i_1i_2\cdots i_{j-2}0}$.
Finally we obtain the vector $\boldu=\boldu(\tilde\run)$ using (\ref{uEq}).
Since the number of leaf-nodes of the truncated Takahashi tree for $r$
is at most $2^{n-1}$, we then have $|\tU(r)|\le 2^{n-1}$.

To illustrate this construction, again consider the case
$p'=223$, $p=69$ and $r=37$. In the case of the leaf-node
$r_{000}$, we have $\vert r_{00}-r_{001}\vert=2=\tkappa_4$,
$\vert r_{0}-r_{01}\vert=9=\tkappa_8$, and $r_1=39=p-\tkappa_{11}$,
and hence $\tilde\sigma_3=4$, $\tilde\sigma_2=8$, and $\tilde\sigma_1=11$.
Also, we have $r_{00}-r_{01}=4=\tkappa_6$ and $r_{0}-r_{1}=9=\tkappa_8$,
and thus $\tilde\tau_3=6$, $\tilde\tau_2=8$ and $\tilde\tau_1=14$.
Then, after noting that $r_1\in\tT'$,
we obtain the run $\tilde\run=\{\{14,8,6\},\{11,8,4\},\{1,-1,-1\}\}$. 
The use of (\ref{uEq}) then yields
$\boldu(\tilde\run)=(0,0,0,1,0,-1,0,0,0,0,1,0,1)$.
After performing similar calculations for the
leaf-nodes $r_{010}$ and $r_{10}$, we obtain the set:
\begin{align*}
\tU(37)=\{  & (0,0,0,1,0,-1,0,0,0,0,1,0,1), \\
            & (0,0,0,1,0,-1,1,-1,0,0,1,0,1), \\
            & (0,0,0,1,0,-1,0,-1,0,0,1,0,0) \}.
\end{align*}

\subsection{The linear term}\label{LinSec}

With $t_0,t_1,t_2,\ldots,t_{n+1}$ and $t$ as defined in Section \ref{ContFSec},
for each $t$-dimensional vector
$\boldu=(u_1,u_2,\ldots,u_t)$, define the $(t-1)$-dimensional vector
$\boldu^{\flat}=(u^{\flat}_1,u^{\flat}_2,\ldots,u^{\flat}_{t-1})$ by:
\begin{equation}
u^{\flat}_j=
\left\{
  \begin{array}{cl}
       0 &\quad\mbox{if } t_k<j\le t_{k+1},\; k\equiv 0\,(\mod 2);\\[0.5mm]
     u_j &\quad\mbox{if } t_k<j\le t_{k+1},\; k\not\equiv 0\,(\mod 2),
  \end{array} \right.
\end{equation}
and the $(t-1)$-dimensional vector
$\boldu^{\sharp}=(u^{\sharp}_1,u^{\sharp}_2,\ldots,u^{\sharp}_{t-1})$ by:
\begin{equation}
u^{\sharp}_j=
\left\{
  \begin{array}{cl}
     u_j &\quad\mbox{if } t_k<j\le t_{k+1},\; k\equiv 0\,(\mod 2);\\[0.5mm]
       0 &\quad\mbox{if } t_k<j\le t_{k+1},\; k\not\equiv 0\,(\mod 2).
  \end{array} \right.
\end{equation}
Then, of course, $(\boldu)_j=(\boldu^{\flat}+\boldu^{\sharp})_j$
for $1\le j<t$.

We also define
$\oboldu=(u_{t_1+2},u_{t_1+3},\ldots,u_t)$,
$\oboldu^{\flat}=(u_{t_1+1}^\flat,u_{t_1+2}^\flat,\ldots,u_{t-1}^\flat)$, and
$\oboldu^{\sharp}=(u_{t_1+1}^\sharp,u_{t_1+2}^\sharp,\ldots,u_{t-1}^\sharp)$.

For convenience, we sometimes write
$\boldu_\flat$, $\boldu_\sharp$, $\oboldu_\flat$ and $\oboldu_\sharp$
for $\boldu^\flat$, $\boldu^\sharp$, $\oboldu^\flat$ and $\oboldu^\sharp$
respectively.

\subsection{The constant term}\label{ConSec}

Here we obtain the constant terms $\gamma(\run^L,\run^R)$
that appear in (\ref{F2Eq}).

First, given the run $\run=\{\tau_j,\sigma_j,\Delta_j\}_{j=1}^d$, define
the $t$-dimensional vector $\boldDelta(\run)$ by:%
\footnote{
The values $\{\Delta_j\}_{j=1}^d$ should not be confused with the
components of $\boldDelta(\run)$.}
\begin{equation}\label{DeltaEq}
\boldDelta(\run)=
\sum_{m=1}^d \Delta_m\boldu_{\sigma_m,\tau_m}
-
\left\{\begin{array}{l@{\quad\quad}l}
0             & \mathrm{if\ } \Delta_1=-1;\\
\bolde_t      & \mathrm{if\ } \Delta_1=1.
 \end{array}\right.
\end{equation}
Note the similarity between (\ref{DeltaEq}) and (\ref{uEq}).
To illustrate this definition, consider the run
$\run_1=\{\{14,8,6,2\},$ $\{10,7,5,1\},\{-1,1,-1,1\}\}$
obtained in Section \ref{FormVSec} above in the case where
$p'=223$, $p=69$ and $a=66$.
Here, (\ref{DeltaEq}) yields:
$$
\boldDelta(\run_1)=(1,-1,0,0,-1,1,1,-1,0,-1,0,0,0).
$$

Using (\ref{DeltaEq}), the runs
$\run^L=\{\tau^L_j,\sigma^L_j,\Delta^L_j\}_{j=1}^{d^L}$
and $\run^R=\{\tau^R_j,\sigma^R_j,\Delta^R_j\}_{j=1}^{d^R}$
give rise to vectors
$\boldDelta^L=\boldDelta(\run^L)$ and
$\boldDelta^R=\boldDelta(\run^R)$.
{}From these, we iteratively generate sequences
$(\beta_{t},\beta_{t-1},\ldots,\beta_0)$,
$(\alpha_{t},\alpha_{t-1},\ldots,\alpha_0)$,
and $(\gamma_{t},$ $\gamma_{t-1},\ldots,\gamma_0)$ as follows.
Let $\alpha_t=\beta_t=\gamma_t=0$.
Now, for $j\in\{t,t-1,\ldots,1\}$, obtain $\alpha_{j-1}$, $\beta_{j-1}$,
and $\gamma_{j-1}$
from $\alpha_j$, $\beta_j$, and $\gamma_j$ in the following three stages.
Firstly, obtain:
\begin{equation}\label{Const1Eq}
(\beta'_{j-1},\gamma'_{j-1})=
(\beta_j+(\boldDelta^L)_j-(\boldDelta^R)_j,\gamma_j+2\alpha_j(\boldDelta^R)_j).
\end{equation}
Then obtain:
\begin{equation}\label{Const2Eq}
(\alpha''_{j-1},\gamma''_{j-1})=
(\alpha_j+\beta'_{j-1},\gamma'_{j-1}-(\beta'_{j-1})^2).
\end{equation}
Finally set:
\begin{equation}\label{Const3Eq}
(\alpha_{j-1},\beta_{j-1},\gamma_{j-1})=
\begin{cases}
(\alpha''_{j-1},\alpha_{j},
 -(\alpha''_{j-1})^2-\gamma''_{j-1})\\
\hskip33mm\text{if $j=t_k+1$ with $1\le k\le n$;}\\
\hbox to 0mm{$\displaystyle(\alpha''_{j-1},\beta'_{j-1},\gamma''_{j-1})$\hss}
\hskip33mm\text{otherwise.}
\end{cases}
\hskip-7mm
\end{equation}

Having thus obtained $\gamma_0$, we set
$\gamma(\run^L,\run^R)=\gamma_0$ whenever $\sigma^R_{d^R}>0$.
This definition suffices for the purposes of writing down
the fermionic expressions for $\chi^{p,p'}_{r,s}$ because,
with $\tT=\{\tkappa_i\}_{i=t_1+1}^{t-1}$, the construction of
Section \ref{TruncSec} ensures that $\sigma^R_{d^R}>0$.

However, for later purposes, values of $\gamma(\run^L,\run^R)$
are required when $\sigma^R_{d^R}=0$.
We specify this and also $\gamma'(\run^L,\run^R)$ as follows.
First obtain $\gamma_0$ as above.
Then, with $\sigma=\sigma^R_{d^R}$ and $\Delta=\Delta^R_{d^R}$:

          if $p'>2p$, $\sigma=0$, $\Delta=1$,
          and either $\left\lfloor\frac{(b+1)p}{p'}\right\rfloor=
               \left\lfloor\frac{bp}{p'}\right\rfloor$ or $b=p'-1$,
          then set
\begin{equation}
\nonumber
\gamma=\gamma_0-2(a-b),\qquad \gamma'=\gamma-2L;
\end{equation}

          if $p'>2p$, $\sigma=0$, $\Delta=-1$,
          and either $\left\lfloor\frac{(b-1)p}{p'}\right\rfloor=
               \left\lfloor\frac{bp}{p'}\right\rfloor$ or $b=1$,
          then set
\begin{equation}
\nonumber
\gamma=\gamma_0+2(a-b),\qquad \gamma'=\gamma-2L;
\end{equation}

          if $p'<2p$, $\sigma=0$, $\Delta=1$,
          and either $\left\lfloor\frac{(b+1)p}{p'}\right\rfloor\ne
               \left\lfloor\frac{bp}{p'}\right\rfloor$ or $b=p'-1$,
          then set
\begin{equation}
\nonumber
\gamma=\gamma_0+2(a-b),\qquad \gamma'=\gamma+2L;
\end{equation}

          if $p'<2p$, $\sigma=0$, $\Delta=-1$,
          and either $\left\lfloor\frac{(b-1)p}{p'}\right\rfloor\ne
               \left\lfloor\frac{bp}{p'}\right\rfloor$ or $b=1$,
          then set
\begin{equation}
\nonumber
\gamma=\gamma_0-2(a-b),\qquad \gamma'=\gamma+2L.
\end{equation}

\noindent Otherwise, we set $\gamma=\gamma'=\gamma_0$.
We then define
$\gamma(\run^L,\run^R)=\gamma$ and
$\gamma'(\run^L,\run^R)=\gamma'$.

Note that although the constant term $\gamma'(\run^L,\run^R)$
is independent of the components of $\boldm$,
in certain cases it does have a dependence on $L$.

\subsection{The quadratic term}\label{QuadSec}

In this section, we define the $(t-t_1-1)\times(t-t_1-1)$
matrices ${\boldCC}$, ${\boldCC}^*$, and the $t_1\times t_1$
matrix $\boldB$ that occurs in (\ref{FermCEq}).

First, with $t_0,t_1,t_2,\ldots,t_{n+1}$ and $t$ as defined in
Section \ref{ContFSec}, define $C_{ji}$ for $0\le i,j\le t$ by, 
when the indices are in this range,
\begin{equation}\label{CijDef}
\begin{array}{cccl}
C_{j,j-1}=-1,
&C_{j,j}=1,
&C_{j,j+1}=\phantom{-}1,
&\hbox{if $j=t_k,\quad k=1,2,\ldots,n$;}\\
C_{j,j-1}=-1,
&C_{j,j}=2,
&C_{j,j+1}=-1,
&\hbox{$0\le j<t$ otherwise,}
\end{array}
\end{equation}
and set $C_{ji}=0$ for $|i-j|>1$.
Also define $B_{ji}=\min\{i,j\}$ for $1\le i,j\le t_1$.
Then define the matrices:
\begin{equation*}
\boldCC=\{C_{ji}\}_{\begin{subarray}{c} t_1+1\le j<t\\
                                              t_1+1\le i<t
                          \end{subarray}},
\qquad
\boldCC^*=\{C_{ji}\}_{\begin{subarray}{c} t_1+2\le j\le t\\
                                              t_1+1\le i<t
                          \end{subarray}},
\qquad
\boldB=\{B_{ji}\}_{\begin{subarray}{c} 1\le j\le t_1\\
                                        1\le i\le t_1
                    \end{subarray}}.
\end{equation*}
Note that $\boldCC$ is tri-diagonal, and $\boldCC^*$ is upper-diagonal.
The matrices $\boldCC$ and $\boldB^{-1}$ may be viewed as minor
generalisations of Cartan matrices of type $A$.

To illustrate these constructions, consider the case $p=9$ and $p'=67$,
where the continued fraction of $p'/p$ is $[7,2,4]$
and $t_1=6$, $t_2=8$, $t_3=12$ and $t=11$. Here, we obtain:
\begin{equation*}
\boldCC=
\left(
\begin{smallmatrix}
\phantom{-}2 &          -1 &\phantom{-}. &\phantom{-}. \\[1mm]
          -1 &\phantom{-}1 &\phantom{-}1 &\phantom{-}. \\[1mm]
\phantom{-}. &-1           &\phantom{-}2 &          -1 \\[1mm]
\phantom{-}. &\phantom{-}. &          -1 &\phantom{-}2
\end{smallmatrix}
\right),
\qquad
\boldCC^*=
\left(
\begin{smallmatrix}
          -1 &\phantom{-}1 &\phantom{-}1 &   \\[1mm]
\phantom{-}. &          -1 &\phantom{-}2 &-1 \\[1mm]
\phantom{-}. &\phantom{-}. &          -1 &\phantom{-}2 \\[1mm]
\phantom{-}. &\phantom{-}. &\phantom{-}. &-1
\end{smallmatrix}
\right),
\qquad
\boldB=
\left(
\begin{smallmatrix}
 1\; & 1\; &\; 1 &\; 1 &\; 1 &\; 1\\[1mm]
 1\; & 2\; &\; 2 &\; 2 &\; 2 &\; 2\\[1mm]
 1\; & 2\; &\; 3 &\; 3 &\; 3 &\; 3\\[1mm]
 1\; & 2\; &\; 3 &\; 4 &\; 4 &\; 4\\[1mm]
 1\; & 2\; &\; 3 &\; 4 &\; 5 &\; 5\\[1mm]
 1\; & 2\; &\; 3 &\; 4 &\; 5 &\; 6
\end{smallmatrix}
\right).
\end{equation*}

The $t\times t$ matrices
${\boldC}$ and ${\boldC}^*$ are conveniently defined here as well.
With $C_{ji}$ as in (\ref{CijDef}) above, set:
\begin{equation*}
\boldC=\{C_{ji}\}_{\begin{subarray}{c} 0\le j<t\\
                                        0\le i<t
                    \end{subarray}},
\qquad
\boldC^*=\{C_{ji}\}_{\begin{subarray}{c} 1\le j\le t\\
                                        0\le i<t
                    \end{subarray}}.
\end{equation*}
Then $\boldC$ is tri-diagonal, $\boldC^*$ is upper-diagonal,
and $\boldCC$ and $\boldCC^*$ are the 
lower left $(t-t_1-1)\times(t-t_1-1)$
submatrices of $\boldC$ and $\boldC^*$ respectively.

To illustrate these constructions, consider the example of $p=9$
and $p'=31$, where the continued fraction
of $p'/p$ is $[3,2,4]$ and $t_1=2$, $t_2=4$, $t_3=8$ and $t=7$. Here:
\begin{equation*}
\boldC=
\left(
\begin{smallmatrix}
\phantom{-}2 &          -1 &\phantom{-}. &\phantom{-}.
             &\phantom{-}. &\phantom{-}. &\phantom{-}. \\[1mm]
          -1 &\phantom{-}2 &          -1 &\phantom{-}.
             &\phantom{-}. &\phantom{-}. &\phantom{-}. \\[1mm]
\phantom{-}. &          -1 &\phantom{-}1 &\phantom{-}1
             &\phantom{-}. &\phantom{-}. &\phantom{-}. \\[1mm]
\phantom{-}. &\phantom{-}. &          -1 &\phantom{-}2
             &          -1 &\phantom{-}. &\phantom{-}. \\[1mm]
\phantom{-}. &\phantom{-}. &\phantom{-}. &          -1
             &\phantom{-}1 &\phantom{-}1 &\phantom{-}. \\[1mm]
\phantom{-}. &\phantom{-}. &\phantom{-}. &\phantom{-}.
             &          -1 &\phantom{-}2 &          -1 \\[1mm]
\phantom{-}. &\phantom{-}. &\phantom{-}. &\phantom{-}.
             &\phantom{-}. &          -1 &\phantom{-}2
\end{smallmatrix}
\right),
\qquad
\boldC^*=
\left(
\begin{smallmatrix}
          -1 &\phantom{-}2 &          -1 &\phantom{-}.
             &\phantom{-}. &\phantom{-}. &\phantom{-}. \\[1mm]
\phantom{-}. &          -1 &\phantom{-}1 &\phantom{-}1
             &\phantom{-}. &\phantom{-}. &\phantom{-}. \\[1mm]
\phantom{-}. &\phantom{-}. &          -1 &\phantom{-}2
             &          -1 &\phantom{-}. &\phantom{-}. \\[1mm]
\phantom{-}. &\phantom{-}. &\phantom{-}. &          -1
             &\phantom{-}1 &\phantom{-}1 &\phantom{-}. \\[1mm]
\phantom{-}. &\phantom{-}. &\phantom{-}. &\phantom{-}.
             &          -1 &\phantom{-}2 &          -1 \\[1mm]
\phantom{-}. &\phantom{-}. &\phantom{-}. &\phantom{-}.
             &\phantom{-}. &          -1 &\phantom{-}2 \\[1mm]
\phantom{-}. &\phantom{-}. &\phantom{-}. &\phantom{-}.
             &\phantom{-}. &\phantom{-}. &          -1
\end{smallmatrix}
\right).
\end{equation*}

\subsection{The $\boldm\boldn$-system and the parity vector}\label{MNsysSec}

The $\boldm\boldn$-system is a set of $t$ linear equations
that depends on $\boldu=(u_1,u_2,\ldots,u_t)=\boldu^L+\boldu^R$,
and which defines an interdependence between two $t$-dimensional vectors
$\boldn=(n_1,n_2,\ldots,n_{t})$
and $\hat{\boldm}=(m_0,m_1,\ldots,m_{t-1})$.
The equations are given by, for $1\le j\le t$:
\begin{align}
m_{j-1}-m_{j+1}&=m_{j}+2n_{j}-u_{j}\quad
&&\text{if $j=t_k,\quad k=1,2,\ldots,n$;}
\label{MNEq1}\\
m_{j-1}+m_{j+1}&=2m_j+2n_j-u_{j}\quad
&&\text{otherwise,}
\label{MNEq2}
\end{align}
where we set $m_{t}=m_{t+1}=0$.

Then (\ref{MNEq1}) and (\ref{MNEq2}) imply that
for $\hat{\boldm}$ and $\boldn$ satisfying the $\boldm\boldn$-system,
\begin{equation}\label{MNmatEq1}
-{\boldC}^*\hat{\boldm}=2{\boldn}-{\boldu}.
\end{equation}
Since ${\boldC}^*$ is upper-triangular, its inverse is readily obtained
to yield:
\begin{equation}\label{MNmatEq2}
\hat{\boldm}=({\boldC^*}){}^{-1}(-2{\boldn}+\boldu).
\end{equation}

We now define $Q_i\in\{0,1\}$ for $0\le i<t$, by:
\begin{equation}\label{ParityDef}
(Q_0,Q_1,Q_2,\ldots,Q_{t-1})\equiv(\boldC^*){}^{-1}\boldu,
\end{equation}
and define the $(t-1)$-dimensional {\em parity vector}
$\boldQ(\boldu)=(Q_1,Q_2,\ldots,Q_{t-1})$
and the $(t-t_1-1)$-dimensional {\em parity vector}
${\boldQQ}(\boldu)=(Q_{t_1+1},Q_{t_1+2},\ldots,Q_{t-1})$.
It is readily seen that we also have
${\boldQQ}(\boldu)\equiv(\boldCC^*){}^{-1}\oboldu$.

The restriction ${\boldm}\equiv{\boldQQ}(\boldu^L+\boldu^R)$
in expression (\ref{F2Eq}) ensures that the quantity
$\frac{1}{2}(\boldCC^*{\boldm}-\oboldu^L-\oboldu^R)_j$
is an integer for $t_1+2\le j\le t$.

\subsection{The extra term}\label{ExtraSec}

Define $\{\xi_\ell\}_{\ell=0}^{2c_n-1}$ and
$\{\tilde\xi_\ell\}_{\ell=0}^{2c_n-1}$ according to
$\xi_0=\tilde\xi_0=0$, $\xi_{2c_n-1}=p'$, $\tilde\xi_{2c_n-1}=p$, and
\begin{displaymath}
\begin{array}{lcllcl}
\xi_{2k-1}&=&ky_n;&\tilde\xi_{2k-1}&=&kz_n;\\
\xi_{2k}&=&ky_n+y_{n-1};\quad&\tilde\xi_{2k}&=&kz_n+z_{n-1},
\end{array}
\end{displaymath}
for $1\le k<c_n$.
For example, when $p'=223$, $p=69$, we obtain:
$$
\begin{array}{ll}
\{\xi_\ell\}_{\ell=0}^{9}&=\quad
\{0,42,55,84,97,126,139,168,181,223\};\\[0.5mm]
\{\tilde\xi_\ell\}_{\ell=0}^{9}&=\quad
\{0,13,17,26,30,39,43,52,56,69\}.
\end{array}
$$
For $1\le s<p'$, define $\eta(s)=\ell$ where $\xi_\ell\le s<\xi_{\ell+1}$.
Similarly, for $1\le r<p$, define $\tilde\eta(r)=\ell$
where $\tilde\xi_\ell\le r<\tilde\xi_{\ell+1}$.
The values of $\hat p'$, $\hat p$, $\hat r$ and $\hat s$ that occur
in (\ref{FermCEq}) are then defined by:
\begin{equation}\label{ExtraEq1}
\begin{split}
\hat p'&=\xi_{\eta(s)+1}-\xi_{\eta(s)},\\
\hat p\phantom{{}'}&=\tilde\xi_{\eta(s)+1}-\tilde\xi_{\eta(s)},\\
\hat r\phantom{{}'}&=r-\tilde\xi_{\eta(s)},\\
\hat s\phantom{{}'}&=s-\xi_{\eta(s)}.
\end{split}
\end{equation}

In fact, it may be shown (see Appendix \ref{TakApp}) that if
$\eta(s)=\tilde\eta(r)$ and $\hat r\ne 0\ne\hat s$ so that
the extra term actually appears in (\ref{FermCEq}), then the
continued fraction of $\hat p'/\hat p$ is given by:
\begin{equation}\label{ExtraCfEq}
\begin{array}{ll}
[c_0,c_1,\ldots,c_{n-2}]
  &\text{ if $\eta(s)$ is odd;}\\
{}[c_0,c_1,\ldots,c_{n-1}]
  &\text{ if $\eta(s)\in\{0,2c_n-2\}$;}\\
{}[c_0,c_1,\ldots,c_{n-1}-1]
  &\text{ if $\eta(s)$ is even, $0<\eta(s)<2c_n-2$ and $c_{n-1}>1$;}
                 \hskip-5mm\\
{}[c_0,c_1,\ldots,c_{n-3}]
  &\text{ if $\eta(s)$ is even, $0<\eta(s)<2c_n-2$ and $c_{n-1}=1$.
                 \hskip-5mm}
\end{array}
\end{equation}
(Here, if a continued fraction $[c_0,\ldots,c_{n'-1},c_{n'}]$ with
$c_{n'}=1$ arises, then it is to be equated with the continued
fraction $[c_0,\ldots,c_{n'-1}+1]$.)

We then see that the height of the continued fraction of
$\hat p'/\hat p$ is at least one less than that of $p'/p$.
Since, we cannot have $\hat p=1$, the recursive process
implied by (\ref{FermCEq}) terminates after at most $n$ steps,
where $n$ is the height of $p'/p$.

\subsection{Finitized fermionic expressions}\label{FinFermSec}

In this and the following sections, we provide fermionic expressions for
the finitized characters $\chi^{p,p'}_{a,b,c}(L)$ where $1\le p<p'$
with $p$ and $p'$ coprime, $1\le a,b<p'$ and $c=b\pm1$ and
$L\equiv a-b\,(\mod2)$. 
In fact, in the proof that occupies the bulk of this paper,
these expressions for $\chi^{p,p'}_{a,b,c}(L)$ are proved first,
and the expressions (\ref{FermCEq}) for $\chi^{p,p'}_{r,s}$
derived therefrom.

The expressions for $\chi^{p,p'}_{a,b,c}(L)$ fall naturally into
two categories. The first deals with those cases for which:%
\footnote{Later in this paper, we refer to values of $b$
which satisfy (\ref{interface}) as {\em interfacial}.}
\begin{equation}\label{interface}
2\le b\le p'-2
\quad\mbox{ and }\quad
\left\lfloor{\frac{(b+1)p}{p'}}\right\rfloor
=\left\lfloor{\frac{(b-1)p}{p'}}\right\rfloor+1.
\end{equation}
Here we may take either $c=b\pm1$ because for these values of $b$,
$\chi^{p,p'}_{a,b,b-1}(L)=\chi^{p,p'}_{a,b,b+1}(L)$.

In these cases, we have the identity:

\begin{equation}\label{FermEq}
\chi^{p,p'}_{a,b,c}(L)=
  \sum_{\begin{subarray}{c} \boldu^L\in\U(a)\\
                            \boldu^R\in\U(b)
        \end{subarray}}
\hskip-2mm
F(\boldu^L,\boldu^R,L)
+ \left\{ \begin{array}{cl}
    \chi^{\hat p,\hat p'}_{\hat a,\hat b,\hat c}(L)
   &\mbox{if } \eta(a)=\eta(b) \mbox{ and }\hat a\ne0\ne\hat b;\\[1.5mm]
   0 &\mbox{otherwise}.
         \end{array} \right.
\end{equation}
where the sets $\U(a)$ and $\U(b)$ are defined
in Section \ref{FormVSec}, $F(\boldu^L,\boldu^R,L)$ is defined below,
and using the definitions of Section \ref{ExtraSec}, we set:
\begin{equation}\label{ExtraEq2}
\begin{split}
\hat p'&=\xi_{\eta(a)+1}-\xi_{\eta(a)},\\
\hat p\phantom{{}'}&=\tilde\xi_{\eta(a)+1}-\tilde\xi_{\eta(a)},\\
\hat a\phantom{{}'}&=a-\xi_{\eta(a)},\\
\hat b\phantom{{}'}&=b-\xi_{\eta(a)},
\end{split}
\end{equation}
and
\begin{equation}\label{ExtraEq3}
\hat c= \begin{cases}
       2&\text{if $c=\xi_{\eta(a)}>0$ and
	  $\left\lfloor\frac{(b+1)p}{p'}\right\rfloor=
 	   \left\lfloor\frac{(b-1)p}{p'}\right\rfloor+1$;}\\[1.5mm]
        \hat p'-2&\text{if $c=\xi_{\eta(a)+1}<p'$ and
 	  $\left\lfloor\frac{(b+1)p}{p'}\right\rfloor=
 	   \left\lfloor\frac{(b-1)p}{p'}\right\rfloor+1$;}\\[1.5mm]
        c-\xi_{\eta(a)} &\text{otherwise}.
         \end{cases}
\end{equation}
When $\eta(a)=\eta(b)$ and $\hat a\ne0\ne\hat b$
(so that the extra term actually appears in (\ref{FermEq})),
the continued fraction of $\hat p'/\hat p$
is again as specified in Section \ref{ExtraSec}, after setting $s=a$.

For the finitized characters $\chi^{p,p'}_{a,b,c}(L)$,
the fundamental fermionic form is defined by:
\begin{gather}\label{FEq}
\hskip-83mm
F(\boldu^L,\boldu^R,L)=\\[0.5mm]
  \sum_{\boldm\equiv\boldQ(\boldu^L+\boldu^R)}
  \hskip-7mm
  q^{\frac{1}{4}\hat{\boldm}^T\boldC\hat{\boldm}-\frac{1}{4} L^2
   -\frac{1}{2}(\boldu^L_\flat+\boldu^R_\sharp)\cdot\boldm
  +\frac{1}{4}\gamma'(\run^L,\run^R)}
  \prod_{j=1}^{t-1}
  \left[\!
  {m_j\!-\!\frac{1}{2} ({\boldC}^*\hat{\boldm}
                       \!-\!\boldu^L\!-\!\boldu^R)_j\!\atop m_j}
  \right]_q\!\!,\notag
\end{gather}
where the sum is over all $(t-1)$-dimensional vectors
$\boldm=(m_1,m_2,\ldots,m_{t-1})$,
each of whose integer components is congruent, modulo 2,
to the corresponding component of the vector $\boldQ(\boldu^L+\boldu^R)$,
and where we set $\hat{\boldm}=(L,m_1,m_2,\ldots,m_{t-1})$.
We also set $\run^L=\run(a,\boldu^L)$ and $\run^R=\run(b,\boldu^R)$,
as defined in Section \ref{FormVSec}.
When $t=1$, $F(\boldu^L,\boldu^R,L)$ is still defined by (\ref{FEq}),
after interpreting the (empty) product as $1$ and omitting the summation.

As will be seen later, the first component of (\ref{MNmatEq2})
gives:
\begin{equation}\label{Particles1Eq}
m_0=\sum_{i=1}^{t} l_i\left(2n_i-u_i\right),
\end{equation}
where $l_1,l_2,\ldots,l_{t}$ are the string lengths defined earlier.
For the practical purposes of evaluating the expression (\ref{FEq}),
we identify $m_0$ with $L$ and write (\ref{Particles1Eq}) in the form:
\begin{equation}\label{Particles2Eq}
\sum_{i=1}^{t} l_in_i=\frac{1}{2}\left(L+\sum_{i=1}^t l_iu_i\right).
\end{equation}
The summands in (\ref{FEq}) then correspond to solutions of
this partition problem with $n_i\in\Z_{\ge0}$ for $1\le i\le t$,
with $\{m_j\}_{j=1}^{t-1}$ obtained using (\ref{MNmatEq2}),
or equivalently, using (\ref{MNEq1}) and (\ref{MNEq2}).

Note that when using (\ref{FermEq}), the expression (\ref{ExtraEq2})
can sometimes yield $\hat b\in\{1,\hat p-1\}$.
In the subsequent iteration, (\ref{interface}) is thus not satisfied,
and the fermionic expression for
$\chi^{\hat p,\hat p'}_{\hat a,\hat b,\hat c}(L)$ must be obtained
using the expressions of the following Section \ref{FermLikeSec}.
Also note that (\ref{ExtraEq3}) might yield $\hat c\in\{0,\hat p'\}$.
Although this is no impediment to using the expressions of
Section \ref{FermLikeSec}, the bosonic expression (\ref{FinRochaEq})
does not immediately apply for $c\in\{0,p'\}$.
In Section \ref{WeightSec} we show how (\ref{FinRochaEq}) may be
extended to deal with $c\in\{0,p'\}$.

\subsection{Fermionic-like expressions}\label{FermLikeSec}

In this section, we provide {\it fermionic-like} expressions
for $\chi^{p,p'}_{a,b,c}(L)$ in the cases for which $b$ does not
satisfy (\ref{interface}).
In contrast to the case where (\ref{interface}) holds, the value of $c$
here does affect the resulting fermionic expression.
Indeed here, $\chi^{p,p'}_{a,b,b-1}(L)\ne\chi^{p,p'}_{a,b,b+1}(L)$.

The statement of these expressions involves sets
$\U^\pm(b)$ where, for $\Delta\in\{\pm1\}$, we define:
\begin{equation*}
\U^{\Delta}(b)=\{\boldu\in\U(b):\Delta(\boldu)=\Delta\}.
\end{equation*}
Therefore, we have the disjoint union $\U(b)=\U^+(b)\cup\U^-(b)$.

For $p'>2$, we then have the following identity:%
\footnote{
In some very particular cases (\ref{Ferm2Eq}) reduces to (\ref{FermEq}).
To describe these, let $t_*=t_1$ if $p'>2p$, and let $t_*=t_2$ if $p'<2p$.
If $1\le b<t_*$ then $\U^-(b)=\U(b)$ and $\U^+(b)=\emptyset$.
If $p'-t_*<b\le p'-1$ then $\U^+(b)=\U(b)$ and $\U^-(b)=\emptyset$.
Thus, if $1\le b<t_*$ and $c=b-1$, or $p'-t_*<b\le p'-1$ and $c=b+1$
then the second summation in (\ref{Ferm2Eq}) is zero,
thus yielding (\ref{FermEq}).}
\begin{align}\label{Ferm2Eq}
&\chi^{p,p'}_{a,b,c}(L)=
  \sum_{\begin{subarray}{c} \boldu^L\in\U(a)\\
                            \boldu^R\in\U^{c-b}(b) \end{subarray}}
\hskip-2mm  F(\boldu^L,\boldu^R,L)
+
  \sum_{\begin{subarray}{c} \boldu^L\in\U(a)\\
                            \boldu^R\in\U^{b-c}(b) \end{subarray}}
\hskip-2mm  \widetilde F(\boldu^L,\boldu^R,L) \\[0.5mm]
\notag
&\hskip35mm
+ \left\{ \begin{array}{cl}
    \chi^{\hat p,\hat p'}_{\hat a,\hat b,\hat c}(L)
   &\mbox{if } \eta(a)=\eta(b) \mbox{ and }\hat a\ne0\ne\hat b;\\[1.5mm]
   0 &\mbox{otherwise},
         \end{array} \right.
\end{align}
where $F(\boldu^L,\boldu^R,L)$ is as defined in (\ref{FEq}),
$\tilde F(\boldu^L,\boldu^R,L)$ is defined below, and the other
notation is as in Section \ref{FinFermSec}.

Let $\boldu\in\U(b)$ and let $\run=\run(b,\boldu)$.
If $\run=\{\tau_j,\sigma_j,\Delta_j\}_{j=1}^d$, we now define runs
$\run^+=\{\tau_j,\sigma^+_j,\Delta_j\}_{j=1}^d$ and
$\run^{++}=\{\tau_j,\sigma^{++}_j,\Delta_j\}_{j=1}^d$,
by setting $\sigma^{++}_j=\sigma^+_j=\sigma_j$ for $1\le j<d$,
and setting $\sigma^{++}_d=\sigma_j+2$ and $\sigma^{+}_d=\sigma_j+1$.
Using (\ref{uEq}), we then define $\boldu^+=\boldu(\run^+)$
and $\boldu^{++}=\boldu(\run^{++})$.

We now define
the fermionic-like terms $\widetilde F(\boldu^L,\boldu,L)$ which appear
in (\ref{Ferm2Eq}) in terms of $F(\boldu^L,\boldu^R,L')$ given by (\ref{FEq}),
with $\boldu^R$ being set to either of $\boldu$, $\boldu^+$ or $\boldu^{++}$.
In these three cases, the $\run^R$ that appears in (\ref{FEq}) should
be set to $\run$, $\run^+$ and $\run^{++}$ respectively,
and $\run^L=\run(a,\boldu^L)$.

Set $\Delta=\Delta_d$ and set:
\begin{equation}\label{tauDef}
\tau=
\left\{ \begin{array}{ll}
             \tau_{d} & \quad\mbox{if } d>1;\\
             t_n & \quad\mbox{if } d=1 \mbox{ and } \sigma_1<t_n;\\
             t   & \quad\mbox{if } d=1 \mbox{ and } \sigma_1\ge t_n.
  \end{array} \right.
\end{equation}
Then if $p'>2p$, define:
\begin{equation}\label{Tilde1Def}
\widetilde F(\boldu^L,\boldu,L)
 =\left\{ \begin{array}{ll}
    q^{\frac12(L+\Delta(a-b))} F(\boldu^L,\boldu,L)
     &\hskip-31mm\mbox{if } \sigma_d=0;\\[1.5mm]
    q^{-\frac12(L-\Delta(a-b))}
       \left( F(\boldu^L,\boldu^+,L+1) - F(\boldu^L,\boldu^{++},L)
       \right)\\[1.5mm]
   &\hskip-31mm\mbox{if } 0<\sigma_d<\tau-1;\\[1.5mm]
    q^{-\frac12(L-\Delta(a-b))}
       \left( F(\boldu^L,\boldu,L) + (q^L-1)F(\boldu^L,\boldu^+,L-1)
       \right)\\[1.5mm]
   &\hskip-31mm\mbox{if } 0<\sigma_d=\tau-1,
         \end{array} \right.
\hskip-4mm
\end{equation}
and if $p'<2p$, define:
\begin{equation}\label{Tilde2Def}
\widetilde F(\boldu^L,\boldu,L)
  =\left\{ \begin{array}{ll}
    q^{-\frac12(L+\Delta(a-b))} F(\boldu^L,\boldu,L)
     &\hskip-25mm\mbox{if } \sigma_d=0;\\[1.5mm]
       F(\boldu^L,\boldu^+,L+1)
         -q^{\frac12(L+2+\Delta(a-b))}
                   F(\boldu^L,\boldu^{++},L)\\[1.5mm]
   &\hskip-25mm\mbox{if } 0<\sigma_d<\tau-1;\\[1.5mm]
    q^{\frac12(L-\Delta(a-b))} F(\boldu^L,\boldu,L)
             +(1-q^L)F(\boldu^L,\boldu^+,L-1)\\[1.5mm]
   &\hskip-25mm\mbox{if } 0<\sigma_d=\tau-1.
         \end{array} \right.
\hskip-4mm
\end{equation}
(c.f. \cite[Eqs.\ (3.35) and (12.9)]{bms}.)
Note that the right side of (\ref{Ferm2Eq}) depends on the value of $c$ 
whereas, as will become clear later in this paper, the right side of
(\ref{FermEq}) is independent of the value of $c$.

The exceptional case of $(p,p')=(1,2)$ is dealt with by:
\begin{equation}\label{FinExceptEq}
\chi^{1,2}_{1,1,0}(L)=\chi^{1,2}_{1,1,2}(L)=\delta_{L,0},
\end{equation}
where $\delta_{i,j}$ is the Kronecker delta.
This expression (\ref{FinExceptEq}) is sometimes required when
iterating (\ref{FermEq}) and (\ref{Ferm2Eq}).

Equating expressions (\ref{FermEq}) and (\ref{Ferm2Eq})
with the corresponding instances of (\ref{FinRochaEq}),
yields polynomial identities.
These may be viewed as finitizations of the $q$-series identities
for the characters $\chi^{p,p'}_{r,s}$
obtained by equating (\ref{FermCEq}) with (\ref{RochaEq}),
because these latter identities result on taking the $L\to\infty$
limit of the above polynomial identities.
In the case of $b$ satisfying (\ref{interface}), the
polynomial identities are of genuine bosonic-fermionic type
because all terms in (\ref{FermEq}) are manifestly positive.
This is not the case for those identities obtained by
equating (\ref{Ferm2Eq}) with (\ref{FinRochaEq})
when $b$ does not satisfy (\ref{interface}),
because the fermionic side consists of linear combinations of
fundamental fermionic forms, not positive sums.
Nonetheless, in the $L\to\infty$ limit,
genuine bosonic-fermionic $q$-series identities result.


\newpage

\setcounter{section}{1}

\section{Path combinatorics}\label{CombinSec}

\subsection{Outline of proof}

In this and the following sections, we prove the fermionic expressions
stated in Section \ref{PrologueSec}.
This proof makes use of the path picture for the (finitized) Virasoro
characters $\chi^{p,p'}_{a,b,c}(L)$ and is combinatorial in the sense
that we obtain relationships between different $\chi^{p,p'}_{a,b,c}(L)$
(or in fact, similar generating functions) by manipulating the paths.
By means of these relationships, we are able to express
$\chi^{p,p'}_{a,b,c}(L)$ in terms of \lq simpler\rq\ such generating
functions. Iterating this process eventually leads to a trivial
generating function and thence to a fermionic expression for
$\chi^{p,p'}_{a,b,c}(L)$. Taking the $L\to\infty$ limit then yields
fermionic expressions for $\chi^{p,p'}_{r,s}$.
The techniques that we use refine and extend those that were developed
in \cite{flpw,foda-welsh,foda-welsh-kyoto}. The account given
here is self-contained.

In Section \ref{BandSec}, we endow the path picture of
Section \ref{IntroSec} with further structure by shading certain regions.
This shaded path picture is referred to as the $(p,p')$-model.
The weighting function that we then define for the paths
$h\in\P^{p,p'}_{a,b,c}(L)$, although apparently unrelated to that of
Forrester and Baxter \cite{forrester-baxter}, turns out to be a
renormalisation thereof
(\cite{flpw} describes how it arises via a bijection between the
Forrester-Baxter paths and the partitions with hook-difference constraints
that appear in \cite{abbbfv}).
Through this weighting function, each vertex
of the path $h$ (the $i$th vertex is the shape of the path at $(i,h_i)$)
is naturally designated as either scoring or non-scoring, with the latter
contributing $0$ to $\owt(h)$.
We actually proceed using a slightly different set of paths
$\P^{p,p'}_{a,b,e,f}(L)$ which have assigned pre- and post-segments
specified by $e,f\in\{0,1\}$. For these paths, we use a
weighting function $\mwt(h)$ that differs subtly from that defined before.
This permits a greater range of consistent combinatorial manipulations.
Thus, we carry out the bulk of our analysis under this weighting function,
only reverting to the original weighting $\owt(h)$ to obtain fermionic
expressions for $\ochi^{p,p'}_{a,b,c}(L)$ later in the paper.

After further refining the set of paths to
$\P^{p,p'}_{a,b,e,f}(L,m)$, each element of which contains exactly
$m$ non-scoring vertices, we proceed in Sections \ref{BTranSec} and
\ref{DTranSec} to define transforms that map paths between different
models.
The $\B$-transform of Section \ref{BTranSec} shows how for specific
$a',b'$, the elements of $\P^{p,p'+p}_{a',b',e,f}(L',L)$ may be
obtained combinatorially from those of $\P^{p,p'}_{a,b,e,f}(L,m)$
for various $m$.
This transform was inspired by \cite{bressoud,agarwal-bressoud}.
It has three stages.
The first stage is known as the $\B_1$-transform and enlarges the features
of a path, so that the resultant path resides in a larger model.
The second stage, referred to as a $\B_2(k)$-transform,
lengthens a path by appending $k$ pairs of segments to the path.
Each of these pairs is known as a particle.
The third stage, the $\B_3(\lambda)$-transform deforms the
path in a particular way. This process may be viewed as the
particles {\em moving} through the path.

The $\D$-transform of Section \ref{DTranSec} notes that the elements
of $\P^{p'-p, p'}_{a, b, 1-e, 1-f}(L)$ may be obtained from those
of $\P^{p, p'}_{a, b, e, f}(L)$ in a combinatorially trivial way.
In fact, it is more convenient to use the $\D$-transform
combined with the $\B$-transform.
The combined $\B\D$-transform is discussed in Section \ref{BDTranSec}.

Up to this point, the development deviates only marginally from that
of \cite{foda-welsh-kyoto}. However here, we require something more general
than the transformations of generating functions that were specified
in Corollaries 3.14 and 4.6 of \cite{foda-welsh-kyoto}.
This is achieved in Section \ref{MazySec}, where first,
for vectors $\boldmu$, $\boldnu$, $\boldmu^*$ and $\boldnu^*$,
we define $\P^{p,p'}_{a,b,e,f}(L,m)\left\{
{\boldmu^{\phantom{*}};\boldnu^{\phantom{*}}\atop\boldmu^*;\boldnu^*}
\right\}$
to be the subset of $\P^{p,p'}_{a,b,e,f}(L,m)$ whose elements
attain certain heights specified by $\boldmu$, $\boldnu$,
$\boldmu^*$ and $\boldnu^*$, in a certain order.
The bijections that the $\B$- and $\D$-transforms yield between different
such sets are specified in Lemmas \ref{MazyBijLem} and \ref{MazyDijLem}.
The relationships between the corresponding generating functions
$\mchi^{p,p'}_{a,b,e,f}(L,m)\left\{
{\boldmu^{\phantom{*}};\boldnu^{\phantom{*}}\atop\boldmu^*;\boldnu^*}
\right\}$
are given in Corollaries \ref{MazyBijCor} and \ref{MazyDijCor}.

In Section \ref{ExTrunSec}, we consider extending or truncating the
paths on either the right or the left.
This yields further relationships between different
$\mchi^{p,p'}_{a,b,e,f}(L,m)\left\{
{\boldmu^{\phantom{*}};\boldnu^{\phantom{*}}\atop\boldmu^*;\boldnu^*}
\right\}$
These relationships are given in Lemmas \ref{ExtGen2Lem},
\ref{ExtGen1Lem}, \ref{AttenGen2Lem} and \ref{AttenGen1Lem}.

We bring together all these results in Section \ref{GenFermSec}.
In Theorem \ref{CoreThrm} we state that
$\mchi^{p,p'}_{a,b,e,f}(L)\left\{
{\boldmu^{\phantom{*}};\boldnu^{\phantom{*}}\atop\boldmu^*;\boldnu^*}
\right\}$
is equal (possibly up to a factor) to the fundamental fermionic form
$F(\boldu^L,\boldu^R,L)$
where $\boldu^L$ and $\boldu^R$ are certain vectors that are related to
$\boldmu$, $\boldnu$, $\boldmu^*$ and $\boldnu^*$.
Sections \ref{CoreSec}, \ref{IndParamSec} and \ref{ProofIndSec}
are dedicated to its proof.
The centrepiece of the proof is Lemma \ref{ProofIndSec}.
The proof of this lemma is quite long and intricate,
requiring an induction argument to show many things simultaneously.

In Section \ref{TransSec}, this result is transferred back to a similar
one concerning the original weighting function.
Thereby, in Theorem \ref{Core2Thrm},
we find $F(\boldu^L,\boldu^R,L)$ to be the generating function
for a subset of $\P^{p,p'}_{a,b,c}(L)$ that attains certain heights
and does so in a certain order. However, the value of $c$ here
sometimes depends on $\boldu^R$. The purpose of Section \ref{OtherSec}
is to obtain the generating function for precisely the same
set of paths, but with the value of $c$ switched.
This generating function turns out to be the {\em fermion-like} term
$\widetilde F(\boldu^L,\boldu^R,L)$.

Section \ref{MN2sysSec} describes the $\boldm\boldn$-system which
aids the evaluation of $F(\boldu^L,\boldu^R,L)$ or
$\widetilde F(\boldu^L,\boldu^R,L)$.

As indicated above, the fermionic forms $F(\boldu^L,\boldu^R,L)$ and
$\widetilde F(\boldu^L,\boldu^R,L)$ are the generating functions
for certain subsets of $\P^{p,p'}_{a,b,c}(L)$.
In Section \ref{CollateSec}, we sum these fermionic forms
over all $\boldu^L\in\U(a)$ and $\boldu^R\in\U(b)$, where $\U(a)$
and $\U(b)$ are obtained from the Takahashi trees for $a$ and $b$
respectively as described in Section \ref{FormVSec}.
As shown in Theorems \ref{AllNode1Thrm} and \ref{AllNode2Thrm},
for many cases, this gives the required $\ochi^{p,p'}_{a,b,c}(L)$.
However, in the other cases, some paths are not accounted for.
Theorems \ref{AllNode1Thrm} and \ref{AllNode2Thrm} show
that these uncounted paths may be viewed as
elements of $\P^{\hat p,\hat p'}_{\hat a,\hat b,\hat c}(L)$ for
certain values $\hat p,\hat p',\hat a,\hat b,\hat c$.
Thus $\ochi^{\hat p,\hat p'}_{\hat a,\hat b,\hat c}(L)$ sometimes
appears in the fermionic expression for $\ochi^{p,p'}_{a,b,c}(L)$,
and in this way, the fermionic expressions that we give are recursive.
Theorem \ref{AllNode1Thrm} proves the expression (\ref{FermEq}), and
Theorem \ref{AllNode2Thrm} proves the expression (\ref{Ferm2Eq}).

In Section \ref{FermCharSec}, we use the fermionic expressions for
$\ochi^{p,p'}_{a,b,c}(L)$ to obtain fermionic expressions
for all characters $\chi^{p,p'}_{r,s}$.
In the first instance, in Section \ref{FermLimSec}, we set $a=s$,
and choose $b$ and $c$ such that (\ref{groundstatelabel}) is satisfied.
Taking the $L\to\infty$ limit of the expressions for
$\ochi^{p,p'}_{a,b,c}(L)$ then yields fermionic expressions
for $\chi^{p,p'}_{r,s}$.
These resulting expressions, stated in Theorems \ref{Limit1Thrm} and
\ref{Limit2Thrm}, still make use of the Takahashi trees for $a$ and $b$.
However, in each of the two cases, certain subsets of $\U(b)$ must
be omitted: they correspond to terms of $\ochi^{p,p'}_{a,b,c}(L)$
that are zero in the $L\to\infty$ limit, and yet their inclusion
would yield an incorrect expression for $\chi^{p,p'}_{r,s}$.

This undesirable feature of the expressions for $\chi^{p,p'}_{r,s}$,
as well as the inconveniences of finding suitable $b$ and $c$, and
considering two separate cases, is ameliorated in
Section \ref{AssimilateTrees}.
Here it is shown that the required subset of $\U(b)$ is
(apart from inconsequential differences) equal to the set $\tU(r)$ that
is obtained from the truncated Takahashi tree for $r$ which
is described in Section \ref{TruncSec}.
The final result, Theorem \ref{FinalThrm},
proves the fermionic expression (\ref{FermCEq}) stated in
Section \ref{CharSec} for an arbitrary character $\chi^{p,p'}_{r,s}$.

We discuss our results in Section \ref{DissSec}.

\subsection{The band structure}\label{BandSec}

Consider again the path picture introduced in Section \ref{IntroSec}.
The regions of the picture between adjacent heights will be known as
{\it bands}. There are $p'-2$ of these. The $h$th band is that
between heights $h$ and $h+1$.

We now assign a parity to each band: the $h$th band is said to be an {\it even}
band if $\lfloor hp/p'\rfloor=\lfloor (h+1)p/p'\rfloor$; and an {\it odd}
band if $\lfloor hp/p'\rfloor\ne\lfloor (h+1)p/p'\rfloor$.
The array of odd and even bands so obtained will be referred to
as the $(p,p')$-model.
It may immediately be deduced that the $(p,p')$-model has $p'-p-1$ even
bands and $p-1$ odd bands.
In addition, it is easily shown that for $1\le r<p$, the band lying
between heights $\lfloor rp'/p\rfloor$ and $\lfloor rp'/p\rfloor+1$
is odd: it will be referred to as the $r$th odd band.

\begin{figure}[ht]
\includegraphics[scale=1.00]{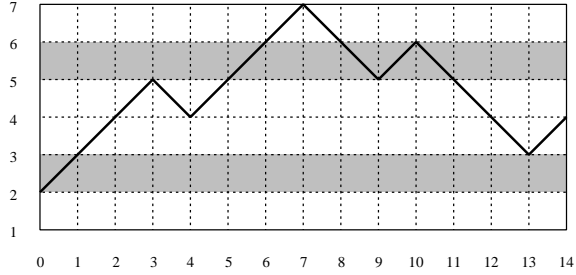}
\caption{Typical path with shaded bands.}
\label{TypicalShadedFig}
\medskip
\end{figure}

When drawing the $(p,p')$-model, we distinguish the bands by shading
the odd bands.
For example, we obtain Fig.~\ref{TypicalShadedFig} when we impose
the shading for the $(3,8)$-model on the path picture of
Fig.~\ref{TypicalBasicFig}.

We note that the band structure is up-down symmetrical.
Furthermore, if $p'>2p$ then both the 1st and $(p'-2)$th
bands are even, and there are no two adjacent odd bands.

\begin{note}\label{ExtraBandsNote}
It will be useful to consider there being a band (the $0$th)
below the bottom edge of the $(p,p')$-model, and a band (the $(p'-1)$th)
above the top edge of the $(p,p')$-model. If $p'>2p$ these two bands
are both designated as even, and if $p'<2p$ they are both designated as odd.
Thus they have the same parity as both the $1$st and $(p'-2)$th bands.
\end{note}

For $2\le a\le p'-2$, we say that $a$ is {\em interfacial\/} if
$\lfloor (a+1)p/p'\rfloor=\lfloor(a-1)p/p'\rfloor+1$.
In addition, $0$ and $p'$ are defined to be always interfacial,
and $1$ and $p'-1$ are defined to be always non-interfacial.
Thus, for $1\le a\le p'-1$,
$a$ is interfacial if and only if $a$ lies between an odd and even band
in the $(p,p')$-model.
Thus for the case of the $(3,8)$-model depicted in
Fig.~\ref{TypicalShadedFig}, $a$ is interfacial for $a=0,2,3,5,6,8$.
We define $\rho^{p,p'}(a)=\lfloor (a+1)p/p'\rfloor$ so that
if $a$ is interfacial with $1\le a<p'$
then the odd band that it borders is the $\rho^{p,p'}(a)$th odd band.
Note that $\rho^{p,p'}(0)=0$ and $\rho^{p,p'}(p')=p$.

It is easily seen that the $(p'-p,p')$-model  
differs from the $(p,p')$-model in that each band has changed parity.
It follows that if $0<a<p'$ and $a$ is interfacial in the $(p,p')$-model
then $a$ is also interfacial in the $(p'-p,p')$-model.

For $2\le a\le p'-2$, we say that $a$ is {\em multifacial\/} if
$\lfloor (a+1)p/p'\rfloor=\lfloor(a-1)p/p'\rfloor+2$.
In addition, $1$ and $p'-1$ are defined to be multifacial if and only
if $p'<2p$.
Then $a$ is multifacial if and only if $a$ lies between two odd bands.
Since there are no two adjacent odd bands when $p'>2p$,
multifacial $a$ may only occur if $p'<2p$.

Some basic results relating to the above notions are derived
in Appendix \ref{AppBSec}.

\subsection{Weighting function}\label{WeightSec}
 
\begin{figure}[ht]
  \psfrag{a}{\footnotesize$i\!-\!1$}
  \psfrag{b}{\footnotesize$i$}
  \psfrag{c}{\footnotesize$i\!+\!1$}
  \psfrag{h}{\footnotesize$h_i$}
\includegraphics[scale=1.00]{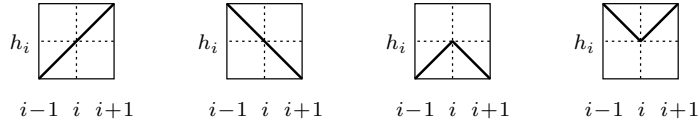}
\caption{Vertex shapes.}
\label{VertexFig}
\medskip
\end{figure}

Given a path $h$ of length $L$, for $1\le i<L$, the values
of $h_{i-1}$, $h_i$ and $h_{i+1}$
determine the shape of the $i$th vertex.
The four possible shapes are given in Fig.~\ref{VertexFig}.

\noindent
The four types of vertices shown in Fig.~\ref{VertexFig} are referred
to as a {\em straight-up} vertex, a {\em straight-down} vertex,
a {\em peak-up} vertex and a {\em peak-down} vertex respectively.
Each vertex is also assigned a parity: the parity of the
$i$th vertex is defined to be the parity of the band
that lies between heights $h_i$ and $h_{i+1}$.

For paths $h\in{\P}^{p,p'}_{a,b,c}(L)$, we define $h_{L+1}=c$,
whereupon the shape and parity of the vertex at $i=L$ is
defined as above.
This also applies if we extend the definition of ${\P}^{p,p'}_{a,b,c}(L)$
to encompass the cases $c=0$ and $c=p'$.

The weight function for the paths is best specified in terms of a
$(x,y)$-coordinate system which is inclined at $45^o$ to the original
$(i,h)$-coordinate
system and whose origin is at the path's initial point at $(i=0,h=a)$.
Specifically,
$$
x=\frac{i-(h-a)}{2},\qquad y=\frac{i+(h-a)}{2}.
$$
Note that at each step in the path, either $x$ or $y$ is incremented
and the other is constant. In this system, the path depicted
in Fig.~\ref{TypicalShadedFig} has its first few vertices at
$(0,1)$, $(0,2)$, $(0,3)$, $(1,3)$, $(1,4)$, $(1,5)$, $(1,6)$,
$(2,6)$, $\ldots$

Now, for $1\le i\le L$, we define the weight 
$c_i= c(h_{i-1},h_i,h_{i+1})$ of the $i$th vertex
according to its shape, its parity and its $(x,y)$-coordinate,
as specified in Table~\ref{WtsTable}.

\begin{table}[ht]
\begin{center}
\begin{tabular}{|c|@{\hspace{3mm}}c@{\hspace{3mm}}|c|@{\hspace{3mm}}
c@{\hspace{3mm}}|}
\hline
Vertex&
${c}_i$&
Vertex&
${c}_i$\\
\hline\hline\wombat
{\includegraphics{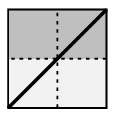}}&
\raisebox{14pt}[0pt]{$x$}&
{\includegraphics{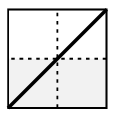}}&
\raisebox{14pt}[0pt]{$0$}\\
\hline\wombat
{\includegraphics{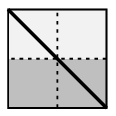}}&
\raisebox{14pt}[0pt]{$y$}&
{\includegraphics{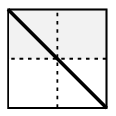}}&
\raisebox{14pt}[0pt]{$0$}\\
\hline\wombat
{\includegraphics{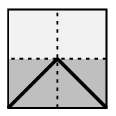}}&
\raisebox{14pt}[0pt]{$0$}&
{\includegraphics{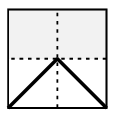}}&
\raisebox{14pt}[0pt]{$x$}\\
\hline\wombat
{\includegraphics{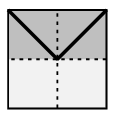}}&
\raisebox{14pt}[0pt]{$0$}&
{\includegraphics{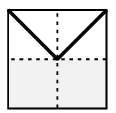}}&
\raisebox{14pt}[0pt]{$y$}\\
\hline
\end{tabular}
\end{center}
\medskip
\caption{Vertex weights.}
\label{WtsTable}
\end{table}

\noindent
In Table~\ref{WtsTable}, the lightly shaded bands can be
either even or odd bands
(including the 0th or $(p'-1)$th bands as specified in
Note \ref{ExtraBandsNote}).
Note that for each vertex shape, only one parity case contributes
non-zero weight in general.
We shall refer to those four vertices, with assigned parity, for
which in general, the weight is non-zero, as {\em scoring} vertices.
The other four vertices will be termed {\em non-scoring}.

We now define:

\begin{equation}\label{WtDef}
\owt(h)=\sum_{i=1}^L c_i.
\end{equation}

To illustrate this procedure, consider again the path
$h\in{\P}^{3,8}_{2,4,3}(14)$ depicted in Fig.~\ref{TypicalShadedFig}.
With $h_{15}=c=3$, Table~\ref{WtsTable} indicates that there are
scoring vertices at $i=3$, $4$, $5$, $7$, $8$, $13$ and $14$.
This leads to:
$$
\owt(h)=0+3+1+1+6+7+6=24.
$$

The generating function $\ochi^{p,p'}_{a,b,c}(L)$ for
the set of paths ${\P}^{p,p'}_{a,b,c}(L)$ is defined to be:
\begin{equation}\label{PathGenDef1}
\ochi^{p,p'}_{a,b,c}(L;q)
=\sum_{h\in{\P}^{p,p'}_{a,b,c}(L)} q^{\owt(h)}.
\end{equation}

\noindent
Often, we drop the base $q$ from the notation so that
$\ochi^{p,p'}_{a,b,c}(L)=\ochi^{p,p'}_{a,b,c}(L;q)$.
The same will be done for other functions without comment.

The expression (\ref{FinExceptEq}) for $\chi^{1,2}_{1,1,0}(L)$ and
$\chi^{1,2}_{1,1,2}(L)$ is immediately obtained from (\ref{PathGenDef1}).

Appendix A proves the bosonic expression for
$\ochi^{p,p'}_{a,b,c}(L)$ given in (\ref{FinRochaEq}).
With the $0$th and $(p'-1)$th bands as specified in Note \ref{ExtraBandsNote},
this result extends to the cases $c=0$ and $c=p'$
(when $b=1$ and $b=p'-1$ respectively), provided that we extend
the definition (\ref{groundstatelabel}) as follows:
\begin{equation}\label{groundstatelabel2}
r= \begin{cases}
       \lfloor pc/p'\rfloor+(b-c+1)/2 &\text{if $1\le c<p'$;}\\
             1 &\text{if $c=0$ and $p'>2p$;}\\
             0 &\text{if $c=0$ and $p'<2p$;}\\
             p-1 &\text{if $c=p'$ and $p'>2p$;}\\
             p &\text{if $c=p'$ and $p'<2p$.}
   \end{cases}
\end{equation}

\subsection{Winged generating functions}\label{WingSec}

For $h\in\P^{p,p'}_{a,b,c}(L)$, the description of the previous
section, in effect, uses the values of $b$ and $c$ to specify
a path {\em post-segment} that extends from $(L,b)$ and $(L+1,c)$.
In this section, we define another set of paths
$\P^{p, p'}_{a, b, e, f}(L)$ for which the values of $e$ and $f$
serve to specify both a {\em pre-segment} and a post-segment.

Let $p$ and $p'$ be positive coprime integers for which $1\le p<p'$.
Then, given $a,b,L\in\Z_{\ge0}$ such that $1\le a,b\le p'-1$,
$L + a - b \equiv 0$ ($\mod2$), and $e,f\in\{0,1\}$,
a path $h \in \P^{p, p'}_{a, b, e, f}(L)$ is a sequence
$h_0, h_1, h_2, \ldots, h_L,$ of integers such that:

\begin{enumerate}
\item $1 \le h_i \le p'-1$  for $0 \le i \le L$, 
\item $h_{i+1} = h_i \pm 1$ for $0 \le i <   L$, 
\item $h_0 = a, 
       h_L = b.$
\end{enumerate}
For $h \in \P^{p, p'}_{a, b, e, f}(L)$, we define $e(h)=e$ and $f(h)=f$.

As with the elements of ${\P}^{p, p'}_{a, b, c}(L)$,
each $h \in \P^{p, p'}_{a, b, e, f}(L)$ can be depicted on
the $(p,p')$-model of length $L$, by connecting the points
$(i,h_i)$ and $(i+1,h_{i+1})$ for $0\le i<L$.
For these values of $i$, the shape and parity of the $i$th vertex
is determined as in Section \ref{WeightSec}.
In addition, $f$ specifies the direction of a path post-segment
which starts at $(L,b)$ and is in the NE (resp.\ SE) direction
if $f=0$ (resp.\ $f=1$).
Similarly, $e$ specifies the direction of a path pre-segment
which ends at $(0,a)$ and is in the SE (resp.\ NE) direction
if $e=0$ (resp.\ $e=1$).
The pre- and post-segments enable a shape and a parity to be assigned
to both the zeroth and the $L$th vertices of $h$.
If a path $h$ is displayed without its pre- and post-segments,
their directions may be inferred from $e(h)$ and $f(h)$.

We now define a weight $\mwt(h)$, for
$h \in \P^{p, p'}_{a, b, e, f}(L)$.
For $1\le i<L$, set $\tilde c_i=c(h_{i-1},h_i,h_{i+1})$ using
Table~\ref{WtsTable} as in Section \ref{WeightSec}.
Then, set
\begin{equation}\label{ModWtDef}
\tilde c_L=
\left\{
  \begin{array}{ll}
x \qquad &\mbox{if } h_L-h_{L-1}=\phantom{-}1 \mbox{ and } f(h)=1;\\[1.5mm]
y \qquad &\mbox{if } h_L-h_{L-1}=-1 \mbox{ and } f(h)=0;\\[1.5mm]
0 \qquad &\mbox{otherwise,}
  \end{array} \right.
\end{equation}
where $(x,y)$ is the coordinate of the $L$th vertex of $h$.
We then designate this vertex as scoring if it is a peak
vertex ($h_L=h_{L-1}-(-1)^{f(h)}$), and as non-scoring otherwise.
We define:
\begin{equation}\label{WtDef2}
\mwt(h)=\sum_{i=1}^L \tilde c_i.
\end{equation}

Consider the corresponding path $h'\in\P^{p,p'}_{a,b,c}(L)$
with $c=b+(-1)^f$, defined by $h_i'=h_i$ for $0\le i\le L$.
{}From Table \ref{WtsTable}, we see that $\mwt(h)=\owt(h')$
if the post-segment of $h$ lies in an even band.

Define the generating function
\begin{equation}\label{PathGenDef2}
\mchi^{p,p'}_{a,b,e,f}(L;q)
=\sum_{h\in{\P}^{p,p'}_{a,b,e,f}(L)} q^{\mwt(h)},
\end{equation}
\noindent
where $\mwt(h)$ is given by (\ref{WtDef2}).
Of course, $\mchi^{p,p'}_{a,b,0,f}(L)=\mchi^{p,p'}_{a,b,1,f}(L)$.

\subsection{Striking sequence of a path}\label{StrikeSec}
 
In this section, for each $h\in{\P}^{p,p'}_{a,b,e,f}(L)$, we define
$\pi(h)$, $d(h)$, $L(h)$, $m(h)$, $\alpha(h)$ and $\beta(h)$.
Define $\pi(h)\in\{0,1\}$ to be the
parity of the band between heights $h_0$ and $h_1$
(if $L(h)=0$, we set $h_1=h_0+(-1)^{f(h)}$).
Thus, for the path $h$ shown in Fig.~\ref{TypicalShadedFig},
we have $\pi(h)=1$.
Next, define $d(h)=0$ when $h_1-h_0=1$ and $d(h)=1$ when $h_1-h_0=-1$.
We then see that if $e(h)+d(h)+\pi(h)\equiv0\,(\mod2)$ then
the $0$th vertex is a scoring vertex,
and if $e(h)+d(h)+\pi(h)\equiv1\,(\mod2)$ then it is a non-scoring vertex.

Now consider each path $h\in{\P}^{p,p'}_{a,b,e,f}(L)$ as a sequence of
straight lines, alternating in direction between NE and SE.
Then, reading from the left, let the lines be of lengths
$w_1$, $w_2$, $w_3,\ldots,w_l,$ for some $l$, with $w_i>0$ for $1\le i\le l$.
Thence $w_1+w_2+\cdots+w_l=L(h)$, where $L(h)=L$ is the length of $h$.

For each of these lines, the last vertex will be considered to be
part of the line but the first will not. Then, the $i$th of these 
lines contains $w_i$ vertices, the first $w_i-1$ of which are
straight vertices. Then write $w_i=a_i+b_i$ so that $b_i$ is the 
number of scoring vertices in the $i$th line. The striking sequence 
of $h$ is then the array:
\begin{displaymath}\label{HseqDef}
\left(\begin{array}{ccccc}
  a_1&a_2&a_3&\cdots&a_l\\ b_1&b_2&b_3&\cdots&b_l
 \end{array}\right)^{(e(h),f(h),d(h))}.
\end{displaymath}

With $\pi=\pi(h)$, $e=e(h)$, $f=f(h)$ and $d=d(h)$, we define
\begin{displaymath}
m(h)=
\left\{
  \begin{array}{ll}
       (e+d+\pi)\,\mod2+\sum_{i=1}^l a_i
          \qquad
          &\mbox{if } L>0;\\[1.5mm]
       (e+f)\,\mod2
          \qquad
          &\mbox{if } L=0,
  \end{array} \right.
\end{displaymath}
whence $m(h)$ is the number of non-scoring vertices possessed
by $h$ (altogether, $h$ has $L(h)+1$ vertices).
We also define $\alpha(h)=(-1)^d((w_1+w_3+\cdots)-(w_2+w_4+\cdots))$
and for $L>0$,
\begin{displaymath}
\beta(h)=
\left\{
  \begin{array}{l}
(-1)^d((b_1+b_3+\cdots)-(b_2+b_4+\cdots))\\
          \hskip40mm \mbox{if } e+d+\pi\equiv0\,(\mod2);\\[1.5mm]
(-1)^d((b_1+b_3+\cdots)-(b_2+b_4+\cdots))+(-1)^e\\
          \hskip40mm \mbox{otherwise.}
  \end{array} \right.
\end{displaymath}
For $L=0$, we set $\beta(h)=f-e$.

For example, for the path shown in Fig.~\ref{TypicalShadedFig} for which
$d(h)=0$ and $\pi(h)=1$, the striking sequence is:
$$
\def\qua{\hskip5pt}
\left(
{2\qua0\qua1\qua1\qua1\qua2\qua0\atop
 1\qua1\qua2\qua1\qua0\qua1\qua1}
\right)^{(e,1,0)}.
$$
In this case, $m(h)=8-e$, $\alpha(h)=2$, and $\beta(h)=2-e$.

We note that given the startpoint $h_0=a$ of the path, the path
can be reconstructed from its striking sequence\footnote{We only
need $w_1,w_2,\ldots,w_l$ together with $d$.}.
In particular, $h_L=b=a+\alpha(h)$. In addition, the nature of the
final vertex may be deduced from $a_l$ and $b_l$%
\footnote{Thus the value of $f$ in the striking sequence
is redundant --- we retain it for convenience.}

\begin{lemma}\label{WtHashLem}
Let the path $h$ have the striking sequence
$\left({a_1 \atop b_1}\,{a_2 \atop b_2}\,{a_3 \atop b_3}\,
 {\cdots\atop\cdots}\,{a_l\atop b_l} \right)^{(e,f,d)}\!,$
with $w_i=a_i+b_i$ for $1\le i\le l$.
Then
$$
\mwt(h)=\sum_{i=1}^l b_i(w_{i-1}+w_{i-3}+\cdots+w_{1+i\bmod2}).
$$
\end{lemma}

\Proof For $L=0$, both sides are clearly $0$.
So assume $L>0$.
First consider $d=0$.
For $i$ odd, the $i$th line is in the NE direction and
its $x$-coordinate is $w_2+w_4+\cdots+w_{i-1}$. By the prescription
of the previous section, and the definition of $b_i$, this line
thus contributes $b_i(w_2+w_4+\cdots+w_{i-1})$ to the weight
$\mwt(h)$ of $h$. Similarly, for $i$ even, the $i$th line is in
the SE direction and contributes $b_i(w_1+w_3+\cdots+w_{i-1})$
to $\mwt(h)$.
The lemma then follows for $d=0$.
The case $d=1$ is similar.
\cqfd
\medskip

\subsection{Path parameters}\label{PparamSec}
 
We make the following definitions:
\begin{displaymath}
\begin{array}{ll}
\alpha^{p,p'}_{a,b}&=\: b-a;\\[0.5mm]
\beta^{p,p'}_{a,b,e,f} &=\:
  \left\lfloor\frac{bp}{p'}\right\rfloor 
       - \left\lfloor\frac{ap}{p'}\right\rfloor + f-e;\\[2.5mm]
\delta^{p,p'}_{a,e} &=\:
  \left\{
    \begin{array}{ll}
       0 \quad &
            \mbox{if }
              \left\lfloor\frac{(a+(-1)^e)p}{p'}\right\rfloor
              =\left\lfloor\frac{ap}{p'}\right\rfloor;\\[3mm]
       1 \quad &
            \mbox{if }
              \left\lfloor\frac{(a+(-1)^e)p}{p'}\right\rfloor
              \ne\left\lfloor\frac{ap}{p'}\right\rfloor.
    \end{array} \right.
\\
\end{array}
\end{displaymath}

\noindent
(The superscripts of $\alpha^{p,p'}_{a,b}$ are superfluous, of course.)
It may be seen that the value of $\delta^{p,p'}_{a,e}$
gives the parity of the band in which the path pre-segment resides.

\begin{lemma}\label{BetaConstLem}
Let $h\in{\P}^{p,p'}_{a,b,e,f}(L)$.
Then $\alpha(h)=\alpha^{p,p'}_{a,b}$ and $\beta(h)=\beta^{p,p'}_{a,b,e,f}$.
\end{lemma}

\Proof That $\alpha(h)=\alpha^{p,p'}_{a,b}$ follows immediately
from the definitions.

The second result is proved by induction on $L$.
If $h\in{\P}^{p,p'}_{a,b,e,f}(0)$ then $a=b$,
whence $\beta^{p,p'}_{a,b,e,f}=f-e=\beta(h)$,
immediately from the definitions.

For $L>0$, let $h\in{\P}^{p,p'}_{a,b,e,f}(L)$ and assume that
the result holds for all $h'\in{\P}^{p,p'}_{a,b',e,f'}(L-1)$.
We consider a particular $h'$ by setting $h_i'=h_i$ for
$0\le i<L$, $b'=h_{L-1}$ and choosing $f'\in\{0,1\}$ so that
$f'=0$ if either $b-b'=1$ and the $L$th segment of $h$ lies in an
even band, or $b-b'=-1$ and the $L$th segment of $h$ lies in an
odd band; and $f'=1$ otherwise.
It may easily be checked that the $(L-1)$th vertex of $h'$ is
scoring if and only if the $(L-1)$th vertex of $h$ is scoring.
Then, from the definition of $\beta(h)$, we see that:
\begin{displaymath}
\beta(h)=
  \left\{
    \begin{array}{ll}
       \beta(h')+1 \quad &
            \mbox{if } b-b'=\phantom{-}1 \mbox{ and } f=1;\\
       \beta(h')-1 \quad &
            \mbox{if } b-b'=-1 \mbox{ and } f=0;\\
       \beta(h') \quad &
            \mbox{otherwise.}
    \end{array} \right.
\end{displaymath}
The induction hypothesis gives
$\beta(h')= \lfloor b'p/p'\rfloor - \lfloor ap/p'\rfloor +f'-e$.
Then when the $L$th segment of $h$ lies in an even band so
that $\lfloor bp/p'\rfloor=\lfloor b'p/p'\rfloor$, consideration
of the four cases of $b-b'=\pm1$ and $f\in\{0,1\}$ shows that
$\beta(h)= \lfloor bp/p'\rfloor - \lfloor ap/p'\rfloor +f-e$.
When the $L$th segment of $h$ lies in an odd band so
that $\lfloor bp/p'\rfloor=\lfloor b'p/p'\rfloor+b-b'$, consideration
of the four cases of $b-b'=\pm1$ and $f\in\{0,1\}$ again shows that
$\beta(h)= \lfloor bp/p'\rfloor - \lfloor ap/p'\rfloor +f-e$.
The result follows by induction.
\cqfd
\medskip

\subsection{Scoring generating functions}
 
We now define a generating function for paths that have a
particular number of non-scoring vertices.
First define the subset ${\P}^{p,p'}_{a,b,e,f}(L,m)$ of
${\P}^{p,p'}_{a,b,e,f}(L)$ to comprise those paths $h$ for which
$m(h)=m$. Then define:
\begin{equation}\label{ResPathGenDef}
\ochi^{p,p'}_{a,b,e,f}(L,m;q)
=\sum_{h\in{\P}^{p,p'}_{a,b,e,f}(L,m)} q^{\mwt(h)}.
\end{equation}

\begin{lemma}\label{ResPathGenLem}
Let $\beta=\beta^{p,p'}_{a,b,e,f}$. Then
\begin{displaymath}
\ochi^{p,p'}_{a,b,e,f}(L)
=\sum_{\scriptstyle m\equiv L+\beta\atop
  \scriptstyle\strut(\mbox{\scriptsize\rm mod}\,2)}
\ochi^{p,p'}_{a,b,e,f}(L,m).
\end{displaymath}
\end{lemma}

\Proof Let $h\in\P^{p,p'}_{a,b,e,f}(L)$.
We claim that $m(h)+L(h)+\beta(h)\equiv0\,(\mod2)$.
This will follow from showing that $L(h)-m(h)+(-1)^{d(h)}\beta(h)$
is even.
If $h$ has striking sequence
$\left({a_1 \atop b_1}\:{a_2 \atop b_2}\:{a_3 \atop b_3}\:
 {\cdots\atop\cdots}\:{a_l\atop b_l} \right)^{(e,f,d)},$
then $L(h)-m(h)=(b_1+b_2+\cdots+b_l)-(e+d+\pi)\,\mod2$,
where $\pi=\pi(h)$.
For $e+d+\pi\equiv0\,(\mod2)$, we immediately obtain
$L(h)-m(h)+(-1)^{d}\beta(h)=2(b_1+b_3+\ldots)$.
For $e+d+\pi\not\equiv0\,(\mod2)$, we obtain
$L(h)-m(h)+(-1)^{d}\beta(h)=2(b_1+b_3+\ldots)-1+(-1)^{d+e}$,
whence the claim is proved in all cases.
The lemma then follows, once it is noted, via Lemma \ref{BetaConstLem},
that $\beta(h)=\beta^{p,p'}_{a,b,e,f}$.
\cqfd
\medskip

\begin{note}
Since each element of ${\P}^{p,p'}_{a,b,e,f}(L,m)$ has $L+1$
vertices, it follows that
$\ochi^{p,p'}_{a,b,e,f}(L,m)$ is non-zero only if $0\le m\le L+1$.
Therefore the sum in Lemma \ref{ResPathGenLem} may be further
restricted to $0\le m\le L+1$.
\end{note}

Later, we shall define generating functions for certain
subsets of ${\P}^{p,p'}_{a,b,e,f}(L,m)$.

\subsection{A seed}
 
The following result provides a seed on which the results of
later sections will act.

\begin{lemma}\label{SeedLem}
If $L\ge0$ is even then:
\begin{displaymath}
\ochi^{1,3}_{1,1,0,0}(L,m)=
\ochi^{1,3}_{2,2,1,1}(L,m)=
\delta_{m,0}q^{\frac14 L^2}.
\end{displaymath}
If $L>0$ is odd then:
\begin{displaymath}
\ochi^{1,3}_{1,2,0,1}(L,m)=
\ochi^{1,3}_{2,1,1,0}(L,m)=
\delta_{m,0}q^{\frac14 (L^2-1)}.
\end{displaymath}
\end{lemma}

\Proof
The $(1,3)$-model comprises one even band.
Thus when $L$ is even, there is precisely one $h\in\P^{1,3}_{1,1,0,0}(L)$.
It has $h_i=1$ for $i$ even, and $h_i=2$ for $i$ odd.
We see that $h$ has striking sequence 
$\left({0 \atop 1}\:{0 \atop 1}\:{0 \atop 1}\:
 {\cdots\atop\cdots}\:{0\atop 1} \right)^{(0,0,0)}$
and $m(h)=0$.
Lemma \ref{WtHashLem} then yields
$\mwt(h)=0+1+1+2+2+3+\cdots+(\frac12L-1)+\frac12L=(L/2)^2$, as required.

The other expressions follow in a similar way.
\cqfd
\medskip

\subsection{Partitions}
 
A partition $\lambda=(\lambda_1,\lambda_2,\ldots,\lambda_k)$ is a
sequence of $k$
integer parts $\lambda_1,\lambda_2,\ldots,\lambda_k,$ satisfying
$\lambda_1\ge\lambda_2\ge\cdots\ge\lambda_k>0$.
For $i>k$, define $\lambda_i=0$.
The weight $\wt(\lambda)$ of $\lambda$ is given
by $\wt(\lambda)=\sum_{i=1}^k\lambda_i$.

We define $\Y(k,m)$ to be the set of all partitions $\lambda$
with at most $k$ parts, and for which $\lambda_1\le m$.
A proof of the following well known result may be found in
\cite{andrews-red-book}.

\begin{lemma}\label{PartitionGenLem}
The generating function,
$$
\sum_{\lambda\in\Y(k,m)} q^{\wt(\lambda)}=
\left[{m+k\atop m}\right]_q.
$$
\end{lemma}

We also require the following easily proved result.
\begin{lemma}\label{PartitionInvLem}
$$
\left[{m+k\atop m}\right]_{q^{-1}}
=q^{-mk}\left[{m+k\atop m}\right]_q.
$$
\end{lemma}

\newpage

\setcounter{section}{2}

\section{The $\B$-transform}\label{BTranSec}

In this section, we introduce the $\B$-transform which maps
paths ${\P}^{p,p'}_{a,b,e,f}(L)$ into
${\P}^{p,p'+p}_{a',b',e,f}(L')$ for certain $a'$, $b'$
and various $L'$.

The band structure of the $(p,p'+p)$-model is easily obtained from
that of the $(p,p')$-model. Indeed, according to Section \ref{BandSec},
for $1\le r<p$, the $r$th odd band of the $(p,p'+p)$-model lies
between heights $\lfloor r(p'+p)/p\rfloor=\lfloor rp'/p\rfloor+r$ and
$\lfloor r(p'+p)/p\rfloor+1=\lfloor rp'/p\rfloor+r+1$.
Thus the height of the $r$th odd band in the $(p,p'+p)$-model
is $r$ greater than that in the $(p,p')$-model.
Therefore, the $(p,p'+p)$-model may be obtained from the $(p,p')$-model
by increasing the distance between neighbouring odd bands by one
unit and appending an extra even band to both the top and the bottom
of the grid. For example, compare the $(3,8)$-model of
Fig.~\ref{TypicalShadedFig} with the $(3,11)$-model of
Fig.~\ref{B_1aFig}.

The $\B$-transform has three components, which will be termed
{\em path-dilation}, {\em particle-insertion}, and {\em particle-motion}.
These three components will also be known as the $\B_1$-, $\B_2$-
and $\B_3$-transforms respectively.
In fact, particle-insertion is dependent on a parameter $k\in\Z_{\ge0}$,
and particle-motion is dependent on a partition $\lambda$
that has certain restrictions.
Consequently, we sometimes refer to particle-insertion and particle-motion
as $\B_2(k)$- and $\B_3(\lambda)$-transforms respectively.
Then, combining the $\B_1$-, $\B_2(k)$- and $\B_3(\lambda)$-transforms
produces the $\B(k,\lambda)$-transform.

\subsection{Path-dilation: the $\B_1$-transform}

The $\B_1$-transform acts on a path $h\in{\P}^{p,p'}_{a,b,e,f}(L)$
to yield a path $h^{(0)}\in{\P}^{p,p'+p}_{a',b',e,f}(L^{(0)})$,
for certain $a'$, $b'$ and $L^{(0)}$.
First, the starting point $a'$ of the new path $h^{(0)}$ is specified
to be:
\begin{equation}\label{startdilateEq}
a'=a+\left\lfloor\frac{ap}{p'}\right\rfloor + e.
\end{equation}
If $r=\lfloor ap/p'\rfloor$ then $r$ is the number of odd bands below
$h=a$ in the $(p,p')$-model.
Since the height of the $r$th odd band in the $(p,p'+p)$-model
is $r$ greater than that in the $(p,p')$-model,
we thus see that under path-dilation, the height of the
startpoint above the next lowermost odd band
(or if there isn't one, the bottom of the grid)
has either increased by one or remained constant.

We define $d(h^{(0)})=d(h)$.
The above definition specifies that $e(h^{(0)})=e(h)$ and
$f(h^{(0)})=f(h)$.

In the case that $L=0$ and $e=f$, we specify $h^{(0)}$ by
setting $L^{(0)}=L(h^{(0)})=0$. When $L=0$ and $e\ne f$, we
leave the action of the $\B_1$-transform on $h$ undefined
(it will not be used in this case).
Thus in Lemmas \ref{BHashLem}, \ref{BWtLem}, \ref{WtShiftLem},
\ref{BresLem}, \ref{BDresLem}, \ref{MazyBijLem}, \ref{MazyDijLem},
and Corollary \ref{EndPtCor}, we implicitly exclude consideration of the
case $L=0$ and $e\ne f$.
However, it must be considered in the proofs of
Corollaries \ref{MazyBijCor} and \ref{MazyDijCor}.

In the case $L>0$ consider, as in Section \ref{StrikeSec}, $h$ to
comprise $l$ straight lines that alternate in direction, the
$i$th of which is of length $w_i$ and possesses $b_i$ scoring
vertices. $h^{(0)}$ is then defined to comprise $l$ straight
lines that alternate in direction (since $d(h^{(0)})=d(h)$, the direction
of the first line in $h^{(0)}$ is the same as that in $h$),
the $i$th of which has length
\begin{displaymath}
w'_i=
  \left\{
    \begin{array}{ll}
       w_i+b_i \quad &
            \mbox{if } i\ge2 \mbox{ or } e(h)+d(h)+\pi(h)\equiv0\,(\mod2);\\
       w_1+b_1+2\pi(h)-1 \quad &
            \mbox{if } i=1 \mbox{ and } e(h)+d(h)+\pi(h)\not\equiv0\,(\mod2).
    \end{array} \right.
\end{displaymath}
In particular, this determines $L^{(0)}=L(h^{(0)})$ and
$b'=h^{(0)}_{L^{(0)}}$.

\begin{figure}[ht]
\includegraphics[scale=1.00]{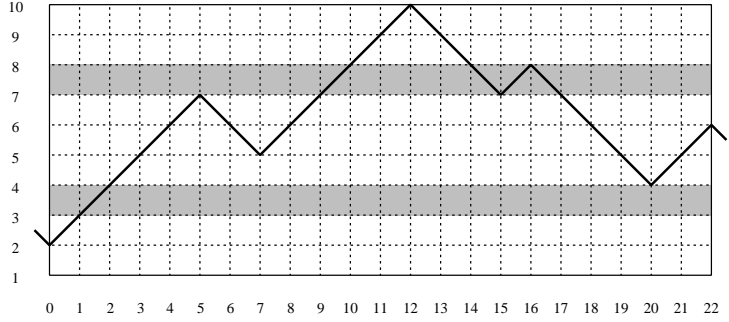}
\caption{Result of $\B_1$-transform when $e=0$ and $f=1$.}
\label{B_1aFig}
\medskip
\end{figure}

As an example, consider the path $h$ shown in Fig.~\ref{TypicalShadedFig}
as an element of $\P^{3,8}_{2, 4, e, 1}(14)$.
Here $d(h)=0$, $\pi(h)=1$ and $\lfloor ap/p'\rfloor=0$.
Thus when $e=0$, the action of path-dilation on $h$ produces the
path given in Fig.~\ref{B_1aFig}.

\noindent
This path is an element of $\P^{3,11}_{2, 6, e, 1}(22)$.

\begin{figure}[ht]
\includegraphics[scale=1.00]{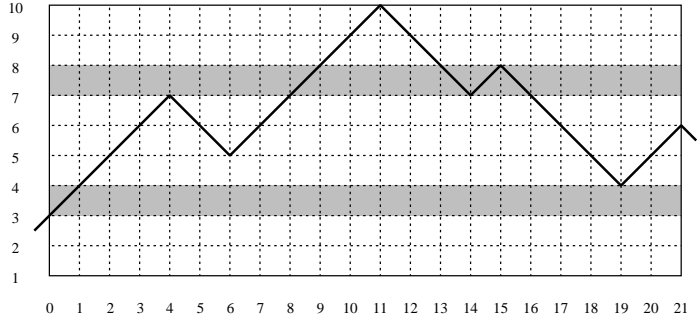}
\caption{Result of $\B_1$-transform when $e=1$ and $f=1$.}
\label{B_1bFig}
\medskip
\end{figure}

In contrast, when $e=1$, the action of path-dilation on $h$ produces the 
element of $\P^{3,11}_{3, 6, e, 1}(21)$ given in Fig.~\ref{B_1bFig}.

The situation at the startpoint may be considered as
falling into one of eight cases, corresponding to
$e(h),d(h),\pi(h)\in\{0,1\}$.\footnote{These cases may be seen
to correspond to the eight cases of vertex type as listed
in Table~\ref{WtsTable}.}
In Table~\ref{BatStartTable}, we illustrate the four cases that arise
when $d(h)=0$ (the four cases for $d(h)=1$ may be obtained
from these by an up-down reflection and changing the value
of $e(h)$).\footnote{The examples here
are such that $w_1\ge3$ and $c_1\ge3$}
\begin{table}[ht]
\bigskip\noindent\quad
\raisebox{20pt}[0pt]{$e(h)=0;\atop\pi(h)=0:$}
$\quad\includegraphics{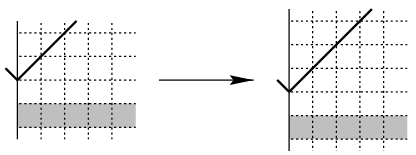}$
\qquad
\raisebox{20pt}[0pt]{$e(h)=1;\atop\pi(h)=0:$}
$\quad\includegraphics{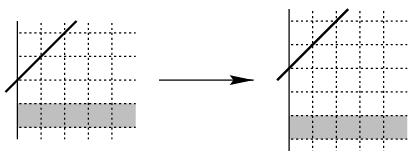}$

\bigskip\noindent\quad
\raisebox{20pt}[0pt]{$e(h)=0;\atop\pi(h)=1:$}
$\quad\includegraphics{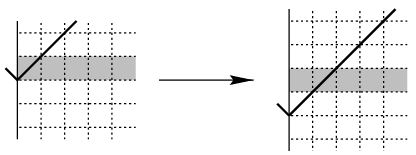}$
\qquad
\raisebox{20pt}[0pt]{$e(h)=1;\atop\pi(h)=1:$}
$\quad\includegraphics{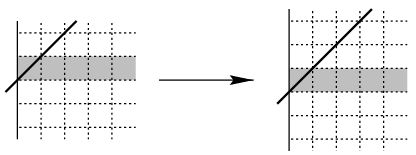}$
\medskip

\caption{$\B_1$-transforms at the startpoint.}
\label{BatStartTable}
\end{table}
{}From Table~\ref{BatStartTable} we see that the pre-segment
of $h^{(0)}$ always lies in an even band. This also follows from
Lemma \ref{StartPtLem}.

\begin{note}\label{Start1Note}
The action of path-dilation on $h\in\P^{p,p'}_{a,b,e,f}(L)$ yields a
path $h^{(0)}\in\P^{p,p'+p}_{a',b',e,f}(L^{(0)})$
that has, including the vertex at $i=0$, no adjacent scoring vertices,
except in the case where $\pi(h)=1$ {\em and} $e(h)=d(h)$, when a single
pair of scoring vertices occurs in $h^{(0)}$ at $i=0$ and $i=1$.

Also note that $\pi(h^{(0)})=\pi(h)$ unless $\pi(h)=1$ {\em and}
$e(h)=d(h)$, in which case $\pi(h^{(0)})=0$.
\end{note}

Now compare the $i$th line of $h^{(0)}$ (which has length $w'_i$)
with the $i$th line of $h$ (which has length $w_i$).
If the lines in question are in the NE direction, we claim that
the height of the final vertex of that in $h^{(0)}$ above
the next lowermost odd band (or if there isn't one, the bottom of
the model at $a=1$) is one greater than that in $h$.

If the lines in question are in the SE direction, we claim that
the height of the final vertex of that in $h^{(0)}$ below
the next highermost odd band (or if there isn't one, the top
of the model) is one greater than that in $h$.
In particular, if either the first or last segment
of the $i$th line is in an odd band, then the corresponding segment
of $h^{(0)}$ lies in the same odd band.

We also claim that if that of $h$ has a straight vertex that passes
into the $k$th odd band in the $(p,p')$-model then
that of $h^{(0)}$ has a straight vertex that passes into the
$k$th odd band in the $(p,p'+p)$-model.

Comparing Figs.~\ref{B_1aFig} and \ref{B_1bFig} with
Fig.~\ref{TypicalShadedFig} demonstrates these claims.

These claims follow because in passing from the $(p,p')$-model
to the $(p,p'+p)$-model, the distance between neighbouring odd bands
has increased by one, and because the length of each line has
increased by one for every scoring vertex with a small
adjustment made to the length of the first line in certain cases.
In effect, a new straight vertex has been inserted
immediately prior to each scoring vertex
and, if $e(h)+d(h)+\pi(h)\not\equiv0\,(\mod2)$,
adjusting the length of the resulting first line by $2\pi(h)-1$.

\begin{lemma}\label{BHashLem}

Let $h \in \P^{p, p'}_{a,b,e,f}(L)$
have striking sequence
$\left({a_1\atop b_1}\:{a_2\atop b_2}\:{a_3\atop b_3}\:
 {\cdots\atop\cdots}\:{a_l\atop b_l} \right)^{(e,f,d)},$
and let
$h^{(0)} \in \P^{p,p'+p}_{a',b',e,f}(L^{(0)})$ be obtained from
the action of the $\B_1$-transform on $h$.
Let $\pi=\pi(h)$. If $e+d+\pi\equiv0\,(\mod2)$ then $h^{(0)}$ has
striking sequence:
\begin{displaymath}
\left(\begin{array}{ccccc}
  a_1+b_1&a_2+b_2&a_3+b_3&\cdots&a_l+b_l\\
  b_1&b_2&b_3&\cdots&b_l
 \end{array}\right)^{(e,f,d)},
\end{displaymath}
and if $e+d+\pi\not\equiv0\,(\mod2)$ then $h^{(0)}$ has
striking sequence:
\begin{displaymath}
\left(\begin{array}{ccccc}
  a_1+b_1+\pi-1&a_2+b_2&a_3+b_3&\cdots&a_l+b_l\\
  b_1+\pi&b_2&b_3&\cdots&b_l
 \end{array}\right)^{(e,f,d)}.
\end{displaymath}
Moreover, if $m=m(h)$:
\begin{itemize}
\item $m(h^{(0)})=L$;
\item $L^{(0)}=
\left\{
  \begin{array}{ll}
2L-m+2 \quad &
          \mbox{if } \pi=1 \mbox{ and } e=d,\\[1.5mm]
2L-m \quad & \mbox{otherwise};
  \end{array} \right. $
\item $\alpha(h^{(0)})=\alpha(h)+\beta(h)$;
\item $\beta(h^{(0)})=\beta(h)$.
\end{itemize}
\end{lemma}

\Proof The form of the striking sequence for $h^{(0)}$ follows
because, for $i>1$, every scoring vertex in the $i$th line of $h$
accounts for an extra non-scoring vertex in that line. The same is true
when $i=1$, except in the case $(e(h)+d(h)+\pi(h))\equiv1$
(throughout this proof, we take all equivalences, modulo 2)
when the length of the new $1$st line becomes $a_1+2b_1+2\pi-1$.
That there are $b_1+\pi$ scoring vertices in this case,
follows from examining Table~\ref{BatStartTable}.

Let $\pi'=\pi(h^{(0)})$.
By definition, $e(h^{(0)})=e$ and $d(h^{(0)})=d$.

If $(e+d+\pi)\equiv0$ then $(e+d+\pi')\equiv0$ by Note \ref{Start1Note}.
Thereupon $m^{(0)}=\sum_{i=1}^l(a_i+b_i)=L$.
Additionally, $L^{(0)}=\sum_{i=1}^l(a_i+2b_i)=2L-\sum_{i=1}^la_i=2L-m$.
That $\beta(h^{(0)})=\beta(h)$ and
$\alpha(h^{(0)})=\alpha(h)+\beta(h)$ both follow immediately
in this case.

On the other hand, if $(e+d+\pi)\not\equiv0$ then
$\pi=0\implies e\ne d$ and $\pi=1\implies e=d$.
In each instance, Note \ref{Start1Note} implies that $\pi'=0$.
Thereupon,
$m^{(0)}=(e+d+\pi')\,\mod2+\pi-1+\sum_{i=1}^l(a_i+b_i)
=\sum_{i=1}^l(a_i+b_i)=L$.
Additionally,
$L^{(0)}=2\pi-1+\sum_{i=1}^l(a_i+2b_i)=2L-(1+\sum_{i=1}^la_i)+2\pi=2L-m+2\pi$.
This is the required value.
Now in this case,
$\beta(h)=(-1)^d((b_1+b_3+\cdots)-(b_2+b_4+\cdots))+(-1)^e$.
When $\pi=0$ so that $(e+d+\pi')\equiv1$ then
$\beta(h^{(0)})=\beta(h)$ follows immediately.
When $\pi=1$, we have
$\beta(h^{(0)})=(-1)^d((b_1+1+b_3+\cdots)-(b_2+b_4+\cdots))$.
$\beta(h^{(0)})=\beta(h)$ now follows in this case because
$(e+d+\pi)\not\equiv0$ implies that $e=d$.
Finally,
$\alpha(h^{(0)})=\alpha(h)+(-1)^d((b_1+b_3+\cdots)-(b_2+b_4+\cdots))
+(-1)^d(2\pi-1)$.
Since $(-1)^d(2\pi-1)=-(-1)^d(-1)^\pi=(-1)^e$,
the lemma then follows.
\cqfd
\medskip

\begin{corollary}\label{EndPtCor}
Let $h \in{\P}^{p,p'}_{a,b,e,f}(L)$
and
$h^{(0)} \in \P^{p,p'+p}_{a',b',e,f}(L^{(0)})$ be the path obtained
by the action of the $\B_1$-transform on $h$.
Then $a'=a+\lfloor ap/p'\rfloor+e$ and $b'=b+\lfloor bp/p'\rfloor+f$.
\end{corollary}
 
\Proof $a'=a+\lfloor ap/p'\rfloor+e$ is by definition.
Lemma \ref{BHashLem} gives $\alpha(h^{(0)})=\alpha(h)+\beta(h)$,
whence Lemma \ref{BetaConstLem} implies that
$\alpha^{p,p'+p}_{a',b'}=\alpha^{p,p'}_{a,b}+\beta^{p,p'}_{a,b,e,f}$.
Expanding this gives
$b'-a'=b-a+\lfloor bp/p'\rfloor-\lfloor ap/p'\rfloor+f-e$,
whence $b'=b+\lfloor bp/p'\rfloor+f$.
\cqfd
\medskip
 
\noindent
The above result implies that the $\B_1$-transform maps
${\P}^{p,p'}_{a,b,e,f}(L)$ into a set of paths whose startpoints are
equal and whose endpoints are equal.
However, the lengths of these paths are not necessarily equal.
We also see that the transformation of the endpoint is analogous
to that which occurs at the startpoint.
In particular, Lemma \ref{StartPtLem} implies that
$\delta^{p,p'+p}_{b',f}=0$ so that the path post-segment of $h^{(0)}$
always resides in an even band.
For the four cases where $h_L=h_{L-1}-1$, the $\B_1$-transform
affects the endpoint as in Table~\ref{BatEndTable}
(the value $\pi'(h)$ is the parity of the band in which the
$L$th segment of $h$ lies).

\begin{table}[ht]
\bigskip\noindent
\hbox{
\raisebox{20pt}[0pt]{$f(h)=0;\atop\pi'(h)=0:$}
$\quad\includegraphics{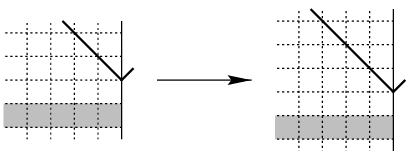}$
\qquad
\raisebox{20pt}[0pt]{$f(h)=1;\atop\pi'(h)=0:$}
$\quad\includegraphics{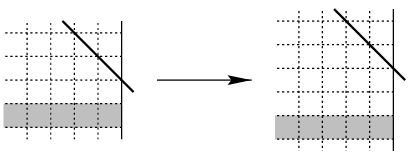}$}

\bigskip\noindent
\hbox{
\raisebox{20pt}[0pt]{$f(h)=0;\atop\pi'(h)=1:$}
$\quad\includegraphics{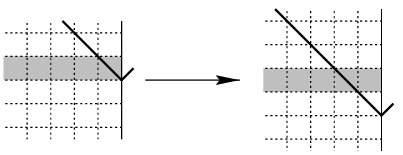}$
\qquad
\raisebox{20pt}[0pt]{$f(h)=1;\atop\pi'(h)=1:$}
$\quad\includegraphics{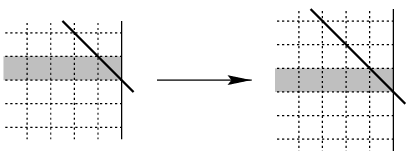}$}
\medskip

\caption{$\B_1$-transforms at the endpoint.}
\label{BatEndTable}
\end{table}

\begin{lemma}\label{ParamLem}
Let $1\le p<p'$, $1\le a,b<p'$, $e,f\in\{0,1\}$,
$a'=a+\lfloor ap/p'\rfloor+e$, and $b'=b+\lfloor bp/p'\rfloor+f$.
Then
$\alpha^{p,p'+p}_{a',b'}=\alpha^{p,p'}_{a,b}+\beta^{p,p'}_{a,b,e,f}$
and
$\beta^{p,p'+p}_{a',b',e,f}=\beta^{p,p'}_{a,b,e,f}$.
\end{lemma}
 
\Proof
Lemma \ref{StartPtLem} implies that
$\lfloor a'p/(p'+p)\rfloor=\lfloor ap/p'\rfloor$,
$\lfloor b'p/(p'+p)\rfloor=\lfloor bp/p'\rfloor$.
The results then follow immediately from the definitions.
\cqfd
\medskip

\begin{lemma}\label{BWtLem}
Let $h \in{\P}^{p,p'}_{a,b,e,f}(L)$ 
and
$h^{(0)} \in \P^{p,p'+p}_{a',b',e,f}(L^{(0)})$ be the path obtained 
by the action of the $\B_1$-transform on $h$.
Then
$$
\mwt(h^{(0)})=\mwt(h)+\frac{1}{4}\left((L^{(0)}-m^{(0)})^2-\beta^2\right),
$$
where $m^{(0)}=m(h^{(0)})$ and $\beta=\beta^{p,p'}_{a,b,e,f}$.
\end{lemma}

\Proof Let $h$ have striking sequence
$\left({a_1\atop b_1}\:{a_2\atop b_2}\:{a_3\atop b_3}\:
 {\cdots\atop\cdots}\:{a_l\atop b_l} \right)^{(e,f,d)}$,
and let $\pi=\pi(h)$.
If $(e+d+\pi)\equiv0\,(\mod2)$, then
Lemmas \ref{BHashLem} and \ref{WtHashLem} show that
$$
\mwt(h^{(0)})-\mwt(h)=(b_1+b_3+b_5+\cdots)(b_2+b_4+b_6+\cdots).
$$
Via Lemma \ref{BHashLem}, we obtain
$L^{(0)}-m^{(0)}=L-m(h)=b_1+b_2+\cdots+b_l$.
Then since $\beta(h)=\pm((b_1+b_3+b_5+\cdots)-(b_2+b_4+b_6+\cdots))$,
it follows that
$$
\mwt(h^{(0)})-\mwt(h)=\frac{1}{4}((L^{(0)}-m^{(0)})^2-\beta(h)^2).
$$

If $(e+d+\pi)\not\equiv0\,(\mod2)$, then
Lemmas \ref{BHashLem} and \ref{WtHashLem} show that
\begin{align*}
\mwt(h^{(0)})-\mwt(h)&=(2\pi-1+b_1+b_3+b_5+\cdots)(b_2+b_4+b_6+\cdots)\\
&=\frac{1}{4}((L^{(0)}-m^{(0)})^2-\beta(h)^2),
\end{align*}
the second equality resulting because
$L^{(0)}-m^{(0)}=L-m(h)+2\pi=b_1+b_2+\cdots+b_l+2\pi-1$ and
\begin{align*}
\beta(h)&=(-1)^d((b_1+b_3+b_5+\cdots)-(b_2+b_4+b_6+\cdots))+(-1)^e\\
&=\pm((2\pi-1+b_1+b_3+b_5+\cdots)-(b_2+b_4+b_6+\cdots)),
\end{align*}
on using $(-1)^{e+d}=-(-1)^\pi=2\pi-1$.

Finally, Lemma \ref{BetaConstLem} gives $\beta(h)=\beta^{p,p'}_{a,b,e,f}=
\beta$.
\cqfd
\medskip

\subsection{Particle insertion: the $\B_2$-transform}\label{InsertSec}

Let $p'>2p$ so that the $(p,p')$-model has no two neighbouring odd bands,
and let $a',e$ be such that $\delta^{p,p'}_{a',e}=0$.
Then if $h^{(0)}\in\P^{p,p'}_{a',b',e,f}(L^{(0)})$,
the pre-segment of $h^{(0)}$ lies in an even band.
By {\em inserting a particle}
into $h^{(0)}$, we mean displacing $h^{(0)}$ two positions to
the right and inserting two segments: the leftmost of these
is in the NE (resp.~SE) direction if $e=0$ (resp.~$e=1$),
and the rightmost is in the opposite direction, which is thus
the direction of the pre-segment of $h^{(0)}$.
In this way, we obtain a path $h^{(1)}$ of length $L^{(0)}+2$.
We assign $e(h^{(1)})=e$ and $f(h^{(1)})=f$.
Note also that $d(h^{(1)})=e$ and $\pi(h^{(1)})=0$.

Thereupon, we may repeat this process of particle insertion.
After inserting $k$ particles into $h^{(0)}$,
we obtain a path $h^{(k)}\in\P^{p,p'}_{a',b',e,f}(L^{(0)}+2k)$.
We say that $h^{(k)}$ has been obtained 
by the action of a $\B_2(k)$-transform on $h^{(0)}$.

\begin{figure}[ht]
\includegraphics[scale=1.00]{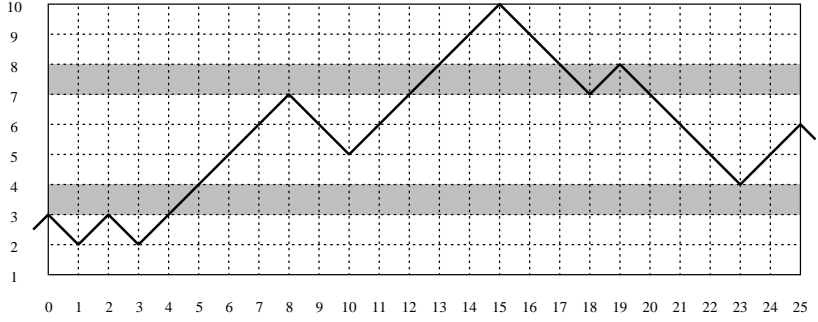}
\caption{Result of $\B_2(k)$-transform.}
\label{B_2fig}
\medskip
\end{figure}

In the case of the element of  $\P^{3,11}_{3, 6, 1, 1}(21)$ shown
in Fig.~\ref{B_1bFig}, the insertion of two particles produces the
element of $\P^{3,11}_{3, 6, 1, 1}(25)$ shown in Fig.~\ref{B_2fig}.
 
\begin{lemma}\label{WtShiftLem}
Let $h \in \P^{p,p'}_{a,b,e,f}(L)$.
Apply a $\B_1$-transform to $h$ to obtain the path
$h^{(0)} \in \P^{p,p'+p}_{a',b',e,f}(L^{(0)})$. 
Then obtain $h^{(k)}\in\P^{p,p'+p}_{a',b',e,f}(L^{(k)})$ by
applying a $\B_2(k)$-transform to $h^{(0)}$.
If $m^{(k)}=m(h^{(k)})$,
then $L^{(k)}=L^{(0)}+2k$, $m^{(k)}=m^{(0)}$ and
$$
\mwt(h^{(k)})=\mwt(h)+\frac{1}{4}((L^{(k)}-m^{(k)})^2-\beta^2),
$$
where $\beta=\beta^{p,p'}_{a,b,e,f}$.
\end{lemma}

\Proof
That $L^{(k)}=L^{(0)}+2k$ follows immediately from the
definition of a $\B_2(k)$-transform.
Lemma \ref{BWtLem} yields:
\begin{equation}\label{Proof1Eq}
\mwt(h^{(0)})=\mwt(h)+\frac{1}{4}\left((L^{(0)}-m(h^{(0)}))^2-\beta^2\right).
\end{equation}
Let the striking sequence of $h^{(0)}$ be
$\left({a_1\atop b_1}\:{a_2\atop b_2}\:
 {\cdots\atop\cdots}\:{a_l\atop b_l} \right)^{(e,f,d)},$
and let $\pi=\pi(h^{(0)})$.
 
If $e=d$, we are restricted to the case $\pi=0$,
since $\delta^{p,p'+p}_{a',e}=0$ by Lemma \ref{StartPtLem}.
The striking sequence of $h^{(1)}$ is then
$\left({0\atop1}\:{0\atop1}\:{a_1\atop b_1}\:{a_2\atop b_2}\:
 {\cdots\atop\cdots}\:{a_l\atop b_l} \right)^{(e,f,e)}$.
Thereupon $m(h^{(1)})=\sum_{i=1}^l a_i=m(h^{(0)})$.
In this case, Lemma \ref{WtHashLem} shows that
$\mwt(h^{(1)})-\mwt(h^{(0)})=1+b_1+b_2+\cdots+b_l=L^{(0)}-m^{(0)}+1$.
 
If $e\ne d$, the striking sequence of $h^{(1)}$ is
$\left({0\atop1}\:{a_1+1-\pi\atop b_1+\pi}\:{a_2\atop b_2}\:
 {\cdots\atop\cdots}\:{a_l\atop b_l} \right)^{(e,f,e)}$.
Then $m(h^{(1)})=1-\pi+\sum_{i=1}^la_i$
which equals $m(h^{(0)})=(e+d+\pi)\,\mod2+\sum_{i=1}^la_i$
for both $\pi=0$ and $\pi=1$.
Here, Lemma \ref{WtHashLem} shows that
$\mwt(h^{(1)})-\mwt(h^{(0)})=\pi+b_1+b_2+\cdots+b_l$.
Since $L^{(0)}-m^{(0)}=-(e+d+\pi)\,\mod2+b_1+b_2+\cdots+b_l$,
we once more have
$\mwt(h^{(1)})-\mwt(h^{(0)})=L^{(0)}-m^{(0)}+1$.
 
Repeated application of these results, yields
$m(h^{(k)})=m(h^{(0)})$ and
$$
\mwt(h^{(k)})=\mwt(h^{(0)})+k\left(L^{(0)}-m(h^{(0)})\right)+k^2.
$$
Then, on using (\ref{Proof1Eq}) and $L^{(k)}=L^{(0)}+2k$, the lemma follows.
\cqfd
\medskip

\subsection{Particle moves}\label{PartMovesSec}

In this section, we once more restrict to the case
$p'>2p$ so that the $(p,p')$-model has no two neighbouring odd bands,
and consider only paths $h\in\P^{p,p'}_{a',b',e,f}(L')$ for which
either $d(h)=e$ or $\pi(h)=0$.

We specify four types of local deformations of a path.
These deformations will be known as {\em particle moves}.
In each of the four cases, 
a particular sequence of contiguous segments of a path is changed to a 
different sequence, the remainder of the path being unchanged. The 
moves are as follows --- the path portion to the left of the arrow is 
changed to that on the right:

\bigskip
\centerline{\includegraphics[scale=1.00]{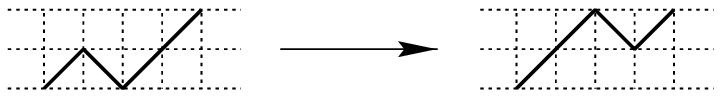}}
\nobreak
\centerline{Move 1.}
\bigskip

\centerline{\includegraphics[scale=1.00]{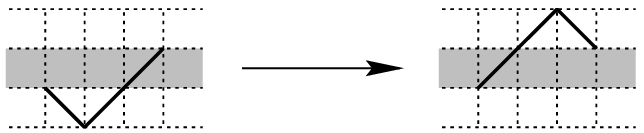}}
\nobreak
\centerline{Move 2.}
\bigskip

\centerline{\includegraphics[scale=1.00]{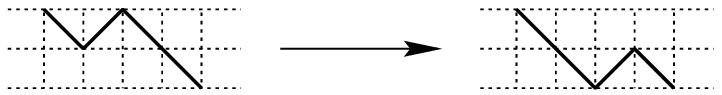}}
\nobreak
\centerline{Move 3.}
\bigskip

\centerline{\includegraphics[scale=1.00]{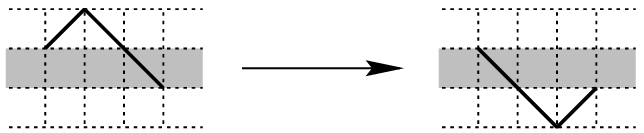}}
\nobreak
\centerline{Move 4.}
\bigskip

\noindent Since $p'>2p$, each odd band is straddled by a pair
of even bands. Thus, there is no impediment to enacting moves
2 and 4 for paths in $\P^{p,p'}_{a',b',e,f}(L)$.
Note that moves 1 and 2 are reflections of moves 3 and 4.

In addition to the four moves described above, we permit certain
additional deformations of a path close to its left and right extremities in
certain circumstances.
Each of these moves will be referred to as an {\em edge-move}.
They, together with their validities, are as follows:
 
\bigskip
\centerline{\includegraphics{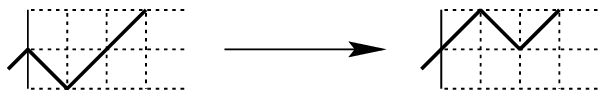}}
\nobreak
\centerline{Edge-move 1: applies only if $e=1$.}
\bigskip

\goodbreak
\centerline{\includegraphics{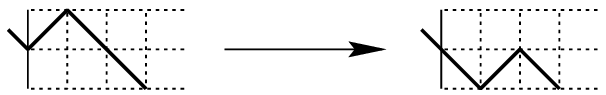}}
\nobreak
\centerline{Edge-move 2: applies only if $e=0$.}
\bigskip

\goodbreak
\centerline{\includegraphics{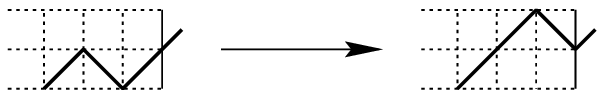}}
\nobreak
\centerline{Edge-move 3: applies only if $f=0$.}
\bigskip

\goodbreak
\centerline{\includegraphics{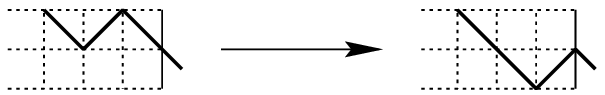}}
\nobreak
\centerline{Edge-move 4: applies only if $f=1$.}
\bigskip

With the pre- and post-segments drawn, we see that
these edge-moves may be considered as instances of
moves 1 and 3 described beforehand.

\begin{lemma}\label{WtChngLem}
Let the path $\hat h$ differ from the path $h$ in that
one change has been made according to one of the
four moves described above, or to one of the four edge-moves
described above (subject to their restrictions).
Then
$$
\mwt(\hat h)=\mwt(h)+1.
$$
Additionally, $L(\hat h)=L(h)$ and $m(\hat h)=m(h)$.
\end{lemma}
\Proof For each of the four moves and four edge-moves, take the
$(x,y)$-coordinate of the leftmost point of the depicted portion
of $h$ to be $(x_0,y_0)$. Now consider 
the contribution to the weight of the three vertices in question 
before and after the move (although the vertex immediately before 
those considered may change, its contribution doesn't). In each 
of the eight cases, the contribution is $x_0+y_0+1$ before the move 
and $x_0+y_0+2$ afterwards. Thus $\mwt(\hat h)=\mwt(h)+1$.
The other statements are immediate on inspecting all eight moves.
\cqfd
\medskip

Now observe that for each of the eight moves specified above, the
sequence of path segments before the move consists of an adjacent
pair of scoring vertices followed by a non-scoring vertex.
The specified move replaces this combination with a non-scoring
vertex followed by two scoring vertices.
As anticipated above, the pair of adjacent scoring vertices is
viewed as a particle. Thus each of the above eight moves describes
a particle moving to the right by one step.
 
\begin{figure}[ht]
\includegraphics[scale=1.00]{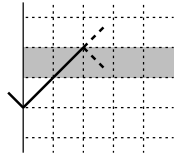}
\caption{Not a particle}
\label{NotPFig}
\medskip
\end{figure}

When $p'>2p$, so that there are no two adjacent odd bands in the
$(p,p')$-model, and noting that $\delta^{p,p'}_{b',f}=0$, we
see that each sequence comprising two scoring vertices followed
by a non-scoring vertex is present amongst the eight configurations
prior to a move, except for the case depicted in Fig.~\ref{NotPFig}
and its up-down reflection.
 
\noindent
Only in these cases, where the 0th and 1st segments are scoring and
the first two segments are in the same direction, do we {\em not}
refer to the adjacent pair of scoring vertices as a particle.
 
Also note that when $p'>2p$ and $\delta^{p,p'}_{a',e}=
\delta^{p,p'}_{b',f}=0$, each
sequence of a non-scoring vertex followed by two scoring vertices
appears amongst the eight configurations that result from a move.
In such cases, the move may thus be reversed.

\subsection{The $\B_3$-transform}

Since in each of the moves described in Section \ref{PartMovesSec},
a pair of scoring vertices shifts to the right by one step,
we see that a succession of such 
moves is possible until the pair is followed by another scoring 
vertex. If this itself is followed by yet another scoring vertex,
we forbid further movement. However, if it is followed by a non-scoring 
vertex, further movement is allowed after considering the latter two 
of the three consecutive scoring vertices to be the particle (instead 
of the first two).

As above, let $h^{(k)}$ be a path resulting from a $\B_2(k)$-transform
acting on a path that itself is the image of a $\B_1$-transform.
We now consider moving the $k$ particles that have been inserted.

\begin{lemma}\label{GaussLem} Let $\delta^{p,p'}_{b',f}=0$.
There is a bijection between the set of paths obtained by moving the
particles in $h^{(k)}$ and $\Y(k,m)$, where $m=m(h^{(k)})$.
This bijection is such that if $\lambda\in\Y(k,m)$ is the
bijective image of a particular $h$ then
$$
\mwt(h)=\mwt(h^{(k)})+\wt(\lambda).
$$
Additionally, $L(h)=L(h^{(k)})$ and $m(h)=m(h^{(k)})$.
\end{lemma}

\Proof
Since each particle moves by traversing a non-scoring vertex,
and there are $m$ of these to the right of the rightmost particle 
in $h^{(k)}$, and there are no consecutive scoring vertices to its 
right, this particle can make $\lambda_1$ moves to the right, with 
$0\le\lambda_1\le m$. Similarly, the next rightmost particle can 
make $\lambda_2$ moves to the right with $0\le\lambda_2\le\lambda_1$.
Here, the upper restriction arises because the two scoring vertices 
would then be adjacent to those of the first particle. Continuing in 
this way, we obtain that all possible final positions of the particles 
are indexed by $\lambda=(\lambda_1,\lambda_2,\ldots,\lambda_k)$ with 
$m\ge\lambda_1\ge\lambda_2\ge\cdots\ge\lambda_k\ge0$, that is, by 
partitions of at most $k$ parts with no part exceeding $m$. Moreover, 
since by Lemma \ref{WtChngLem} the weight increases by one for each 
move, the weight increase after the sequence of moves specified by 
a particular $\lambda$ is equal to $\wt(\lambda)$.
The final statement also follows from Lemma \ref{WtChngLem}.
\cqfd
\medskip

\noindent We say that a path obtained by moving the particles
in $h^{(k)}$ according to the partition $\lambda$ has been obtained
by the action of a $\B_3(\lambda)$-transform.

Having defined $\B_1$, $\B_2(k)$ for $k\ge0$ and $\B_3(\lambda)$
for $\lambda$ a partition with at most $k$ parts, we now define a 
$\B(k,\lambda)$-transform as the composition 
$\B(k,\lambda)=\B_3(\lambda)\circ\B_2(k)\circ\B_1$.

\begin{lemma}\label{BresLem}
Let $h' \in \P^{p,p'+p}_{a',b',e,f}(L')$ be obtained from
$h \in \P^{p,p'}_{a,b,e,f}(L)$ by the action of the $\B(k,\lambda)$-transform.
If $\pi=\pi(h)$ and $m=m(h)$ then:
\begin{displaymath}
\begin{array}{l}
\bullet\quad L'=
  \left\{
    \begin{array}{ll}
  2L-m+2k+2 \quad &
            \mbox{if } \pi=1 \mbox{ and } e=d,\\[1.5mm]
  2L-m+2k \quad & \mbox{otherwise};
    \end{array} \right.
\\[1mm]
\bullet\quad m(h')=L;\\[1mm]
\bullet\quad \mwt(h')=\mwt(h) + \frac{1}{4}\left( (L'-L)^2
 - \beta^2 \right) + \wt(\lambda),
\end{array}
\end{displaymath}
where $\beta=\beta^{p,p'}_{a,b,e,f}$.
\end{lemma}

\Proof These results follow immediately from
Lemmas \ref{BHashLem}, \ref{WtShiftLem} and \ref{GaussLem}.
\cqfd
\medskip

\begin{note}\label{Start2Note}
Since particle insertion and the particle moves don't change the
startpoint, endpoint or value $e(h)$ or $f(h)$ of a path $h$,
then in view of Lemma \ref{StartPtLem} and Corollary \ref{EndPtCor},
we see that the action of a $\B$-transform on
$h\in\P^{p,p'}_{a,b,e,f}(L)$
yields a path $h'\in\P^{p,p'+p}_{a',b',e,f}(L')$,
where
$a'=a+\lfloor ap/p'\rfloor+e$, $b'=b+\lfloor bp/p'\rfloor+f$,
and $\delta^{p,p'+p}_{a',e}=\delta^{p,p'+p}_{b',f}=0$.
\end{note}

\subsection{Particle content of a path}\label{ParticleSec}

We again restrict to the case
$p'>2p$ so that the $(p,p')$-model has no two neighbouring odd bands,
and let $h'\in\P^{p,p'}_{a',b',e,f}(L')$.
In the following lemma, we once more restrict to the cases for which
$\delta^{p,p'}_{a',e}=\delta^{p,p'}_{b',f}=0$,
and thus only consider the cases for which
the pre-segment and the post-segment of $h'$ lie in even bands.

\begin{lemma}\label{UniqueLem}
For $1\le p<p'$ with $p'>2p$, let $1\le a',b'<p'$ and $e,f\in\{0,1\}$,
with $\delta^{p,p'}_{a',e}=\delta^{p,p'}_{b',f}=0$.
Set $a=a'-\lfloor a'p/p'\rfloor-e$
and $b=b'-\lfloor b'p/p'\rfloor-f$.
If $h'\in\P^{p,p'}_{a',b',e,f}(L')$,
then there is a unique triple $(h,k,\lambda)$
where $h\in\P^{p,p'-p}_{a,b,e,f}(L)$ for some $L$, such that
the action of a $\B(k,\lambda)$-transform on $h$ results in $h'$.
\end{lemma}

\Proof This is proved by reversing the constructions described
in the previous sections. Locate the leftmost pair of consecutive
scoring vertices in $h'$, and move them leftward by reversing the
particle moves, until they occupy the $0$th and $1$st positions.
This is possible in all cases when
$\delta^{p,p'}_{a',e}=\delta^{p,p'}_{b',f}=0$.
Now ignoring these two vertices, do the same with the next leftmost
pair of consecutive scoring vertices, moving them leftward until
they occupy the second and third positions. Continue in this way
until all consecutive scoring vertices occupy the leftmost positions
of the path.
Denote this path by $h^{(\cdot)}$.
At the leftmost end of $h^{(\cdot)}$, there will be a number of even
segments (possibly zero) alternating in direction.
Let this number be $2k$ or $2k+1$ according to whether is it even or odd.
Clearly $h'$ results from $h^{(\cdot)}$
by a $\B_3(\lambda)$-transform for a particular $\lambda$ with at most
$k$ parts.
 
Removing the first $2k$ segments of $h^{(\cdot)}$ yields a path
$h^{(0)}\in\P^{p,p'}_{a',b',e,f}$. This path thus has no
two consecutive scoring vertices, except possibly at the $0$th
and $1$st positions, and then only if the first vertex is a straight
vertex (as in Fig.~\ref{NotPFig}).
Moreover, $h^{(k)}$ arises by the action of
a $\B_2(k)$-transform on $h^{(0)}$.
 
Ignoring for the moment the case where there are scoring vertices at the
$0$th and $1$st positions, $h^{(0)}$ has by construction no pair of
consecutive scoring vertices.
Therefore, beyond the $0$th vertex, we may remove a non-scoring vertex
before every scoring vertex (and increment the length of the first line if
$\pi(h')=1$ and $e(h')\ne f(h')$)
to obtain a path $h\in\P^{p,p'-p}_{a,b,e,f}(L)$ for some $L$.
Since $\lfloor a'p/p'\rfloor=\lfloor a'p/(p'-p)\rfloor$, and
$\lfloor b'p/p'\rfloor=\lfloor b'p/(p'-p)\rfloor$ by Lemma \ref{StartPtLem},
and thus $a'=a+\lfloor ap/(p'-p)\rfloor+e$ and
$b'=b+\lfloor bp/(p'-p)\rfloor+f$,
it follows that $h^{(0)}$ arises by the action of a
$\B_1$-transform on $h$.
 
On examining the third case depicted in Table~\ref{BatStartTable},
we see that the case where $h^{(0)}$ has a pair of scoring vertices
at the $0$th and $1$st positions, arises similarly from a particular
$h\in\P^{p,p'-p}_{a,b,e,f}(L)$ for some $L$.
The lemma is then proved.
\cqfd
\medskip

\noindent
The value of $k$ obtained above will be referred to as the
particle content of $h'$.

We could now proceed as in Lemma 3.13 and Corollary 3.14 of
\cite{foda-welsh-kyoto} to, for appropriate parameters, give a
weight-respecting bijection between $\bigcup_{k}
\P^{p,p'}_{a,b,e,f}(m_1,2k+2m_1-m_0)\times\Y(k,m_1)$
and $\P^{p,p'+p}_{a',b',e,f}(m_0,m_1)$, and state the corresponding
relationship between the generating functions.
However, something more general is required in this paper.
This will be derived in Lemma \ref{MazyBijLem} and Corollary \ref{MazyBijCor}
for certain restricted paths once the notation is developed in
Section \ref{MazySec}.

\newpage

\setcounter{section}{3}

\section{The $\D$-transform}\label{DTranSec}

\subsection{The pure $\D$-transform}

The {\it $\D$-transform}
is defined to act on each $h \in \P^{p,p'}_{a,b,e,f}(L)$ to yield
a path $\hat h\in\P^{p'-p,p'}_{a,b,1-e,1-f}(L)$ with exactly the same
sequence of integer heights, i.e., $\hat h_i=h_i$ for $0\le i\le L$.
Note that, by definition, $e(\hat h)=1-e(h)$ and $f(\hat h)=1-f(h)$.

Since the band structure of
the $(p'-p,p')$-model is obtained from that of the $(p,p')$-model
simply by replacing odd bands by even bands and vice-versa,
then, ignoring the vertex at $i=0$,
each scoring vertex maps to a non-scoring vertex and vice-versa.
That $e(h)$ and $e(\hat h)$ differ implies that the vertex
at $i=0$ is both scoring or both non-scoring in $h$ and $\hat h$.

\begin{lemma}\label{DresLem}
Let $\hat h \in \P^{p'-p,p'}_{a,b,1-e,1-f}(L)$ be obtained from
$h \in \P^{p,p'}_{a,b,e,f}(L)$ by the action of the $\D$-transform.
Then $\pi(\hat h)=1-\pi(h)$.
Moreover, if $m=m(h)$ then:
\begin{displaymath}
\begin{array}{l}
\bullet\quad L(\hat h)=L;\\[1mm]
\bullet\quad m(\hat h)=
  \left\{
    \begin{array}{ll}
  L-m \quad &
            \mbox{if } e+d+\pi(h)\equiv0\,(\mod2),\\[1.5mm]
  L-m+2 \quad & \mbox{if } e+d+\pi(h)\not\equiv0\,(\mod2);
    \end{array} \right.
\\[1mm]
\bullet\quad \mwt(\hat h)=\frac{1}{4}
   \left( L^2 - \alpha(h)^2 \right) - \mwt(h).
\end{array}
\end{displaymath}
\end{lemma}

\Proof
Let $h$ have striking sequence
$\left({a_1\atop b_1}\:{a_2\atop b_2}\:{a_3\atop b_3}\:
 {\cdots\atop\cdots}\:{a_l\atop b_l} \right)^{(e,f,d)}$.
Since, beyond the zeroth vertex, the $\D$-transform
exchanges scoring vertices for non-scoring vertices and vice-versa,
it follows that the striking sequence for $\hat h$ is
$\left({b_1\atop a_1}\:{b_2\atop a_2}\:{b_3\atop a_3}\:
 {\cdots\atop\cdots}\:{b_l\atop a_l} \right)^{(1-e,1-f,d)}$.
It is immediate that $L(\hat h)=L$, $\pi(\hat h)=1-\pi(h)$,
$e(\hat h)=1-e(h)$ and $d(\hat h)=d(h)$.
Then
$m(\hat h)=(e(\hat h)+d(\hat h)+\pi(\hat h))\,\mod2+\sum_{i=1}^l b_i
=(e+d+\pi(h))\,\mod2+L-\sum_{i=1}^l a_i
=2((e+d+\pi(h))\,\mod2)+L-m(h)$.
 
Now let $w_i=a_i+b_i$ for $1\le i\le l$.
Then, using Lemma \ref{WtHashLem}, we obtain
\begin{align*}
\mwt(h)+\mwt(\hat h)&=
\sum_{i=1}^l b_i(w_{i-1}+w_{i-3}+\cdots+w_{1+i\bmod2})\\
&\qquad
+ \sum_{i=1}^l a_i(w_{i-1}+w_{i-3}+\cdots+w_{1+i\bmod2})\\[0.5mm]
&=
\sum_{i=1}^l w_i(w_{i-1}+w_{i-3}+\cdots+w_{1+i\bmod2})\\[0.5mm]
&=
(w_1+w_3+w_5+\cdots)(w_2+w_4+w_6+\cdots).
\end{align*}
The lemma then follows because
$(w_1+w_3+w_5+\cdots)+(w_2+w_4+w_6+\cdots)=L$ and
$(w_1+w_3+w_5+\cdots)-(w_2+w_4+w_6+\cdots)=\pm\alpha(h)$.
\cqfd
\medskip

\begin{lemma}\label{DPathParamLem}
If $1\le p<p'$ with $p$ coprime to $p'$,
$1\le a,b<p'$ and $e,f\in\{0,1\}$ then
$\alpha^{p'-p,p'}_{a,b}=\alpha^{p,p'}_{a,b}$ and
$\beta^{p'-p,p'}_{a,b,1-e,1-f}+\beta^{p,p'}_{a,b,e,f}=\alpha^{p,p'}_{a,b}$.
\end{lemma}

\Proof Lemma \ref{DParamLem} gives
$\lfloor ap/p'\rfloor+\lfloor a(p'-p)/p'\rfloor=a-1$ and
likewise, $\lfloor bp/p'\rfloor+\lfloor b(p'-p)/p'\rfloor=b-1$.
The required results then follow immediately.
\cqfd
\medskip

\subsection{The $\B\D$-pair}\label{BDTranSec}

It will often be convenient to consider the combined action
of a $\D$-transform followed immediately by a $\B$-transform.
Such a pair will naturally be referred to as a $\B\D$-transform
and maps a path $h\in\P^{p'-p,p'}_{a,b,e,f}(L)$ to a path
$h'\in\P^{p,p'+p}_{a',b',1-e,1-f}(L')$, where $a',b',L'$
are determined by our previous results.

In what follows, the $\B\D$-transform will always follow a $\B$-transform.
Thus we restrict consideration to where $2(p'-p)<p'$.

\begin{lemma}\label{BDresLem}
With $p'<2p$, let $h\in{\P}^{p'-p,p'}_{a,b,e,f}(L)$.
Let $h'\in\P^{p,p'+p}_{a',b',1-e,1-f}(L')$ result from the action
of a $\D$-transform on $h$, followed by a $\B(k,\lambda)$-transform.
Then:
\begin{displaymath}
\begin{array}{l}
\bullet\quad L'=
  \left\{
    \begin{array}{ll}
  L+m(h)+2k-2 \quad &
            \mbox{if } \pi(h)=1 \mbox{ and } e=d(h),\\[1.5mm]
  L+m(h)+2k \quad & \mbox{otherwise};
    \end{array} \right.
\\[1mm]
\bullet\quad m(h')=L;\\[1mm]
\bullet\quad \mwt(h')=\frac{1}{4}\left( L^2 + (L'-L)^2
- \alpha^2 - \beta^2 \right) + \wt(\lambda) - \mwt(h),
\end{array}
\end{displaymath}
where $\alpha=\alpha^{p,p'}_{a,b}$ and $\beta=\beta^{p,p'}_{a,b,1-e,1-f}$.
\end{lemma}

\Proof Let $\hat h$ result from the action of the $\D$-transform
on $h$, and let $d=d(h)$, $\pi=\pi(h)$, $\hat e=e(\hat h)$
$\hat d=d(\hat h)$, $\hat\pi=\pi(\hat h)$. Then we immediately
have $\hat d=d$, $\hat e=1-e$, and $\hat\pi=1-\pi$.

In the case where $\pi=0$ and $e\ne d$, we then have,
using Lemmas \ref{BresLem} and \ref{DresLem},
$L'=2L(\hat h)-m(\hat h)+2k+2=2L-(L-m(h)+2)+2k+2=L+m(h)+2k$.

In the case where $\pi=1$ and $e=d$, we then have,
using Lemmas \ref{BresLem} and \ref{DresLem},
$L'=2L(\hat h)-m(\hat h)+2k+2=2L-(L-m(h)+2)+2k=L+m(h)+2k-2$.

In the other cases, $e+d+\pi\equiv0\,(\mod2)$ and
so $\hat e+\hat d+\hat\pi\equiv0\,(\mod2)$.
Lemmas \ref{BresLem} and \ref{DresLem} yield
$L'=2L(\hat h)-m(\hat h)+2k=2L-(L-m(h))+2k=L+m(h)+2k$.

The expressions for $m(h')$
and $\mwt(h')$ follow immediately from Lemmas \ref{BresLem},
\ref{DresLem} and \ref{BetaConstLem}.
\cqfd
\medskip

We could now proceed as in Lemma 4.5 and
Corollary 4.6 of \cite{foda-welsh-kyoto} to show that the $\B\D$-transform
implies, for certain parameters, a weight respecting bijection between
$\bigcup_{k}
\P^{p'-p,p'}_{a,b,e,f}(m_1,m_0-m_1-2k)\times\Y(k,m_1)$
and $\P^{p,p'+p}_{a',b',1-e,1-f}(m_0,m_1)$,
and consequently a relationship between the corresponding generating
functions.
However again, something more general is required in the current paper.
This is developed in Section \ref{MazySec}, resulting in
Lemma \ref{MazyDijLem} and Corollary \ref{MazyDijCor}.

We finish this section by examining how the $\B\D$-transform affects
some of the parameters.

\begin{lemma}\label{BDPathParamLem}
Let $1\le p<p'<2p$ with $p$ coprime to $p'$,
$1\le a,b<p'$ and $e,f\in\{0,1\}$, and set
$a'=a+1-e+\lfloor ap/p'\rfloor$ and  $b'=b+1-f+\lfloor bp/p'\rfloor$.
Then
$\alpha^{p,p'+p}_{a',b'}=2\alpha^{p'-p,p'}_{a,b}-\beta^{p'-p,p'}_{a,b,e,f}$
and
$\beta^{p,p'+p}_{a',b',1-e,1-f}
  =\alpha^{p'-p,p'}_{a,b}-\beta^{p'-p,p'}_{a,b,e,f}$.
\end{lemma}

\Proof By using first Lemma \ref{ParamLem} and then
Lemma \ref{DPathParamLem}, we obtain:
$\alpha^{p,p'+p}_{a',b'}=
\alpha^{p,p'}_{a,b}+\beta^{p,p'}_{a,b,1-e,1-f}=
2\alpha^{p'-p,p'}_{a,b}-\beta^{p'-p,p'}_{a,b,e,f}$.
By a similar route, we obtain:
$\beta^{p,p'+p}_{a',b',1-e,1-f}=
\beta^{p,p'}_{a,b,1-e,1-f}=
\alpha^{p'-p,p'}_{a,b}-\beta^{p'-p,p'}_{a,b,e,f}$.
\cqfd
\medskip

\newpage

\setcounter{section}{4}

\section{Mazy runs}\label{MazySec}

In this section, we examine paths that are constrained to attain
certain heights in a certain order. We will find that the subset
of $\P^{p,p'}_{a,b,e,f}(L)$ that corresponds to an individual
term $F_{a,b}(\boldu^L,\boldu^R,L)$ as defined in Section \ref{FinFermSec},
may be characterised in such a way.
 
For $1\le a<p'$, let $a$ be interfacial in the $(p,p')$-model,
and let $r=\rho^{p,p'}(a)$ so that $a$ borders the $r$th odd band.
Then $a=\lfloor p'r/p\rfloor$ if and only if
$r=\lfloor pa/p'\rfloor+1$, and
$a=\lfloor p'r/p\rfloor+1$ if and only if
$r=\lfloor pa/p'\rfloor$.
In the former case, we write $\omega^{p,p'}(a)=r^-$,
and in the latter case we write $\omega^{p,p'}(a)=r^+$.
It will also be convenient to define $\omega^{p,p'}(a)=\infty$
when $a$ is not interfacial
(simply so that certain statements make sense).
When $p$ and $p'$ are implicit, we write $\omega(a)=\omega^{p,p'}(a)$.

\subsection{The passing sequence}
For each path $h\in\P^{p,p'}_{a,b,e,f}(L)$, the passing sequence $\omega(h)$
of $h$ is a word that indicates how the path $h$ meanders between
the interfacial values. It is a word in the elements of
$\R^{p,p'}=\{1^-,1^+,2^-,2^+,\ldots,(p-1)^-,(p-1)^+\}$
and is obtained as follows.
First start with an empty word.
For $i=0,1,\ldots,L$, if $h_i$ is interfacial,
append $\omega(h_i)$ to the word.
Then replace each subsequence of consecutive identical symbols within
this word by a single instance of that symbol.
We define $\omega(h)$ to be this resulting word.

For example, in the case of the path $h\in\P^{3,8}_{2,4,e,1}(14)$
depicted in Fig.~\ref{TypicalShadedFig}, we obtain the passing sequence
$\omega(h)=1^-1^+2^-2^+2^-2^+2^-1^+$.

In the case $p'<2p$, it will be useful to note that not all the elements
of $\R^{p,p'}$ can occur in a passing sequence, because some
of them correspond to multifacial heights. The set of elements that
cannot occur in this $p'<2p$ case is thus given by:
\begin{equation*}
\overline{\R}{}^{p,p'}=\{(r-1)^+,r^-:1<r<p',
\lfloor(r-1)p'/p\rfloor=\lfloor rp'/p\rfloor-1\}
\cup\{1^-,(p-1)^+\}.
\end{equation*}

The utility of the passing function lies in its near invariance
under the $\B$-transform, and its simple transformation under
the $\B\D$-transform.

\begin{lemma}\label{PassBLem1}
Let $p'>2p$ and let $h'\in{\P}^{p,p'+p}_{a',b',e,f}(L')$ be obtained from
$h\in{\P}^{p,p'}_{a,b,e,f}(L)$ through the action of a
$\B(k,\lambda)$-transform.
Then either $\omega(h')=\omega(h)$,
$\omega(a)\omega(h')=\omega(h)$,
$\omega(h')\omega(b)=\omega(h)$,
or $\omega(a)\omega(h')\omega(b)=\omega(h)$.

Moreover, if either $\omega(a)\omega(h')=\omega(h)$ or
$\omega(a)\omega(h')\omega(b)=\omega(h)$ then
$$\delta^{p,p'}_{a,e}=1 \text{ and } h_1=a-(-1)^e,$$
and if either $\omega(h')=\omega(h)\omega(b)$ or
$\omega(h')=\omega(a)\omega(h)\omega(b)$ then
$$\delta^{p,p'}_{b,f}=1 \text{ and } h_{L-1}=b-(-1)^f.$$
\end{lemma}

\Proof Let $h^{(0)}$ result from the action of the $\B_1$-transform on $h$.
The discussion of the $\B_1$-transform immediately following
Note \ref{Start1Note} together with Tables \ref{BatStartTable}
and \ref{BatEndTable}
show that either $\omega(h^{(0)})=\omega(h)$,
$\omega(a)\omega(h^{(0)})=\omega(h)$,
$\omega(h^{(0)})\omega(b)=\omega(h)$,
or $\omega(a)\omega(h^{(0)})\omega(b)=\omega(h)$.
In addition (see Table \ref{BatStartTable})
$\omega(a)\omega(h^{(0)})=\omega(h)$
or $\omega(a)\omega(h^{(0)})\omega(b)=\omega(h)$,
only if $\delta^{p,p'}_{a,e}=1$ and $h_1=a-(-1)^e$.
Similarly (see Table \ref{BatEndTable})
$\omega(h^{(0)})\omega(b)=\omega(h)$
or $\omega(a)\omega(h^{(0)})\omega(b)=\omega(h)$,
only if $\delta^{p,p'}_{b,f}=1$ and $h_1=b-(-1)^f$.

Now let $h^{(k)}$ be the result of applying the $\B_2(k)$-transform
of Section \ref{InsertSec} to $h^{(0)}$.
Noting that because, when $p'>2p$, neighbouring odd bands in the
$(p,p'+p)$-model are separated by at least two even bands, we see
that inserting particles into $h^{(0)}$ can change the passing
sequence only if $a'+(-1)^{e}$ is interfacial and
$h^{(0)}_1=a'-(-1)^{e}$.
This occurs only in that case where $\delta^{p,p'}_{a,e}=1$ and $h_1=a-(-1)^e$.
In this case, $\omega(h^{(k)})=\omega(a)\omega(h^{(0)})=\omega(h)$.
Otherwise $\omega(h^{(k)})=\omega(h^{(0)})$.

Now let $h'$ be the result of applying the $\B_3(\lambda)$-transform
to $h^{(k)}$,
By inspection, again noting that neighbouring odd bands are separated by at
least two even bands, each of the moves from Section \ref{PartMovesSec}
does not change the passing sequence.
This is not necessarily true of the edge-moves.
Edge-moves 1 and 2 can change the passing sequence
only if $a'+(-1)^e$ is interfacial, $\omega(h^{(k)})$
begins with $\omega(a)$, and $h_1=a-(-1)^e$.
In this case, the effect of the edge-move is to remove the $\omega(a)$.
Edge-moves 3 and 4 can change the passing sequence
only if $b'+(-1)^f$ is interfacial, $\omega(h^{(k)})$
does not end with $\omega(b)$, and $h_1=b-(-1)^f$.
In this case, the effect of the edge-move is to append $\omega(b)$.
Because $a'+(-1)^e$ is interfacial implies that $\delta^{p,p'}_{a,e}=1$,
and $b'+(-1)^f$ is interfacial implies that $\delta^{p,p'}_{b,f}=1$,
the lemma then follows.
\cqfd
\medskip


For the case $p'<2p$, the above result doesn't hold and
$\omega(h')$ can be very different from $\omega(h)$.
However, in this case, we can make do with the following result.
For the sake of clarification, by a subword $\omega'$ of
$\omega=\omega_1\omega_2\cdots\omega_k$, we mean a sequence
$\omega'=\omega_{i_1}\omega_{i_2}\cdots \omega_{i_{k'}}$ where
$1\le i_1<i_2<\cdots<i_{k'}\le k$.

\begin{lemma}\label{PassBLem2}
Let $p'<2p$ and let
$h'\in{\P}^{p,p'+p}_{a',b',e,f}(L')$ be obtained from
the path $h\in{\P}^{p,p'}_{a,b,e,f}(L)$ through the action of a
$\B(k,\lambda)$-transform.
Let $\omega'=\omega'_1\omega'_2\cdots \omega'_{k'}$ be the longest subword
of $\omega(h')$ for which $\omega'_i\ne \omega'_{i+1}$ for $1\le i<k'$,
that comprises the symbols in
$\R^{p,p'}\backslash\overline{\R}{}^{p,p'}$.
Then either $\omega(h')=\omega(h)$,
$\omega(a)\omega(h')=\omega(h)$,
$\omega(h')\omega(b)=\omega(h)$,
or $\omega(a)\omega(h')\omega(b)=\omega(h)$.

Moreover, if either $\omega(a)\omega(h')=\omega(h)$ or
$\omega(a)\omega(h')\omega(b)=\omega(h)$ then
$\delta^{p,p'}_{a,e}=1$ and $h_1=a-(-1)^e,$
and if either $\omega(h')\omega(b)=\omega(h)$ or
$\omega(a)\omega(h')\omega(b)=\omega(h)$ then
$\delta^{p,p'}_{b,f}=1$ and $h_{L-1}=b-(-1)^f$.
\end{lemma}

\Proof Note first that the heights in the $(p,p'+p)$-model that
correspond to the symbols from
$\R^{p,p'}\backslash\overline{\R}{}^{p,p'}$
are either separated by a single odd band or at least two even bands.
We then compare the actions of the
$\B_1$-, $\B_2(k)$-, and $\B_3(k,\lambda)$-transforms to those described
in the proof of Lemma \ref{PassBLem1}.
As there, $\omega(a)$ might get removed from the start of the passing
sequence if $a$ is interfacial in the $(p,p')$-model,
$\delta^{p,p'}_{a,e}=1$ and $h_1=a-(-1)^e$.
Similarly $\omega(b)$ might get removed from the end of of the passing
sequence if $b$ is interfacial in the $(p,p')$-model,
$\delta^{p,p'}_{b,f}=1$ and $h_{L-1}=b-(-1)^f$.
Apart from these changes,
the action of the $\B_1$-transform can only insert elements from
$\overline{\R}{}^{p,p'}$ into the passing sequence.
The same is true of the $\B_2(k)$-transform.
In the case of the $\B_3(k)$-transform, we see that only those
symbols that correspond to heights separated by a single
even band might get inserted or removed.
The lemma then follows.
\cqfd
\medskip

Since heights that are interfacial in the $(p'-p,p')$-model are
also interfacial in the $(p,p')$-model and vice-versa,
we see that the effect of a $\D$-transform on a passing
sequence is just to change the name of each symbol.
So as to be able to describe this change, for $1\le r<(p'-p)$ set
$R=\lfloor rp/(p'-p)\rfloor$ and define:
\begin{displaymath}
\begin{array}{l}
\sigma(r^-)=R^+;\\[1.5mm]
\sigma(r^+)=(R+1)^-,
\end{array}
\end{displaymath}
Then, if $\omega(h)=\omega_1\omega_2\cdots \omega_k$, define
$\sigma(\omega(h))=\sigma(\omega_1)\sigma(\omega_2)\cdots\sigma(\omega_k)$.

\begin{lemma}\label{PassBDLem}
Let $p'<2p$ and let
$h'\in{\P}^{p,p'+p}_{a',b',1-e,1-f}(L')$ be obtained from
the path $h\in{\P}^{p'-p,p'}_{a,b,e,f}(L)$ through the action of a
$\B\D$-transform.
If $\omega'=\omega'_1\omega'_2\cdots \omega'_{k'}$ is
the longest subword of $\omega(h')$ for which
$\omega'_i\ne \omega'_{i+1}$ for $1\le i<k'$, and which comprises
the symbols in $\{\sigma(x):x\in\R^{p'-p,p'}\}$,
then either $\omega'=\sigma(\omega(h))$,
$\omega(a)\omega'=\sigma(\omega(h))$,
$\omega'\omega(b)=\sigma(\omega(h))$,
or $\omega(a)\omega'\omega(b)=\sigma(\omega(h))$.

Moreover, if either $\omega(a)\omega'=\sigma(\omega(h))$ or
$\omega(a)\omega'\omega(b)=\sigma(\omega(h))$ then
$a$ is interfacial in the $(p,p')$-model,
$\delta^{p,p'}_{a,e}=1$ and $h_1=a+(-1)^e,$
and if either $\omega'\omega(b)=\sigma(\omega(h))$ or
$\omega(a)\omega'\omega(b)=\sigma(\omega(h))$ then
$b$ is interfacial in the $(p,p')$-model,
$\delta^{p,p'}_{b,f}=1$ and $h_{L-1}=b+(-1)^f$.
\end{lemma}

\Proof
For $1\le r<(p'-p)$, the lowermost edge of the $r$th
odd band in the $(p'-p,p')$-model is at height
$a=\lfloor rp'/(p'-p)\rfloor$.
Since $p'>2(p'-p)$, no two odd bands in the $(p'-p,p')$-model
are adjacent and so $a$ is the height of the
uppermost edge of the $R$th even band since
$R=\lfloor rp'/(p'-p)\rfloor-r$.
Thus, in the $(p,p')$-model, the uppermost edge of the
$R$th odd band is at height $a$.
We also see that the uppermost edge of the $r$th odd band
in the $(p'-p,p')$-model is at height $a+1$ which is also the height
of the lowermost edge of the $(R+1)$th odd band of the $(p,p')$-model.

Thus if $\hat h\in{\P}^{p,p'}_{a,b,1-e,1-f}(L)$ results from
a $\D$-transform acting on $h$, then we immediately obtain
$\omega(\hat h)=\sigma(\omega(h))$.
The lemma then follows on applying Lemma \ref{PassBLem2} to $\hat h$.
\cqfd
\medskip

\subsection{The passing function}\label{PassingSec}

Let $p'>2p$, $1\le a,b<p'$, $e,f\in\{0,1\}$ and $L\ge0$.
We now define certain subsets of $\P^{p,p'}_{a,b,e,f}(L)$
that depend on $d^L$-dimensional vectors
$\boldmu=(\mu_1,\mu_2,\ldots,\mu_{d^L})$ and
$\boldmu^*=(\mu^*_1,\mu^*_2,\ldots,\mu^*_{d^L})$,
and $d^R$-dimensional vectors
$\boldnu=(\nu_1,\nu_2,\ldots,\nu_{d^R})$ and
$\boldnu^*=(\nu^*_1,\nu^*_2,\ldots,\nu^*_{d^R})$,
where $d^L,d^R\ge0$.
We say that a path $h\in\P^{p,p'}_{a,b,e,f}(L)$ is {\em mazy-compliant}
(with $\boldmu,\boldmu^*,\boldnu,\boldnu^*$) if the following five
conditions are satisfied:
\begin{enumerate}
\item for $1\le j\le d^L$, there exists $i$ with $0\le i\le L$ such that
$h_i=\mu^*_j$;
\item for $1\le j\le d^L$, if there exists $i'$ with $0\le i'\le L$
such that $h_{i'}=\mu_j$, then there exists $i$ with $0\le i<i'$ such that
$h_i=\mu^*_j$;
\item for $1\le j\le d^R$, there exists $i$ with $0\le i\le L$ such that
$h_i=\nu^*_j$;
\item for $1\le j\le d^R$, if there exists $i'$ with $0\le i'\le L$
such that $h_{i'}=\nu_j$, then there exists $i$ with $i'<i\le L$ such that
$h_i=\nu^*_j$;
\item if $d^L,d^R>0$ and $i$ is the smallest value such that
$0\le i\le L$ and $h_i=\mu^*_1$ and $i'$ is the largest value such that
$0\le i'\le L$ and $h_{i'}=\nu^*_1$, then $i\le i'$.
\end{enumerate}
Loosely speaking, these are the paths which, for $1\le j\le d^L$,
attain $\mu^*_j$, and do so before they attain any $\mu_j$;
which, for $1\le k\le d^R$,
attain $\nu^*_k$, and do so after they attain any $\nu_k$;
and which, if $d^L,d^R>0$, attain $\nu^*_1$ after they attain
the first $\mu^*_1$.

The set
$\P^{p,p'}_{a,b,e,f}(L)\left\{
\begin{smallmatrix} \boldmu^{\phantom{*}}\!&;&\boldnu^{\phantom{*}}\\
                    \boldmu^*\!&;&\boldnu^* \end{smallmatrix} \right\}$
is defined to be the subset of $\P^{p,p'}_{a,b,e,f}(L)$
comprising those paths $h$ which are mazy-compliant with
$\boldmu,\boldmu^*,\boldnu,\boldnu^*$.
The generating function for these paths is defined to be:
\begin{equation}\label{MazyPathGenDef}
\mchi^{p,p'}_{a,b,e,f}(L;q)\left\{
\begin{matrix} \boldmu^{\phantom{*}};\boldnu^{\phantom{*}}\\
                    \boldmu^*;\boldnu^* \end{matrix} \right\}
=\sum_{h\in\P^{p,p'}_{a,b,e,f}(L)\left\{
\genfrac{}{}{0pt}{3}{\boldmu^{\phantom{*}}}{\boldmu^*}\!
\genfrac{}{}{0pt}{3}{;}{;}
\genfrac{}{}{0pt}{3}{\boldnu^{\phantom{*}}}{\boldnu^*} \right\}}
\hskip-5mm
q^{\mwt(h)}.
\end{equation}

We also define
$\P^{p,p'}_{a,b,e,f}(L,m)\left\{
\begin{smallmatrix} \boldmu^{\phantom{*}}\!&;&\boldnu^{\phantom{*}}\\
                    \boldmu^*\!&;&\boldnu^* \end{smallmatrix} \right\}$
to be the subset of $\P^{p,p'}_{a,b,e,f}(L,m)$
comprising those paths $h$ which are mazy-compliant with
$\boldmu,\boldmu^*,\boldnu,\boldnu^*$.
We then define the generating function:
\begin{equation*}
\mchi^{p,p'}_{a,b,e,f}(L,m;q)\left\{
\begin{matrix} \boldmu^{\phantom{*}};\boldnu^{\phantom{*}}\\
                    \boldmu^*;\boldnu^* \end{matrix} \right\}
=\sum_{h\in\P^{p,p'}_{a,b,e,f}(L,m)\left\{
\genfrac{}{}{0pt}{3}{\boldmu^{\phantom{*}}}{\boldmu^*}\!
\genfrac{}{}{0pt}{3}{;}{;}
\genfrac{}{}{0pt}{3}{\boldnu^{\phantom{*}}}{\boldnu^*} \right\}}
\hskip-5mm
q^{\mwt(h)}.
\end{equation*}
The following result is proved in exactly the
same way as Lemma \ref{ResPathGenLem}:
\begin{lemma}\label{MazyResPathGenLem}
Let $\beta=\beta^{p,p'}_{a,b,e,f}$. Then
\begin{equation*}
\mchi^{p,p'}_{a,b,e,f}(L;q)\left\{
\begin{matrix} \boldmu^{\phantom{*}};\boldnu^{\phantom{*}}\\
                    \boldmu^*;\boldnu^* \end{matrix} \right\}
=\sum_{\scriptstyle m\equiv L+\beta\atop
  \scriptstyle\strut(\mbox{\scriptsize\rm mod}\,2)}
\mchi^{p,p'}_{a,b,e,f}(L,m;q)\left\{
\begin{matrix} \boldmu^{\phantom{*}};\boldnu^{\phantom{*}}\\
                    \boldmu^*;\boldnu^* \end{matrix} \right\}.
\end{equation*}
\end{lemma}

Our vectors $\boldmu$, $\boldmu^*$, $\boldnu$, $\boldnu^*$,
will satisfy certain constraints.
We say that $\boldmu$ and $\boldmu^*$ are a {\em mazy-pair} in the
$(p,p')$-model sandwiching $a$ if they satisfy the following three
conditions:
\begin{enumerate}
\item $0\le\mu_j,\mu^*_j\le p'$ for $1\le j\le d^L$
(we permit $\mu^*_j=0$, $\mu^*_j=p'$ although the set of paths
constrained by such parameters will be empty);
\item $\mu_j\ne\mu^*_j$ for $1\le j\le d^L$;
\item $a$ is strictly between\footnote{By $x$ between $y$ and $z$, we mean
that either $y\le x\le z$ or $z\le x\le y$.  By $x$ strictly between
$y$ and $z$, we mean that either $y<x<z$ or $z<x<y$.}
$\mu_{d^L}$ and $\mu^*_{d^L}$, and for $1\le j<d^L$,
both $\mu_{j+1}$ and $\mu_{j+1}^*$ are between $\mu_j$ and $\mu^*_j$.
\end{enumerate}
If in addition,
\begin{enumerate}
\item[4.] $\mu_j$ and $\mu^*_j$ are interfacial in the $(p,p')$-model
for $1\le j\le d^L$,
\end{enumerate}
we say that $\boldmu$ and $\boldmu^*$ are an {\em interfacial mazy-pair} in
the $(p,p')$-model sandwiching $a$.
If $\boldmu$ and $\boldmu^*$ are a mazy-pair (resp.\ interfacial mazy-pair)
in the $(p,p')$-model sandwiching $a$, and $\boldnu$ and $\boldnu^*$ are
a mazy-pair (resp.\ interfacial mazy-pair) in the $(p,p')$-model
sandwiching $b$, then we say that
$\boldmu$, $\boldmu^*$, $\boldnu$, $\boldnu^*$, are a
{\em mazy-four} (resp.\ {\em interfacial mazy-four}) in the $(p,p')$-model
sandwiching $(a,b)$.

\begin{lemma}\label{CutParam1Lem}
Let $\boldmu,\boldmu^*,\boldnu,\boldnu^*$ be a
mazy-four in the $(p,p')$-model sandwiching $(a,b)$.
If $d^L\ge2$ and $\mu_{d^L}=\mu_{d^L-1}$, then
\begin{equation*}
\mchi^{p,p'}_{a,b,e,f}(L,m)\left\{
\begin{matrix}
         \mu_1,\ldots,\mu_{d^L-1},\mu_{d^L};\boldnu^{\phantom{*}}\\
         \mu^*_1,\ldots,\mu^*_{d^L-1},\mu^*_{d^L};\boldnu^*
\end{matrix} \right\}
=\mchi^{p,p'}_{a,b,e,f}(L,m)\left\{
\begin{matrix}
         \mu_1,\ldots,\mu_{d^L-1};\boldnu^{\phantom{*}}\\
         \mu^*_1,\ldots,\mu^*_{d^L-1};\boldnu^*
\end{matrix} \right\}.
\end{equation*}
If $d^R\ge2$ and $\nu_{d^R}=\nu_{d^R-1}$, then
\begin{equation*}
\mchi^{p,p'}_{a,b,e,f}(L,m)\left\{
\begin{matrix}
\boldmu^{\phantom{*}};\nu_1,\ldots,\nu_{d^R-1},\nu_{d^R}\\
 \boldmu^*;\nu^*_1,\ldots,\nu^*_{d^R-1},\nu^*_{d^R}
\end{matrix} \right\}
=\mchi^{p,p'}_{a,b,e,f}(L,m)\left\{
\begin{matrix}
\boldmu^{\phantom{*}};\nu_1,\ldots,\nu_{d^R-1}\\
 \boldmu^*;\nu^*_1,\ldots,\nu^*_{d^R-1}
\end{matrix} \right\}.
\end{equation*}
\end{lemma}

\Proof
Immediately from the definition,
$
\P^{p,p'}_{a,b,e,f}(L,m)\left\{
\begin{smallmatrix}
         \mu_1,\ldots,\mu_{d^L-1},\mu_{d^L};\boldnu^{\phantom{*}}\\
         \mu^*_1,\ldots,\mu^*_{d^L-1},\mu^*_{d^L};\boldnu^*
\end{smallmatrix} \right\}
\subseteq
$
\break
$
\P^{p,p'}_{a,b,e,f}(L,m)\left\{
\begin{smallmatrix}
         \mu_1,\ldots,\mu_{d^L-1};\boldnu^{\phantom{*}}\\
         \mu^*_1,\ldots,\mu^*_{d^L-1};\boldnu^*
\end{smallmatrix} \right\}
$.
Now consider a path $h$ belonging to the second of these sets.
Since $a$ is between $\mu_{d^L}=\mu_{d^L-1}$ and
$\mu^*_{d^L}$, and $\mu^*_{d^L}$
is between $\mu_{d^L}=\mu_{d^L-1}$ and $\mu^*_{d^L-1}$,
then $\mu^*_{d^L}$ is between $a$ and $\mu^*_{d^L-1}$.
Thus if $h$ attains $\mu^*_{d^L-1}$ then it necessarily attains
$\mu^*_{d^L}$, and if it attains $\mu^*_{d^L-1}$ before it attains
$\mu_{d^L}=\mu_{d^L-1}$, then it necessarily attains $\mu^*_{d^L}$
before it attains $\mu_{d^L}=\mu_{d^L-1}$.
Thus $h$ belongs to the first of the above sets, whence the first
part of the lemma follows. The second expression is proved
in an analogous way.
\cqfd
\medskip

\begin{lemma}\label{CutParam2Lem}
Let $\boldmu,\boldmu^*,\boldnu,\boldnu^*$ be a
mazy-four in the $(p,p')$-model sandwiching $(a,b)$.
If $1\le j<d^L$ and $\mu^*_j=\mu^*_{j+1}$, then
\begin{multline*}
\mchi^{p,p'}_{a,b,e,f}(L)\left\{
\begin{matrix}
  \mu_1,\ldots,\mu_{j-1},\mu_{j},\mu_{j+1},\ldots,\mu_{d^L};
                                                  \boldnu^{\phantom{*}}\\
  \mu^*_1,\ldots,\mu^*_{j-1},\mu^*_{j},\mu^*_{j+1},\ldots,\mu^*_{d^L};
                                                  \boldnu^*
\end{matrix} \right\}\\
=
\mchi^{p,p'}_{a,b,e,f}(L)\left\{
\begin{matrix}
  \mu_1,\ldots,\mu_{j-1},\mu_{j+1},\ldots,\mu_{d^L};
                                                  \boldnu^{\phantom{*}}\\
  \mu^*_1,\ldots,\mu^*_{j-1},\mu^*_{j+1},\ldots,\mu^*_{d^L};
                                                  \boldnu^*
\end{matrix} \right\}.
\end{multline*}
If $1\le j<d^R$ and $\nu^*_j=\nu^*_{j+1}$, then
\begin{multline*}
\mchi^{p,p'}_{a,b,e,f}(L)\left\{
\begin{matrix}
  \boldmu^{\phantom{*}};
        \nu_1,\ldots,\nu_{j-1},\nu_{j},\nu_{j+1},\ldots,\nu_{d^R}\\
  \boldmu^*;
        \nu^*_1,\ldots,\nu^*_{j-1},\nu^*_{j},\nu^*_{j+1},\ldots,\nu^*_{d^R}
\end{matrix} \right\}\\
=
\mchi^{p,p'}_{a,b,e,f}(L)\left\{
\begin{matrix}
  \boldmu^{\phantom{*}};
        \nu_1,\ldots,\nu_{j-1},\nu_{j+1},\ldots,\nu_{d^R}\\
  \boldmu^*;
        \nu^*_1,\ldots,\nu^*_{j-1},\nu^*_{j+1},\ldots,\nu^*_{d^R}
\end{matrix} \right\}.
\end{multline*}
\end{lemma}

\Proof Immediately from the definition,
$
\P^{p,p'}_{a,b,e,f}(L)\left\{
\begin{smallmatrix}
        \mu_1,\ldots,\mu_{j-1},\mu_{j},\mu_{j+1},\ldots,\mu_{d^R};
  \boldnu^{\phantom{*}}\\
        \mu^*_1,\ldots,\mu^*_{j-1},\mu^*_{j},\mu^*_{j+1},\ldots,\mu^*_{d^R};
  \boldnu^*
\end{smallmatrix} \right\}
\subseteq
$
\break
$
\P^{p,p'}_{a,b,e,f}(L)\left\{
\begin{smallmatrix}
        \mu_1,\ldots,\mu_{j-1},\mu_{j+1},\ldots,\mu_{d^R};
  \boldnu^{\phantom{*}}\\
        \mu^*_1,\ldots,\mu^*_{j-1},\mu^*_{j+1},\ldots,\mu^*_{d^R};
  \boldnu^*
\end{smallmatrix} \right\}
$.
Now consider a path $h$ belonging to the second of these sets.
Since $\mu_{j+1}$ is between $\mu^*_j=\mu^*_{j+1}$ and $\mu_j$,
then if $h$ attains $\mu^*_{j+1}$ and does so before it attains
any $\mu_{j+1}$, then it necessarily attains $\mu^*_j$ and does so
before it attains any $\mu_j$.
Thus $h$ belongs to the first of the above sets, whence the first
part of the lemma follows.
The second expression is proved in an analogous way.
\cqfd
\medskip

\subsection{Transforming the passing function}\label{PassingGenSec}

In this section, we determine how the generating functions defined
in the previous section behave under the $\B$-transform and the
$\B\D$-transform.

\begin{lemma}\label{MazyBijLem}
For $1\le p<p'$ with $p'>2p$, let $1\le a,b<p'$, $e,f\in\{0,1\}$,
and $m_0,m_1\ge0$.
Let $\boldmu,\boldmu^*,\boldnu,\boldnu^*$ be an
interfacial mazy-four in the $(p,p')$-model sandwiching $(a,b)$.
If $\delta^{p,p'}_{a,e}=1$ we restrict to $d^L>0$,
and likewise if $\delta^{p,p'}_{b,f}=1$ we restrict to $d^R>0$.
Set $a'=a+e+\lfloor ap/p'\rfloor$ and $b'=b+f+\lfloor bp/p'\rfloor$.
Define the vectors 
$\boldmu',\boldmu^{*\prime},\boldnu',\boldnu^{*\prime}$,
by setting $\mu_j'=\mu_j+\rho^{p,p'}(\mu_j)$ and
$\mu^{*\prime}_j=\mu^*_j+\rho^{p,p'}(\mu^*_j)$ for $1\le j\le d^L$;
and setting $\nu_j'=\nu_j+\rho^{p,p'}(\nu_j)$ and
$\nu^{*\prime}_j=\nu^*_j+\rho^{p,p'}(\nu^*_j)$ for $1\le j\le d^R$.

If $\delta^{p,p'}_{a,e}=0$,
then the map $(h,k,\lambda)\mapsto h'$ effected by the action of
a $\B(k,\lambda)$-transform on $h$, is a bijection between
$$
\bigcup_{k}
\P^{p,p'}_{a,b,e,f}(m_1,2k+2m_1-m_0)\left\{
\begin{matrix} \boldmu^{\phantom{*}};\boldnu^{\phantom{*}}\\
           \boldmu^*;\boldnu^* \end{matrix} \right\}
\times\Y(k,m_1)
$$
and
$\P^{p,p'+p}_{a',b',e,f}(m_0,m_1)\left\{
\begin{smallmatrix} \boldmu^{\prime\phantom{*}};\boldnu^{\prime\phantom{*}}\\
           \boldmu^{*\prime};\boldnu^{*\prime} \end{smallmatrix} \right\}$.

If $\delta^{p,p'}_{a,e}=1$
and $\mu_{d^L}=a-(-1)^e$,
then the map $(h,k,\lambda)\mapsto h'$ effected by the action of
a $\B(k,\lambda)$-transform on $h$, is a bijection between
$$
\bigcup_{k}
\P^{p,p'}_{a,b,e,f}(m_1,2k+2m_1-m_0+2)\left\{
\begin{matrix} \boldmu^{\phantom{*}};\boldnu^{\phantom{*}}\\
           \boldmu^*;\boldnu^* \end{matrix} \right\}
\times\Y(k,m_1)
$$
and
$\P^{p,p'+p}_{a',b',e,f}(m_0,m_1)\left\{
\begin{smallmatrix} \boldmu^{\prime\phantom{*}};\boldnu^{\prime\phantom{*}}\\
           \boldmu^{*\prime};\boldnu^{*\prime} \end{smallmatrix} \right\}$.

In both cases,
\begin{displaymath}
\mwt(h')=\mwt(h)+\frac{1}{4}\left( (m_0-m_1)^2 - \beta^2 \right)+\wt(\lambda),
\end{displaymath}
where $\beta=\beta^{p,p'}_{a,b,e,f}$.

Additionally, $\boldmu',\boldmu^{*\prime},\boldnu',\boldnu^{*\prime}$ are an
interfacial mazy-four in the $(p,p'+p)$-model sandwiching $(a',b')$.
\end{lemma}

\Proof It is immediate that
$\boldmu',\boldmu^{*\prime},\boldnu',\boldnu^{*\prime}$
are a mazy-four in the $(p,p'+p)$-model sandwiching $(a',b')$.
Let $1\le j\le {d^L}$. Since $\mu^*_j$ is interfacial
in the $(p,p')$-model, Lemma \ref{InterLem}(1) implies that
$\mu^{*\prime}_j=\mu^*_j+\rho^{p,p'}(\mu^*_j)$
is interfacial in the $(p,p'+p)$-model,
The analogous results hold for $\mu_j$ for $1\le j\le d^L$ and for
$\nu^*_j$ and $\nu_j$ for $1\le j\le d^R$.
It follows that $\boldmu',\boldmu^{*\prime},\boldnu',\boldnu^{*\prime}$
are an interfacial mazy-four in the $(p,p'+p)$-model sandwiching $(a',b')$.

Lemma \ref{InterLem}(1) also shows that
$\omega^{p,p'+p}(\mu^{*\prime}_j)
=\omega^{p,p'+p}(\mu^*_j+\rho^{p,p'}(\mu^*_j))
=\omega^{p,p'}(\mu^*_j)$,
with again, the analogous results holding for $\mu_j$ for $1\le j\le d^L$
and for $\nu^*_j$ and $\nu_j$ for $1\le j\le d^R$.

Let $h\in
\P^{p,p'}_{a,b,e,f}(m_1,m)\left\{
\begin{matrix} \boldmu^{\phantom{*}};\boldnu^{\phantom{*}}\\
           \boldmu^*;\boldnu^* \end{matrix} \right\}$
and let $h'$ result from the action of the $\B(k,\lambda)$-transform on $h$.
Lemma \ref{PassBLem1} implies that either
$\omega(h')=\omega(h)$,
$\omega(a)\omega(h')=\omega(h)$,
$\omega(h')\omega(b)=\omega(h)$,
or $\omega(a)\omega(h')\omega(b)=\omega(h)$.
In the case $\omega(h')=\omega(h)$, it follows immediately that $h'$ is
mazy-compliant with $\boldmu',\boldmu^{*\prime},\boldnu',\boldnu^{*\prime}$.

In the case $\omega(a)\omega(h')=\omega(h)$, then necessarily
$a$ is interfacial in the $(p,p')$-model.
Since $\boldmu,\boldmu^*$ are a mazy-pair sandwiching $a$,
$a\ne\mu^*_j$ for $1\le j\le d^L$.
It follows that $h'$ satisfies conditions 1 and 2 for being
mazy-compliant with $\boldmu',\boldmu^{*\prime},\boldnu',\boldnu^{*\prime}$.
For the other conditions, note first that Lemma \ref{PassBLem1}
implies that $\delta^{p,p'}_{a,e}=1$ so that $d^L>0$, whereupon if
$d^R>0$, there exists $i'\ge i>0$ such that
$h_i=\mu^*_1$ and $h_{i'}=\nu^*_1$.
Since $\boldnu,\boldnu^{*}$ are a mazy-pair sandwiching $b$,
it follows that if $1\le j\le d^R$, then there exists $i''\ge i'$
such that $h_{i''}=\nu^*_j$.
Then from $\omega(a)\omega(h')=\omega(h)$, it follows that $h'$
satisfies the remaining criteria for being mazy-compliant with 
$\boldmu',\boldmu^{*\prime},\boldnu',\boldnu^{*\prime}$.

The cases $\omega(h')\omega(b)=\omega(h)$ and
$\omega(a)\omega(h')\omega(b)=\omega(h)$ are dealt with in a similar way.
Thereupon, making use of Lemma \ref{BresLem}, $h'\in
\P^{p,p'+p}_{a',b',e,f}(L',m_1)\left\{
\begin{smallmatrix} \boldmu^{\prime\phantom{*}};\boldnu^{\prime\phantom{*}}\\
           \boldmu^{*\prime};\boldnu^{*\prime} \end{smallmatrix} \right\}$,
for some $L'$.

In the case $\delta^{p,p'}_{a,e}=0$,
if $\pi(h)=1$ then $d(h)\ne e$.
Lemma \ref{BresLem} thus gives $L'=2m_1-m+2k$.

In the case $\delta^{p,p'}_{a,e}=1$, we have $\mu_{d^L}=a-(-1)^e$
which implies that $h_1=a+(-1)^e$.
Thence, $\pi(h)=1$ and $d(h)=e$.
Lemma \ref{BresLem} thus gives $L'=2m_1-m+2k+2$.

Now consider
$h'\in\P^{p,p'+p}_{a',b',e,f}(L',m_1)\left\{
\begin{smallmatrix} \boldmu^{\prime\phantom{*}};\boldnu^{\prime\phantom{*}}\\
           \boldmu^{*\prime};\boldnu^{*\prime} \end{smallmatrix} \right\}$.
Lemma \ref{StartPtLem} shows that
$\delta^{p,p'+p}_{a',e}=\delta^{p,p'+p}_{b',f}=0$,
whereupon Lemma \ref{UniqueLem} shows that there is a unique triple
$(h,k,\lambda)$ with $h\in\P^{p,p'}_{a,b,e,f}(m_1,m)$ 
and $\lambda$ is a partition having at most $k$ parts,
the greatest of which does not exceed $m_1$,
such that $h'$ arises from the action of a $\B(k,\lambda)$-transform on $h$.
Now employ Lemma \ref{PassBLem1}.
Whether $\omega(h')=\omega(h)$,
$\omega(a)\omega(h')=\omega(h)$,
$\omega(h')\omega(b)=\omega(h)$,
or $\omega(a)\omega(h')\omega(b)=\omega(h)$,
it is readily seen that $h$ is mazy-compliant with
$\boldmu,\boldmu^{*},\boldnu,\boldnu^{*}$.
The bijection follows on taking $m_0=L'$.

The expression for $\mwt(h)$ also results from Lemma \ref{BresLem}.
\cqfd
\medskip

\begin{corollary}\label{MazyBijCor} Let all parameters be as in
the first paragraph of Lemma \ref{MazyBijLem}.
If either $\delta^{p,p'}_{a,e}=0$, or
$\delta^{p,p'}_{a,e}=1$ and $\mu_{d^L}=a-(-1)^e$, then:
\begin{multline*}
\mchi^{p,p'+p}_{a',b',e,f}(m_0,m_1)\left\{
\begin{matrix} \boldmu^{\prime\phantom{*}};\boldnu^{\prime\phantom{*}}\\
           \boldmu^{*\prime};\boldnu^{*\prime} \end{matrix} \right\} \\
=
q^{\frac{1}{4}\left( (m_0-m_1)^2 - \beta^2 \right)}\!\!
\sum_{\scriptstyle m\equiv m_0\atop
  \scriptstyle\strut(\mbox{\scriptsize\rm mod}\,2)}
\left[{\frac{1}{2}(m_0+m)\atop m_1}\right]_q
\mchi^{p,p'}_{a,b,e,f}(m_1,m+2\delta^{p,p'}_{a,e})\left\{
\begin{matrix} \boldmu^{\phantom{*}};\boldnu^{\phantom{*}}\\
           \boldmu^*;\boldnu^* \end{matrix} \right\}\!,
\end{multline*}
where $\beta=\beta^{p,p'}_{a,b,e,f}$.

In addition,
$\alpha^{p,p'+p}_{a',b'}=\alpha^{p,p'}_{a,b}+\beta^{p,p'}_{a,b,e,f}$
and
$\beta^{p,p'+p}_{a',b',e,f}=\beta^{p,p'}_{a,b,e,f}$.
\end{corollary}

\Proof Apart from the case for which $m_1=0$ and $e\ne f$,
the first statement follows immediately from Lemma \ref{MazyBijLem}
on setting $m=2k+2m_1-m_0$, once it is noted, via Lemma \ref{PartitionGenLem},
that $\left[{k+m_1\atop m_1}\right]_q$ is the
generating function for $\Y(k,m_1)$.

In the case $m_1=0$ and $e\ne f$, after noting that $\delta^{p,p'}_{b',f}=0$,
it is readily seen that if a path $h'$ is to contribute to the left side,
then $m_0$ is odd, $\vert a'-b'\vert=1$, the band between $a'$ and $b'$ is
even and $h'$ alternates between heights $a'$ and $b'$.
Also note that $a'-e=b'-f$ from which we obtain $a=b$.
There is only one path $h'$ satisfying the above.
Via the same calculation as in the proof of \ref{SeedLem},
$\mwt(h')=\frac14(m_0^2-1)$.
Then, since $\mchi^{p,p'}_{a,b,e,f}(0,m)\{\}=\delta_{m,1}$,
the expression certainly holds if $d^L=d^R=0$.
If either $d^L>0$ or $d^R>0$ then the right side is zero since each
contributing path must attain $\mu^*_{d^L}\ne a$ or
$\nu^*_{d^R}\ne b$ respectively.
Correspondingly, the left side is zero in the case $d^L>0$ since because
$\mu^*_{d^L}$ is interfacial and
$\mu^{*\prime}_{d^L}=\mu^*_{d^L}+\rho^{p,p'}(\mu^*_{d^L})$,
it follows that $\mu^{*\prime}_{d^L}\ne a'$ and $\mu^{*\prime}_{d^L}\ne b'$.
That the left side is zero in the $d^R>0$ case follows similarly.

By the above reasoning, the left side is zero if either $m_0$ is even
or $a\ne b$. It is easily checked that the right side is also zero in
these cases.

The second statement is Lemma \ref{ParamLem}.
\cqfd
\medskip


We also require $\B\D$-transform analogues of Lemma \ref{MazyBijLem}
and Corollary \ref{MazyBijCor}. However, the restrictions that
applied in the cases where $\delta^{p,p'}_{a,e}=1$ or
$\delta^{p,p'}_{b,f}=1$ are not sufficient for what is required later.

\begin{lemma}\label{MazyDijLem}
For $1\le p<p'$ with $p'<2p$, let $1\le a,b<p'$, $e,f\in\{0,1\}$,
and $m_0,m_1\ge0$.
Let $\boldmu,\boldmu^*,\boldnu,\boldnu^*$ be an
interfacial mazy-four in the $(p,p')$-model sandwiching $(a,b)$.
If $\delta^{p'-p,p'}_{a,e}=1$ we restrict to $d^L>0$,
and likewise if $\delta^{p'-p,p'}_{b,f}=1$ we restrict to $d^R>0$.
Set $a'=a+1-e+\lfloor ap/p'\rfloor$ and $b'=b+1-f+\lfloor bp/p'\rfloor$.
Define the vectors 
$\boldmu',\boldmu^{*\prime},\boldnu',\boldnu^{*\prime}$,
by setting $\mu_j'=\mu_j+\rho^{p,p'}(\mu_j)$ and
$\mu^{*\prime}_j=\mu^*_j+\rho^{p,p'}(\mu^*_j)$ for $1\le j\le d^L$;
and setting $\nu_j'=\nu_j+\rho^{p,p'}(\nu_j)$ and
$\nu^{*\prime}_j=\nu^*_j+\rho^{p,p'}(\nu^*_j)$ for $1\le j\le d^R$.

If $\delta^{p'-p,p'}_{a,e}=0$,
then the map $(h,k,\lambda)\mapsto h'$ effected by the action of
a $\D$-transform on $h$ followed by a $\B(k,\lambda)$-transform,
is a bijection between
$$
\bigcup_{k}
\P^{p'-p,p'}_{a,b,e,f}(m_1,m_0-m_1-2k)\left\{
\begin{matrix} \boldmu^{\phantom{*}};\boldnu^{\phantom{*}}\\
           \boldmu^*;\boldnu^* \end{matrix} \right\}
\times\Y(k,m_1)
$$
and
$\P^{p,p'+p}_{a',b',1-e,1-f}(m_0,m_1)\left\{
\begin{smallmatrix} \boldmu^{\prime\phantom{*}};\boldnu^{\prime\phantom{*}}\\
           \boldmu^{*\prime};\boldnu^{*\prime} \end{smallmatrix} \right\}$.

If $\delta^{p'-p,p'}_{a,e}=1$ and $\mu_{d^L}=a-(-1)^e$,
then the map $(h,k,\lambda)\mapsto h'$ effected by the action of
a $\B(k,\lambda)$-transform on $h$, is a bijection between
$$
\bigcup_{k}
\P^{p'-p,p'}_{a,b,e,f}(m_1,m_0-m_1-2k+2)\left\{
\begin{matrix} \boldmu^{\phantom{*}};\boldnu^{\phantom{*}}\\
           \boldmu^*;\boldnu^* \end{matrix} \right\}
\times\Y(k,m_1)
$$
and
$\P^{p,p'+p}_{a',b',1-e,1-f}(m_0,m_1)\left\{
\begin{smallmatrix} \boldmu^{\prime\phantom{*}};\boldnu^{\prime\phantom{*}}\\
           \boldmu^{*\prime};\boldnu^{*\prime} \end{smallmatrix} \right\}$.

In both cases,
\begin{displaymath}
\mwt(h')=
\frac{1}{4}\left( m_1^2 + (m_0-m_1)^2 - \alpha^2 
- \beta^2 \right)+\wt(\lambda)-\mwt(h),
\end{displaymath}
where $\alpha=\alpha^{p,p'}_{a,b}$ and $\beta=\beta^{p,p'}_{a,b,1-e,1-f}$.

Additionally, $\boldmu',\boldmu^{*\prime},\boldnu',\boldnu^{*\prime}$ are an
interfacial mazy-four in the $(p,p'+p)$-model sandwiching $(a',b')$.
\end{lemma}

\Proof It is immediate that
$\boldmu',\boldmu^{*\prime},\boldnu',\boldnu^{*\prime}$
are a mazy-four in the $(p,p'+p)$-model sandwiching $(a',b')$.
Let $1\le j\le d^L$. Since $\mu^*_j$ 
is interfacial in the $(p'-p,p')$-model,
Lemma \ref{InterLem} implies that $\mu^*_j$
is interfacial in the $(p,p')$-model and then also that
$\mu^{*\prime}_j=\mu^*_j+\rho^{p,p'}(\mu^*_j)$ is interfacial in
the $(p,p'+p)$-model.
The analogous results hold for $\mu_j$ for $1\le j\le d^L$ and for
$\nu^*_j$ and $\nu_j$ for $1\le j\le d^R$.
It follows that $\boldmu',\boldmu^{*\prime},\boldnu',\boldnu^{*\prime}$
are an interfacial mazy-four in the $(p,p'+p)$-model sandwiching $(a',b')$.

Lemma \ref{InterLem} also shows that if $R=\lfloor rp/(p'-p)\rfloor$ then
$\omega^{p'-p,p'}(\mu^*_j)=r^-\implies
\omega^{p,p'}(\mu^*_j)=R^+\implies
\omega^{p,p'+p}(\mu^{*\prime}_j)=R^+$, and
$\omega^{p'-p,p'}(\mu^*_j)=r^+\implies
\omega^{p,p'}(\mu^*_j)=(R+1)^-\implies
\omega^{p,p'+p}(\mu^{*\prime}_j)=(R+1)^-$.
The analogous results hold for $\mu_j$ for $1\le j\le d^L$,
and for $\nu^*_j$ and $\nu_j$ for $1\le j\le d^R$.

Let
$h\in\P^{p'-p,p'}_{a,b,e,f}(m_1,m)\left\{
\begin{smallmatrix} \boldmu^{\phantom{*}};\boldnu^{\phantom{*}}\\
           \boldmu^*;\boldnu^* \end{smallmatrix} \right\}$,
and let $\hat h$ result from the action of a $\D$-transform on $h$,
and let $h'$ result from the action of a
$\B(k,\lambda)$-transform on $\hat h$.
Lemma \ref{PassBDLem} then implies that the longest subword
$\omega'=\omega_1'\omega_2'\cdots\omega_{k'}'$ of $\omega(h')$
that comprises the symbols in $\{\sigma(x):x\in\R^{p'-p,p'}\}$,
and for which $\omega_i'\ne\omega_{i+1}'$ for $1\le i<k'$,
is such that either $\omega'=\sigma(\omega(h))$,
$\omega(a)\omega'=\sigma(\omega(h))$,
$\omega'\omega(b)=\sigma(\omega(h))$, or
$\omega(a)\omega'\omega(b)=\sigma(\omega(h))$.
In the case $\omega'=\omega(h)$, it follows immediately that $h'$ is
mazy-compliant with $\boldmu',\boldmu^{*\prime},\boldnu',\boldnu^{*\prime}$.
The other cases are dealt with as in the proof of Lemma \ref{MazyBijLem}.
Thereupon, making use of Lemma \ref{BDresLem}, $h'\in
\P^{p,p'+p}_{a',b',1-e,1-f}(L',m_1)\left\{
\begin{smallmatrix} \boldmu^{\prime\phantom{*}};\boldnu^{\prime\phantom{*}}\\
           \boldmu^{*\prime};\boldnu^{*\prime} \end{smallmatrix} \right\}$,
for some $L'$.


To determine $L'$, first consider $\delta^{p'-p,p'}_{a,e}=0$. 
Here, if $\pi(h)=1$ then $d(h)\ne e$.
Lemma \ref{BDresLem} thus gives $L'=m+m_1+2k$.

In the case $\delta^{p'-p,p'}_{a,e}=1$, we have $\mu_{d^L}=a-(-1)^e$
which implies that $h_1=a+(-1)^e$.
Thence, $\pi(h)=1$ and $d(h)=e$.
Lemma \ref{BDresLem} thus gives $L'=m+m_1+2k+2$.

Now consider
$h'\in\P^{p,p'+p}_{a',b',e,f}(L',m_1)\left\{
\begin{smallmatrix} \boldmu^{\prime\phantom{*}};\boldnu^{\prime\phantom{*}}\\
           \boldmu^{*\prime};\boldnu^{*\prime} \end{smallmatrix} \right\}$
and let $\omega'$ be the longest subword
$\omega'=\omega_1'\omega_2'\cdots\omega_{k'}'$ of $\omega(h')$
that comprises the symbols in $\{\sigma(x):x\in\R^{p'-p,p'}\}$,
and for which $\omega_i'\ne\omega_{i+1}'$ for $1\le i<k'$.
Lemma \ref{StartPtLem} shows that
$\delta^{p,p'+p}_{a',e}=\delta^{p,p'+p}_{b',f}=0$,
whereupon Lemma \ref{UniqueLem} shows that there is a unique triple
$(\hat h,k,\lambda)$ with $\hat h\in\P^{p,p'}_{a,b,e,f}(m_1,m)$ 
and $\lambda$ is a partition having at most $k$ parts,
the greatest of which does not exceed $m_1$,
such that $h'$ arises from the action of a $\B(k,\lambda)$-transform
on $\hat h$.
The $\D$-transform maps $\hat h$ to a unique
$h\in\P^{p'-p,p'}_{a,b,e,f}(m_1,m)$, so that $h'$ arises from a
unique $h$ by a $\B\D$-transform.
Lemma \ref{PassBDLem} implies that
$\omega'=\sigma(\omega(h))$,
$\omega(a)\omega'=\sigma(\omega(h))$,
$\omega'\omega(b)=\sigma(\omega(h))$, or
$\omega(a)\omega'\omega(b)=\sigma(\omega(h))$.
In either case, it follows that $h$ is mazy-compliant with
$\boldmu,\boldmu^{*},\boldnu,\boldnu^{*}$.
Therefore,
$h\in\P^{p'-p,p'}_{a,b,e,f}(m_1,m)\left\{
\begin{smallmatrix} \boldmu^{\phantom{*}};\boldnu^{\phantom{*}}\\
           \boldmu^*;\boldnu^* \end{smallmatrix} \right\}$,
and the bijection follows on taking $m_0=L'$.

Lemma \ref{BDresLem} also yields the expression for $\mwt(h)$.
\cqfd
\medskip

\begin{corollary}\label{MazyDijCor} Let all parameters be as in
the premise of Lemma \ref{MazyDijLem}.
If either $\delta^{p'-p,p'}_{a,e}=0$, or
$\delta^{p'-p,p'}_{a,e}=1$ and $\mu_{d^L}=a-(-1)^e$,
then
\begin{multline*}
\mchi^{p,p'+p}_{a',b',1-e,1-f}(m_0,m_1;q)\left\{
\begin{matrix} \boldmu^{\prime\phantom{*}};\boldnu^{\prime\phantom{*}}\\
           \boldmu^{*\prime};\boldnu^{*\prime} \end{matrix} \right\}\\
=
q^{\frac{1}{4}\left( m_1^2 + (m_0-m_1)^2 - \alpha^2 - \beta^2 \right)}\!\!
\sum_{\scriptstyle m\equiv m_0-m_1\atop
  \scriptstyle\strut(\mbox{\scriptsize\rm mod}\,2)}
\left(
\hskip2mm
\left[{\frac{1}{2}(m_0+m_1-m)\atop m_1}\right]_q
\right.\\
\left.
\times
\hskip2mm
\mchi^{p'-p,p'}_{a,b,e,f}(m_1,m+2\delta^{p'-p,p'}_{a,e};q^{-1})\left\{
\begin{matrix} \boldmu^{\phantom{*}};\boldnu^{\phantom{*}}\\
           \boldmu^*;\boldnu^* \end{matrix} \right\}
\hskip2mm
\right) ,
\end{multline*}
where $\alpha=\alpha^{p,p'}_{a,b}$ and $\beta=\beta^{p,p'}_{a,b,1-e,1-f}$.

In addition,
$\alpha^{p,p'+p}_{a',b'}=2\alpha^{p'-p,p'}_{a,b}-\beta^{p'-p,p'}_{a,b,e,f}$
and
$\beta^{p,p'+p}_{a',b',1-e,1-f}
  =\alpha^{p'-p,p'}_{a,b}-\beta^{p'-p,p'}_{a,b,e,f}$.
\end{corollary}

\Proof Apart from the case for which $m_1=0$ and $e\ne f$,
the first statement follows immediately from Lemma \ref{MazyDijLem}
on setting $m=m_0-m_1-2k$, once it is noted, via Lemma \ref{PartitionGenLem},
that $\left[{k+m_1\atop m_1}\right]_q$ is the
generating function for $\Y(k,m_1)$.
The case $m_1=0$ and $e\ne f$ is dealt with exactly as in the proof
of Corollary \ref{MazyBijCor}.

The second statement is Lemma \ref{BDPathParamLem}.
\cqfd
\medskip

\subsection{Mazy runs in the original weighting}\label{MazyOSec}

Although not required until Section \ref{TransSec}, here we define
generating functions for paths restricted in the same way as in
Section \ref{PassingSec}, but which have the original weight
function (\ref{WtDef}) applied to them.

With $\boldmu=(\mu_1,\mu_2,\ldots,\mu_{d^L})$,
$\boldmu^*=(\mu^*_1,\mu^*_2,\ldots,\mu^*_{d^L})$,
$\boldnu=(\nu_1,\nu_2,\ldots,\nu_{d^R})$,
and $\boldnu^*=(\nu^*_1,\nu^*_2,\ldots,\nu^*_{d^R})$
for $d^L,d^R\ge0$,
we say that a path $h\in\P^{p,p'}_{a,b,c}(L)$ is {\em mazy-compliant}
(with $\boldmu,\boldmu^*,\boldnu,\boldnu^*$) if the five
conditions of Section \ref{PassingSec} hold.
The set
$\P^{p,p'}_{a,b,c}(L)\left\{
\begin{smallmatrix} \boldmu^{\phantom{*}}\!&;&\boldnu^{\phantom{*}}\\
                    \boldmu^*\!&;&\boldnu^* \end{smallmatrix} \right\}$
is defined to be the subset of $\P^{p,p'}_{a,b,c}(L)$
comprising those paths $h$ which are mazy-compliant with
$\boldmu,\boldmu^*,\boldnu,\boldnu^*$.
The generating function for these paths is defined to be:
\begin{equation}\label{MazyPathOGenDef}
\ochi^{p,p'}_{a,b,c}(L;q)\left\{
\begin{matrix} \boldmu^{\phantom{*}};\boldnu^{\phantom{*}}\\
                    \boldmu^*;\boldnu^* \end{matrix} \right\}
=\sum_{h\in\P^{p,p'}_{a,b,c}(L)\left\{
\genfrac{}{}{0pt}{3}{\boldmu^{\phantom{*}}}{\boldmu^*}\!
\genfrac{}{}{0pt}{3}{;}{;}
\genfrac{}{}{0pt}{3}{\boldnu^{\phantom{*}}}{\boldnu^*} \right\}}
\hskip-5mm
q^{\owt(h)}.
\end{equation}

\begin{lemma}\label{RecurseLem}
Let $1\le a,b<p'$ with $b$ non-interfacial in the
$(p,p')$-model, and let $\Delta\in\{\pm1\}$.
Let $\boldmu,\boldmu^*,\boldnu,\boldnu^*$ be a mazy-four in
the $(p,p')$-model sandwiching $(a,b)$.

If $d^R=0$ or $b-\Delta\ne\nu^*_{d^R}$ set $\boldnu^+=\boldnu$ and
$\boldnu^{+*}=\boldnu^*$.
Otherwise, if $d^R>0$ and $b-\Delta=\nu^*_{d^R}$,
set $d^+\ge0$ to be the smallest integer such that
$b-\Delta=\nu^*_{d^++1}$, and
set $\boldnu^+=(\nu_1,\nu_2,\ldots,\nu_{d^+})$ and
$\boldnu^{+*}=(\nu^*_1,\nu^*_2,\ldots,\nu^*_{d^+})$.

1) If $p'>2p$ and either $d^R=0$ or $b+\Delta\ne\nu_{d^R}$ then:
\begin{equation*}
\begin{split}
\ochi^{p,p'}_{a,b,b-\Delta}(L)\left\{
\begin{matrix} \boldmu^{\phantom{*}};\boldnu^{\phantom{*}}\\
                    \boldmu^*;\boldnu^* \end{matrix} \right\}
&=
q^{\frac12(L+\Delta(a-b))}
\ochi^{p,p'}_{a,b-\Delta,b}(L-1)\left\{
\begin{matrix} \boldmu^{\phantom{*}};\boldnu^{+\phantom{*}}\\
                    \boldmu^*;\boldnu^{+*} \end{matrix} \right\}\\
&\hskip10mm
+
\ochi^{p,p'}_{a,b+\Delta,b}(L-1)\left\{
\begin{matrix} \boldmu^{\phantom{*}};\boldnu^{\phantom{*}}\\
                    \boldmu^*;\boldnu^* \end{matrix} \right\},\\
\ochi^{p,p'}_{a,b,b+\Delta}(L)\left\{
\begin{matrix} \boldmu^{\phantom{*}};\boldnu^{\phantom{*}}\\
                    \boldmu^*;\boldnu^* \end{matrix} \right\}
&=
\ochi^{p,p'}_{a,b-\Delta,b}(L-1)\left\{
\begin{matrix} \boldmu^{\phantom{*}};\boldnu^{+\phantom{*}}\\
                    \boldmu^*;\boldnu^{+*} \end{matrix} \right\}\\
&\hskip10mm
+ q^{\frac12(L-\Delta(a-b))}
\ochi^{p,p'}_{a,b+\Delta,b}(L-1)\left\{
\begin{matrix} \boldmu^{\phantom{*}};\boldnu^{\phantom{*}}\\
                    \boldmu^*;\boldnu^* \end{matrix} \right\}.
\end{split}
\end{equation*}

2) If $p'>2p$ and either both $d^R>0$ and $b+\Delta=\nu_{d^R}$, or
$b+\Delta\in\{0,p'\}$ then:
\begin{equation*}
\begin{split}
\ochi^{p,p'}_{a,b,b-\Delta}(L)\left\{
\begin{matrix} \boldmu^{\phantom{*}};\boldnu^{\phantom{*}}\\
                    \boldmu^*;\boldnu^* \end{matrix} \right\}
&=
q^{\frac12(L+\Delta(a-b))}
\ochi^{p,p'}_{a,b,b+\Delta}(L)\left\{
\begin{matrix} \boldmu^{\phantom{*}};\boldnu^{\phantom{*}}\\
                    \boldmu^*;\boldnu^{*} \end{matrix} \right\}\\
&=
q^{\frac12(L+\Delta(a-b))}
\ochi^{p,p'}_{a,b-\Delta,b}(L-1)\left\{
\begin{matrix} \boldmu^{\phantom{*}};\boldnu^{+\phantom{*}}\\
                    \boldmu^*;\boldnu^{+*} \end{matrix} \right\}.
\end{split}
\end{equation*}

3) If $p'<2p$ and either $d^R=0$ or $b+\Delta\ne\nu_{d^R}$ then:
\begin{equation*}
\begin{split}
\ochi^{p,p'}_{a,b,b-\Delta}(L)\left\{
\begin{matrix} \boldmu^{\phantom{*}};\boldnu^{\phantom{*}}\\
                    \boldmu^*;\boldnu^* \end{matrix} \right\}
&=
\ochi^{p,p'}_{a,b-\Delta,b}(L-1)\left\{
\begin{matrix} \boldmu^{\phantom{*}};\boldnu^{+\phantom{*}}\\
                    \boldmu^*;\boldnu^{+*} \end{matrix} \right\}\\
&\hskip10mm
+
q^{\frac12(L-\Delta(a-b))}
\ochi^{p,p'}_{a,b+\Delta,b}(L-1)\left\{
\begin{matrix} \boldmu^{\phantom{*}};\boldnu^{\phantom{*}}\\
                    \boldmu^*;\boldnu^* \end{matrix} \right\},\\
\ochi^{p,p'}_{a,b,b+\Delta}(L)\left\{
\begin{matrix} \boldmu^{\phantom{*}};\boldnu^{\phantom{*}}\\
                    \boldmu^*;\boldnu^* \end{matrix} \right\}
&= q^{\frac12(L+\Delta(a-b))}
\ochi^{p,p'}_{a,b-\Delta,b}(L-1)\left\{
\begin{matrix} \boldmu^{\phantom{*}};\boldnu^{+\phantom{*}}\\
                    \boldmu^*;\boldnu^{+*} \end{matrix} \right\}\\
&\hskip10mm
+
\ochi^{p,p'}_{a,b+\Delta,b}(L-1)\left\{
\begin{matrix} \boldmu^{\phantom{*}};\boldnu^{\phantom{*}}\\
                    \boldmu^*;\boldnu^* \end{matrix} \right\}.
\end{split}
\end{equation*}

4) If $p'<2p$ and either both $d^R>0$ and $b+\Delta=\nu_{d^R}$, or
$b+\Delta\in\{0,p'\}$ then:
\begin{equation*}
\begin{split}
\ochi^{p,p'}_{a,b,b-\Delta}(L)\left\{
\begin{matrix} \boldmu^{\phantom{*}};\boldnu^{\phantom{*}}\\
                    \boldmu^*;\boldnu^* \end{matrix} \right\}
&=
q^{-\frac12(L+\Delta(a-b))}
\ochi^{p,p'}_{a,b,b+\Delta}(L)\left\{
\begin{matrix} \boldmu^{\phantom{*}};\boldnu^{\phantom{*}}\\
                    \boldmu^*;\boldnu^{*} \end{matrix} \right\}\\
&=
\ochi^{p,p'}_{a,b-\Delta,b}(L-1)\left\{
\begin{matrix} \boldmu^{\phantom{*}};\boldnu^{+\phantom{*}}\\
                    \boldmu^*;\boldnu^{+*} \end{matrix} \right\}.
\end{split}
\end{equation*}

\end{lemma}

\Proof Consider
$h\in\P^{p,p'}_{a,b,b-\Delta}(L)\left\{
\begin{smallmatrix} \boldmu^{\phantom{*}}\!&;&\boldnu^{\phantom{*}}\\
                    \boldmu^*\!&;&\boldnu^* \end{smallmatrix} \right\}$.
If either $d^R=0$ or $b+\Delta\ne\nu_{d^R}$, the first $(L-1)$ segments
of $h$ constitute a path which is either a member of
$\P^{p,p'}_{a,b-\Delta,b}(L-1)\left\{
\begin{smallmatrix} \boldmu^{\phantom{*}}\!&;&\boldnu^{+\phantom{*}}\\
                    \boldmu^*\!&;&\boldnu^{+*} \end{smallmatrix} \right\}$
or a member of
$\P^{p,p'}_{a,b+\Delta,b}(L-1)\left\{
\begin{smallmatrix} \boldmu^{\phantom{*}}\!&;&\boldnu^{\phantom{*}}\\
                    \boldmu^*\!&;&\boldnu^* \end{smallmatrix} \right\}$.
(In the first of these, we use $\boldnu^{+}$ and $\boldnu^{+*}$
instead of $\boldnu$ and $\boldnu^{*}$ when $b-\Delta=\nu^*_{d^R}$
because then $\boldnu$ and $\boldnu^{*}$ are not a mazy-pair
sandwiching $b-\Delta$.)
In either case $p'>2p$ or $p'<2p$,
consideration of the weight of the $L$th vertex of $h$,
then yields the first identity between generating functions in 1) and 3).
The second identity in these cases follows in a similar way.

Now consider $h$ as above when either both $d^R>0$ and $b+\Delta=\nu_{d^R}$,
or $b+\Delta\in\{0,p'\}$.
Necessarily $h_{L-1}=b-\Delta$.
Thus when $p'>2p$ the $L$th vertex is scoring, and when $p'<2p$
the $L$th vertex is non-scoring.
Changing its direction yields the first identity in both 2) and 4).
The first $(L-1)$ segments of $h$ yields a path which is a member of
$\P^{p,p'}_{a,b-\Delta,b}(L-1)\left\{
\begin{smallmatrix} \boldmu^{\phantom{*}}\!&;&\boldnu^{+\phantom{*}}\\
                    \boldmu^*\!&;&\boldnu^{+*} \end{smallmatrix} \right\}$.
The second identity in 2) and 4) follows as in the first paragraph above.
\cqfd

\newpage

\setcounter{section}{5}

\section{Extending and truncating paths}\label{ExTrunSec}

\subsection{Extending paths}\label{ExtendSec}

In this section, we specify a process by which a path
$h\in\P^{p,p'}_{a,b,e,f}(L)$ may be extended by
a single unit to its left, or by a single unit to its right.
Consequently, the new path $h'$ is of length $L'=L+1$.
An extension on the right may follow one on the left to yield
a path of length $L+2$.

Throughout this section,
$\boldmu,\boldmu^*,\boldnu,\boldnu^*$ are a
mazy-four in the $(p,p')$-model sandwiching $(a,b)$, with
$\boldmu$, $\boldmu^*$ each of dimension $d^L\ge0$, and
$\boldnu$, $\boldnu^*$ each of dimension $d^R\ge0$.

We restrict path extension on the right to the cases where
$\delta^{p,p'}_{b,f}=0$
so that the post-segment of $h$ lies in the even band.
The extended path $h'$ has endpoint $b'=b+\Delta$,
where if $b$ is interfacial we permit both $\Delta=\pm1$,
and if $b$ is not interfacial we permit only $\Delta=(-1)^f$
(in this latter case, the $(L+1)$th segment of $h'$ lies in the
same direction as the post-segment of $h$).
We specify that the $(L+1)$th vertex of $h'$ is a scoring (peak) vertex.
Thus, on setting $f'=0$ if $b'=b-1$, and $f'=1$ if $b'=b+1$,
we have $h'\in\P^{p,p'}_{a,b',e,f'}(L+1)$ given by
$h'_{i}=h_{i}$ for $0\le i\le L$ and $h'_{L+1}=b'$.

\begin{figure}[ht]
  \psfrag{b}{\footnotesize$b$}
  \psfrag{c}{\footnotesize$b'$}
\includegraphics[scale=1.00]{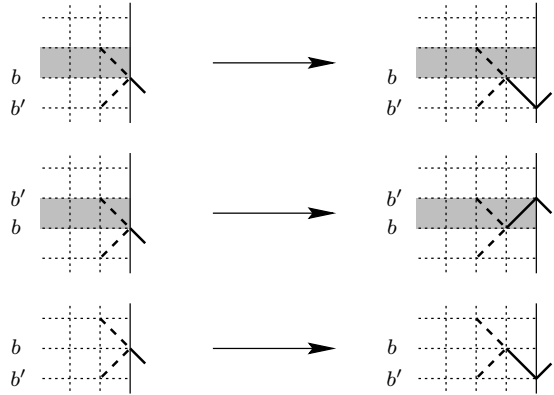}
\caption{Path extension on the right.}
\label{Extend2Fig}
\medskip
\end{figure}

The three diagrams in Fig.~\ref{Extend2Fig} depict the extending process
when $f=1$.

\noindent The $f=0$ cases may be obtained by reflecting these diagrams
in a horizontal axis.

\begin{lemma}\label{Extend2Lem}

Let $1\le a,b<p'$ and $e,f\in\{0,1\}$ with $\delta^{p,p'}_{b,f}=0$,
and let $L\ge0$.
If $b$ is interfacial let $\Delta=\pm1$, and otherwise set $\Delta=(-1)^f$.
Then set $b'=b+\Delta$ and let $f'\in\{0,1\}$ be such that $\Delta=-(-1)^{f'}$.
If $h'\in\P^{p,p'}_{a,b',e,f'}(L')$ is obtained from
$h\in\P^{p,p'}_{a,b,e,f}(L)$ by the above process of path extension, then:
\begin{displaymath}
\begin{array}{l}
\bullet\quad L'=L+1;\\[1mm]
\bullet\quad m(h')=m(h);\\[1mm]
\bullet\quad \mwt(h')=\mwt(h)+ \frac{1}{2}(L-\Delta\alpha(h)).
\end{array}
\end{displaymath}
In addition, if $\lfloor bp/p'\rfloor=\lfloor b'p/p'\rfloor$ then
$f'=1-f$, and if $\lfloor bp/p'\rfloor\ne\lfloor b'p/p'\rfloor$
then $f'=f$.
Furthermore, $\alpha^{p,p'}_{a,b'}=\alpha^{p,p'}_{a,b}+\Delta$
and $\beta^{p,p'}_{a,b',e,f'}=\beta^{p,p'}_{a,b,e,f}+\Delta$.
\end{lemma}

\Proof
Let $h$ have striking sequence
$\left({a_1\atop b_1}\:{a_2\atop b_2}\:{a_3\atop b_3}\:
 {\cdots\atop\cdots}\:{a_l\atop b_l} \right)^{(e,f,d)}$
and assume first that $L>0$.
With $\delta^{p,p'}_{b,f}=0$,
it is readily checked that the $L$th vertex of $h'$ is scoring
if and only if the $L$th vertex of $h$ is scoring.
This holds even if $L=0$.
Since the $(L+1)$th vertex of $h'$ is scoring, it immediately
follows that $m(h')=m(h)$.

If $L=0$, then necessarily $\alpha(h)=0$ whereupon the expression for
$\mwt(h')$ holds.
Now consider $L>0$.
If the extending segment is in the same direction as the $L$th
segment, $h'$ has striking sequence
$\left({a_1\atop b_1}\:{a_2\atop b_2}\:{a_3\atop b_3}\:
 {\cdots\atop\cdots}\:{a_l\atop b_l+1} \right)^{(e,f',d)}$
and $\Delta=-(-1)^{d+l}$.
If the extending segment is in the direction opposite to that of
the $L$th segment, $h'$ has striking sequence
$\left({a_1\atop b_1}\:{a_2\atop b_2}\:
 {\cdots\atop\cdots}\:{a_l\atop b_l}\:{0\atop 1} \right)^{(e,f',d)}$
and $\Delta=(-1)^{d+l}$.


For $1\le i\le l$, let $w_i=a_i+b_i$.
We find $\alpha(h)=-(-1)^{d+l}((w_l+w_{l-2}\cdots)-(w_{l-1}+w_{l-3}+\cdots))$.
In the first case above, Lemma \ref{WtHashLem} gives
$\mwt(h')=\mwt(h)+(w_{l-1}+w_{l-3}+w_{l-5}+\cdots)$, whereupon
we obtain $\mwt(h')=\mwt(h)+\frac12(L(h)-\Delta\alpha(h))$.
In the second case above, Lemma \ref{WtHashLem} gives
$\mwt(h')=\mwt(h)+(w_{l}+w_{l-2}+w_{l-4}+\cdots)$, and we again
obtain $\mwt(h')=\mwt(h)+\frac12(L(h)-\Delta\alpha(h))$.

In the $L=0$ case, let $d'=d(h')$.
It is easily seen that $d'=f$ if and only if $\pi(h')=0$.
It then follows that $h'$ has striking sequence
$\left( {0\atop 1} \right)^{(e,f',d')}$ and that
$m(h')=m(h)$, $\alpha(h)=0$ and $\mwt(h')=\mwt(h)=0$ thus
verifying the first statement in the $L=0$ case.

The penultimate statement follows from the definitions, noting that
$\lfloor b'p/p'\rfloor\ne\lfloor bp/p'\rfloor$ only if
$b$ is interfacial in the $(p,p')$-model and $\Delta=-(-1)^f$.

That $\alpha^{p,p'}_{a,b'}=\alpha^{p,p'}_{a,b}+\Delta$ is immediate.
$\beta^{p,p'}_{a,b',e,f'}=\beta^{p,p'}_{a,b,e,f}+\Delta$ follows
in the $\lfloor bp/p'\rfloor\ne\lfloor b'p/p'\rfloor$ case because then
$\lfloor b'p/p'\rfloor=\lfloor bp/p'\rfloor+\Delta$ and $f'=f$,
and in the $\lfloor bp/p'\rfloor=\lfloor b'p/p'\rfloor$ case
because $f=1-f'$ and $\Delta=-(-1)^{f'}=2f'-1$.
\cqfd
\medskip

\begin{lemma}\label{ExtGen2Lem}
Let $1\le a,b<p'$ and $e,f\in\{0,1\}$ with $\delta^{p,p'}_{b,f}=0$.
If $b$ is interfacial let $\Delta=\pm1$, and otherwise set $\Delta=(-1)^f$.
Then set $b'=b+\Delta$ and let $f'\in\{0,1\}$
be such that $\Delta=-(-1)^{f'}$.
Define $\boldnu'=(\nu_1,\ldots,\nu_{d^R},b'+\Delta)$ and
$\boldnu^{*\prime}=(\nu^*_1,\ldots,\nu^*_{d^R},b'-\Delta)$.
If either $d^R=0$ or both $b'\ne\nu_{d^R}$ and $b'\ne\nu^*_{d^R}$, then
$\boldmu,\boldmu^*,\boldnu',\boldnu^{*\prime}$ are a mazy-four
sandwiching $(a,b')$ and
\begin{equation*}
\mchi^{p,p'}_{a,b',e,f'}(L,m)\!\left\{
\begin{matrix}
         \boldmu^{\phantom{*}};\boldnu^{\prime\phantom{*}}\\
         \boldmu^*;\boldnu^{*\prime}
\end{matrix} \!\right\}
= q^{\frac{1}{2}(L-1-\Delta\alpha)}
\mchi^{p,p'}_{a,b,e,f}(L-1,m)\!\left\{
\begin{matrix} \boldmu^{\phantom{*}};\boldnu^{\phantom{*}}\\
           \boldmu^*;\boldnu^* \end{matrix} \!\right\},
\end{equation*}
where $\alpha=\alpha^{p,p'}_{a,b}$.

In addition, $\alpha^{p,p'}_{a,b'}=\alpha^{p,p'}_{a,b}+\Delta$,
$\beta^{p,p'}_{a,b',e,f'}=\beta^{p,p'}_{a,b,e,f}+\Delta$,
and if $f'=1-f$ then $\lfloor b'p/p'\rfloor=\lfloor bp/p'\rfloor$,
and if $f'=f$ then $\lfloor b'p/p'\rfloor=\lfloor bp/p'\rfloor+\Delta$.

\end{lemma}

\Proof
That $\boldmu,\boldmu^*,\boldnu',\boldnu^{*\prime}$ are a mazy-four
sandwiching $(a,b')$ follows immediately if $d^R=0$.
If $d^R>0$, it follows after noting that
since $b$ is strictly between $\nu_{d^R}$ and $\nu^*_{d^R}$,
and $b'\ne\nu_{d^R}$ and $b'\ne\nu^*_{d^R}$, then $b'+\Delta$
and $b'-\Delta$ are both between $\nu_{d^R}$ and $\nu^*_{d^R}$.

Let $L>0$ and $h\in\P^{p,p'}_{a,b,e,f}(L-1,m)\left\{
\begin{smallmatrix} \boldmu^{\phantom{*}};\boldnu^{\phantom{*}}\\
           \boldmu^*;\boldnu^* \end{smallmatrix} \right\}$.
Extend this path on the right to obtain $h'$ with $h'_L=b'=b+\Delta$.
Clearly, $h'$ attains $b=b'-\Delta$ and does so after it attains
any $b'+\Delta$. It follows that $h'$ is mazy-compliant
with $\boldmu,\boldmu^*,\boldnu',\boldnu^{*\prime}$.
Then, via Lemma \ref{Extend2Lem},
$h'\in\P^{p,p'}_{a,b',e,f'}(L,m)\left\{
\begin{smallmatrix}
         \boldmu^{\phantom{*}};\boldnu^{\prime\phantom{*}}\\
         \boldmu^*;\boldnu^{*\prime}
\end{smallmatrix} \right\}$.
Conversely, any such $h'$ arises from some
$h\in\P^{p,p'}_{a,b,e,f}(L-1,m)\left\{
\begin{smallmatrix} \boldmu^{\phantom{*}};\boldnu^{\phantom{*}}\\
           \boldmu^*;\boldnu^* \end{smallmatrix} \right\}$,
in this way.
For $L>0$, the required result then follows from the expression for $\mwt(h')$
given in Lemma \ref{Extend2Lem}, and $\alpha(h)=\alpha^{p,p'}_{a,b}$
from Lemma \ref{BetaConstLem}.
For $L\le0$, both sides are clearly equal to $0$.

The final statement follows from Lemma \ref{Extend2Lem}.
\cqfd
\medskip


For $h\in\P^{p,p'}_{a,b,e,f}(L)$, we now define path
extension to the left in a similar way.
Here, we restrict path extension to the cases where
$\delta^{p,p'}_{a,e}=0$
so that the pre-segment of $h$ lies in the even band.
If $L=0$, we also restrict to the case $\delta^{p,p'}_{b,f}=0$.
The extended path $h'$ has startpoint $a'=a+\Delta$,
where if $a$ is interfacial we permit both $\Delta=\pm1$,
and if $a$ is not interfacial we permit only $\Delta=(-1)^e$
(in this latter case, the $0$th segment of $h'$ lies in the
same direction as the pre-segment of $h$).
We specify that the $0$th vertex of $h'$ is a scoring (peak) vertex.
Thus, on setting $e'=0$ if $a'=a-1$, and $e'=1$ if $a'=a+1$,
we have $h'\in\P^{p,p'}_{a',b,e',f}(L+1)$ given by
$h'_{i}=h_{i-1}$ for $1\le i\le L+1$ and $h'_{0}=b'$.

\begin{figure}[ht]
  \psfrag{a}{\footnotesize$a$}
  \psfrag{c}{\footnotesize$a'$}
\includegraphics[scale=1.00]{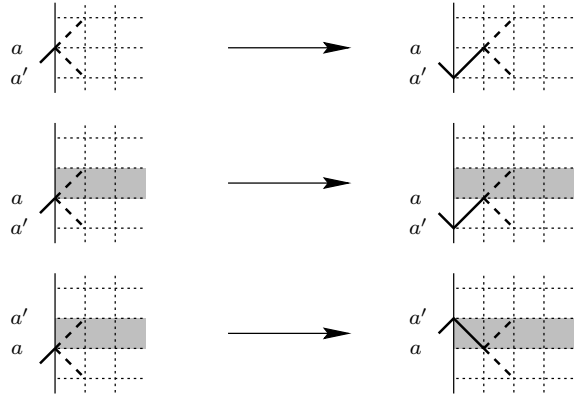}
\caption{Path extension on the left.}
\label{Extend1Fig}
\medskip
\end{figure}

The three diagrams in Fig.~\ref{Extend1Fig} depict the extending process
when $e=1$.

\noindent The cases where $e=0$
may be obtained by reflecting these diagrams in a horizontal axis.

\begin{lemma}\label{Extend1Lem}
Let $1\le a,b<p'$ and $e,f\in\{0,1\}$ with $\delta^{p,p'}_{a,e}=0$.
Let $L\ge0$, but if $\delta^{p,p'}_{b,f}=1$ restrict to $L>0$.
If $a$ is interfacial let $\Delta=\pm1$, and otherwise set $\Delta=(-1)^e$.
Then set $a'=a+\Delta$ and let $e'\in\{0,1\}$ be such that $\Delta=-(-1)^{e'}$.
If $h'\in\P^{p,p'}_{a',b,e',f}(L')$ is obtained from
$h\in\P^{p,p'}_{a,b,e,f}(L)$ by the above process of path extension, then:
\begin{displaymath}
\begin{array}{l}
\bullet\quad L'=L+1;\\[2mm]
\bullet\quad m(h')=
\left\{
  \begin{array}{ll}
m(h) \quad &
          \mbox{if } \lfloor ap/p'\rfloor=\lfloor a'p/p'\rfloor,\\[1mm]
m(h)+2 \quad & \mbox{if } \lfloor ap/p'\rfloor\ne\lfloor a'p/p'\rfloor;
  \end{array} \right.
\\[5mm]
\bullet\quad \mwt(h')=\mwt(h)+ \frac{1}{2}(L-m(h)+\Delta\beta(h)).
\end{array}
\end{displaymath}
In addition, if $\lfloor ap/p'\rfloor=\lfloor a'p/p'\rfloor$ then
$e'=1-e$, and if $\lfloor ap/p'\rfloor\ne\lfloor a'p/p'\rfloor$ then $e'=e$.
Furthermore, $\alpha^{p,p'}_{a',b}=\alpha^{p,p'}_{a,b}-\Delta$
and $\beta^{p,p'}_{a',b,e',f}=\beta^{p,p'}_{a,b,e,f}-\Delta$.
\end{lemma}

\Proof
Let $h$ have striking sequence
$\left({a_1\atop b_1}\:{a_2\atop b_2}\:{a_3\atop b_3}\:
 {\cdots\atop\cdots}\:{a_l\atop b_l} \right)^{(e,f,d)}$,
let $\pi=\pi(h)$, and assume first that $L>0$.
Consider the case $\lfloor ap/p'\rfloor=\lfloor a'p/p'\rfloor$,
where necessarily $\Delta=(-1)^e$
(in the case of interfacial $a$, this follows from $\delta^{p,p'}_{a,e}=0$).
It follows immediately from the definition that
$d(h')=e'=1-e$ and $\pi(h')=0$.

In the subcase where $e=d$ then $h'$ has striking sequence
$\left({0\atop1}\:{a_1\atop b_1}\:{a_2\atop b_2}\:
 {\cdots\atop\cdots}\:{a_l\atop b_l} \right)^{(e',f,e')}$ and $\pi=0$.
We see that $m(h')=m(h)$ and
$\mwt(h')=\mwt(h)+(b_1+b_3+\cdots)$.
Since $\beta(h)=\Delta((b_1+b_3+\cdots)-(b_2+b_4+\cdots))$
and $m(h)=(a_1+a_2+a_3+\cdots)$, we obtain
$\mwt(h')=\mwt(h)+(L-m(h)+\Delta\beta(h))/2$.

In the subcase where $e\ne d$ and $\pi=1$ then
$\left({a_1\atop b_1+1}\:{a_2\atop b_2}\:
 {\cdots\atop\cdots}\:{a_l\atop b_l} \right)^{(e',f,e')}$
is the striking sequence of $h'$.
We see that $m(h')=m(h)$ and
$\mwt(h')=\mwt(h)+(b_2+b_4+\cdots)$.
Since $\beta(h)=-\Delta((b_1+b_3+\cdots)-(b_2+b_4+\cdots))$
and $m(h)=(a_1+a_2+a_3+\cdots)$, we obtain
$\mwt(h')=\mwt(h)+(L-m(h)+\Delta\beta(h))/2$.

In the subcase where $e\ne d$ and $\pi=0$ then
$\left({a_1+1\atop b_1}\:{a_2\atop b_2}\:
 {\cdots\atop\cdots}\:{a_l\atop b_l} \right)^{(e',f,e')}$
is the striking sequence of $h'$.
We see that $m(h')=m(h)$ and
$\mwt(h')=\mwt(h)+(b_2+b_4+\cdots)$.
Since $\beta(h)=-\Delta((b_1+b_3+\cdots)-(1+b_2+b_4+\cdots))$
and $m(h)=(1+a_1+a_2+a_3+\cdots)$, we obtain
$\mwt(h')=\mwt(h)+(L-m(h)+\Delta\beta(h))/2$.

Now consider the case $\lfloor ap/p'\rfloor\ne\lfloor a'p/p'\rfloor$,
where necessarily $a$ is interfacial and $\Delta=-(-1)^e$.
It follows immediately from the definition that
$d(h')=e'=e$ and $\pi(h')=1$.

In the subcase where $e=d$ then $h'$ has striking sequence
$\left({a_1+1\atop b_1}\:{a_2\atop b_2}\:
 {\cdots\atop\cdots}\:{a_l\atop b_l} \right)^{(e',f,e')}$ and $\pi=0$.
We see that $m(h')=m(h)+2$ and
$\mwt(h')=\mwt(h)+(b_2+b_4+\cdots)$.
Since $\beta(h)=-\Delta((b_1+b_3+\cdots)-(b_2+b_4+\cdots))$
and $m(h)=(a_1+a_2+a_3+\cdots)$, we obtain
$\mwt(h')=\mwt(h)+(L-m(h)+\Delta\beta(h))/2$.

In the subcase where $e\ne d$ then $h'$ has striking sequence
$\left({1\atop0}\:{a_1\atop b_1}\:{a_2\atop b_2}\:
 {\cdots\atop\cdots}\:{a_l\atop b_l} \right)^{(e',f,e')}$ and $\pi=1$.
We see that $m(h')=m(h)+2$ and
$\mwt(h')=\mwt(h)+(b_1+b_3+\cdots)$.
Since $\beta(h)=\Delta((b_1+b_3+\cdots)-(b_2+b_4+\cdots))$
and $m(h)=(a_1+a_2+a_3+\cdots)$, we obtain
$\mwt(h')=\mwt(h)+(L-m(h)+\Delta\beta(h))/2$.

Now consider $L=0$ when $\mwt(h)=0$. If $e=f$ then $m(h)=\beta(h)=0$.
If in addition $\pi(h')=0$, then $h'$ has striking sequence
$\left({0\atop1}\right)^{(1-e,f,1-e)}$ so that $m(h')=\mwt(h')=0$,
as required.
Otherwise, if in addition $\pi(h')=1$, then $h'$ has striking sequence
$\left({1\atop0}\right)^{(e,f,e)}$ so that $m(h')=2$ and $\mwt(h')=0$,
as required.
For $e\ne f$, we have $m(h)=1$ and $\beta(h)=f-e$.
The restriction that $\delta^{p,p'}_{b,f}=0$ forces $\pi(h')=0$,
whereupon $h'$ has striking sequence
$\left({1\atop0}\right)^{(1-e,f,1-e)}$, so that
$\Delta=f-e=\beta(h)$, $m(h')=1$ and $\mwt(h')=0$,
as required.

That $\alpha^{p,p'}_{a',b}=\alpha^{p,p'}_{a,b}-\Delta$ is immediate.
$\beta^{p,p'}_{a',b,e',f}=\beta^{p,p'}_{a,b,e,f}-\Delta$ follows
in the $\lfloor ap/p'\rfloor\ne\lfloor a'p/p'\rfloor$ case because then
$\lfloor a'p/p'\rfloor=\lfloor ap/p'\rfloor+\Delta$ and $e'=e$,
and in the $\lfloor ap/p'\rfloor=\lfloor a'p/p'\rfloor$ case
because $e=1-e'$ and $\Delta=-(-1)^{e'}=2e'-1$.
\cqfd
\medskip

\begin{lemma}\label{ExtGen1Lem}
Let $1\le a,b<p'$ and $e,f\in\{0,1\}$ with $\delta^{p,p'}_{a,e}=0$.
In addition, if $\delta^{p,p'}_{b,f}=1$, restrict to $d^R>0$.
If $a$ is interfacial let $\Delta=\pm1$, and otherwise set $\Delta=(-1)^e$.
Then set $a'=a+\Delta$ and let $e'\in\{0,1\}$ be such that $\Delta=-(-1)^{e'}$.
Define $\boldmu'=(\mu_1,\ldots,\mu_{d^L},a'+\Delta)$ and
$\boldmu^{*\prime}=(\mu^*_1,\ldots,\mu^*_{d^L},a'-\Delta)$.
If either $d^L=0$ or both $a'\ne\nu_{d^L}$ and $a'\ne\nu^*_{d^L}$, then
$\boldmu',\boldmu^{*\prime},\boldnu,\boldnu^*$ are a mazy-four
sandwiching $(a',b)$ and
\begin{equation*}
\mchi^{p,p'}_{a',b,e',f}(L,m+2\delta^{p,p'}_{a',e'})\left\{
\begin{matrix}
         \boldmu^{\prime\phantom{*}};\boldnu^{\phantom{*}}\\
         \boldmu^{*\prime};\boldnu^{*}
\end{matrix} \right\}
=
q^{\frac{1}{2}(L-1-m+\Delta\beta)}
\mchi^{p,p'}_{a,b,e,f}(L-1,m)\left\{
\begin{matrix} \boldmu^{\phantom{*}};\boldnu^{\phantom{*}}\\
           \boldmu^*;\boldnu^* \end{matrix} \right\},
\end{equation*}
where $\beta=\beta^{p,p'}_{a,b,e,f}$.

In addition, $\alpha^{p,p'}_{a',b}=\alpha^{p,p'}_{a,b}-\Delta$,
$\beta^{p,p'}_{a',b,e',f}=\beta^{p,p'}_{a,b,e,f}-\Delta$,
and if $e'=1-e$ then $\lfloor a'p/p'\rfloor=\lfloor ap/p'\rfloor$,
and if $e'=e$ then $\lfloor a'p/p'\rfloor=\lfloor ap/p'\rfloor+\Delta$.
\end{lemma}

\Proof
That $\boldmu',\boldmu^{*\prime},\boldnu,\boldnu^{*}$ are a mazy-four
sandwiching $(a',b)$ follows immediately if $d^L=0$.
If $d^L>0$, it follows after noting that since $a$ is strictly between
$\mu_{d^L}$ and $\mu^*_{d^L}$,
and $a'\ne\mu_{d^L}$ and $a'\ne\mu^*_{d^L}$, then $a'+\Delta$
and $a'-\Delta$ are both between $\mu_{d^L}$ and $\mu^*_{d^L}$.

Let $L>0$ and $h\in\P^{p,p'}_{a,b,e,f}(L-1,m)\left\{
\begin{smallmatrix} \boldmu^{\phantom{*}};\boldnu^{\phantom{*}}\\
           \boldmu^*;\boldnu^* \end{smallmatrix} \right\}$.
Extend $h$ on the left to obtain $h'$ with $h'_0=a'=a+\Delta$.
Clearly, $h'$ attains $a=a'-\Delta$ and does so before it attains
any $a'+\Delta$. It follows that $h'$ is mazy-compliant
with $\boldmu',\boldmu^{*\prime},\boldnu,\boldnu^{*}$.
Then, via Lemma \ref{Extend1Lem},
$h'\in\P^{p,p'}_{a',b,e',f}(L,m)\left\{
\begin{smallmatrix}
         \boldmu^{\prime\phantom{*}};\boldnu^{\phantom{*}}\\
         \boldmu^{*\prime};\boldnu^{*}
\end{smallmatrix} \right\}$.
Conversely, any such $h'$ arises from
some $h\in\P^{p,p'}_{a,b,e,f}(L-1,m)\left\{
\begin{smallmatrix} \boldmu^{\phantom{*}};\boldnu^{\phantom{*}}\\
           \boldmu^*;\boldnu^* \end{smallmatrix} \right\}$
in this way.
For $L>0$, the required result then follows from the expression for $\mwt(h')$
given in Lemma \ref{Extend1Lem}, and $\beta(h)=\beta^{p,p'}_{a,b,e,f}$
from Lemma \ref{BetaConstLem}.
For $L\le0$, both sides are clearly equal to $0$.

The final statement also follows from Lemma \ref{Extend1Lem}.
\cqfd
\medskip

\subsection{Truncating paths}\label{AttenSec}
In this section, we specify a process by which a path
$h\in\P^{p,p'}_{a,b,e,f}(L)$, where $L>0$, may be shortened by removing
just the first segment, or by removing just the $L$th segment.
Consequently, the new path $h'$ is of length $L'=L-1$.
A shortening on the right may follow one on the left to yield
a path of length $L-2$.

Throughout this section,
$\boldmu,\boldmu^*,\boldnu,\boldnu^*$ are a
mazy-four in the $(p,p')$-model sandwiching $(a,b)$, with
$\boldmu$, $\boldmu^*$ each of dimension $d^L$, and
$\boldnu$, $\boldnu^*$ each of dimension $d^R$.

Removing the $L$th segment of $h\in\P^{p,p'}_{a,b,e,f}(L)$ is permitted
only if $h_{L-1}=b+(-1)^f$ and $\delta^{p,p'}_{b,f}=0$ so that the $L$th
segment of $h$ lies in an even band and the $L$th vertex is scoring.
We specify that $f(h')=1-f(h)$ so that the post-segment of $h'$ is in
the same direction as the $L$th segment of $h$.
If $b'=h_{L-1}$ and $f'=1-f$ then
$h'\in\P^{p,p'}_{a,b',e,f'}(L-1)$ with $h'_{i}=h_{i}$ for $0\le i\le L-1$.

\begin{lemma}\label{Atten2Lem}
Let $1\le a,b<p'$ and $e,f\in\{0,1\}$ with $\delta^{p,p'}_{b,f}=0$.
For $L>0$, let $h\in\P^{p,p'}_{a,b,e,f}(L)$ be such that $h_{L-1}=b+(-1)^f$
and let $h'\in\P^{p,p'}_{a,b',e,f'}(L')$ be obtained from $h$ by
the above process of path truncation.
If $\Delta=b-b'$ then $\Delta=(-1)^{f'}$, $f=1-f'$,
and
\begin{displaymath}
\begin{array}{l}
\bullet\quad L'=L-1;\\[1mm]
\bullet\quad m(h')=m(h);\\[1mm]
\bullet\quad \mwt(h')=\mwt(h)- \frac{1}{2}(L-\Delta\alpha(h)).
\end{array}
\end{displaymath}
In addition, $\alpha^{p,p'}_{a,b'}=\alpha^{p,p'}_{a,b}-\Delta$
and $\beta^{p,p'}_{a,b',e,f'}=\beta^{p,p'}_{a,b,e,f}-\Delta$.
\end{lemma}

\Proof That $\Delta=(-1)^{f'}$ and $f=1-f'$ and $L'=L-1$ are
immediate from the definition.
Let $h$ have striking sequence
$\left({a_1\atop b_1}\:{a_2\atop b_2}\:{a_3\atop b_3}\:
 {\cdots\atop\cdots}\:{a_l\atop b_l} \right)^{(e,f,d)}$,
whence $\Delta=-(-1)^{d+l}$.
Since $h_{L-1}=b+(-1)^f$, the $L$th vertex of $h$ is scoring
and thus $b_l\ge1$.
Since $\delta^{p,p'}_{b,f}=0$ and $f'=1-f$, it follows that
the $(L-1)$th vertex of $h'$ is scoring if and only if the
$(L-1)$th vertex of $h$ is scoring.
This holds even if $L=1$. Thus $m(h')=m(h)$.

When $a_l+b_l>1$, $h'$ has striking sequence
$\left({a_1\atop b_1}\:{a_2\atop b_2}\:{a_3\atop b_3}\:
 {\cdots\atop\cdots}\:{a_l\atop b_l-1} \right)^{(e,f',d)}$
and when $a_l+b_l>1$, $h'$ has striking sequence
$\left({a_1\atop b_1}\:{a_2\atop b_2}\:{a_3\atop b_3}\:
 {\cdots\atop\cdots}\:{a_{l-1}\atop b_{l-1}} \right)^{(e,f',d)}$.
It follows that $L'=L-1$ and $m(h')=m(h)$.

For $1\le i\le l$, let $w_i=a_i+b_i$.
Lemma \ref{WtHashLem} gives
$\mwt(h')=\mwt(h)-(w_{l-1}+w_{l-3}+w_{l-5}+\cdots)$, whereupon,
since $\alpha(h)=\Delta((w_l+w_{l-2}\cdots)-(w_{l-1}+w_{l-3}+\cdots))$,
we obtain $\mwt(h')=\mwt(h)-\frac12(L(h)-\Delta\alpha(h))$.

That $\alpha^{p,p'}_{a,b'}=\alpha^{p,p'}_{a,b}-\Delta$ is immediate.
Since $\delta^{p,p'}_{b,f}=0$,
we have $\lfloor b'p/p'\rfloor
=\lfloor (b+(-1)^f)p/p'\rfloor=\lfloor bp/p'\rfloor$.
Then $\beta^{p,p'}_{a,b',e,f'}=\lfloor b'p/p'\rfloor-\lfloor ap/p'\rfloor
+f'-e=\lfloor bp/p'\rfloor-\lfloor ap/p'\rfloor+f-(1-2f')-e=
\beta^{p,p'}_{a,b,e,f}-\Delta$, as required.
\cqfd
\medskip

\begin{lemma}\label{AttenGen2Lem}
Let $1\le a,b,b'<p'$ and $e,f\in\{0,1\}$ with
$b-b'=-(-1)^f$ and $\delta^{p,p'}_{b,f}=0$.
Set $\Delta=b-b'$ and $f'=1-f$.
If $d^R>0$, let $\nu_{d^R}^*\ne b'$ and $\nu_{d^R}=b+\Delta$.
If $d^R=0$, restrict to $b\in\{1,p'-1\}$.
Then $\boldmu,\boldmu^*,\boldnu,\boldnu^{*}$ are a mazy-four
sandwiching $(a,b')$, and if $L\ge0$ then:
\begin{equation*}
\mchi^{p,p'}_{a,b',e,f'}(L,m)\left\{
\begin{matrix} \boldmu^{\phantom{*}};\boldnu^{\phantom{*}}\\
           \boldmu^*;\boldnu^* \end{matrix} \right\}
= q^{-\frac{1}{2}(L+1-\Delta\alpha)}
\mchi^{p,p'}_{a,b,e,f}(L+1,m)\left\{
\begin{matrix} \boldmu^{\phantom{*}};\boldnu^{\phantom{*}}\\
           \boldmu^*;\boldnu^* \end{matrix} \right\},
\end{equation*}
where $\alpha=\alpha^{p,p'}_{a,b}$.

In addition, $\alpha^{p,p'}_{a,b'}=\alpha^{p,p'}_{a,b}-\Delta$,
$\beta^{p,p'}_{a,b',e,f'}=\beta^{p,p'}_{a,b,e,f}-\Delta$
and $\lfloor b'p/p'\rfloor=\lfloor bp/p'\rfloor$.
\end{lemma}

\Proof If $d^R>0$ then since $b$ is strictly between $\nu_{d^R}$ and
$\nu_{d^R}^*$, and $b'\ne\nu_{d^R}$ and $b'\ne\nu_{d^R}^*$, it follows
that $b'$ is strictly between $\nu_{d^R}$ and $\nu_{d^R}^*$.
Thus $\boldmu,\boldmu^*,\boldnu,\boldnu^{*}$ are a mazy-four
sandwiching $(a,b')$.

Let $h\in\P^{p,p'}_{a,b,e,f}(L+1,m)\left\{
\begin{smallmatrix} \boldmu^{\phantom{*}};\boldnu^{\phantom{*}}\\
           \boldmu^*;\boldnu^* \end{smallmatrix} \right\}$.
First consider $d^R>0$.
Since $\nu_{d^R}=b+\Delta$ and $\Delta=\pm1$, and since
$h$ attains $\nu^*_{d^R}$ after any $\nu_{d^R}$,
then necessarily $h_{L-1}=b-\Delta=b'$.
Then, on removing the final segment of $h$, we obtain a path $h'$
for which, via Lemma \ref{Atten2Lem},
$h'\in\P^{p,p'}_{a,b',e,f'}(L,m)\left\{
\begin{smallmatrix} \boldmu^{\phantom{*}};\boldnu^{\phantom{*}}\\
           \boldmu^*;\boldnu^* \end{smallmatrix} \right\}$.
Clearly, each $h'$ in this latter set arises this way.
If $d^R=0$, the same argument may be used after replacing
$\nu_{d^R}$ by 0 (resp.~$p'$) when $b=1$ (resp.~$p'-1$) and
necessarily $\Delta=-1$ (resp.~${}+1$).
The lemma then follows on using the expression for
$\mwt(h')$ given in Lemma \ref{Atten2Lem},
and $\alpha(h)=\alpha^{p,p'}_{a,b}$ from Lemma \ref{BetaConstLem}.

The final statement follows from Lemma \ref{Atten2Lem} and noting
that $\delta^{p,p'}_{b,f}=0$ and $b'=b+(-1)^f$ implies
$\lfloor b'p/p'\rfloor=\lfloor bp/p'\rfloor$.
\cqfd
\medskip


Removing the first segment of $h$ is permitted only when
$\pi(h)=0$ and when $d(h)=e$ so that the 0th vertex of $h$ is scoring.
We specify that $e(h')=1-e$ so that the pre-segment of $h'$ is in
the same direction as the first segment of $h$.
Let $a'=h_1$.
Then, on setting $e'=1-e$,
we obtain $h'\in\P^{p,p'}_{a',b,e',f}(L-1)$ with
$h'_{i}=h_{i+1}$ for $0\le i\le L-1$.

\begin{lemma}\label{Atten1Lem}
Let $1\le a,b<p'$ and $e,f\in\{0,1\}$.
For $L>0$, let $h\in\P^{p,p'}_{a,b,e,f}(L)$ be such that $d(h)=e$ and
$\pi(h)=0$,
and let $h'\in\P^{p,p'}_{a',b,e',f}(L')$ be obtained from $h$ by
the above process of path truncation.
If $\Delta=a-a'$ then $\Delta=(-1)^{e'}$, $e=1-e'$,
and
\begin{displaymath}
\begin{array}{l}
\bullet\quad L'=L-1;\\[1mm]
\bullet\quad m(h')=m(h);\\[1mm]
\bullet\quad \mwt(h')=\mwt(h) - \frac{1}{2}(L-m(h)+\Delta\beta(h)).
\end{array}
\end{displaymath}
In addition, $\alpha^{p,p'}_{a',b}=\alpha^{p,p'}_{a,b}+\Delta$
and $\beta^{p,p'}_{a',b,e',f}=\beta^{p,p'}_{a,b,e,f}+\Delta$.
\end{lemma}

\Proof That $\Delta=(-1)^{e'}$ and $e=1-e'$ and $L'=L-1$ are
immediate from the definition.
Let $h$ have striking sequence
$\left({a_1\atop b_1}\:{a_2\atop b_2}\:{a_3\atop b_3}\:
 {\cdots\atop\cdots}\:{a_l\atop b_l} \right)^{(e,f,d)}$.
If $a_1+b_1\ge2$ and the first vertex of $h$ is non-scoring, then
$h'$ has striking sequence
$\left({a_1-1\atop b_1}\:{a_2\atop b_2}\:{a_3\atop b_3}\:
 {\cdots\atop\cdots}\:{a_l\atop b_l} \right)^{(1-e,f,d)}$
and $\pi(h')=0$.
If $a_1+b_1\ge2$ and the first vertex of $h$ is scoring, then
$h'$ has striking sequence
$\left({a_1\atop b_1-1}\:{a_2\atop b_2}\:{a_3\atop b_3}\:
 {\cdots\atop\cdots}\:{a_l\atop b_l} \right)^{(1-e,f,d)}$
and $\pi(h')=1$.
If $a_1+b_1=1$, then 
$h'$ has striking sequence
$\left({a_2\atop b_2}\:{a_3\atop b_3}\:
 {\cdots\atop\cdots}\:{a_l\atop b_l} \right)^{(1-e,f,1-d)}$
and $\pi(h')=0$.
In each case,
Lemma \ref{WtHashLem} gives $\mwt(h')=\mwt(h)-(b_2+b_4+\cdots)$,
whereupon, since $\beta(h)=-\Delta((b_1+b_3+\cdots)-(b_2+b_4+\cdots))$,
we obtain $\mwt(h')=\mwt(h)-\frac12(L(h)-m(h)+\Delta\beta(h))$.
Also when $L>1$, we immediately see that $m(h')=m(h)$ in each case.
For $L=1$, it is readily verified that if $e\ne f$ then
$m(h)=m(h')=0$, and if $e=f$ then $m(h)=m(h')=1$.


That $\alpha^{p,p'}_{a',b}=\alpha^{p,p'}_{a,b}+\Delta$ is immediate.
Since $\pi(h)=0$, we have $\lfloor a'p/p'\rfloor=\lfloor ap/p'\rfloor$.
Then $\beta^{p,p'}_{a',b,e',f}=\lfloor bp/p'\rfloor-\lfloor a'p/p'\rfloor
+f-e'=\lfloor bp/p'\rfloor-\lfloor ap/p'\rfloor+f-e+(1-2e')=
\beta^{p,p'}_{a,b,e,f}+\Delta$, as required.
\cqfd
\medskip

\begin{lemma}\label{AttenGen1Lem}
Let $1\le a,a',b<p'$ and $e,f\in\{0,1\}$ with
$a-a'=-(-1)^e$ and $\delta^{p,p'}_{a,e}=0$.
Set $\Delta=a-a'$ and $e'=1-e$.
If ${d^L}>0$, let $\mu_{d^L}^*\ne a'$ and $\mu_{d^L}=a+\Delta$.
If ${d^L}=0$, restrict to $a\in\{1,p'-1\}$.
Then $\boldmu,\boldmu^*,\boldnu,\boldnu^{*}$ are a mazy-four
sandwiching $(a',b)$, and if $L\ge0$ then:
\begin{equation*}
\mchi^{p,p'}_{a',b,e',f}(L,m)\left\{
\begin{matrix} \boldmu^{\phantom{*}};\boldnu^{\phantom{*}}\\
           \boldmu^*;\boldnu^* \end{matrix} \right\}
=q^{-\frac{1}{2}(L+1-m+\Delta\beta)}
\mchi^{p,p'}_{a,b,e,f}(L+1,m)\left\{
\begin{matrix} \boldmu^{\phantom{*}};\boldnu^{\phantom{*}}\\
           \boldmu^*;\boldnu^* \end{matrix} \right\},
\end{equation*}
where $\beta=\beta^{p,p'}_{a,b,e,f}$.

In addition, $\alpha^{p,p'}_{a',b}=\alpha^{p,p'}_{a,b}+\Delta$,
$\beta^{p,p'}_{a',b,e',f}=\beta^{p,p'}_{a,b,e,f}+\Delta$
and $\lfloor a'p/p'\rfloor=\lfloor ap/p'\rfloor$.
\end{lemma}

\Proof
If $d^L>0$ then since $a$ is strictly between
$\mu_{d^L}$ and $\mu_{d^L}^*$, and $a'\ne\mu_{d^L}$ and $a'\ne\mu_{d^L}^*$,
it follows that $a'$ is strictly between $\mu_{d^L}$ and $\mu_{d^L}^*$.
Thus $\boldmu,\boldmu^*,\boldnu,\boldnu^{*}$ are a mazy-four
sandwiching $(a',b)$.

Let $h\in\P^{p,p'}_{a,b,e,f}(L+1,m)\left\{
\begin{smallmatrix} \boldmu^{\phantom{*}};\boldnu^{\phantom{*}}\\
           \boldmu^*;\boldnu^* \end{smallmatrix} \right\}$.
First consider $d^L>0$.
Since $\mu_{d^L}=a+\Delta$ and $\Delta=\pm1$, and since
$h$ attains $\mu^*_{d^L}$ before any $\mu_{d^L}$,
then necessarily $h_1=a-\Delta=a'$.
Then, on removing the first segment of $h$, we obtain a path $h'$
for which, via Lemma \ref{Atten1Lem},
$h'\in\P^{p,p'}_{a',b,e',f}(L,m)\left\{
\begin{smallmatrix} \boldmu^{\phantom{*}};\boldnu^{\phantom{*}}\\
           \boldmu^*;\boldnu^* \end{smallmatrix} \right\}$.
Clearly, each $h'$ in this latter set arises this way.
If $d^L=0$, the same argument may be used after replacing
$\mu_{d^L}$ by 0 (resp.~$p'$) when $a=1$ (resp.~$p'-1$) and
necessarily $\Delta=-1$ (resp.~${}+1$).
The lemma then follows on using the expression for
$\mwt(h')$ given in Lemma \ref{Atten1Lem},
and $\beta(h)=\beta^{p,p'}_{a,b,e,f}$ from Lemma \ref{BetaConstLem}.

The final statement follows from Lemma \ref{Atten1Lem}, and noting
that $\delta^{p,p'}_{a,e}=0$ and $a'=a+(-1)^e$ implies that
$\lfloor a'p/p'\rfloor=\lfloor ap/p'\rfloor$.
\cqfd
\medskip

\newpage

\setcounter{section}{6}

\section{Generating the fermionic expressions}\label{GenFermSec}

In this section, we prove the core result of this work: namely,
for certain vectors $\boldu^{L}$ and $\boldu^{R}$, we
precisely specify a subset of the set of paths $\P^{p,p'}_{a,b,c}(m_0)$
for which $F(\boldu^{L},\boldu^{R},m_0)$ is the generating function.
Throughout this section, we fix a pair $p,p'$ of coprime integers with
$1\le p<p'$, and employ the notation of Sections \ref{ContFSec},
\ref{TakSec}, \ref{LinSec}, \ref{QuadSec} and \ref{MNsysSec}.
We assume that $t>0$ so that (only) the case $(p,p')=(1,2)$ is
not considered hereafter.

\subsection{Runs and the core result}\label{CoreSec}

For $d\ge1$, we refer to a set
$\{\tau_j,\sigma_j,\Delta_j\}_{j=1}^{d}$
as a {\em run} if $\tau_1=t+1$,
each $\Delta_j=\pm1$,
\begin{equation}\label{Seq1Eq}
0\le
\sigma_{d}<\tau_{d}
<\sigma_{d-1}\le\tau_{d-1}
<\cdots
<\sigma_2\le\tau_2
<\sigma_1<t,
\end{equation}
and
\begin{equation}
\tau_j=\sigma_j\quad\implies\quad
\Delta_j=\Delta_{j+1}
\mbox{ and }
2\le j<d.
\end{equation}
For $0\le i\le t$, define $\eta(i)$ to be such that
$\tau_{\eta(i)+1}\le i<\tau_{\eta(i)}$,
where we set $\tau_{d+1}=0$.
In addition, we define
\begin{equation}\label{D0Eq}
d_0=
\left\{
  \begin{array}{ll}
  \displaystyle 0 \quad &\mbox{if }\sigma_1<t_n;\\[1mm]
  \displaystyle 1 \quad &\mbox{if }\sigma_1\ge t_n.
  \end{array}
\right.
\end{equation}

In what follows, superscripts $L$ or $R$
(to designate {\em left} or {\em right}) may be appended to the
quantities defined above.

Given a run
$\run^L=\{\tau^{L}_j,\sigma^{L}_j,\Delta^{L}_j\}_{j=1}^{d^{L}}$,
we define values $a$, $e$, $t$-dimensional vectors $\boldu^L$,
$\boldDelta^L$,
and $(d^L-d^L_0)$-dimensional vectors $\boldmu,\boldmu^*$.

Define:
\begin{equation*}
a=\sum_{m=2}^{d^{L}}
        \Delta^{L}_m (\kappa_{\tau_m^{L}}-\kappa_{\sigma_m^{L}})
+\left\{
  \begin{array}{ll}
  \displaystyle
  \kappa_{\sigma^{L}_1} \quad &\mbox{if }\Delta^{L}_1=-1;\\[1mm]
  \displaystyle
  p'-\kappa_{\sigma^{L}_1}\quad &\mbox{if }\Delta^{L}_1=+1.
  \end{array}
\right.
\end{equation*}

Define $e\in\{0,1\}$ to be such that:
\begin{equation*}
\Delta^{L}_{d^{L}}=
\left\{
  \begin{array}{ll}
  \displaystyle
  -(-1)^e\quad &\mbox{if }\sigma^{L}_{d^{L}}=0;\\[1mm]
  \displaystyle
  (-1)^{e+k^{L}}\quad &\mbox{if }\sigma^{L}_{d^{L}}>0,
  \end{array}
\right.
\end{equation*}
where $k^{L}=\zeta(\sigma^{L}_{d^{L}})$.


Define $t$-dimensional vectors $\boldu^L$ and $\boldDelta^L$ as follows
(c.f.\ (\ref{uEq}) and (\ref{DeltaEq})):
\begin{align*}
\boldu^{L}&=
\sum_{m=1}^{d^{L}} \boldu_{\sigma^{L}_m,\tau^{L}_m}
&&\hskip-25mm+\left\{
  \begin{array}{ll}
  \displaystyle
  0\quad
  &\mbox{if }\Delta^{L}_1=-1;\\[1mm]
  \displaystyle
  \bolde_t\quad
  &\mbox{if }\Delta^{L}_1=+1;
  \end{array}
\right.\\[1mm]
\boldDelta^{L}&=
\sum_{m=1}^{d^{L}}
\Delta^{L}_m\boldu_{\sigma^{L}_m,\tau^{L}_m}
&&\hskip-25mm-\left\{
  \begin{array}{ll}
  \displaystyle
  0\quad
  &\mbox{if }\Delta^{L}_1=-1;\\[1mm]
  \displaystyle
  \bolde_t\quad
  &\mbox{if }\Delta^{L}_1=+1.
  \end{array}
\right.
\end{align*}
Note that $\Delta^L_m$ are not the components of $\boldDelta^L$.

Define values $\{\mu_j,\mu^*_j\}_{j=0}^{d^L-1}$ as follows:
\begin{align*}
\mu^*_0&=
\left\{
  \begin{array}{ll}
  \displaystyle
  \kappa_{t_n} \quad &\mbox{if }\Delta^{L}_1=-1;\\[1mm]
  \displaystyle
  p'-\kappa_{t_n}\quad &\mbox{if }\Delta^{L}_1=+1;
  \end{array}
\right.\\[1mm]
\mu_0&=
\left\{
  \begin{array}{ll}
  \displaystyle
  0 \quad &\mbox{if }\Delta^{L}_1=-1;\\[1mm]
  \displaystyle
  p'\phantom{\displaystyle{}-\kappa_{t_n}}\quad &\mbox{if }\Delta^{L}_1=+1;
  \end{array}
\right.\displaybreak[0]\\[1mm]
\mu^*_j&=
 \sum_{m=2}^{j}
        \Delta^{L}_m (\kappa_{\tau_m^{L}}-\kappa_{\sigma_m^{L}})
+\left\{
  \begin{array}{ll}
  \displaystyle
  \kappa_{\sigma^{L}_1} \quad &\mbox{if }\Delta^{L}_1=-1;\\[1mm]
  \displaystyle
  p'-\kappa_{\sigma^{L}_1}\quad &\mbox{if }\Delta^{L}_1=+1;
  \end{array}
\right.\\[1mm]
\mu_j&=
\mu^*_j+\Delta^{L}_{j+1}\kappa_{\tau^{L}_{j+1}},
\end{align*}
for $1\le j<d^L$.
The $(d^L-d^L_0)$-dimensional vectors $\boldmu$ and $\boldmu^*$
are defined by $\boldmu=(\mu_{d^L_0},\ldots,\mu_{d^L-1})$
and $\boldmu^*=(\mu^*_{d^L_0},\ldots,\mu^*_{d^L-1})$.

In a totally analogous way, a run
$\run^R=\{\tau^{R}_j,\sigma^{R}_j,\Delta^{R}_j\}_{j=1}^{d^{R}}$,
is used to define values $b$, $f$, $\{\nu_j,\nu^*_j\}_{j=0}^{d^R-1}$,
$t$-dimensional vectors $\boldu^R$, $\boldDelta^R$,
and $(d^R-d^R_0)$-dimensional vectors $\boldnu,\boldnu^*$.

In view of (\ref{Seq1Eq}), we see that $\boldmu,\boldmu^*,\boldnu,\boldnu^*$
satisfy the three criteria for being a
mazy-four in the $(p,p')$-model sandwiching $(a,b)$.
That they are in fact interfacial will be established later
Their neighbouring odd bands will be seen to be specified by the values
$\{\tilde\mu_j,\tilde\mu^*_j\}_{j=0}^{d^L-1}$ and
$\{\tilde\nu_j,\tilde\nu^*_j\}_{j=0}^{d^R-1}$, where we define:
\begin{align*}
\tilde\mu^*_0&=
\left\{
  \begin{array}{ll}
  \displaystyle
  \tkappa_{t_n} \quad &\mbox{if }\Delta^{L}_1=-1;\\[1mm]
  \displaystyle
  p-\tkappa_{t_n}\quad &\mbox{if }\Delta^{L}_1=+1;
  \end{array}
\right.\\[1mm]
\tilde\mu_0&=
\left\{
  \begin{array}{ll}
  \displaystyle
  0 \quad &\mbox{if }\Delta^{L}_1=-1;\\[1mm]
  \displaystyle
  p\phantom{\displaystyle{}-\tkappa_{t_n}}\quad &\mbox{if }\Delta^{L}_1=+1;
  \end{array}
\right.\\[1mm]
\tilde\mu^*_j&=
 \sum_{m=2}^{j}
        \Delta^{L}_m (\tkappa_{\tau_m^{L}}-\tkappa_{\sigma_m^{L}})
+\left\{
  \begin{array}{ll}
  \displaystyle
  \tkappa_{\sigma^{L}_1} \quad &\mbox{if }\Delta^{L}_1=-1;\\[1mm]
  \displaystyle
  p-\tkappa_{\sigma^{L}_1}\quad &\mbox{if }\Delta^{L}_1=+1;
  \end{array}
\right.\\[1mm]
\tilde\mu_j&=
\tilde\mu^*_j+\Delta^{L}_{j+1}\tkappa_{\tau^{L}_{j+1}},
\end{align*}
for $1\le j<d^L$, with
$\{\tilde\nu_j,\tilde\nu^*_j\}_{j=0}^{d^R-1}$ defined analogously.

Finally, note that $\boldDelta^L=\boldDelta(\run^L)$
and $\boldDelta^R=\boldDelta(\run^R)$, and use the procedure given
in Section \ref{ConSec} to calculate $\gamma'$ and well as
$\{\alpha_j,\alpha''_j,\beta_j,\beta_j'',\gamma_j,\gamma_j',\gamma_j''
 \}_{j=0}^{t}$.

A run $\{\tau_j,\sigma_j,\Delta_j\}_{j=1}^{d}$
is said to be {\em naive} when
if either $1<j<d$, or $1=j<d$ and $\sigma_1\le t_n$, then
\begin{displaymath}
\begin{array}{lll}
\Delta_{j+1}=\Delta_j\quad
&\implies&\quad\tau_{j+1}=t_{\zeta(\sigma_j-1)},\\
\Delta_{j+1}\ne\Delta_j\quad
&\implies&\quad\tau_{j+1}=t_{\zeta(\sigma_j)};
\end{array}
\end{displaymath}
and if $d>1$ and $\sigma_1>t_n$ then
\begin{displaymath}
\begin{array}{lll}
\Delta_2=\Delta_1\quad
&\implies&\quad\tau_{2}=t_{n-1},\\
\Delta_2\ne\Delta_1\mbox{ and }c_{n-1}>1\quad
&\implies&\quad\tau_{2}=t_{n}-1,\\
\Delta_2\ne\Delta_1\mbox{ and }c_{n-1}=1\quad
&\implies&\quad\tau_{2}=t_{n-2}.
\end{array}
\end{displaymath}

\noindent
Note that if $\{\tau_j,\sigma_j,\Delta_j\}_{j=1}^{d}$
is a naive run then
for $1\le j<d$, there exists (at least one value) $k$ such
that $\tau_{j+1}\le t_k<\sigma_j$.

Later, in Section \ref{CollateSec}, we show that the values
$\{\tau_j,\sigma_j,\Delta_j\}_{j=1}^{d}$ that are obtained from each
leaf node of a Takahashi tree as described in Section \ref{FormVSec},
provide a naive run.

The core result is the following:

\begin{theorem}\label{CoreThrm}
Let $p'>2p$, and
let both $\run^L=\{\tau^{L}_j,\sigma^{L}_j,\Delta^{L}_j\}_{j=1}^{d^{L}}$ and
$\run^R=\{\tau^{R}_j,\sigma^{R}_j,\Delta^{R}_j\}_{j=1}^{d^{R}}$
be naive runs.
If $m_0\equiv b-a\;(\mod2)$ then:
\begin{displaymath}
\mchi^{p,p'}_{a,b,e,f}(m_0)
\left\{
\begin{matrix} \boldmu^{\phantom{*}};\boldnu^{\phantom{*}}\\
                    \boldmu^*;\boldnu^* \end{matrix} \right\}
=q^{\frac14(\gamma_0-\gamma')}
F(\boldu^{L},\boldu^{R},m_0).
\end{displaymath}
\end{theorem}

We prove this result in the following sections.
The strategy for the proof is simple: we begin with
the trivial generating functions for the $(1,3)$-model stated in
Lemma \ref{SeedLem},\footnote{It is actually possible to start
with $\mchi^{1,2}_{1,1,0,0}(0,0)$: an application of a single
$\B$-transform yields the results of Lemma \ref{SeedLem}.
We don't do this in order to avoid a further increase in
the notational complexities in the main induction proof.}
and for $i=t-1,t-2,\ldots,1$, apply a $\B$-transform if $i\ne t_k$ for
all $k$, and apply a $\B\D$-transform if $i=t_k$ for some $k$.
If $(\boldu^{L})_i=-1$ (resp.~$+1$), we follow the transform
with path extension (resp.~truncation) on the left.
If $(\boldu^{R})_i=-1$ (resp.~$+1$), we follow
with path extension (resp.~truncation) on the right.
This process is carried out in Lemma \ref{CoreIndLem} to yield,
in Corollary \ref{CoreCor}, an expression for
$\displaystyle
\mchi^{p,p'}_{a,b,e,f}(m_0,m_1)
\left\{
\begin{smallmatrix} \boldmu^{\phantom{*}};\boldnu^{\phantom{*}}\\
                    \boldmu^*;\boldnu^* \end{smallmatrix} \right\},
$
for certain values of $m_1$.
Then, at the end of Section \ref{ProofIndSec}, a final sum over $m_1$
is performed to prove Theorem \ref{CoreThrm}.

In Section \ref{TransSec}, we transfer our result across to generating
functions in terms of the original weighting function of
(\ref{WtDef2}).
In this way, for a specific $c$, we will have identified a subset
of $\P^{p,p'}_{a,b,c}(L)$ for which the generating function is
precisely $F(\boldu^{L},\boldu^{R},L)$.

We actually first prove an analogue of
Theorem \ref{CoreThrm} involving {\em reduced} runs instead of
naive runs, where a run
$\{\tau_j,\sigma_j,\Delta_j\}_{j=1}^{d}$
is said to be reduced if:
\begin{equation}\label{Seq2Eq}
0\le
\sigma_{d}<\tau_{d}
<\sigma_{d-1}<\tau_{d-1}
<\cdots
<\sigma_2<\tau_2
<\sigma_1<t,
\end{equation}
and if $d>1$ then
\begin{displaymath}
\begin{array}{l}
1.\quad \sigma_1>t_n \mbox{ and } \Delta_2=\Delta_1
      \quad
      \implies\quad
      \tau_2=t_{k} \mbox{ for }  1\le k\le n-1;\\[0.8mm]
2.\quad \sigma_1>t_n \mbox{ and } \Delta_2\ne\Delta_1
                       \mbox{ and } c_{n-1}>1 \\
 \hskip35mm
      \implies\quad
      \tau_2=t_n-1 \mbox{ or } \tau_2=t_{k}
                      \mbox{ for }  1\le k\le n-1;\\[0.8mm]
3.\quad \sigma_1>t_n \mbox{ and } \Delta_2\ne\Delta_1
                       \mbox{ and } c_{n-1}=1 \\
 \hskip35mm
      \implies\quad
      \tau_2=t_{k} \mbox{ for }  1\le k\le n-2;\\[0.8mm]
4.\quad \sigma_1\le t_n \mbox{ and } \Delta_2=\Delta_1
      \quad
      \implies\quad
      \tau_2=t_{k} \mbox{ for }  1\le k\le \zeta(\sigma_1-1);\\[0.8mm]
5.\quad \sigma_1\le t_n \mbox{ and } \Delta_2\ne\Delta_1
      \quad
      \implies\quad
      \tau_2=t_{k} \mbox{ for }  1\le k\le \zeta(\sigma_1),\\[0.8mm]
\end{array}
\end{displaymath}
and if $1<j<d$ then
\begin{displaymath}
\begin{array}{l}
6.\quad \Delta_{j+1}=\Delta_j
      \quad
      \implies\quad
      \tau_{j+1}=t_{k} \mbox{ for }  1\le k\le \zeta(\sigma_j-1);\\[0.8mm]
7.\quad \Delta_{j+1}\ne\Delta_j
      \quad
      \implies\quad
      \tau_{j+1}=t_{k} \mbox{ for }  1\le k\le \zeta(\sigma_j).
\end{array}
\end{displaymath}
 
\noindent
The following result is immediately obtained:

\begin{lemma}\label{RunLem}
Let $\{\tau_j,\sigma_j,\Delta_j\}_{j=1}^{d}$
be a run, and let the run
$\{\tau^{\prime}_j,\sigma^{\prime}_j,
                       \Delta^{\prime}_j\}_{j=1}^{d^{\prime}}$
be obtained from it by removing each triple
$\{\tau_j,\sigma_j,\Delta_j\}$ whenever
$\tau_j=\sigma_j$. Naturally $d^{\prime}\le d$.
If $\{\tau_j,\sigma_j,\Delta_j\}_{j=1}^{d}$
is a naive run then
$\{\tau^{\prime}_j,\sigma^{\prime}_j,
                       \Delta^{\prime}_j\}_{j=1}^{d^{\prime}}$
is a reduced run.
\end{lemma}

\Proof A routine verification.
\cqfd
\medskip

\noindent
The reduced run
$\{\tau^{\prime}_j,\sigma^{\prime}_j,
                        \Delta^{\prime}_j\}_{j=1}^{d^{\prime}}$
identified in the above lemma will be referred to as the
reduced run {\em corresponding} to the
run $\{\tau_j,\sigma_j,\Delta_j\}_{j=1}^{d}$.

\begin{lemma}\label{BigHashLem}
Let $\{\tau_j,\sigma_j,\Delta_j\}_{j=1}^{d}$
be a naive run and
$\{\tau^{\prime}_j,\sigma^{\prime}_j,
                        \Delta^{\prime}_j\}_{j=1}^{d^{\prime}}$
the corresponding reduced run.
Let $d_0<j<d$. If $\tau^\prime_{j+1}=t_k$ for some $k$ and
$c_k=1$ then:
\begin{displaymath}
\begin{array}{l}
1.\quad \sigma^\prime_j=t_{k+1}\phantom{{}+1}\quad\implies\quad
           \Delta^\prime_{j+1}\ne\Delta^\prime_j;\\
2.\quad \sigma^\prime_j=t_{k+1}+1 \quad\implies\quad
           \Delta^\prime_{j+1}=\Delta^\prime_j;\\
3.\quad \sigma^\prime_j>t_{k+1}+1 \quad\implies\quad
           \sigma^\prime_j>t_{k+2}.
\end{array}
\end{displaymath}
\end{lemma}

\Proof
We have $t_{k+1}=t_k+1$ and
$\tau^\prime_{j+1}<\sigma^\prime_j<\tau^\prime_j$.
Let $J$ and $J'$ be such that $\tau_{J+1}=\tau^\prime_{j+1}$ and
$\sigma_{J'}=\sigma^\prime_j$.
Then $J'\le J$ and $\sigma_i=\tau_i$ for $J'<i\le J$.
Moreover, $\Delta_{J'+1}=\Delta_{J'+2}=\cdots=\Delta_{J}=\Delta_{J+1}$.
First note that
$\sigma_{J}=\tau_{J}=t_{k+1}$ is forbidden here since
$\Delta_{J}=\Delta_{J+1}$ would then dictate that
$\tau_{J+1}=t_{\zeta(\sigma_J-1)}<t_k$ thereby contradicting
the specified value of $\tau_{J+1}=\tau^\prime_{j+1}=t_k$.

Now consider the case $\sigma^\prime_j=t_{k+1}$.
Necessarily $J'=J$. Then $\Delta^\prime_{j+1}=\Delta^\prime_j$
is excluded since this would also imply that
$\tau_{J+1}=t_{\zeta(\sigma_J-1)}<t_k$.

Now consider the case $\sigma^\prime_j=t_{k+1}+1$.
Since $\sigma_{J}=\tau_{J}=t_{k+1}$ is forbidden then we again
have $J'=J$. Here we exclude $\Delta^\prime_{j+1}\ne\Delta^\prime_j$
since this would imply that $\Delta_{J+1}\ne\Delta_J$ and
then $\tau_{J+1}=t_{\zeta(\sigma_J)}=t_{\zeta(\sigma^\prime_j)}=t_{k+1}$.

Finally consider the case $\sigma^\prime_j>t_{k+1}+1$.
Here $J'<J$ because if $J'=J$ then
$\tau_{J+1}=t_{\zeta(\sigma^\prime_j)}\ge t_{k+1}$
or $\tau_{J+1}=t_{\zeta(\sigma^\prime_j-1)}\ge t_{k+1}$.
Therefore $\sigma_{J}=\tau_{J}$ and $\Delta_{J}=\Delta_{J+1}$.
If $\sigma_J>t_{k+1}+1$ then
$\tau^\prime_{j+1}=\tau_{J+1}=t_{\zeta(\sigma_J-1)}\ge t_{k+1}$.
Therefore, since $\sigma_J=t_{k+1}$ is forbidden, we conclude that
$\sigma_J=t_{k+1}+1$.
So $\tau_J=t_{k+1}+1$.
But either $\tau_J=t_{k'}$ for some $k'$, or
$\tau_J=t_n-1$, $J=2$ and $\sigma_1>t_n$.
In the former case, necessarily $k'=k+2$ (and thus $c_{k+1}=1$),
whereupon $\sigma^\prime_j=\sigma_{J'}>\sigma_J=\tau_J=t_{k+2}$,
as required.
In the latter case, necessarily $n\ge k+2$ and $J'=j=1$.
Then $\sigma^\prime_1=\sigma_1>t_n\ge t_{k+2}$ as required.
\cqfd
\medskip

Now note that the quantities $a$, $e$, $\boldu^L$, and $\boldDelta^L$
(resp.\ $b$, $f$ ,$\boldu^R$ and $\boldDelta^R$) obtained from
a run $\{\tau^{L}_j,\sigma^{L}_j,\Delta^{L}_j\}_{j=1}^{d^{L}}$
(resp.\ $\{\tau^{R}_j,\sigma^{R}_j,\Delta^{R}_j\}_{j=1}^{d^{R}}$)
are equal to those obtained from the
corresponding reduced run
$\{\tau^{L\prime}_j,\sigma^{L\prime}_j,
                        \Delta^{L\prime}_j\}_{j=1}^{d^{L\prime}}$
(resp.\ $\{\tau^{R\prime}_j,\sigma^{R\prime}_j,
                        \Delta^{R\prime}_j\}_{j=1}^{d^{R\prime}}$).

The following lemma shows that the generating function is
unchanged on replacing naive runs with their corresponding
reduced runs.

\begin{lemma}\label{ReducedLem}
Let $\{\tau^{L}_j,\sigma^{L}_j,\Delta^{L}_j\}_{j=1}^{d^{L}}$
and $\{\tau^{R}_j,\sigma^{R}_j,\Delta^{R}_j\}_{j=1}^{d^{R}}$
be naive runs, and let
$\{\mu_j, \mu^*_j\}_{j=0}^{d^L-1}$ and 
$\{\nu_j, \nu^*_j\}_{j=0}^{d^R-1}$ be obtained from these as above.

Then let
$\{\tau^{L\prime}_j,\sigma^{L\prime}_j,
                        \Delta^{L\prime}_j\}_{j=1}^{d^{L\prime}}$
and $\{\tau^{R\prime}_j,\sigma^{R\prime}_j,
                        \Delta^{R\prime}_j\}_{j=1}^{d^{R\prime}}$
be the corresponding reduced runs, and let
$\{\mu^{\prime}_j, \mu^{*\prime}_j\}_{j=0}^{d^{L\prime}-1}$ and 
$\{\nu^{\prime}_j, \nu^{*\prime}_j\}_{j=0}^{d^{R\prime}-1}$ be
obtained from these as above.
Then:
\begin{equation*}
\begin{split}
\mchi^{p,p'}_{a,b,e,f}(m_0)
&\left\{
{\mu_{d^L_0},\atop \mu^*_{d^L_0},}
{\mu_1,\atop \mu^*_1,}
{\ldots,\atop\ldots,}
{\mu_{d^L-1};\atop \mu^*_{d^L-1};}
{\nu_{d^R_0},\atop \nu^*_{d^R_0},}
{\nu_1,\atop \nu^*_1,}
{\ldots,\atop\ldots,}
{\nu_{d^R-1}\atop \nu^*_{d^R-1}}
\right\}\\
&\quad=
\mchi^{p,p'}_{a,b,e,f}(m_0)
\left\{
{\mu^{\prime}_{d^L_0},\atop \mu^{\prime*}_{d^L_0},}
{\mu^{\prime}_1,\atop \mu^{\prime*}_1,}
{\ldots,\atop\ldots,}
{\mu^{\prime}_{d^{L\prime}-1};\atop \mu^{\prime*}_{d^{L\prime}-1};}
{\nu^{\prime}_{d^R_0},\atop \nu^{\prime*}_{d^R_0},}
{\nu^{\prime}_1,\atop \nu^{\prime*}_1,}
{\ldots,\atop\ldots,}
{\nu^{\prime}_{d^{R\prime}-1}\atop \nu^{\prime*}_{d^{R\prime}-1}}
\right\}.
\end{split}
\end{equation*}
\end{lemma}

\Proof
Let $j$ be the smallest value for which $\tau^L_j=\sigma^L_j$.
If such a value exists, then necessarily
$2\le j<d^L$ and $\Delta^L_j=\Delta^L_{j+1}$.
{}From the definitions, we see that
$\mu^*_j=\mu^*_{j-1}$.
Thus, by Lemma \ref{CutParam2Lem},
\begin{equation*}
\begin{split}
\mchi^{p,p'}_{a,b,e,f}(m_0)
&\left\{
{\ldots,\atop\ldots,}
{\mu_{j-2},\atop \mu^*_{j-2},}
{\mu_{j-1},\atop \mu^*_{j-1},}
{\mu_{j},\atop \mu^*_{j},}
{\ldots;\atop\ldots;}
{\nu_{d^R_0},\atop \nu^*_{d^R_0},}
{\nu_1,\atop \nu^*_1,}
{\ldots,\atop\ldots,}
{\nu_{d^R-1}\atop \nu^*_{d^R-1}}
\right\}\\
&\quad=
\mchi^{p,p'}_{a,b,e,f}(m_0)
\left\{
{\ldots,\atop\ldots,}
{\mu_{j-2},\atop \mu^*_{j-2},}
{\mu_{j},\atop \mu^*_{j},}
{\ldots;\atop\ldots;}
{\nu_{d^R_0},\atop \nu^*_{d^R_0},}
{\nu_1,\atop \nu^*_1,}
{\ldots,\atop\ldots,}
{\nu_{d^R-1}\atop \nu^*_{d^R-1}}
\right\}.
\end{split}
\end{equation*}
By recursively applying this process, and doing likewise for any $j$
for which $\tau^R_j=\sigma^R_j$, we obtain the required result.
\cqfd
\medskip

Thus, to prove Theorem \ref{CoreThrm} for a pair of arbitrary naive runs,
it is only necessary to prove it for the corresponding reduced runs.
This is what will be done.

\subsection{Induction parameters}\label{IndParamSec}

In this and the following Section \ref{ProofIndSec},
we assume that $p'>2p$.
We also fix a pair 
$\{\tau^{L}_j,\sigma^{L}_j,\Delta^{L}_j\}_{j=1}^{d^{L}}$,
$\{\tau^{R}_j,\sigma^{R}_j,\Delta^{R}_j\}_{j=1}^{d^{R}}$
of reduced runs.
In our main induction (Lemma \ref{CoreIndLem}), we will for
$i=t-1,t-2,\ldots,0$, step through a sequence of
$(p^{(i)},p^{(i)\prime})$-models, constructing the generating function
for a certain set of paths at each stage.
In this section, we define various sets of parameters that pertain
to those sets of paths, and identify some basic relationships
between them.

Firstly, for $0\le i<t$,
let $k(i)$ be such that $t_{k(i)}\le i<t_{k(i)+1}$ (i.e.\ $k(i)=\zeta(i+1)$),
and define $p^{(i)}$ and $p^{(i)\prime}$ to be the positive coprime
integers for which $p^{(i)\prime}/p^{(i)}$ has continued
fraction $[t_{k(i)+1}+1-i,c_{k(i)+1},\ldots,c_n]$.
Thus $p^{(i)\prime}/p^{(i)}$ has rank $t^{(i)}$ where we set $t^{(i)}=t-i$.
As in Section \ref{TakSec}, we obtain
Takahashi lengths $\{\kappa^{(i)}_j\}_{j=0}^{t^{(i)}}$ and
truncated Takahashi lengths $\{\tilde\kappa^{(i)}_j\}_{j=0}^{t^{(i)}}$
for $p^{(i)\prime}/p^{(i)}$.

\begin{lemma}\label{IndTakLem} Let $1\le i<t$. If $i\ne t_{k(i)}$
then:
\begin{displaymath}
\begin{array}{lll}
p^{(i-1)\prime}
 &\hskip-2.5mm=p^{(i)\prime}+p^{(i)};\\[0.5mm]
p^{(i-1)}
 &\hskip-2.5mm=p^{(i)};\\[0.5mm]
\kappa^{(i-1)}_j
 &\hskip-2.5mm=\kappa^{(i)}_{j-1}+\tkappa^{(i)}_{j-1}\qquad
 &(1\le j\le t^{(i-1)});\\[0.5mm]
\tkappa^{(i-1)}_j
 &\hskip-2.5mm=\tkappa^{(i)}_{j-1}\qquad
 &(1\le j\le t^{(i-1)}).
\end{array}
\end{displaymath}
If $i=t_{k(i)}$ then:
\begin{displaymath}
\begin{array}{lll}
p^{(i-1)\prime}
 &\hskip-2.5mm=2p^{(i)\prime}-p^{(i)};\\[0.5mm]
p^{(i-1)}
 &\hskip-2.5mm=p^{(i)\prime}-p^{(i)};\\[0.5mm]
\kappa^{(i-1)}_j
 &\hskip-2.5mm=2\kappa^{(i)}_{j-1}-\tkappa^{(i)}_{j-1}\qquad
 &(2\le j\le t^{(i-1)});\\[0.5mm]
\tkappa^{(i-1)}_j
 &\hskip-2.5mm=\tkappa^{(i)}_{j-1}-\tkappa^{(i)}_{j-1}\qquad
 &(2\le j\le t^{(i-1)}).
\end{array}
\end{displaymath}
\end{lemma}

\Proof If $i\ne t_{k(i)}$ then $k(i-1)=k(i)$.
Then $p^{(i)\prime}/p^{(i)}$ and
$p^{(i-1)\prime}/p^{(i-1)}$ have continued fractions
$[t_{k(i)}+1-i,c_{k(i)+1},\ldots,c_n]$ and
$[t_{k(i)}+2-i,c_{k(i)+1},\ldots,c_n]$ respectively.
That $p^{(i-1)\prime}=p^{(i)\prime}+p^{(i)}$ and $p^{(i-1)}=p^{(i)}$
follows immediately.
The expressions for $\kappa^{(i-1)}_j$ and $\tkappa^{(i-1)}_j$
then follow from Lemma \ref{BmodelLem}.

If $i=t_{k(i)}$ then $k(i-1)=k(i)-1$.
Then $p^{(i)\prime}/p^{(i)}$ and
$p^{(i-1)\prime}/p^{(i-1)}$ have continued fractions
$[c_{k(i)}+1,c_{k(i)+1},\ldots,c_n]$ and
$[2,c_{k(i)},c_{k(i)+1},\ldots,c_n]$ respectively.
That $p^{(i-1)\prime}=2p^{(i)\prime}-p^{(i)}$
and $p^{(i-1)}=p^{(i)\prime}-p^{(i)}$
follows immediately.
The expressions for $\kappa^{(i-1)}_j$ and $\tkappa^{(i-1)}_j$
then follow from combining Lemma \ref{DmodelLem} with Lemma \ref{BmodelLem}.
\cqfd
\medskip


For $0\le i< t$, we now define $e^{(i)},f^{(i)}\in\{0,1\}$.
With $k^{L}(i)$ such that
$t_{k^{L}(i)}<\sigma^{L}_{\eta^{L}(i)}\le t_{k^{L}(i)+1}$,
define $e^{(i)}$ such that:
\begin{equation}\label{DefEiEq}
\Delta^{L}_{\eta^{L}(i)}=
\left\{
  \begin{array}{ll}
  -(-1)^{e^{(i)}}\quad& \mbox{if } i\ge\sigma^{L}_{\eta^{L}(i)};\\
  (-1)^{e^{(i)}+k(i)-k^{L}(i)}\quad& \mbox{if } i<\sigma^{L}_{\eta^{L}(i)},
  \end{array}
\right.
\end{equation}
and with $k^{R}(i)$ such that
$t_{k^{R}(i)}<\sigma^{R}_{\eta^{R}(i)}\le t_{k^{R}(i)+1}$,
define $f^{(i)}$ such that:
\begin{equation}\label{DefFiEq}
\Delta^{R}_{\eta^{R}(i)}=
\left\{
  \begin{array}{ll}
  -(-1)^{f^{(i)}}\quad& \mbox{if } i\ge\sigma^{R}_{\eta^{R}(i)};\\
  (-1)^{f^{(i)}+k(i)-k^{R}(i)}\quad& \mbox{if } i<\sigma^{R}_{\eta^{R}(i)}.
  \end{array}
\right.
\end{equation}


We now define what will be the starting and ending points of the
paths that we consider at the $i$th induction step.
For $0\le i<t$, define:
\begin{equation*}
a^{(i)}=
\left\{
  \begin{split}
   &-\Delta^{L}_{\eta^{L}(i)}
      +\sum_{m=2}^{\eta^{L}(i)}
                          \Delta^{L}_m \kappa^{(i)}_{\tau_m^{L}-i}
      -\sum_{m=1}^{\eta^{L}(i)-1}
                          \Delta^{L}_m \kappa^{(i)}_{\sigma_m^{L}-i}\\
   &\hskip65mm
        \mbox{if }\Delta^{L}_1=-1\mbox{ and } i\ge\sigma^{L}_{\eta^{L}(i)};\\
   &   \sum_{m=2}^{\eta^{L}(i)}
                          \Delta^{L}_m \kappa^{(i)}_{\tau_m^{L}-i}
      -\sum_{m=1}^{\eta^{L}(i)}
                          \Delta^{L}_m \kappa^{(i)}_{\sigma_m^{L}-i}\\
   &\hskip65mm
        \mbox{if }\Delta^{L}_1=-1\mbox{ and } i\le\sigma^{L}_{\eta^{L}(i)};\\
   &p^{(i)\prime} -\Delta^{L}_{\eta^{L}(i)}
      +\sum_{m=2}^{\eta^{L}(i)}
                          \Delta^{L}_m \kappa^{(i)}_{\tau_m^{L}-i}
      -\sum_{m=1}^{\eta^{L}(i)-1}
                          \Delta^{L}_m \kappa^{(i)}_{\sigma_m^{L}-i}\\
   &\hskip65mm
        \mbox{if }\Delta^{L}_1=+1\mbox{ and } i\ge\sigma^{L}_{\eta^{L}(i)};\\
   &p^{(i)\prime}
      +\sum_{m=2}^{\eta^{L}(i)}
                          \Delta^{L}_m \kappa^{(i)}_{\tau_m^{L}-i}
      -\sum_{m=1}^{\eta^{L}(i)}
                          \Delta^{L}_m \kappa^{(i)}_{\sigma_m^{L}-i}\\
   &\hskip65mm
        \mbox{if }\Delta^{L}_1=+1\mbox{ and } i\le\sigma^{L}_{\eta^{L}(i)}.
  \end{split}
\right.
\end{equation*}
Similarly, for $0\le i<t$ define:
\begin{equation*}
b^{(i)}=
\left\{
  \begin{split}
   &-\Delta^{R}_{\eta^{R}(i)}
      +\sum_{m=2}^{\eta^{R}(i)}
                          \Delta^{R}_m \kappa^{(i)}_{\tau_m^{R}-i}
      -\sum_{m=1}^{\eta^{R}(i)-1}
                          \Delta^{R}_m \kappa^{(i)}_{\sigma_m^{R}-i}\\
   &\hskip65mm
        \mbox{if }\Delta^{R}_1=-1\mbox{ and } i\ge\sigma^{R}_{\eta^{R}(i)};\\
   &   \sum_{m=2}^{\eta^{R}(i)}
                          \Delta^{R}_m \kappa^{(i)}_{\tau_m^{R}-i}
      -\sum_{m=1}^{\eta^{R}(i)}
                          \Delta^{R}_m \kappa^{(i)}_{\sigma_m^{R}-i}\\
   &\hskip65mm
        \mbox{if }\Delta^{R}_1=-1\mbox{ and } i\le\sigma^{R}_{\eta^{R}(i)};\\
   &p^{(i)\prime} -\Delta^{R}_{\eta^{R}(i)}
      +\sum_{m=2}^{\eta^{R}(i)}
                          \Delta^{R}_m \kappa^{(i)}_{\tau_m^{R}-i}
      -\sum_{m=1}^{\eta^{R}(i)-1}
                          \Delta^{R}_m \kappa^{(i)}_{\sigma_m^{R}-i}\\
   &\hskip65mm
        \mbox{if }\Delta^{R}_1=+1\mbox{ and } i\ge\sigma^{R}_{\eta^{R}(i)};\\
   &p^{(i)\prime}
      +\sum_{m=2}^{\eta^{R}(i)}
                          \Delta^{R}_m \kappa^{(i)}_{\tau_m^{R}-i}
      -\sum_{m=1}^{\eta^{R}(i)}
                          \Delta^{R}_m \kappa^{(i)}_{\sigma_m^{R}-i}\\
   &\hskip65mm
        \mbox{if }\Delta^{R}_1=+1\mbox{ and } i\le\sigma^{R}_{\eta^{R}(i)}.
  \end{split}
\right.
\end{equation*}


\noindent
We now define what will be the starting and ending points of certain
intermediate paths that we consider at the $i$th induction step.
For $1\le i<t$, define:
\begin{equation*}
a^{(i)\prime}=
\left\{
  \begin{split}
   &-\Delta^{L}_{\eta^{L}(i)}(1+\delta_{i,t_{k(i)}})
      +\sum_{m=2}^{\eta^{L}(i)}
                          \Delta^{L}_m \kappa^{(i-1)}_{\tau_m^{L}-i+1}
      -\sum_{m=1}^{\eta^{L}(i)-1}
                          \Delta^{L}_m \kappa^{(i-1)}_{\sigma_m^{L}-i+1}\\
   &\hskip65mm
     \mbox{if }\Delta^{L}_1=-1\mbox{ and } i\ge\sigma^{L}_{\eta^{L}(i)};\\
   &\sum_{m=2}^{\eta^{L}(i)}
                          \Delta^{L}_m \kappa^{(i-1)}_{\tau_m^{L}-i+1}
      -\sum_{m=1}^{\eta^{L}(i)}
                          \Delta^{L}_m \kappa^{(i-1)}_{\sigma_m^{L}-i+1}\\
   &\hskip65mm
     \mbox{if }\Delta^{L}_1=-1\mbox{ and } i<\sigma^{L}_{\eta^{L}(i)};\\
   &p^{(i-1)\prime}-\Delta^{L}_{\eta^{L}(i)}(1+\delta_{i,t_{k(i)}})
      +\sum_{m=2}^{\eta^{L}(i)}
                          \Delta^{L}_m \kappa^{(i-1)}_{\tau_m^{L}-i+1}
      -\sum_{m=1}^{\eta^{L}(i)-1}
                          \Delta^{L}_m \kappa^{(i-1)}_{\sigma_m^{L}-i+1}\\
   &\hskip65mm
     \mbox{if }\Delta^{L}_1=+1\mbox{ and } i\ge\sigma^{L}_{\eta^{L}(i)};\\
   &p^{(i-1)\prime}
      +\sum_{m=2}^{\eta^{L}(i)}
                          \Delta^{L}_m \kappa^{(i-1)}_{\tau_m^{L}-i+1}
      -\sum_{m=1}^{\eta^{L}(i)}
                          \Delta^{L}_m \kappa^{(i-1)}_{\sigma_m^{L}-i+1}\\
   &\hskip65mm
     \mbox{if }\Delta^{L}_1=+1\mbox{ and } i<\sigma^{L}_{\eta^{L}(i)}.
  \end{split}
\right.
\end{equation*}
Similarly, for $1\le i<t$ define:
\begin{equation*}
b^{(i)\prime}=
\left\{
  \begin{split}
   &-\Delta^{R}_{\eta^{R}(i)}(1+\delta_{i,t_{k(i)}})
      +\sum_{m=2}^{\eta^{R}(i)}
                          \Delta^{R}_m \kappa^{(i-1)}_{\tau_m^{R}-i+1}
      -\sum_{m=1}^{\eta^{R}(i)-1}
                          \Delta^{R}_m \kappa^{(i-1)}_{\sigma_m^{R}-i+1}\\
   &\hskip65mm
     \mbox{if }\Delta^{R}_1=-1\mbox{ and } i\ge\sigma^{R}_{\eta^{R}(i)};\\
   &\sum_{m=2}^{\eta^{R}(i)}
                          \Delta^{R}_m \kappa^{(i-1)}_{\tau_m^{R}-i+1}
      -\sum_{m=1}^{\eta^{R}(i)}
                          \Delta^{R}_m \kappa^{(i-1)}_{\sigma_m^{R}-i+1}\\
   &\hskip65mm
     \mbox{if }\Delta^{R}_1=-1\mbox{ and } i<\sigma^{R}_{\eta^{R}(i)};\\
   &p^{(i-1)\prime}-\Delta^{R}_{\eta^{R}(i)}(1+\delta_{i,t_{k(i)}})
      +\sum_{m=2}^{\eta^{R}(i)}
                          \Delta^{R}_m \kappa^{(i-1)}_{\tau_m^{R}-i+1}
      -\sum_{m=1}^{\eta^{R}(i)-1}
                          \Delta^{R}_m \kappa^{(i-1)}_{\sigma_m^{R}-i+1}\\
   &\hskip65mm
     \mbox{if }\Delta^{R}_1=+1\mbox{ and } i\ge\sigma^{R}_{\eta^{R}(i)};\\
   &p^{(i-1)\prime}
      +\sum_{m=2}^{\eta^{R}(i)}
                          \Delta^{R}_m \kappa^{(i-1)}_{\tau_m^{R}-i+1}
      -\sum_{m=1}^{\eta^{R}(i)}
                          \Delta^{R}_m \kappa^{(i-1)}_{\sigma_m^{R}-i+1}\\
   &\hskip65mm
     \mbox{if }\Delta^{R}_1=+1\mbox{ and } i<\sigma^{R}_{\eta^{R}(i)}.
  \end{split}
\right.
\end{equation*}


\noindent
We now define values that will turn out to be
$\lfloor a^{(i)}p^{(i)}/p^{(i)\prime}\rfloor$ and
$\lfloor b^{(i)}p^{(i)}/p^{(i)\prime}\rfloor$.
For $0\le i<t$, define:
\begin{equation*}
\tilde a^{(i)}=
\left\{
  \begin{split}
    &-e^{(i)}
      +\sum_{m=2}^{\eta^{L}(i)}
                          \Delta^{L}_m \tkappa^{(i)}_{\tau_m^{L}-i}
      -\sum_{m=1}^{\eta^{L}(i)-1}
                          \Delta^{L}_m \tkappa^{(i)}_{\sigma_m^{L}-i}\\
   &\hskip65mm
     \mbox{if }\Delta^{L}_1=-1\mbox{ and } i\ge\sigma^{L}_{\eta^{L}(i)};\\
    &-e^{(i)}
      +\sum_{m=2}^{\eta^{L}(i)}
                          \Delta^{L}_m \tkappa^{(i)}_{\tau_m^{L}-i}
      -\sum_{m=1}^{\eta^{L}(i)}
                          \Delta^{L}_m \tkappa^{(i)}_{\sigma_m^{L}-i}\\
   &\hskip65mm
     \mbox{if }\Delta^{L}_1=-1\mbox{ and } i<\sigma^{L}_{\eta^{L}(i)};\\
    &p^{(i)}-e^{(i)}
      +\sum_{m=2}^{\eta^{L}(i)}
                          \Delta^{L}_m \tkappa^{(i)}_{\tau_m^{L}-i}
      -\sum_{m=1}^{\eta^{L}(i)-1}
                          \Delta^{L}_m \tkappa^{(i)}_{\sigma_m^{L}-i}\\
   &\hskip65mm
     \mbox{if }\Delta^{L}_1=+1\mbox{ and } i\ge\sigma^{L}_{\eta^{L}(i)};\\
    &p^{(i)}-e^{(i)}
      +\sum_{m=2}^{\eta^{L}(i)}
                          \Delta^{L}_m \tkappa^{(i)}_{\tau_m^{L}-i}
      -\sum_{m=1}^{\eta^{L}(i)}
                          \Delta^{L}_m \tkappa^{(i)}_{\sigma_m^{L}-i}\\
   &\hskip65mm
     \mbox{if }\Delta^{L}_1=+1\mbox{ and } i<\sigma^{L}_{\eta^{L}(i)}.
  \end{split}
\right.
\end{equation*}
Similarly, for $0\le i<t$ define:
\begin{equation*}
\tilde b^{(i)}=
\left\{
  \begin{split}
    &-f^{(i)}
      +\sum_{m=2}^{\eta^{R}(i)}
                          \Delta^{R}_m \tkappa^{(i)}_{\tau_m^{R}-i}
      -\sum_{m=1}^{\eta^{R}(i)-1}
                          \Delta^{R}_m \tkappa^{(i)}_{\sigma_m^{R}-i}\\
   &\hskip65mm
     \mbox{if }\Delta^{R}_1=-1\mbox{ and } i\ge\sigma^{R}_{\eta^{R}(i)};\\
    &-f^{(i)}
      +\sum_{m=2}^{\eta^{R}(i)}
                          \Delta^{R}_m \tkappa^{(i)}_{\tau_m^{R}-i}
      -\sum_{m=1}^{\eta^{R}(i)}
                          \Delta^{R}_m \tkappa^{(i)}_{\sigma_m^{R}-i}\\
   &\hskip65mm
     \mbox{if }\Delta^{R}_1=-1\mbox{ and } i<\sigma^{R}_{\eta^{R}(i)};\\
    &p^{(i)}-f^{(i)}
      +\sum_{m=2}^{\eta^{R}(i)}
                          \Delta^{R}_m \tkappa^{(i)}_{\tau_m^{R}-i}
      -\sum_{m=1}^{\eta^{R}(i)-1}
                          \Delta^{R}_m \tkappa^{(i)}_{\sigma_m^{R}-i}\\
   &\hskip65mm
     \mbox{if }\Delta^{R}_1=+1\mbox{ and } i\ge\sigma^{R}_{\eta^{R}(i)};\\
    &p^{(i)}-f^{(i)}
      +\sum_{m=2}^{\eta^{R}(i)}
                          \Delta^{R}_m \tkappa^{(i)}_{\tau_m^{R}-i}
      -\sum_{m=1}^{\eta^{R}(i)}
                          \Delta^{R}_m \tkappa^{(i)}_{\sigma_m^{R}-i}\\
   &\hskip65mm
     \mbox{if }\Delta^{R}_1=+1\mbox{ and } i<\sigma^{R}_{\eta^{R}(i)}.
  \end{split}
\right.
\end{equation*}

\begin{lemma}\label{A_AprimeLem}
Let $1\le i<t$.
Then:
\begin{equation*}
a^{(i)\prime}=
\left\{
  \begin{array}{ll}
     \displaystyle a^{(i)}+\tilde a^{(i)}+e^{(i)}
       & \mbox{if } i\ne t_{k(i)};  \\[1mm]
     \displaystyle 2a^{(i)}-\tilde a^{(i)}-e^{(i)}
       & \mbox{if } i= t_{k(i)},
  \end{array}
\right.
\end{equation*}
and
\begin{equation*}
b^{(i)\prime}=
\left\{
  \begin{array}{ll}
     \displaystyle b^{(i)}+\tilde b^{(i)}+f^{(i)}
       & \mbox{if } i\ne t_{k(i)};  \\[1mm]
     \displaystyle 2b^{(i)}-\tilde b^{(i)}-f^{(i)}
       & \mbox{if } i= t_{k(i)}.
  \end{array}
\right.
\end{equation*}
\end{lemma}

\Proof
These results follow immediately from Lemma \ref{IndTakLem}.
\cqfd
\medskip

\begin{lemma}\label{Aprime_ALem}
Let $1\le i<t$, $\eta^L=\eta^L(i)$ and $\eta^R=\eta^R(i)$.
Then:
\begin{equation*}
a^{(i-1)}=
\left\{
  \begin{array}{ll}
     \displaystyle a^{(i)\prime}+\Delta^L_{\eta^L+1}
       & \mbox{if } i=\tau^L_{\eta^L+1};  \\[1mm]
     \displaystyle a^{(i)\prime}+\Delta^L_{\eta^L}
       & \mbox{if } i>\sigma^L_{\eta^L} \mbox{ and } i=t_{k(i)}; \\[1mm]
     \displaystyle a^{(i)\prime}-\Delta^L_{\eta^L}
       & \mbox{if } i=\sigma^L_{\eta^L} \mbox{ and } i\ne t_{k(i)}; \\[1mm]
     \displaystyle a^{(i)\prime}
       & \mbox{otherwise,}
  \end{array}
\right.
\end{equation*}
and     
\begin{equation*}
b^{(i-1)}=
\left\{
  \begin{array}{ll}
     \displaystyle b^{(i)\prime}+\Delta^R_{\eta^R+1}
       & \mbox{if } i=\tau^R_{\eta^R+1};  \\[1mm]
     \displaystyle b^{(i)\prime}+\Delta^R_{\eta^R}
       & \mbox{if } i>\sigma^R_{\eta^R} \mbox{ and } i=t_{k(i)}; \\[1mm]
     \displaystyle b^{(i)\prime}-\Delta^R_{\eta^R}
       & \mbox{if } i=\sigma^R_{\eta^R} \mbox{ and } i\ne t_{k(i)}; \\[1mm]
     \displaystyle b^{(i)\prime}
       & \mbox{otherwise.}
  \end{array}
\right.
\end{equation*}
\end{lemma}

\Proof
In the case $i=\tau^L_{\eta^L+1}$, we have $i<\sigma^L_{\eta^L}$,
$\sigma^L_{\eta^L+1}\le i-1<\tau^L_{\eta^L+1}$ and $\eta^L(i-1)=\eta^L+1$.
Then, since $\kappa^{(i-1)}_1=2$, we obtain
$a^{(i-1)}=a^{(i)\prime}+\Delta^L_{\eta^L+1}$.

For $\sigma^L_{\eta^L}<i<\tau^L_{\eta^L}$, then
manifestly $a^{(i-1)}=a^{(i)\prime}+\Delta^L_{\eta^L}\delta_{i,t_{k(i)}}$.

In the case $i=\sigma^L_{\eta^L}$, we have
$\tau^L_{\eta^L+1}\le i-1<\sigma^L_{\eta^L}$.
Since $\kappa^{(i-1)}_1=2$, we immediately obtain
$a^{(i-1)}=a^{(i)\prime}+\Delta^L_{\eta^L}(-1+\delta_{i,t_{k(i)}})$.

Finally, for $\tau^L_{\eta^L+1}<i<\sigma^L_{\eta^L}$, we immediately
see that $a^{(i-1)}=a^{(i)\prime}$.

The expressions for $b^{(i-1)}$ are obtained in precisely the same way.
\cqfd
\medskip

\begin{lemma}\label{Aprime_AprimeLem}
Let $1\le i<t$, $\eta^L=\eta^L(i)$ and $\eta^R=\eta^R(i)$.
Then:
\begin{displaymath}
\tilde a^{(i-1)}=
\left\{
  \begin{array}{ll}
     \displaystyle \tilde a^{(i)}+\Delta^L_{\eta^L+1}
       & \mbox{if } i\ne t_{k(i)}, i=\tau^L_{\eta^L+1} \mbox{ and }
                                   e^{(i-1)}= e^{(i)} ; \\[1mm]
     \displaystyle \tilde a^{(i)}
       & \mbox{if } i\ne t_{k(i)} \mbox{ otherwise;} \\[1mm]
     \displaystyle a^{(i)}-1-\tilde a^{(i)}+\Delta^L_{\eta^L+1}
       & \mbox{if } i= t_{k(i)}, i=\tau^L_{\eta^L+1} \mbox{ and }
                                   e^{(i-1)}\ne e^{(i)} ; \\[1mm]
     \displaystyle a^{(i)}-1-\tilde a^{(i)}
       & \mbox{if } i= t_{k(i)} \mbox{ otherwise,}
  \end{array}
\right.
\end{displaymath}
and     
\begin{displaymath}
\tilde b^{(i-1)}=
\left\{
  \begin{array}{ll}
     \displaystyle \tilde b^{(i)}+\Delta^R_{\eta^R+1}
       & \mbox{if } i\ne t_{k(i)}, i=\tau^R_{\eta^R+1} \mbox{ and }
                                   f^{(i-1)}= f^{(i)} ; \\[1mm]
     \displaystyle \tilde b^{(i)}
       & \mbox{if } i\ne t_{k(i)} \mbox{ otherwise;} \\[1mm]
     \displaystyle b^{(i)}-1-\tilde b^{(i)}+\Delta^R_{\eta^R+1}
       & \mbox{if } i= t_{k(i)}, i=\tau^R_{\eta^R+1} \mbox{ and }
                                   f^{(i-1)}\ne f^{(i)} ; \\[1mm]
     \displaystyle b^{(i)}-1-\tilde b^{(i)}
       & \mbox{if } i= t_{k(i)} \mbox{ otherwise.}
  \end{array}
\right.
\end{displaymath}
\end{lemma}

\Proof
First consider $i\ne t_{k(i)}$, so that
$\tkappa^{(i-1)}_j=\tkappa^{(i)}_{j-1}$ for $1\le j\le t^{(i-1)}$
by Lemma \ref{IndTakLem}.

In the subcase $i=\tau^L_{\eta^L+1}$, we have $i<\sigma^L_{\eta^L}$
and $\sigma^L_{\eta^L+1}\le i-1<\tau^L_{\eta^L+1}$
and $\eta^L(i-1)=\eta^L+1$.
Then, since $\tkappa^{(i-1)}_1=1$, we have
$\tilde a^{(i-1)}-\tilde a^{(i)}=-e^{(i-1)}+\Delta^L_{\eta^L+1}+
e^{(i)}$. The $e^{(i)}=e^{(i-1)}$ case follows immediately
and the $e^{(i)}=1-e^{(i-1)}$ case follows on using
$\Delta^L_{\eta^L+1}=2e^{(i-1)}-1$, which is verified immediately from
the definition of $e^{(i-1)}$.

In the subcase $i=\sigma^L_{\eta^L}$, since $\tkappa^{(i-1)}_1=1$, we have
$\tilde a^{(i-1)}-\tilde a^{(i)}=-e^{(i-1)}-\Delta^L_{\eta^L}+
e^{(i)}=0$ because $e^{(i-1)}=1-e^{(i)}$ and
$\Delta^L_{\eta^L}=2e^{(i)}-1$.

For $\sigma^L_{\eta^L}<i<\tau^L_{\eta^L}$
and $\tau^L_{\eta^L+1}<i<\sigma^L_{\eta^L}$,
we see that $e^{(i-1)}=e^{(i)}$,
whence $\tilde a^{(i-1)}=\tilde a^{(i)}$.

Now when $i=t_{k(i)}$, Lemma \ref{IndTakLem} gives
$\tkappa^{(i-1)}_j=\kappa^{(i)}_{j-1}-\tkappa^{(i)}_{j-1}$
for $2\le j\le t^{(i-1)}$.

In the subcase $i=\tau^L_{\eta^L+1}$,
since $\tkappa^{(i-1)}_1=1$, we have
$\tilde a^{(i-1)}-a^{(i)}+1+\tilde a^{(i)}
=-e^{(i-1)}+1+\Delta^L_{\eta^L+1}-e^{(i)}$.
The $e^{(i)}=1-e^{(i-1)}$ case follows immediately
and the $e^{(i)}=e^{(i-1)}$ case follows on using
$\Delta^L_{\eta^L+1}=2e^{(i-1)}-1$.

In the subcase $i=\sigma^L_{\eta^L}$, since $\tkappa^{(i-1)}_1=1$, we have
$\tilde a^{(i-1)}-a^{(i)}+1+\tilde a^{(i)}
=-e^{(i-1)}+\Delta^L_{\eta^L}+1-e^{(i)}-\Delta^L_{\eta^L}=0$
because $e^{(i-1)}=1-e^{(i)}$.

For $\tau^L_{\eta^L+1}<i<\sigma^L_{\eta^L}$, we obtain
$\tilde a^{(i-1)}-a^{(i)}+1+\tilde a^{(i)}
=-e^{(i-1)}-e^{(i)}+1=0$, since here $k(i-1)=k(i)-1$
implies that $e^{(i-1)}=1-e^{(i)}$.

For $\sigma^L_{\eta^L}<i<\tau^L_{\eta^L}$, we obtain
$\tilde a^{(i-1)}-a^{(i)}+1+\tilde a^{(i)}
=-e^{(i-1)}+\Delta^L_{\eta^L}-e^{(i)}+1=0$,
since here $e^{(i-1)}=e^{(i)}$ and $\Delta^L_{\eta^L}=2e^{(i)}-1$.

The expressions for $\tilde b^{(i-1)}$ are obtained in
precisely the same way.
\cqfd
\medskip

 
\noindent
For $1\le j<d^L$ and $0\le i<\tau^L_{j+1}$ (so that $j<\eta^L(i)$), define:
\begin{align*}
\mu^{(i)*}_j&
  =\sum_{m=2}^{j}
        \Delta^{L}_m (\kappa^{(i)}_{\tau_m^{L}-i}
                           -\kappa^{(i)}_{\sigma_m^{L}-i})
+\left\{
  \begin{array}{ll}
  \displaystyle
  \kappa^{(i)}_{\sigma^{L}_1-i} \quad &\mbox{if }\Delta^{L}_1=-1;\\[1mm]
  \displaystyle
  p^{(i)\prime}-\kappa^{(i)}_{\sigma^{L}_1-i}\quad &\mbox{if }\Delta^{L}_1=+1;
  \end{array}
\right.
\\[1mm]
\mu^{(i)\phantom{*}}_j&
  =\mu^{(i)*}_j+\Delta^{L}_{j+1}\kappa^{(i)}_{\tau^{L}_{j+1}-i},
\end{align*}
and for $0\le i<t_n$, define:
\begin{align*}
\mu^{(i)*}_0&
  =\left\{
  \begin{array}{ll}
  \displaystyle
  \kappa^{(i)}_{t_n-i} \quad &\mbox{if }\Delta^{L}_1=-1;\\[1mm]
  \displaystyle
  p^{(i)\prime}-\kappa^{(i)}_{t_n-i}\quad &\mbox{if }\Delta^{L}_1=+1,
  \end{array}
\right.
\\[1mm]
\mu^{(i)\phantom{*}}_0&
  =\left\{
  \begin{array}{ll}
  \displaystyle
  0 \quad &\mbox{if }\Delta^{L}_1=-1;\\[1mm]
  \displaystyle
  p^{(i)\prime}\phantom{\displaystyle{}-\kappa^{(i)}_{t_n-i}}\quad
                  &\mbox{if }\Delta^{L}_1=+1.
  \end{array}
\right.
\end{align*}
We also set $d^L_0(i)=0$ if both $i<t_n$ and $\sigma^L_1<t_n$,
and $d^L_0(i)=1$ otherwise.
For $0\le i<t$, we define the $(\eta^L(i)-d^L_0(i))$-dimensional
vectors $\muIndUp{i}\!\!$ and $\muIndDn{i}$,
and the $(\eta^L(i+1)-d^L_0(i+1))$-dimensional
vectors $\muIndUpA{i}\!\!$ and $\muIndDnA{i}$ by:
\begin{align*}
\muIndUp{i}\!&=(\mu^{(i)}_{d^L_0(i)},\ldots,\mu^{(i)}_{\eta^L(i)-1}),&
\muIndDn{i}&=(\mu^{(i)*}_{d^L_0(i)},\ldots,\mu^{(i)*}_{\eta^L(i)-1}),\\
\muIndUpA{i}\!&=(\mu^{(i)}_{d^L_0(i+1)},\ldots,\mu^{(i)}_{\eta^L(i+1)-1}),&
\muIndDnA{i}&=(\mu^{(i)*}_{d^L_0(i+1)},\ldots,\mu^{(i)*}_{\eta^L(i+1)-1}).
\end{align*}

\noindent
For $1\le j<d^R$ and $0\le i<\tau^R_{j+1}$ (so that $j<\eta^R(i)$), define:
\begin{align*}
\nu^{(i)*}_j&
  =\sum_{m=2}^{j}
        \Delta^{R}_m (\kappa^{(i)}_{\tau_m^{R}-i}
                           -\kappa^{(i)}_{\sigma_m^{R}-i})
+\left\{
  \begin{array}{ll}
  \displaystyle
  \kappa^{(i)}_{\sigma^{R}_1-i} \quad &\mbox{if }\Delta^{R}_1=-1;\\[1mm]
  \displaystyle
  p^{(i)\prime}-\kappa^{(i)}_{\sigma^{R}_1-i}\quad &\mbox{if }\Delta^{R}_1=+1;
  \end{array}
\right.\\[1mm]
\nu^{(i)\phantom{*}}_j&
  =\nu^{(i)*}_j+\Delta^{R}_{j+1}\kappa^{(i)}_{\tau^{R}_{j+1}-i},
\\[1mm]
\end{align*}
and for $0\le i<t_n$, define:
\begin{align*}
\nu^{(i)*}_0&
  =\left\{
  \begin{array}{ll}
  \displaystyle
  \kappa^{(i)}_{t_n-i} \quad &\mbox{if }\Delta^{R}_1=-1;\\[1mm]
  \displaystyle
  p^{(i)\prime}-\kappa^{(i)}_{t_n-i}\quad &\mbox{if }\Delta^{R}_1=+1,
  \end{array}
\right.
\\[1mm]
\nu^{(i)\phantom{*}}_0&
  =\left\{
  \begin{array}{ll}
  \displaystyle
  0 \quad &\mbox{if }\Delta^{R}_1=-1;\\[1mm]
  \displaystyle
  p^{(i)\prime}\phantom{\displaystyle{}-\kappa^{(i)}_{t_n-i}}\quad
                  &\mbox{if }\Delta^{R}_1=+1.
  \end{array}
\right.
\end{align*}
We also set $d^R_0(i)=0$ if both $i<t_n$ and $\sigma^R_1<t_n$,
and $d^R_0(i)=1$ otherwise.
For $0\le i<t$, we define the $(\eta^R(i)-d^R_0(i))$-dimensional
vectors $\nuIndUp{i}\!\!$ and $\nuIndDn{i}$,
and the $(\eta^R(i+1)-d^R_0(i+1))$-dimensional
vectors $\nuIndUpA{i}\!\!$ and $\nuIndDnA{i}$ by:
\begin{align*}
\nuIndUp{i}\!&=(\nu^{(i)}_{d^R_0(i)},\ldots,\nu^{(i)}_{\eta^R(i)-1}),&
\nuIndDn{i}&=(\nu^{(i)*}_{d^R_0(i)},\ldots,\nu^{(i)*}_{\eta^R(i)-1}),\\
\nuIndUpA{i}\!&=(\nu^{(i)}_{d^R_0(i+1)},\ldots,\nu^{(i)}_{\eta^R(i+1)-1}),&
\nuIndDnA{i}&=(\nu^{(i)*}_{d^R_0(i+1)},\ldots,\nu^{(i)*}_{\eta^R(i+1)-1}).
\end{align*}
(The parentheses that delimit $\muIndUp{i}$, $\muIndDn{i}$,
$\muIndUpA{i}$, $\muIndDnA{i}$ $\nuIndUp{i}$, $\nuIndDn{i}$,
$\nuIndUpA{i}$, $\nuIndDnA{i}$ will be dropped when these symbols
are incorporated into the notation for the path generating functions.)

\begin{lemma}\label{Mu_ALem}
Let $1\le j<d^L$. Then:
\begin{align*}
\mu^{(i-1)*}_j &=a^{(i-1)}-\Delta^L_{j+1}
&& \text{if } i=\tau^L_{j+1};\\*
\mu^{(i-1)*}_j &=a^{(i-1)}-\Delta^L_{j+1}(\kappa^{(i-1)}_{\tau^L_{j+1}-i+1}-1)
&& \text{if } \sigma^L_{j+1}<i\le\tau^L_{j+1};\\*
\mu^{(i-1)*}_j &=a^{(i-1)}-\Delta^L_{j+1}(\kappa^{(i-1)}_{\tau^L_{j+1}-i+1}-2)
&& \text{if } i=\sigma^L_{j+1};\\*
\mu^{(i-1)\phantom{*}}_j &=a^{(i-1)}+\Delta^L_{j+1}
&& \text{if } \sigma^L_{j+1}<i\le\tau^L_{j+1};\\*
\mu^{(i-1)\phantom{*}}_j &=a^{(i)\prime}+(1+\delta_{i,t_{k(i)}})\Delta^L_{j+1}
&& \text{if } i=\sigma^L_{j+1}.\\
\intertext{If $\sigma^L_1<t_n$ then:}
\mu^{(i-1)*}_0 &=a^{(i-1)}-\Delta^L_{1}
&& \text{if } i=t_n;\\*
\mu^{(i-1)*}_0 &=a^{(i-1)}-\Delta^L_{1}(\kappa^{(i-1)}_{t_n-i+1}-1)
&& \text{if } \sigma^L_{1}< i\le t_n;\\*
\mu^{(i-1)*}_0 &=a^{(i-1)}-\Delta^L_{1}(\kappa^{(i-1)}_{t_n-i+1}-2)
&& \text{if } i=\sigma^L_{1};\\*
\mu^{(i-1)\phantom{*}}_0 &=a^{(i-1)}+\Delta^L_{1}
&& \text{if } \sigma^L_{1}< i\le t_n;\\*
\mu^{(i-1)\phantom{*}}_0 &=a^{(i)\prime}+(1+\delta_{i,t_{k(i)}})\Delta^L_{1}
&& \text{if } i=\sigma^L_{1}.
\end{align*}

Now let $1\le j<d^R$. Then:
\begin{align*}
\nu^{(i-1)*}_j &=b^{(i-1)}-\Delta^R_{j+1}
&& \text{if } i=\tau^R_{j+1};\\*
\nu^{(i-1)*}_j &=b^{(i-1)}-\Delta^R_{j+1}(\kappa^{(i-1)}_{\tau^R_{j+1}-i+1}-1)
&& \text{if } \sigma^R_{j+1}<i\le\tau^R_{j+1};\\*
\nu^{(i-1)*}_j &=b^{(i-1)}-\Delta^R_{j+1}(\kappa^{(i-1)}_{\tau^R_{j+1}-i+1}-2)
&& \text{if } i=\sigma^R_{j+1};\\*
\nu^{(i-1)\phantom{*}}_j &=b^{(i-1)}+\Delta^R_{j+1}
&& \text{if } \sigma^R_{j+1}<i\le\tau^R_{j+1};\\*
\nu^{(i-1)\phantom{*}}_j &=b^{(i)\prime}+(1+\delta_{i,t_{k(i)}})\Delta^R_{j+1}
&& \text{if } i=\sigma^R_{j+1}.\\
\intertext{If $\sigma^R_1<t_n$ then:}
\nu^{(i-1)*}_0 &=b^{(i-1)}-\Delta^R_{1}
&& \text{if } i=t_n;\\*
\nu^{(i-1)*}_0 &=b^{(i-1)}-\Delta^R_{1}(\kappa^{(i-1)}_{t_n-i+1}-1)
&& \text{if } \sigma^R_{1}< i\le t_n;\\*
\nu^{(i-1)*}_0 &=b^{(i-1)}-\Delta^R_{1}(\kappa^{(i-1)}_{t_n-i+1}-2)
&& \text{if } i=\sigma^R_{1};\\*
\nu^{(i-1)\phantom{*}}_0 &=b^{(i-1)}+\Delta^R_{1}
&& \text{if } \sigma^R_{1}< i\le t_n;\\*
\nu^{(i-1)\phantom{*}}_0 &=b^{(i)\prime}+(1+\delta_{i,t_{k(i)}})\Delta^R_{1}
&& \text{if } i=\sigma^R_{1}.
\end{align*}
\end{lemma}

\Proof
The first four expressions follow immediately from the definitions
after noting that $\sigma^L_{j+1}\le i\le\tau^L_{j+1}$ implies that
$\eta^L(i-1)=j+1$ (and using $\kappa^{(i-1)}_1=2$). The fifth is
immediate. The next five follow from the definitions after noting
that $\sigma^L_1\le i\le t_n$ implies that $\eta^L(i-1)=1$.
The other expressions follow similarly.
\cqfd
\medskip

Using (\ref{Seq2Eq}), it is easily verified that
$\boldmu^{(i)},\boldmu^{(i)*},\boldnu^{(i)},\boldnu^{(i)*}$
satisfy the first three criteria for being a mazy-four in the
$(p^{(i)},p^{(i)\prime})$-model sandwiching $(a^{(i)},b^{(i)})$.
That they satisfy the fourth criterion is
most readily determined during the main induction.
However, the burden will be lessened by,
for $1\le j<d^L$ and $i<\tau^L_{j+1}$, defining:
\begin{align*}
\tilde\mu^{(i)*}_j
  &=\sum_{m=2}^{j}
        \Delta^{L}_m (\tkappa^{(i)}_{\tau_m^{L}-i}
                           -\tkappa^{(i)}_{\sigma_m^{L}-i})
+\left\{
  \begin{array}{ll}
  \displaystyle
  \tkappa^{(i)}_{\sigma^{L}_1-i} \quad &\mbox{if }\Delta^{L}_1=-1;\\[1mm]
  \displaystyle
  p^{(i)}-\tkappa^{(i)}_{\sigma^{L}_1-i}\quad &\mbox{if }\Delta^{L}_1=+1;
  \end{array}
\right.\\[1mm]
\tilde\mu^{(i)\phantom{*}}_j
  &=\tilde\mu^{(i)*}_j+\Delta^{L}_{j+1}\tkappa^{(i)}_{\tau^{L}_{j+1}-i},
\end{align*}
for $1\le j<d^R$ and $i<\tau^R_{j+1}$, defining:
\begin{align*}
\tilde\nu^{(i)*}_j
  &=\sum_{m=2}^{j}
        \Delta^{R}_m (\tkappa^{(i)}_{\tau_m^{R}-i}
                           -\tkappa^{(i)}_{\sigma_m^{R}-i})
+\left\{
  \begin{array}{ll}
  \displaystyle
  \tkappa^{(i)}_{\sigma^{R}_1-i} \quad &\mbox{if }\Delta^{R}_1=-1;\\[1mm]
  \displaystyle
  p^{(i)}-\tkappa^{(i)}_{\sigma^{R}_1-i}\quad &\mbox{if }\Delta^{R}_1=+1;
  \end{array}
\right.\\[1mm]
\tilde\nu^{(i)\phantom{*}}_j
&=\tilde\nu^{(i)*}_j+\Delta^{R}_{j+1}\tkappa^{(i)}_{\tau^{R}_{j+1}-i}.
\end{align*}
and for $0\le i<t_n$, defining:
\begin{align*}
\tilde\mu^{(i)*}_0&
  =\left\{
  \begin{array}{ll}
  \displaystyle
  \tilde\kappa^{(i)}_{t_n-i} \quad &\mbox{if }\Delta^{L}_1=-1;\\[1mm]
  \displaystyle
  p^{(i)}-\tilde\kappa^{(i)}_{t_n-i}\quad &\mbox{if }\Delta^{L}_1=+1;
  \end{array}
\right.
\\[1mm]
\tilde\mu^{(i)\phantom{*}}_0&
  =\left\{
  \begin{array}{ll}
  \displaystyle
  0 \quad &\mbox{if }\Delta^{L}_1=-1;\\[1mm]
  \displaystyle
  p^{(i)}\phantom{\displaystyle{}-\kappa^{(i)}_{t_n-i}}\quad
                  &\mbox{if }\Delta^{L}_1=+1,
  \end{array}
\right.
\\[1mm]
\tilde\nu^{(i)*}_0&
  =\left\{
  \begin{array}{ll}
  \displaystyle
  \tilde\kappa^{(i)}_{t_n-i} \quad &\mbox{if }\Delta^{R}_1=-1;\\[1mm]
  \displaystyle
  p^{(i)}-\tilde\kappa^{(i)}_{t_n-i}\quad &\mbox{if }\Delta^{R}_1=+1;
  \end{array}
\right.
\\[1mm]
\tilde\nu^{(i)\phantom{*}}_0&
  =\left\{
  \begin{array}{ll}
  \displaystyle
  0 \quad &\mbox{if }\Delta^{R}_1=-1;\\[1mm]
  \displaystyle
  p^{(i)}\phantom{\displaystyle{}-\kappa^{(i)}_{t_n-i}}\quad
                  &\mbox{if }\Delta^{R}_1=+1,
  \end{array}
\right.
\end{align*}

\begin{lemma}\label{MuInterLem}
Let $0\le j<d^L$ with $0<i<\tau^L_{j+1}$ if $j>0$, and $0<i<t_n$ if $j=0$.
If $i\ne t_{k(i)}$ then:
\begin{align*}
\mu^{(i-1)\phantom{*}}_j&=\mu^{(i)}_j+\tilde\mu^{(i)}_j;&
\tilde\mu^{(i-1)\phantom{*}}_j&=\tilde\mu^{(i)}_j;\\
\mu^{(i-1)*}_j&=\mu^{(i)*}_j+\tilde\mu^{(i)*}_j;&
\tilde\mu^{(i-1)*}_j&=\tilde\mu^{(i)*}_j,\\
\intertext{and if $i=t_{k(i)}$ then:}
\mu^{(i-1)\phantom{*}}_j&=2\mu^{(i)}_j-\tilde\mu^{(i)}_j;&
\tilde\mu^{(i-1)\phantom{*}}_j&=\mu^{(i)}_j-\tilde\mu^{(i)}_j;\\
\mu^{(i-1)*}_j&=2\mu^{(i)*}_j-\tilde\mu^{(i)*}_j;&
\tilde\mu^{(i-1)*}_j&=\mu^{(i)*}_j-\tilde\mu^{(i)*}_j.
\end{align*}

Let $0\le j<d^R$ with $0<i<\tau^R_{j+1}$ if $j>0$, and $0<i<t_n$ if $j=0$.
If $i\ne t_{k(i)}$ then:
\begin{align*}
\nu^{(i-1)\phantom{*}}_j&=\nu^{(i)}_j+\tilde\nu^{(i)}_j;&
\tilde\nu^{(i-1)\phantom{*}}_j&=\tilde\nu^{(i)}_j;\\
\nu^{(i-1)*}_j&=\nu^{(i)*}_j+\tilde\nu^{(i)*}_j;&
\tilde\nu^{(i-1)*}_j&=\tilde\nu^{(i)*}_j,\\
\intertext{and if $i=t_{k(i)}$ then:}
\nu^{(i-1)\phantom{*}}_j&=2\nu^{(i)}_j-\tilde\nu^{(i)}_j;&
\tilde\nu^{(i-1)\phantom{*}}_j&=\nu^{(i)}_j-\tilde\nu^{(i)}_j;\\
\nu^{(i-1)*}_j&=2\nu^{(i)*}_j-\tilde\nu^{(i)*}_j;&
\tilde\nu^{(i-1)*}_j&=\nu^{(i)*}_j-\tilde\nu^{(i)*}_j.
\end{align*}
\end{lemma}

\Proof
These results follow readily from Lemma \ref{IndTakLem}.
\cqfd
\medskip


For each $t$-dimensional vector
$\boldu=(u_1,u_2,\ldots,u_t)$, define the $(t-1)$-dimensional vector
$\boldu^{(\flat,k)}=(u^{(\flat,k)}_1,u^{(\flat,k)}_2,\ldots,
u^{(\flat,k)}_{t-1})$ by
\begin{equation}
u^{(\flat,k)}_i=
\left\{
  \begin{array}{cl}
       0 &\quad\mbox{if } t_{k'}<i\le t_{k'+1},
                                \; k'\equiv k\;(\mod 2);\\[0.5mm]
      u_i &\quad\mbox{if } t_{k'}<i\le t_{k'+1},
                                \; k'\not\equiv k\;(\mod 2),
  \end{array} \right.
\end{equation}
and the $(t-1)$-dimensional vector
$\boldu^{(\sharp,k)}=(u^{(\sharp,k)}_1,u^{(\sharp,k)}_2,\ldots,
u^{(\sharp,k)}_{t-1})$ by
\begin{equation}
u^{(\sharp,k)}_i=
\left\{
  \begin{array}{cl}
     u_i &\quad\mbox{if } t_{k'}<i\le t_{k'+1},
                                \; k'\equiv k\;(\mod 2);\\[0.5mm]
      0 &\quad\mbox{if } t_{k'}<i\le t_{k'+1},
                                \; k'\not\equiv k\;(\mod 2),
  \end{array} \right.
\end{equation}
For convenience, we sometimes write
$\boldu_{(\flat,k)}$ instead of $\boldu^{(\flat,k)}$, and
$\boldu_{(\sharp,k)}$ instead of $\boldu^{(\sharp,k)}$.


Now for $0\le i\le t-3$, define:
\begin{equation}\label{Ind1Eq}
\begin{split}
&F^{(i)}(\boldu^{L},\boldu^{R},m_i,m_{i+1};q)=\\[3.5mm]
&\qquad
\sum
  q^{\frac{1}{4}\hat{\boldm}^{(i+1)T}\boldC\hat{\boldm}^{(i+1)}
   +\frac{1}{4} m_i^2
   -\frac{1}{2} m_im_{i+1}
   -\frac{1}{2}(\boldu^{L}_{(\flat,k(i))}+\boldu^{R}_{(\sharp,k(i))})\cdot
   \boldm^{(i)}
  +\frac{1}{4}\gamma''_i}\\
&\hskip45mm
\times\quad
  \prod_{j=i+1}^{t-1}
  \left[
  {m_j-\frac{1}{2}(\boldC^*\hat{\boldm}^{(i)}
                       \!-\!\boldu^{L}\!-\!\boldu^{R})_j\atop m_j}
  \right]_q.
\end{split}
\end{equation}
With $(Q_1,Q_2,\ldots,Q_{t-1})=\boldQ(\boldu^{L}+\boldu^{R})$
as defined in Section \ref{MNsysSec},
the sum here is to be taken over all
$(m_{i+2},m_{i+3},\ldots,m_{t-1})\equiv(Q_{i+2},Q_{i+3},\ldots,Q_{t-1})$.
The $(t-1)$-dimensional
$\boldm^{(i)}=(0,\ldots,0,m_{i+1},m_{i+2},m_{i+3},\ldots,m_{t-1})$
has its first $i$ components equal to zero.
The $t$-dimensional
$\hat{\boldm}^{(i)}=(0,\ldots,0,m_{i},m_{i+1},m_{i+2},\ldots,m_{t-1})$
has its first $i$ components equal to zero.
The matrix $\boldC$ is as defined in Section \ref{QuadSec}.
We define $F^{(t-2)}(\boldu^{L},\boldu^{R},m_{t-2},m_{t-1};q)$
using (\ref{Ind1Eq}), but with the summation omitted from the right side.

We also define:
\begin{equation}\label{Ind2Eq}
F^{(t-1)}(\boldu^{L},\boldu^{R},m_{t-1},m_{t};q)
=q^{\frac{1}{4}m_{t-1}^2+\frac{1}{4}\gamma''_{t-1}}
\delta_{m_t,0}.
\end{equation}
For convenience, we also set $Q_t=0$.

By Lemma \ref{PartitionInvLem},
$\left[{m+n\atop m}\right]_{q^{-1}}=q^{-mn}\left[{m+n\atop m}\right]_q$.
Then since $C_{i+1,i}=-1$, it follows that for $0\le i\le t-2$:
\begin{equation}\label{Invq1Eq}
\begin{split}
&F^{(i)}(\boldu^{L},\boldu^{R},m_i,m_{i+1};q^{-1})=\\[3.5mm]
&\qquad
\sum
  q^{\frac{1}{4}\hat{\boldm}^{(i+1)T}\boldC\hat{\boldm}^{(i+1)}
   -\frac{1}{4} m_i^2
   -\frac{1}{2}(\boldu^{L}_{(\flat,k(i)-1)}+\boldu^{R}_{(\sharp,k(i)-1)})
   \cdot \boldm^{(i)}
  -\frac{1}{4}\gamma''_i}\\
&\hskip45mm
\times\quad
  \prod_{j=i+1}^{t-1}
  \left[
  {m_j-\frac{1}{2}(\boldC^*\hat{\boldm}^{(i)}
                       \!-\!\boldu^{L}\!-\!\boldu^{R})_j\atop m_j}
  \right]_q,
\end{split}
\end{equation}
where, as above, the sum is taken over all
$(m_{i+2},m_{i+3},\ldots,m_{t-1})\equiv(Q_{i+2},Q_{i+3},$ $\ldots,Q_{t-1})$,
except in the $i=t-2$ case where the summation is omitted.
Of course, we also have:
\begin{equation}\label{Invq2Eq}
F^{(t-1)}(\boldu^{L},\boldu^{R},m_{t-1},m_{t};q^{-1})
=q^{-\frac{1}{4}m_{t-1}^2-\frac{1}{4}\gamma''_{t-1}}
\delta_{m_t,0}.
\end{equation}


\newpage

\subsection{The induction}\label{ProofIndSec}

In this section, we bring together the results on
the $\B$-transform (Corollary \ref{MazyBijCor}),
the $\D$-transform (Corollary \ref{MazyDijCor}), 
path extension (Lemmas \ref{ExtGen1Lem} and \ref{ExtGen2Lem}),
and path truncation (Lemmas \ref{AttenGen1Lem} and \ref{AttenGen2Lem}),
to obtain Theorem \ref{CoreThrm} by means of a huge induction proof.

Throughout this section, we restrict consideration to the cases for
which $p'>2p$.
We fix runs 
$\run^L=\{\tau^{L}_j,\sigma^{L}_j,\Delta^{L}_j\}_{j=1}^{d^{L}}$ and
$\run^R=\{\tau^{R}_j,\sigma^{R}_j,\Delta^{R}_j\}_{j=1}^{d^{R}}$
and make use of the definitions of Section \ref{CoreSec}.
Except in the final Corollary \ref{StarInterCor}, these runs will
be assumed to be reduced and the definitions of Section \ref{IndParamSec}
will be used.
We also make use of the results of Appendices \ref{AppBSec} and \ref{DireApp}.
For typographical convenience, we abbreviate
$\rho^{p^{(i)},p^{(i)\prime}}$ (defined in Section \ref{BandSec})
to $\rho^{(i)}$.

\begin{lemma}\label{CoreIndLem}
Let $0\le i<t$.
If $m_i\equiv Q_i$ and $m_{i+1}\equiv Q_{i+1}$ then:
\begin{equation*}
\mchi^{p^{(i)},p^{(i)\prime}}_{a^{(i)},b^{(i)},e^{(i)},f^{(i)}}
              (m_i,m_{i+1}+2\delta^{p^{(i)},p^{(i)\prime}}_{a^{(i)},e^{(i)}})
\left\{
\begin{matrix}
\muIndUp{i};\nuIndUp{i}\\
\muIndDn{i};\nuIndDn{i}
\end{matrix} \right\}
= F^{(i)}(\boldu^{L},\boldu^{R},m_i,m_{i+1}).
\end{equation*}
In addition,
\begin{displaymath}
\begin{array}{l}
\bullet\quad 
  \alpha^{p^{(i)},p^{(i)\prime}}_{a^{(i)},b^{(i)}}=\alpha_i'';\\[1mm]
\bullet\quad
  \beta^{p^{(i)},p^{(i)\prime}}_{a^{(i)},b^{(i)},e^{(i)},f^{(i)}}=\beta_i';
                                                              \\[1mm]
\bullet\quad
  \lfloor a^{(i)}p^{(i)}/p^{(i)\prime}\rfloor=\tilde a^{(i)};\\[1mm]
\bullet\quad
  \lfloor b^{(i)}p^{(i)}/p^{(i)\prime}\rfloor=\tilde b^{(i)};\\[1mm]
\bullet\quad
     \mbox{if }\delta^{p^{(i)},p^{(i)\prime}}_{a^{(i)},e^{(i)}}=1
     \mbox{ then } i+1=\tau^L_{j} \mbox{ for some } j
     \mbox{ and } \mu^{(i)}_{j-1}=a^{(i)}-(-1)^{e^{(i)}};\\[1mm]
\bullet\quad
     \mbox{if }\delta^{p^{(i)},p^{(i)\prime}}_{b^{(i)},f^{(i)}}=1
     \mbox{ then } i+1=\tau^R_{j} \mbox{ for some } j
     \mbox{ and } \nu^{(i)}_{j-1}=b^{(i)}-(-1)^{f^{(i)}};\\[1mm]
\bullet\quad
     \mbox{if there exists $j$ and $k$ for which }
     t_k \le i \mbox{ and }
     \tau^L_{j+1}\le i < t_{k+1} < \sigma^L_{j} \mbox{ then}\\[1mm]
\qquad
     a^{(i)} \mbox{ is interfacial in the }
     (p^{(i)},p^{(i)\prime})\mbox{-model.}\\[1mm]
\qquad
     \mbox{If } t_k \le i \mbox{ and }
     \tau^L_{j+1}\le i < t_{k+1} \le \sigma^L_{j}
     \mbox{ and } a^{(i)} \mbox{ is not interfacial in the }\\[1mm]
\qquad
     (p^{(i)},p^{(i)\prime})\mbox{-model then }
     \tau^L_{j}=\sigma^L_j+1;\\[1mm]
\bullet\quad
     \mbox{if there exists $j$ and $k$ for which }
     t_k \le i \mbox{ and }
     \tau^R_{j+1}\le i < t_{k+1} < \sigma^R_{j} \mbox{ then}\\[1mm]
\qquad
     b^{(i)} \mbox{ is interfacial in the }
     (p^{(i)},p^{(i)\prime})\mbox{-model.}\\[1mm]
\qquad
     \mbox{If } t_k \le i \mbox{ and }
     \tau^R_{j+1}\le i < t_{k+1} \le \sigma^R_{j}
     \mbox{ and } b^{(i)} \mbox{ is not interfacial in the }\\[1mm]
\qquad
     (p^{(i)},p^{(i)\prime})\mbox{-model then }
     \tau^R_{j}=\sigma^R_j+1;\\[1mm]
\bullet\quad
     \muIndUp{i}, \muIndDn{i}, \nuIndUp{i}, \nuIndDn{i}
     \mbox{ are an interfacial mazy-four in the}\\[1mm]
\qquad
     (p^{(i)},p^{(i)\prime})\mbox{-model sandwiching }
     (a^{(i)},b^{(i)});\\[1mm]
\bullet\quad
     \mbox{If } d^L_0(i)\le j<\eta^L(i) \mbox{ then }
     \rho^{(i)}(\mu^{(i)*}_j)=\tilde\mu^{(i)*}_j
     \mbox{ and }
     \rho^{(i)}(\mu^{(i)}_j)=\tilde\mu^{(i)}_j;\\[1mm]
\qquad
     \mbox{if } d^R_0(i)\le j<\eta^R(i) \mbox{ then }
     \rho^{(i)}(\nu^{(i)*}_j)=\tilde\nu^{(i)*}_j
     \mbox{ and }
     \rho^{(i)}(\nu^{(i)}_j)=\tilde\nu^{(i)}_j.
%
%
\end{array}
\end{displaymath}
\end{lemma}

\Proof This is proved by downward induction.
For $i=t-1$,
we have $i\ge t_n$, $i\ge\sigma^L_1$ and $i\ge\sigma^R_1$,
and consequently, $d^L_0(i)=d^R_0(i)=1$ and $\eta^L(i)=\eta^R(i)=1$.
The definitions of Section \ref{IndParamSec} then
yield $p^{(i)\prime}=3$, $p^{(i)}=1$,
and if $\Delta^L_1=-1$ then
$a^{(i)}=1$, $e^{(i)}=0$ and $(\boldDelta^L)_t=0$;
and if $\Delta^L_1=+1$ then
$a^{(i)}=2$, $e^{(i)}=1$ and $(\boldDelta^L)_t=-1$.
Similarly,
if $\Delta^R_1=-1$ then
$b^{(i)}=1$, $f^{(i)}=0$ and $(\boldDelta^R)_t=0$;
and if $\Delta^R_1=+1$ then
$b^{(i)}=2$, $f^{(i)}=1$ and $(\boldDelta^R)_t=-1$.
In particular, we now immediately obtain
$\delta^{p^{(i)},p^{(i)\prime}}_{a^{(i)},e^{(i)}}=0$.
Via Section \ref{ConSec}, we obtain
$\alpha_{t-1}''=\beta_{t-1}'=(\boldDelta^L)_t-(\boldDelta^R)_t$
and
$\gamma_{t-1}''=-((\boldDelta^L)_t-(\boldDelta^R)_t)^2$.
For $i=t-1$, the first statement of our induction proposition
is now seen to hold via  Lemma \ref{SeedLem}.
Each of the bulleted items follows readily.

Now assume the result holds for a particular $i$ with $1\le i<t$.
Let $k=k(i)$ so that $t_{k}\le i<t_{k+1}$, and let
$\eta^L=\eta^L(i)$ and $\eta^R=\eta^R(i)$ so that
$\tau_{\eta^L+1}\le i<\tau_{\eta^L}$ and
$\tau_{\eta^R+1}\le i<\tau_{\eta^R}$.
We also set $\eta^{L\prime}=\eta^L(i-1)$ and $\eta^{R\prime}=\eta^R(i-1)$.
Then $\eta^{L\prime}=\eta^{L}$ unless $i=\tau^L_{\eta^L+1}$
in which case $\eta^{L\prime}=\eta^{L}+1$; and
$\eta^{R\prime}=\eta^{R}$ unless $i=\tau^R_{\eta^R+1}$
in which case $\eta^{R\prime}=\eta^{R}+1$.

First consider the case $i\ne t_{k}$.
Equation (\ref{Const3Eq}) gives $\alpha_i=\alpha_i''$,
$\beta_i=\beta_i'$ and $\gamma_i=\gamma_i''$.
Let $m_{i-1}\equiv Q_{i-1}$.
On setting $M=m_{i-1}+u^L_i+u^R_i$, equations (\ref{MNEq2}),
(\ref{MNmatEq2}) and (\ref{ParityDef}) imply that $M\equiv Q_{i+1}$.
Then, use of the induction hypothesis and Lemmas \ref{MazyBijCor},
\ref{IndTakLem}, \ref{A_AprimeLem} and \ref{MuInterLem} yields:
\begin{equation*}
\begin{split}
&\mchi^{p^{(i-1)},p^{(i-1)\prime}}_{a^{(i)\prime},b^{(i)\prime},
                       e^{(i)},f^{(i)}} (M,m_{i})
\left\{
\begin{matrix}
\muIndUpA{i-1};\nuIndUpA{i-1}\\
\muIndDnA{i-1};\nuIndDnA{i-1}
\end{matrix}
\right\}\\[2mm]
&\hskip15mm
=\sum_{m_{i+1}\equiv Q_{i+1}}
q^{\frac{1}{4}(M-m_i)^2-\frac{1}{4}\beta^{2}_i}
\left[{\frac{1}{2}(M+m_{i+1})}\atop m_i\right]_q
F^{(i)}(\boldu^{L},\boldu^{R},m_i,m_{i+1}).
\end{split}
\end{equation*}
Lemma \ref{MazyBijCor} also gives
$\alpha^{p^{(i-1)},p^{(i-1)\prime}}_{a^{(i)\prime},b^{(i)\prime}}
=\alpha_i+\beta_i$, and
$\beta^{p^{(i-1)},p^{(i-1)\prime}}_{a^{(i)\prime},b^{(i)\prime},
e^{(i)},f^{(i)}}=\beta_i$.
Lemma \ref{StartPtLem} implies, via Lemma \ref{A_AprimeLem}, that
$\lfloor a^{(i)\prime}p^{(i-1)}/p^{(i-1)\prime}\rfloor=\tilde a^{(i)}$
and
$\lfloor b^{(i)\prime}p^{(i-1)}/p^{(i-1)\prime}\rfloor=\tilde b^{(i)}$.
Lemma \ref{StartPtLem} also implies that
$\delta^{p^{(i-1)},p^{(i-1)\prime}}_{a^{(i)\prime},e^{(i)}}=
\delta^{p^{(i-1)},p^{(i-1)\prime}}_{b^{(i)\prime},f^{(i)}}=0$.
Lemma \ref{MazyBijCor} also shows that
$\muIndUpA{i-1}$, $\muIndDnA{i-1}$, $\nuIndUpA{i-1}$,
$\nuIndDnA{i-1}$ are an interfacial mazy-four in the
$(p^{(i-1)},p^{(i-1)\prime})$-model sandwiching
$(a^{(i)\prime},b^{(i)\prime})$.
Together with Lemma \ref{InterLem}(1), Lemma \ref{MuInterLem} also
shows that $\rho^{(i-1)}(\mu^{(i-1)*}_j)=\tilde\mu^{(i-1)*}_j$
and $\rho^{(i-1)}(\mu^{(i-1)}_j)=\tilde\mu^{(i-1)}_j$
for $d^L_0(i)\le j<\eta^L(i)$,
and that $\rho^{(i-1)}(\nu^{(i-1)*}_j)=\tilde\nu^{(i-1)*}_j$
and $\rho^{(i-1)}(\nu^{(i-1)}_j)=\tilde\nu^{(i-1)}_j$
for $d^R_0(i)\le j<\eta^R(i)$.

Since $M=m_{i-1}+u^L_i+u^R_i$, on noting that $t_k<i<t_{k+1}$, we have:
\begin{displaymath}
M+m_{i+1}=2m_i-(\boldC^*\hat{\boldm}^{(i-1)}-\boldu^{L}-\boldu^{R})_i,
\end{displaymath}
and
\begin{displaymath}
\begin{array}{l}
\displaystyle
\hat{\boldm}^{(i+1)T}\boldC\hat{\boldm}^{(i+1)}+m_i^2-2m_im_{i+1}+(M-m_i)^2\\
\qquad=
\hat{\boldm}^{(i)T}\boldC\hat{\boldm}^{(i)}+M^2-2Mm_i\\
\qquad=
\hat{\boldm}^{(i)T}\boldC\hat{\boldm}^{(i)}+m_{i-1}^2-2m_im_{i-1}
+2(m_{i-1}-m_i)(u^L_i+u^R_i)+(u^L_i+u^R_i)^2.
\end{array}
\end{displaymath}
(In the case $i=t-1$, we require this expression after substituting
$m_t=0$.) 
Thence,
\begin{equation*}
\begin{split}
&\mchi^{p^{(i-1)},p^{(i-1)\prime}}_{a^{(i)\prime},b^{(i)\prime},
                       e^{(i)},f^{(i)}} (m_{i-1}+u^L_i+u^R_i,m_{i})
\left\{
\begin{matrix}
\muIndUpA{i-1};\nuIndUpA{i-1}\\
\muIndDnA{i-1};\nuIndDnA{i-1}
\end{matrix}
\right\}\\[1.5mm]
&\hskip5mm
=\sum
q^{\frac{1}{4}\hat{\boldm}^{(i)T}\boldC\hat{\boldm}^{(i)}
  +\frac{1}{4}m_{i-1}^2
  -\frac{1}{2}m_im_{i-1}
  +\frac{1}{2}(m_{i-1}-m_i)(u^L_i+u^R_i)
  +\frac{1}{4}(u^L_i+u^R_i)^2
  +\frac{1}{4}\gamma_i-\frac{1}{4}\beta_i^2}\\[1.5mm]
&\hskip20mm
\times\quad
  q^{-\frac{1}{2}(\boldu^{L}_{(\flat,k)}+\boldu^{R}_{(\sharp,k)})\cdot
    \boldm^{(i)}}
  \prod_{j=i}^{t-1}
  \left[
  {m_j-\frac{1}{2}(\boldC^*\hat{\boldm}^{(i-1)}
                       \!-\!\boldu^{L}\!-\!\boldu^{R})_j\atop m_j}
  \right]_q,
\end{split}
\end{equation*}
where the sum is over all
$(m_{i+1},m_{i+2},\ldots,m_{t-1})\equiv(Q_{i+1},Q_{i+2},\ldots,Q_{t-1})$.

We now check that we can apply Lemma \ref{AttenGen2Lem}
(for path truncation on the right) if $i=\sigma^R_{\eta^R}$,
and Lemma \ref{ExtGen2Lem} (for path extension on the right)
if $i=\tau^R_{\eta^R+1}$.

If $i=\sigma^R_{\eta^R}$ then $u^R_i=1$, and we see that
$\Delta^R_{\eta^R}=-(-1)^{f^{(i)}}=(-1)^{f^{(i-1)}}$ from (\ref{DefFiEq}),
and $b^{(i-1)}=b^{(i)\prime}-\Delta^R_{\eta^R}$ from Lemma \ref{Aprime_ALem}.
Also note that $\eta^{R\prime}=\eta^R$.
If $\eta^R>1$ or both $\eta^R=1$ and $\sigma^R_1<t_n$ then
Lemma \ref{Mu_ALem} implies that
$\nu^{(i-1)}_{\eta^R-1}=b^{(i)\prime}+\Delta^R_{\eta^R}$
and $\nu^{(i-1)*}_{\eta^R-1}\ne b^{(i-1)}$.
Otherwise, when $i=\sigma^R_1>t_n$ then
$b^{(i)\prime}\in\{1,p^{(i-1)\prime}-1\}$ from the definition.

If $i=\tau^R_{\eta^R+1}$ (when necessarily $\eta^R<d^R$)
then $u^R_i=-1$ and $\eta^{R\prime}=\eta^R+1$.
We see that $\Delta^R_{\eta^{R\prime}}=-(-1)^{f(i-1)}$ from (\ref{DefFiEq})
and $b^{(i-1)}=b^{(i)\prime}+\Delta^R_{\eta^{R\prime}}$
via Lemma \ref{Aprime_ALem}.
In addition, Lemma \ref{Mu_ALem} gives
$\nu^{(i-1)*}_{\eta^R}=b^{(i-1)}-\Delta^R_{\eta^{R\prime}}$ and
$\nu^{(i-1)}_{\eta^R}=b^{(i-1)}+\Delta^R_{\eta^{R\prime}}$.
The definition of a reduced run implies that this case, $i=\tau^R_{\eta^R+1}$
and $i\ne t_k$, only occurs when $t_{n-1}\le t_n-1=\tau^R_2=i<t_n<\sigma^R_1$.
The induction hypothesis then implies that $b^{(i)}$ is interfacial
in the $(p^{(i)},p^{(i)\prime})$-model and
$\delta^{p^{(i)},p^{(i)\prime}}_{b^{(i)},f^{(i)}}=0$.
Thereupon, by Lemma \ref{StartPtLem}, $b^{(i)\prime}$ is interfacial
in the $(p^{(i-1)},p^{(i-1)\prime})$-model.
Also note that $d^R_0(i)=1$ and so
$\nuIndUpA{i-1}=\nuIndDnA{i-1}=()$ here.

Then, using Lemma \ref{AttenGen2Lem} if $i=\sigma^R_{\eta^R}$,
or Lemma \ref{ExtGen2Lem} if $i=\tau^R_{\eta^R+1}$, yields:

\begin{equation*}
\begin{split}
&\mchi^{p^{(i-1)},p^{(i-1)\prime}}_{a^{(i)\prime},b^{(i-1)},e^{(i)},f^{(i-1)}}
                   (m_{i-1}+u^L_i,m_{i})
\left\{
\begin{matrix}
\muIndUpA{i-1};\nuIndUp{i-1}\\
\muIndDnA{i-1};\nuIndDn{i-1}
\end{matrix}
\right\}\\[1.5mm]
&\hskip15mm
=
q^{-\frac12u^R_i(m_{i-1}+u^L_i+u^R_i
   -\Delta^R_{\eta^{R\prime}}(\alpha_i+\beta_i))}
\\[1.5mm]
&\hskip25mm
\times\quad
\mchi^{p^{(i-1)},p^{(i-1)\prime}}_{a^{(i)\prime},b^{(i)\prime},e^{(i)},f^{(i)}}
                   (m_{i-1}+u^L_i+u^R_i,m_{i})
\left\{
\begin{matrix}
\muIndUpA{i-1};\nuIndUpA{i-1}\\
\muIndDnA{i-1};\nuIndDnA{i-1}
\end{matrix}
\right\},
\end{split}
\end{equation*}
noting that $d^R_0(i)=d^R_0(i-1)$ because $i\ne t_n$.
If neither $i=\sigma^R_{\eta^R}$ nor $i=\tau^R_{\eta^R+1}$ then
(noting that $i\ne t_{k}$) $u^R_i=0$, $f^{(i-1)}=f^{(i)}$ from
(\ref{DefFiEq}) and, via Lemma \ref{Aprime_ALem}, $b^{(i-1)}=b^{(i)\prime}$.
The preceding expression thus also holds (trivially) in this case.

Note that
$\delta^{p^{(i-1)},p^{(i-1)\prime}}_{b^{(i-1)},f^{(i-1)}}=1$ only in
the case $i=\tau^R_{\eta^R+1}$, when
$\nu^{(i-1)}_{\eta^{R\prime}-1}=b^{(i-1)}-(-1)^{f^{(i-1)}}$.

Lemmas \ref{AttenGen2Lem} and Lemma \ref{ExtGen2Lem} also imply that:
\begin{displaymath}
\begin{array}{ll}
\alpha^{p^{(i-1)},p^{(i-1)\prime}}_{a^{(i)\prime},b^{(i-1)}}
&=\alpha_i+\beta_i-u_i^{R}\Delta^R_{\eta^{R\prime}};\\[2mm]
\beta^{p^{(i-1)},p^{(i-1)\prime}}_{a^{(i)\prime},b^{(i-1)},e^{(i)},f^{(i-1)}}
&=\beta_i-u_i^{R}\Delta^R_{\eta^{R\prime}};\\[2mm]
\left\lfloor\frac{b^{(i-1)}p^{(i-1)}}{p^{(i-1)\prime}}\right\rfloor
&=\left\{
     \begin{array}{ll}
	\tilde b^{(i)}+\Delta^R_{\eta^{R\prime}}
                     & \mbox{if }i=\tau^R_{\eta^R+1}
			\mbox{ and } f^{(i)}=f^{(i-1)};\\[1mm]
	\tilde b^{(i)}
                     & \mbox{otherwise.}
     \end{array}
     \right.
\end{array}
\end{displaymath}
They also imply that
$\muIndUpA{i-1}$, $\muIndDnA{i-1}$, $\nuIndUp{i-1}$,
$\nuIndDn{i-1}$ are a mazy-four sandwiching
$(a^{(i)\prime},b^{(i-1)})$.
That they are actually interfacial follows if
in the $i=\tau^R_{\eta^R+1}$ case, we show that
$\nu^{(i-1)*}_{\eta^R}$ and $\nu^{(i-1)}_{\eta^R}$ are interfacial.
{}From above, $\nu^{(i-1)*}_{\eta^R}=b^{(i)\prime}$ is interfacial.
That $\nu^{(i-1)}_{\eta^R}$ is interfacial will be established below.

We now check that we can apply Lemma \ref{AttenGen1Lem}
(for path truncation on the left) if $i=\sigma^L_{\eta^L}$,
and Lemma \ref{ExtGen1Lem} (for path extension on the left)
if $i=\tau^L_{\eta^L+1}$.

If $i=\sigma^L_{\eta^L}$ then $u^L_i=1$, and we see that
$\Delta^L_{\eta^L}=-(-1)^{e^{(i)}}=(-1)^{e^{(i-1)}}$ from (\ref{DefEiEq}),
and $a^{(i-1)}=a^{(i)\prime}-\Delta^L_{\eta^L}$ from Lemma \ref{Aprime_ALem}.
Also note that $\eta^{L\prime}=\eta^L$.
If $\eta^L>1$ or both $\eta^L=1$ and $\sigma^L_1<t_n$ then
Lemma \ref{Mu_ALem} implies that
$\mu^{(i-1)}_{\eta^L-1}=a^{(i)\prime}+\Delta^L_{\eta^L}$ and
$\mu^{(i-1)*}_{\eta^L-1}\ne a^{(i-1)}$.
Otherwise, when $i=\sigma^L_1>t_n$, then
$a^{(i)\prime}\in\{1,p^{(i-1)\prime}-1\}$ from the definition.

If $i=\tau^L_{\eta^L+1}$ (when necessarily $\eta^L<d^L$)
then $u^L_i=-1$ and $\eta^{L\prime}=\eta^L+1$.
We see that $\Delta^L_{\eta^L(i-1)}=-(-1)^{e(i-1)}$ from (\ref{DefEiEq})
and, via Lemma \ref{Aprime_ALem},
$a^{(i-1)}=a^{(i)\prime}+\Delta^L_{\eta^{L\prime}}$.
In addition, Lemma \ref{Mu_ALem} gives
$\mu^{(i-1)*}_{\eta^L}=a^{(i-1)}-\Delta^L_{\eta^{L\prime}}$ and
$\mu^{(i-1)}_{\eta^L}=a^{(i-1)}+\Delta^L_{\eta^{L\prime}}$.
The definition of a reduced run implies that this case, $i=\tau^L_{\eta^L+1}$
and $i\ne t_k$, only occurs when $t_{n-1}\le t_n-1=\tau^L_2=i<t_n<\sigma^L_1$.
The induction hypothesis then implies that $a^{(i)}$ is interfacial
in the $(p^{(i)},p^{(i)\prime})$-model and
$\delta^{p^{(i)},p^{(i)\prime}}_{a^{(i)},e^{(i)}}=0$.
Thereupon, by Lemma \ref{StartPtLem}, $a^{(i)\prime}$ is interfacial
in the $(p^{(i-1)},p^{(i-1)\prime})$-model.
Also note that $d^L_0(i)=1$ and so
$\muIndUpA{i-1}=\muIndDnA{i-1}=()$ here.

Then, on using Lemma \ref{AttenGen1Lem} if $i=\sigma^L_{\eta^L}$,
or Lemma \ref{ExtGen1Lem} if $i=\tau^L_{\eta^L+1}$ yields:

\begin{equation*}
\begin{split}
&\mchi^{p^{(i-1)},p^{(i-1)\prime}}_{a^{(i-1)},b^{(i-1)},e^{(i-1)},f^{(i-1)}}
    (m_{i-1},m_{i}+2\delta^{p^{(i-1)},p^{(i-1)\prime}}_{a^{(i-1)},e^{(i-1)}})
\left\{
\begin{matrix}
\muIndUp{i-1};\nuIndUp{i-1}\\
\muIndDn{i-1};\nuIndDn{i-1}
\end{matrix}
\right\}\\[1.5mm]
&\hskip15mm
=
q^{-\frac12u^L_i(m_{i-1}-m_i+u^L_i
      +\Delta^L_{\eta^{L\prime}}(\beta_i-u_i^{R}\Delta^R_{\eta^{R\prime}}) )}
\\[1.5mm]
&\hskip25mm
\times\quad
\mchi^{p^{(i-1)},p^{(i-1)\prime}}_{a^{(i)\prime},b^{(i-1)},e^{(i)},f^{(i-1)}}
                   (m_{i-1}+u^L_i,m_{i})
\left\{
\begin{matrix}
\muIndUpA{i-1};\nuIndUp{i-1}\\
\muIndDnA{i-1};\nuIndDn{i-1}
\end{matrix}
\right\},
\end{split}
\end{equation*}
noting that $d^L_0(i)=d^L_0(i-1)$ because $i\ne t_n$.
If neither $i=\sigma^L_{\eta^L}$ nor $i=\tau^L_{\eta^L+1}$ then
(noting that $i\ne t_{k}$) $u^L_i=0$, $e^{(i-1)}=e^{(i)}$ from (\ref{DefEiEq})
and, via Lemma \ref{Aprime_ALem}, $a^{(i-1)}=a^{(i)\prime}$.
The preceding expression thus also holds (trivially) in this case.

Note that
$\delta^{p^{(i-1)},p^{(i-1)\prime}}_{a^{(i-1)},e^{(i-1)}}=1$ only in
the case $i=\tau^L_{\eta^L+1}$, so that
$\mu^{(i-1)}_{\eta^{L\prime}-1}=a^{(i-1)}-(-1)^{e^{(i-1)}}$.

Lemmas \ref{AttenGen1Lem} and Lemma \ref{ExtGen1Lem} also imply that:
\begin{displaymath}
\begin{array}{ll}
\alpha^{p^{(i-1)},p^{(i-1)\prime}}_{a^{(i-1)},b^{(i-1)}}
&=\alpha_i+\beta_i
          -u_i^{R}\Delta^R_{\eta^{R\prime}}
          +u_i^{L}\Delta^L_{\eta^{L\prime}};\\[2mm]
\beta^{p^{(i-1)},p^{(i-1)\prime}}_{a^{(i-1)},b^{(i-1)},e^{(i-1)},f^{(i-1)}}
&=\beta_i-u_i^{R}\Delta^R_{\eta^{R\prime}}
          +u_i^{L}\Delta^L_{\eta^{L\prime}};\\[2mm]
\left\lfloor\frac{a^{(i-1)}p^{(i-1)}}{p^{(i-1)\prime}}\right\rfloor
&=\left\{
     \begin{array}{ll}
        \tilde a^{(i)}+\Delta^L_{\eta^{L\prime}}
                     & \mbox{if }i=\tau^L_{\eta^L+1}
                       \mbox{ and } e^{(i)}=e^{(i-1)};\\[1mm]
        \tilde a^{(i)}
                     & \mbox{otherwise.}
     \end{array}
     \right.
\end{array}
\end{displaymath}
They also imply that
$\muIndUp{i-1}$, $\muIndDn{i-1}$, $\nuIndUp{i-1}$,
$\nuIndDn{i-1}$ are a mazy-four sandwiching
$(a^{(i-1)},b^{(i-1)})$.
That they are actually interfacial follows if
in the $i=\tau^L_{\eta^L+1}$ case, we show that
$\mu^{(i-1)*}_{\eta^L}$ and $\mu^{(i-1)}_{\eta^L}$ are interfacial.
{}From above, $\mu^{(i-1)*}_{\eta^L}=a^{(i)\prime}$ is interfacial.
That $\mu^{(i-1)}_{\eta^L}$ is interfacial will be established below.

Combining all the above, and using the expression for $\gamma_{i-1}''$
given by (\ref{Const1Eq}) and (\ref{Const2Eq}), yields:
\begin{equation*}
\begin{split}
&\mchi^{p^{(i-1)},p^{(i-1)\prime}}_{a^{(i-1)},b^{(i-1)},e^{(i-1)},f^{(i-1)}}
    (m_{i-1},m_{i}+2\delta^{p^{(i-1)},p^{(i-1)\prime}}_{a^{(i-1)},e^{(i-1)}})
\left\{
\begin{matrix}
\muIndUp{i-1};\nuIndUp{i-1}\\
\muIndDn{i-1};\nuIndDn{i-1}
\end{matrix}
\right\}\\[1.5mm]
&\hskip15mm
=\sum
  q^{\frac{1}{4}\hat{\boldm}^{(i)T}\boldC\hat{\boldm}^{(i)}
   +\frac{1}{4} m_{i-1}^2
   -\frac{1}{2} m_{i-1}m_{i}
   -\frac{1}{2}(\boldu^{L}_{(\flat,k)}+\boldu^{R}_{(\sharp,k)})
   \cdot \boldm^{(i-1)}
  +\frac{1}{4}\gamma''_{i-1}}\\[1.5mm]
&\hskip55mm
\times\quad
  \prod_{j=i}^{t-1}
  \left[
  {m_j-\frac{1}{2}(\boldC^*\hat{\boldm}^{(i-1)}
                       \!-\!\boldu^{L}\!-\!\boldu^{R})_j\atop m_j}
  \right]_q\\[1.5mm]
&\hskip15mm
=F^{(i-1)}(\boldu^L,\boldu^R,m_{i-1},m_i),
\end{split}
\end{equation*}
which is the required result when $i\ne t_k$, since $k=k(i)=k(i-1)$.

In this $i\ne t_k$ case, making use of (\ref{Const1Eq}),
(\ref{Const2Eq}), and Lemma \ref{Aprime_AprimeLem}, we also immediately
obtain:
\begin{displaymath}
\begin{array}{ll}
\alpha^{p^{(i-1)},p^{(i-1)\prime}}_{a^{(i-1)},b^{(i-1)}}
&=\alpha_{i-1}'';\\[2mm]
\beta^{p^{(i-1)},p^{(i-1)\prime}}_{a^{(i-1)},b^{(i-1)},e^{(i-1)},f^{(i-1)}}
&=\beta_{i-1}';\\[2mm]
\left\lfloor\frac{a^{(i-1)}p^{(i-1)}}{p^{(i-1)\prime}}\right\rfloor
&=\tilde a^{(i-1)};\\[2mm]
\left\lfloor\frac{b^{(i-1)}p^{(i-1)}}{p^{(i-1)\prime}}\right\rfloor
&=\tilde b^{(i-1)}.
\end{array}
\end{displaymath}


Now consider the case for which $i=t_k$.
Equation (\ref{Const3Eq}) gives $\alpha_i=\alpha_i''$,
$\beta_i=\alpha_i''-\beta_i'$ and $\gamma_i=-\alpha_i^2-\gamma_i''$.
Then Lemma \ref{DPathParamLem} gives
$\alpha^{p^{(i)\prime}-p^{(i)},p^{(i)\prime}}_{a^{(i)},b^{(i)}}
=\alpha^{\prime\prime}_i=\alpha_i$ and
$\beta^{p^{(i)\prime}-p^{(i)},p^{(i)\prime}}_%
{a^{(i)},b^{(i)},1-e^{(i)},1-f^{(i)}}
=\alpha^{\prime\prime}_i-\beta^\prime_i=\beta_i$.
Let $m_{i-1}\equiv Q_{i-1}$.
On setting $M=m_{i-1}+u^L_i+u^R_i$, equations (\ref{MNEq1}),
(\ref{MNmatEq2}) and (\ref{ParityDef}) imply that $M-m_i\equiv Q_{i+1}$.
Then, use of the induction hypothesis and Lemmas \ref{MazyDijCor},
\ref{IndTakLem}, \ref{A_AprimeLem} and \ref{MuInterLem} yields:
\begin{equation*}
\begin{split}
&\mchi^{p^{(i-1)},p^{(i-1)\prime}}_{a^{(i)\prime},b^{(i)\prime},
                       1-e^{(i)},1-f^{(i)}} (M,m_{i};q)
\left\{
\begin{matrix}
\muIndUpA{i-1};\nuIndUpA{i-1}\\
\muIndDnA{i-1};\nuIndDnA{i-1}
\end{matrix}
\right\}\\[4mm]
&\hskip15mm
=\sum_{m_{i+1}\equiv Q_{i+1}}
q^{\frac{1}{4}(m_i^2+(M-m_i)^2-\alpha_i^{2}-\beta_i^{2})}
\left[{\frac{1}{2}(M+m_i-m_{i+1})}\atop m_i\right]_q\\[4mm]
&\hskip65mm
\times\quad
F^{(i)}(\boldu^{L},\boldu^{R},m_i,m_{i+1};q^{-1}).
\end{split}
\end{equation*}
Lemma \ref{MazyDijCor} also gives
$\alpha^{p^{(i-1)},p^{(i-1)\prime}}_{a^{(i)\prime},b^{(i)\prime}}
=2\alpha^{\prime\prime}_i-\beta^\prime_i=\alpha_i+\beta_i$, and
$\beta^{p^{(i-1)},p^{(i-1)\prime}}_{a^{(i)\prime},b^{(i)\prime},
1-e^{(i)},1-f^{(i)}}=\alpha^{\prime\prime}_i-\beta^\prime_i=\beta_i$, and
Lemmas \ref{A_AprimeLem} and \ref{BDParamLem} imply that
$\lfloor a^{(i)\prime}p^{(i-1)}/p^{(i-1)\prime}\rfloor=
a^{(i)}-1-\tilde a^{(i)}$
and
$\lfloor b^{(i)\prime}p^{(i-1)}/p^{(i-1)\prime}\rfloor=
b^{(i)}-1-\tilde b^{(i)}$.
Lemma \ref{BDParamLem} also implies that
$\delta^{p^{(i-1)},p^{(i-1)\prime}}_{a^{(i)\prime},1-e^{(i)}}=
\delta^{p^{(i-1)},p^{(i-1)\prime}}_{b^{(i)\prime},1-f^{(i)}}=0$.
Lemma \ref{MazyDijCor} also implies that
$\muIndUpA{i-1}$, $\muIndDnA{i-1}$, $\nuIndUpA{i-1}$,
$\nuIndDnA{i-1}$ are an interfacial mazy-four in the
$(p^{(i-1)},p^{(i-1)\prime})$-model sandwiching
$(a^{(i)\prime},b^{(i)\prime})$.
Together with Lemma \ref{InterLem}, Lemma \ref{MuInterLem} also
shows that $\rho^{(i-1)}(\mu^{(i-1)*}_j)=\tilde\mu^{(i-1)*}_j$
and $\rho^{(i-1)}(\mu^{(i-1)}_j)=\tilde\mu^{(i-1)}_j$
for $d^L_0(i)\le j<\eta^L(i)$,
and that $\rho^{(i-1)}(\nu^{(i-1)*}_j)=\tilde\nu^{(i-1)*}_j$
and $\rho^{(i-1)}(\nu^{(i-1)}_j)=\tilde\nu^{(i-1)}_j$
for $d^R_0(i)\le j<\eta^R(i)$.

Since $M=m_{i-1}+u^L_i+u^R_i$, on noting that $i=t_k$, we have:
\begin{displaymath}
M+m_i-m_{i+1}=2m_i-(\boldC^*\hat{\boldm}^{(i-1)}-\boldu^{L}-\boldu^{R})_i
\end{displaymath}
(in the case $i=t-1$, we require this expression after substituting
$m_t=0$),
and
\begin{displaymath}
\begin{array}{l}
\displaystyle
\hat{\boldm}^{(i+1)T}\boldC\hat{\boldm}^{(i+1)}+(M-m_i)^2\\
\qquad=
\hat{\boldm}^{(i)T}\boldC\hat{\boldm}^{(i)}+M^2-2Mm_i\\
\qquad=
\hat{\boldm}^{(i)T}\boldC\hat{\boldm}^{(i)}+m_{i-1}^2-2m_im_{i-1}
+2(m_{i-1}-m_i)(u^L_i+u^R_i)+(u^L_i+u^R_i)^2.
\end{array}
\end{displaymath}
Use of expression (\ref{Invq1Eq}) or (\ref{Invq2Eq}) for
$F^{(i)}_{a,b}(\boldu^{L},\boldu^{R},m_i,m_{i+1};q^{-1})$ then gives:
\begin{equation*}
\begin{split}
&\mchi^{p^{(i-1)},p^{(i-1)\prime}}_{a^{(i)\prime},b^{(i)\prime},
                       1-e^{(i)},1-f^{(i)}} (m_{i-1}+u^L_i+u^R_i,m_{i})
\left\{
\begin{matrix}
\muIndUpA{i-1};\nuIndUpA{i-1}\\
\muIndDnA{i-1};\nuIndDnA{i-1}
\end{matrix}
\right\}\\[1.5mm]
&\hskip5mm
=\sum
q^{\frac{1}{4}\hat{\boldm}^{(i)T}\boldC\hat{\boldm}^{(i)}
  +\frac{1}{4}m_{i-1}^2
  -\frac{1}{2}m_im_{i-1}
  +\frac{1}{2}(m_{i-1}-m_i)(u^L_i+u^R_i)
  +\frac{1}{4}(u^L_i+u^R_i)^2
  +\frac{1}{4}\gamma_i-\frac{1}{4}\beta_i^2}\\[1.5mm]
&\hskip20mm
\times\quad
  q^{-\frac{1}{2}(\boldu^{L}_{(\flat,k-1)}+\boldu^{R}_{(\sharp,k-1)})
    \cdot \boldm^{(i)}}
  \prod_{j=i}^{t-1}
  \left[
  {m_j-\frac{1}{2}(\boldC^*\hat{\boldm}^{(i-1)}
                       \!-\!\boldu^{L}\!-\!\boldu^{R})_j\atop m_j}
  \right]_q,
\end{split}
\end{equation*}
where the sum is over all
$(m_{i+1},m_{i+2},\ldots,m_{t-1})\equiv(Q_{i+1},Q_{i+2},\ldots,Q_{t-1})$.

Since $i=t_k$, it follows that $u^R_i=-1$ if
$i=\tau^R_{\eta^R+1}$ or
$\sigma^R_{\eta^R}< i<\tau^R_{\eta^R}$.
We will apply Lemma \ref{ExtGen2Lem} (for path extension on the right)
in these cases.

In the case $i=\tau^R_{\eta^R+1}$, we have $i+1\ne\tau^R_{\eta^R}$ which
implies via the induction hypothesis that
$\delta^{p^{(i)},p^{(i)\prime}}_{b^{(i)},f^{(i)}}=0$.
Then Lemma \ref{BDParamLem} implies that
$b^{(i)\prime}$ is interfacial
in the $(p^{(i-1)},p^{(i-1)\prime})$-model.
Note that if $\eta^R=2$ and $\sigma^R_1\ge t_n$ then
$\nuIndUpA{i-1}=\nuIndDnA{i-1}=()$;
if $\eta^R=2$ and $\sigma^R_1<t_n$ (necessarily $k<n$) then
$\nuIndUpA{i-1}=(\nu^{(i-1)}_0)$ and $\nuIndDnA{i-1}=(\nu^{(i-1)*}_0)$,
with $\nu^{(i-1)}_0\ne b^{(i-1)}\ne\nu^{(i-1)*}_0$ direct from the
definitions after noting that $i<\sigma^R_1<t_n$;
and if $\eta^R=j+1$ for $j>1$ then
$\nuIndUpA{i-1}=(\ldots,\nu^{(i-1)}_{j-1})$
and $\nuIndDnA{i-1}=(\ldots,\nu^{(i-1)*}_{j-1})$,
with $\nu^{(i-1)}_{j-1}\ne b^{(i-1)}\ne\nu^{(i-1)*}_{j-1}$ direct from the
definitions after noting that $i<\sigma^R_j<\tau^R_j$.
We now proceed via Lemma \ref{ExtGen2Lem} as in the $i\ne t_k$ case,
to obtain:
\begin{equation*}
\begin{split}
&\mchi^{p^{(i-1)},p^{(i-1)\prime}}_{a^{(i)\prime},b^{(i-1)},
                   1-e^{(i)},f^{(i-1)}} (m_{i-1}+u^L_i,m_{i})
\left\{
\begin{matrix}
\muIndUpA{i-1};\nuIndUp{i-1}\\
\muIndDnA{i-1};\nuIndDn{i-1}
\end{matrix}
\right\}\\[1.5mm]
&\hskip10mm
=
q^{-\frac12u^R_i(m_{i-1}+u^L_i+u^R_i
   -\Delta^R_{\eta^{R\prime}}(\alpha_i+\beta_i))}
\\[1.5mm]
&\hskip20mm
\times\quad
\mchi^{p^{(i-1)},p^{(i-1)\prime}}_{a^{(i)\prime},b^{(i)\prime},
                   1-e^{(i)},1-f^{(i)}} (m_{i-1}+u^L_i+u^R_i,m_{i})
\left\{
\begin{matrix}
\muIndUpA{i-1};\nuIndUpA{i-1}\\
\muIndDnA{i-1};\nuIndDnA{i-1}
\end{matrix}
\right\},
\end{split}
\end{equation*}
having noted that $d^R_0(i)=d^R_0(i-1)$.
As in the $i\ne t_k$ case, we obtain that
$\muIndUpA{i-1}$, $\muIndDnA{i-1}$, $\nuIndUp{i-1}$,
$\nuIndDn{i-1}$ are a mazy-four sandwiching
$(a^{(i)\prime},b^{(i-1)})$, and that they are actually interfacial
once it is established that $\nu^{(i-1)}_{\eta^R}$ is interfacial.
This is done below.

In the case $\sigma^R_{\eta^R}< i<\tau^R_{\eta^R}$,
we have $\eta^{R\prime}=\eta^R$ and
$\Delta^R_{\eta^{R}}=-(-1)^{f^{(i-1)}}=(-1)^{1-f^{(i)}}$
from (\ref{DefFiEq}) and, via Lemma \ref{Aprime_ALem},
$b^{(i-1)}=b^{(i)\prime}+\Delta^R_{\eta^{R}}$.
Note that if $k=n$ (necessarily $\eta^R=1$) then
$\nuIndUpA{i-1}=\nuIndDnA{i-1}=()$;
if $k<n$ and $\eta^R=1$ then
$\nuIndUpA{i-1}=(\nu^{(i-1)}_0)$ and $\nuIndDnA{i-1}=(\nu^{(i-1)*}_0)$,
with $\nu^{(i-1)}_0\ne b^{(i-1)}\ne\nu^{(i-1)*}_0$ direct from the
definitions after noting that $i<t_n$;
and if $k<n$ and $\eta^R=j$ for $j>1$ then
$\nuIndUpA{i-1}=(\ldots,\nu^{(i-1)}_{j-1})$
and $\nuIndDnA{i-1}=(\ldots,\nu^{(i-1)*}_{j-1})$,
with $\nu^{(i-1)}_{j-1}\ne b^{(i-1)}\ne\nu^{(i-1)*}_{j-1}$ direct from the
definitions after noting that $i<\tau^R_j$.
Then Lemma \ref{ExtGen2Lem} yields:
\begin{equation*}
\begin{split}
&\mchi^{p^{(i-1)},p^{(i-1)\prime}}_{a^{(i)\prime},b^{(i-1)},
                   1-e^{(i)},f^{(i-1)}} (m_{i-1}+u^L_i,m_{i})
\left\{
\begin{matrix}
\muIndUpA{i-1};\nuIndUpA{i-1},b^{(i-1)}+\Delta^R_{\eta^R}\\
\muIndDnA{i-1};\nuIndDnA{i-1},b^{(i-1)}-\Delta^R_{\eta^R}
\end{matrix}
\right\}\\[1.5mm]
&\hskip5mm
=
q^{-\frac12u^R_i(m_{i-1}+u^L_i+u^R_i
   -\Delta^R_{\eta^{R\prime}}(\alpha_i+\beta_i))}
\\[1.5mm]
&\hskip15mm
\times\quad
\mchi^{p^{(i-1)},p^{(i-1)\prime}}_{a^{(i)\prime},b^{(i)\prime},
                     1-e^{(i)},1-f^{(i)}} (m_{i-1}+u^L_i+u^R_i,m_{i})
\left\{
\begin{matrix}
\muIndUpA{i-1};\nuIndUpA{i-1}\\
\muIndDnA{i-1};\nuIndDnA{i-1}
\end{matrix}
\right\}.
\end{split}
\end{equation*}
Here, if $i=t_n$ and $\sigma^R_1<t_n$, we have
$d^R_0(i)=1$ and $d^R_0(i-1)=0$.
Also note that $\eta^{R\prime}=\eta^R=1$.
Then $\nu^{(i-1)}_0=b^{(i-1)}+\Delta^R_1$ and
$\nu^{(i-1)*}_0=b^{(i-1)}-\Delta^R_1$ by Lemma \ref{Mu_ALem},
whereupon we obtain precisely the expression obtained above in
the $i=\tau^R_{\eta^R+1}$ case.
Also note that $i+1\ne\tau^R_{1}$ again implies,
via $\delta^{p^{(i)},p^{(i)\prime}}_{b^{(i)},f^{(i)}}=0$
and Lemma \ref{BDParamLem}, that $b^{(i)\prime}$ is interfacial
in the $(p^{(i-1)},p^{(i-1)\prime})$-model.
Lemma \ref{ExtGen2Lem} shows that
$\muIndUpA{i-1}$, $\muIndDnA{i-1}$, $\nuIndUp{i-1}$,
$\nuIndDn{i-1}$ are a mazy-four sandwiching
$(a^{(i)\prime},b^{(i-1)})$.
Since $\nu^{(i-1)*}_0=b^{(i)\prime}$ is interfacial,
to show that $\muIndUpA{i-1}$, $\muIndDnA{i-1}$, $\nuIndUp{i-1}$,
$\nuIndDn{i-1}$ are an interfacial mazy-four requires only that it be
established that $\nu^{(i-1)}_{0}$ is interfacial.
This is done below.

Otherwise for $\sigma^R_{\eta^R}< i<\tau^R_{\eta^R}$, we have
$\nu^{(i-1)}_{\eta^R-1}=b^{(i-1)}+\Delta^R_{\eta^R}$ by Lemma \ref{Mu_ALem},
and $\eta^{R\prime}=\eta^R$,
whence use of Lemma \ref{CutParam1Lem} again yields precisely
the expression obtained above in the $i=\tau^R_{\eta^R+1}$ case.

In addition, the same expression clearly also holds in the case
$\tau^R_{\eta^R+1}< i\le\sigma^R_{\eta^R}$, for which
$u^R_i=0$ and $\eta^{R\prime}=\eta^R$
(in the $i=\sigma^R_{\eta^R}$ case, note that $k(i-1)=k-1=k^R(i-1)$
and consequently $f^{(i-1)}=1-f^{(i)}$ by (\ref{DefFiEq})).

Note that in the case
$\sigma^R_{\eta^R}< i<\tau^R_{\eta^R}$,
we have $f^{(i-1)}=f^{(i)}$ from (\ref{DefFiEq}).
Lemma \ref{ExtGen2Lem} then implies that
$\delta^{p^{(i-1)},p^{(i-1)\prime}}_{b^{(i-1)},f^{(i-1)}}=0$ in
this case.
Consequently,
$\delta^{p^{(i-1)},p^{(i-1)\prime}}_{b^{(i-1)},f^{(i-1)}}=1$ only in
the case $i=\tau^R_{\eta^R+1}$, and then
only if $\nu^{(i-1)}_{\eta^R-1}=b^{(i-1)}-(-1)^{f^{(i-1)}}$.

Lemma \ref{ExtGen2Lem} also implies that:
\begin{equation*}
\begin{split}
&\alpha^{p^{(i-1)},p^{(i-1)\prime}}_{a^{(i)\prime},b^{(i-1)}
\phantom{,e^{(i)},f^{(i-1)}}}
=\alpha_i+\beta_i-u_i^{R}\Delta^R_{\eta^{R\prime}};\\[2mm]
&\beta^{p^{(i-1)},p^{(i-1)\prime}}_{a^{(i)\prime},b^{(i-1)},e^{(i)},f^{(i-1)}}
=\beta_i-u_i^{R}\Delta^R_{\eta^{R\prime}};\\[2mm]
&\left\lfloor\frac{b^{(i-1)}p^{(i-1)}}{p^{(i-1)\prime}}\right\rfloor
=\left\{
     \begin{array}{ll}
        b^{(i)}-1-\tilde b^{(i)}+\Delta^R_{\eta^{R\prime}}
                     & \mbox{if }i=\tau^R_{\eta^R+1}
                       \mbox{ and } f^{(i)}\ne f^{(i-1)};\\[1mm]
        b^{(i)}-1-\tilde b^{(i)}
                     & \mbox{otherwise.}
     \end{array}
     \right.
\end{split}
\end{equation*}

Since $i=t_k$, it follows that $u^L_i=-1$ if
$i=\tau^L_{\eta^L+1}$ or
$\sigma^L_{\eta^L}< i<\tau^L_{\eta^L}$.
We will apply Lemma \ref{ExtGen1Lem} (for path extension on the left)
in these cases.

In the case $i=\tau^L_{\eta^L+1}$, we have $i+1\ne\tau^L_{\eta^L}$ which
implies via the induction hypothesis that
$\delta^{p^{(i)},p^{(i)\prime}}_{a^{(i)},e^{(i)}}=0$.
Then Lemma \ref{BDParamLem} implies that
$a^{(i)\prime}$ is interfacial
in the $(p^{(i-1)},p^{(i-1)\prime})$-model.
Note that if $\eta^L=2$ and $\sigma^L_1\ge t_n$ then
$\muIndUpA{i-1}=\muIndDnA{i-1}=()$;
if $\eta^L=2$ and $\sigma^L_1<t_n$ (necessarily $k<n$) then
$\muIndUpA{i-1}=(\mu^{(i-1)}_0)$ and $\muIndDnA{i-1}=(\mu^{(i-1)*}_0)$,
with $\mu^{(i-1)}_0\ne a^{(i-1)}\ne\mu^{(i-1)*}_0$ direct from the
definitions after noting that $i<\sigma^L_1<t_n$;
and if $\eta^L=j+1$ for $j>1$ then
$\muIndUpA{i-1}=(\ldots,\mu^{(i-1)}_{j-1})$
and $\muIndDnA{i-1}=(\ldots,\mu^{(i-1)*}_{j-1})$,
with $\mu^{(i-1)}_{j-1}\ne a^{(i-1)}\ne\mu^{(i-1)*}_{j-1}$ direct
from the definitions after noting that $i<\sigma^L_j<\tau^L_j$.
We now proceed via Lemma \ref{ExtGen1Lem} as in the $i\ne t_k$ case,
to obtain:
\begin{equation*}
\begin{split}
&\mchi^{p^{(i-1)},p^{(i-1)\prime}}_{a^{(i-1)},b^{(i-1)},e^{(i-1)},f^{(i-1)}}
    (m_{i-1},m_{i}+2\delta^{p^{(i-1)},p^{(i-1)\prime}}_{a^{(i-1)},e^{(i-1)}})
\left\{
\begin{matrix}
\muIndUp{i-1};\nuIndUp{i-1}\\
\muIndDn{i-1};\nuIndDn{i-1}
\end{matrix}
\right\}\\[1.5mm]
&\hskip15mm
=
q^{-\frac12u^L_i(m_{i-1}-m_i+u^L_i
      +\Delta^L_{\eta^{L\prime}}(\beta_i-u_i^{R}\Delta^R_{\eta^{R\prime}}) )}
\\[1.5mm]
&\hskip25mm
\times\quad
\mchi^{p^{(i-1)},p^{(i-1)\prime}}_{a^{(i)\prime},b^{(i-1)},
                   1-e^{(i)},f^{(i-1)}} (m_{i-1}+u^L_i,m_{i})
\left\{
\begin{matrix}
\muIndUpA{i-1};\nuIndUp{i-1}\\
\muIndDnA{i-1};\nuIndDn{i-1}
\end{matrix}
\right\},
\end{split}
\end{equation*}
after noting that $d^L_0(i)=d^L_0(i-1)$.
As in the $i\ne t_k$ case, we obtain that
$\muIndUp{i-1}$, $\muIndDn{i-1}$, $\nuIndUp{i-1}$,
$\nuIndDn{i-1}$ are an interfacial mazy-four sandwiching
$(a^{(i-1)},b^{(i-1)})$, and that they are actually interfacial
once it is established that $\mu^{(i-1)}_{\eta^L}$ is interfacial.
This is done below.

In the case $\sigma^L_{\eta^L}< i<\tau^L_{\eta^L}$,
we have $\eta^{L\prime}=\eta^L$ and
$\Delta^L_{\eta^{L}}=-(-1)^{e^{(i-1)}}=(-1)^{1-e^{(i)}}$
from (\ref{DefEiEq}), and via Lemma \ref{Aprime_ALem},
$a^{(i-1)}=a^{(i)\prime}+\Delta^L_{\eta^{L}}$.
Note that if $k=n$ (necessarily $\eta^L=1$) then
$\muIndUpA{i-1}=\muIndDnA{i-1}=()$;
if $k<n$ and $\eta^L=1$ then
$\muIndUpA{i-1}=(\mu^{(i-1)}_0)$ and $\muIndDnA{i-1}=(\mu^{(i-1)*}_0)$,
with $\mu^{(i-1)}_0\ne a^{(i-1)}\ne\mu^{(i-1)*}_0$ direct from the
definitions after noting that $i<t_n$;
and if $k<n$ and $\eta^L=j$ for $j>1$ then
$\muIndUpA{i-1}=(\ldots,\mu^{(i-1)}_{j-1})$
and $\muIndDnA{i-1}=(\ldots,\mu^{(i-1)*}_{j-1})$,
with $\mu^{(i-1)}_{j-1}\ne a^{(i-1)}\ne\mu^{(i-1)*}_{j-1}$ direct from the
definitions after noting that $i<\tau^L_j$.
Then Lemma \ref{ExtGen1Lem} yields:
\begin{equation*}
\begin{split}
&\mchi^{p^{(i-1)},p^{(i-1)\prime}}_{a^{(i-1)},b^{(i-1)},e^{(i-1)},f^{(i-1)}}
    (m_{i-1},m_{i}+2\delta^{p^{(i-1)},p^{(i-1)\prime}}_{a^{(i-1)},e^{(i-1)}})
\left\{
\begin{matrix}
\muIndUpA{i-1},a^{(i-1)}\!+\!\Delta^L_{\eta^L};\nuIndUp{i-1}\\
\muIndDnA{i-1},a^{(i-1)}\!-\!\Delta^L_{\eta^L};\nuIndDn{i-1}
\end{matrix}
\right\}\\[1.5mm]
&\hskip15mm
=
q^{-\frac12u^L_i(m_{i-1}-m_i+u^L_i
      +\Delta^L_{\eta^{L\prime}}(\beta_i-u_i^{R}\Delta^R_{\eta^{R\prime}}) )}
\\[1.5mm]
&\hskip25mm
\times\quad
\mchi^{p^{(i-1)},p^{(i-1)\prime}}_{a^{(i)},b^{(i-1)},1-e^{(i)},f^{(i-1)}}
                   (m_{i-1}+u^L_i,m_{i})
\left\{
\begin{matrix}
\muIndUpA{i-1};\nuIndUp{i-1}\\
\muIndDnA{i-1};\nuIndDn{i-1}
\end{matrix}
\right\}.
\end{split}
\end{equation*}
Here, if $i=t_n$ and $\sigma^L_1<t_n$, we have
$d^L_0(i)=1$ and $d^L_0(i-1)=0$.
Also note that $\eta^{L\prime}=\eta^L=1$.
Then $\mu^{(i-1)}_0=a^{(i-1)}+\Delta^L_1$ and
$\mu^{(i-1)*}_0=a^{(i-1)}-\Delta^L_1$ by Lemma \ref{Mu_ALem},
whereupon we obtain precisely the expression 
obtained above in the $i=\tau^L_{\eta^L+1}$ case.
Also note that $i+1\ne\tau^L_{1}$ again implies,
via $\delta^{p^{(i)},p^{(i)\prime}}_{a^{(i)},e^{(i)}}=0$
and Lemma \ref{BDParamLem}, that $a^{(i)\prime}$ is interfacial
in the $(p^{(i-1)},p^{(i-1)\prime})$-model.
Lemma \ref{ExtGen1Lem} shows that
$\muIndUp{i-1}$, $\muIndDn{i-1}$, $\nuIndUp{i-1}$,
$\nuIndDn{i-1}$ are a mazy-four sandwiching
$(a^{(i-1)},b^{(i-1)})$.
Since $\mu^{(i-1)*}_0=a^{(i)\prime}$ is interfacial,
to show that $\muIndUp{i-1}$, $\muIndDn{i-1}$, $\nuIndUp{i-1}$,
$\nuIndDn{i-1}$ are an interfacial mazy-four requires only that it be
established that $\mu^{(i-1)}_{0}$ is interfacial.
This is done below.

Otherwise for $\sigma^L_{\eta^L}< i<\tau^L_{\eta^L}$, we have
$\mu^{(i-1)}_{\eta^L-1}=a^{(i-1)}+\Delta^L_{\eta^L}$ by Lemma
\ref{Mu_ALem}, and $\eta^{L\prime}=\eta^L$,
whence use of Lemma \ref{CutParam1Lem} again yields precisely
the expression obtained above in the $i=\tau^L_{\eta^L+1}$ case.

In addition, the same expression clearly also holds in the case
$\tau^L_{\eta^L+1}< i\le\sigma^L_{\eta^L}$, for which
$u^L_i=0$ and $\eta^{L\prime}=\eta^L$
(in the $i=\sigma^L_{\eta^L}$ case, note that $k(i-1)=k-1=k^L(i-1)$
and consequently $e^{(i-1)}=1-e^{(i)}$ by (\ref{DefEiEq})).

Note that in the case $\sigma^L_{\eta^L}< i<\tau^L_{\eta^L}$,
we have $e^{(i-1)}=e^{(i)}$ from (\ref{DefEiEq}).
Lemma \ref{ExtGen1Lem} then implies that
$\delta^{p^{(i-1)},p^{(i-1)\prime}}_{a^{(i-1)},e^{(i-1)}}=0$ in this case.
Consequently,
$\delta^{p^{(i-1)},p^{(i-1)\prime}}_{a^{(i-1)},e^{(i-1)}}=1$ only in
the case $i=\tau^L_{\eta^L+1}$, and then
only if $\mu^{(i-1)}_{\eta^L-1}=a^{(i-1)}-(-1)^{e^{(i-1)}}$.

Lemma \ref{ExtGen1Lem} also implies that:
\begin{equation*}
\begin{split}
&\alpha^{p^{(i-1)},p^{(i-1)\prime}}_{a^{(i-1)},b^{(i-1)}
\phantom{,e^{(i-1)},f^{(i-1)}}}
=\alpha_i+\beta_i
          -u_i^{R}\Delta^R_{\eta^{R\prime}}
          +u_i^{L}\Delta^L_{\eta^{L\prime}};\\[2mm]
&\beta^{p^{(i-1)},p^{(i-1)\prime}}_{a^{(i-1)},b^{(i-1)},e^{(i-1)},f^{(i-1)}}
=\beta_i-u_i^{R}\Delta^R_{\eta^{R\prime}}
          +u_i^{L}\Delta^L_{\eta^{L\prime}};\\[2mm]
&\left\lfloor\frac{a^{(i-1)}p^{(i-1)}}{p^{(i-1)\prime}}\right\rfloor
=\left\{
     \begin{array}{ll}
        a^{(i)}-1-\tilde a^{(i)}+\Delta^L_{\eta^{L\prime}}
                     & \mbox{if }i=\tau^L_{\eta^L+1}
                       \mbox{ and } e^{(i)}\ne e^{(i-1)};\\[1mm]
        a^{(i)}-1-\tilde a^{(i)}
                     & \mbox{otherwise.}
     \end{array}
     \right.
\end{split}
\end{equation*}

Combining all the above cases for $i=t_k$ gives:

\begin{equation*}
\begin{split}
&\mchi^{p^{(i-1)},p^{(i-1)\prime}}_{a^{(i-1)},b^{(i-1)},e^{(i-1)},f^{(i-1)}}
    (m_{i-1},m_{i}+2\delta^{p^{(i-1)},p^{(i-1)\prime}}_{a^{(i-1)},e^{(i-1)}})
\left\{
\begin{matrix}
\muIndUp{i-1};\nuIndUp{i-1}\\
\muIndDn{i-1};\nuIndDn{i-1}
\end{matrix}
\right\}\\[5mm]
&\hskip15mm
=\sum
  q^{\frac{1}{4}\hat{\boldm}^{(i)T}\boldC\hat{\boldm}^{(i)}
   +\frac{1}{4} m_{i-1}^2
   -\frac{1}{2} m_{i-1}m_{i}
   -\frac{1}{2}(\boldu^{L}_{(\flat,k-1)}+\boldu^{R}_{(\sharp,k-1)})
   \cdot\boldm^{(i-1)}
  +\frac{1}{4}\gamma''_{i-1}}\\[1.5mm]
&\hskip55mm
\times\quad
  \prod_{j=i}^{t-1}
  \left[
  {m_j-\frac{1}{2}(\boldC^*\hat{\boldm}^{(i-1)}
                       \!-\!\boldu^{L}\!-\!\boldu^{R})_j\atop m_j}
  \right]_q\\[1.5mm]
&\hskip15mm
=F^{(i-1)}(\boldu^L,\boldu^R,m_{i-1},m_i),
\end{split}
\end{equation*}
which is the required result when $i=t_k$, since $k(i-1)=k-1$.

In this $i=t_k$ case, making use of (\ref{Const1Eq}),
(\ref{Const2Eq}), and Lemma \ref{Aprime_AprimeLem}, we also immediately
obtain:
\begin{displaymath}
\begin{array}{ll}
\alpha^{p^{(i-1)},p^{(i-1)\prime}}_{a^{(i-1)},b^{(i-1)}}
&=\alpha_{i-1}'';\\[2mm]
\beta^{p^{(i-1)},p^{(i-1)\prime}}_{a^{(i-1)},b^{(i-1)},e^{(i-1)},f^{(i-1)}}
&=\beta_{i-1}';\\[2mm]
\left\lfloor\frac{a^{(i-1)}p^{(i-1)}}{p^{(i-1)\prime}}\right\rfloor
&=\tilde a^{(i-1)};\\[2mm]
\left\lfloor\frac{b^{(i-1)}p^{(i-1)}}{p^{(i-1)\prime}}\right\rfloor
&=\tilde b^{(i-1)}.
\end{array}
\end{displaymath}


It remains to consider the four final bulleted items of the
induction statement. 
First consider $t_k\le i-1$ and $\tau^L_{j+1}\le i-1<t_{k+1}<\sigma^L_j$.
If $i<t_{k+1}$ then $a^{(i)}$ is interfacial in
the $(p^{(i)},p^{(i)\prime})$-model by the induction hypothesis.
Then Lemma \ref{StartPtLem} implies that $a^{(i)\prime}$ is
interfacial in the $(p^{(i-1)},p^{(i-1)\prime})$-model.
If $i=t_{k+1}$ first note that since $i\ne\tau^L_j$, the
induction hypothesis yields
$\delta^{p^{(i)},p^{(i)\prime}}_{a^{(i)},e^{(i)}}=0$.
Then Lemma \ref{BDParamLem} implies that $a^{(i)\prime}$
is interfacial in the $(p^{(i-1)},p^{(i-1)\prime})$-model.
In both cases, $a^{(i-1)}=a^{(i)\prime}$ by Lemma \ref{Aprime_ALem}
and thus $a^{(i-1)}$ is interfacial in the
$(p^{(i-1)},p^{(i-1)\prime})$-model as required.

Now consider $t_k\le i-1$ and $\tau^L_{j+1}\le i-1<t_{k+1}=\sigma^L_j$
with $a^{(i-1)}$ not interfacial in the $(p^{(i-1)},p^{(i-1)\prime})$-model.
$a^{(i-1)}=a^{(i)\prime}$ by Lemma \ref{Aprime_ALem}.
In the case $i<t_{k+1}$ note that $i+1\ne\tau^L_j$,
whereupon the induction hypothesis yields
$\delta^{p^{(i)},p^{(i)\prime}}_{a^{(i)},e^{(i)}}=0$,
and then Lemma \ref{StartPtLem} implies that $a^{(i)}$ is not
interfacial in the $(p^{(i)},p^{(i)\prime})$-model.
The induction hypothesis then implies that $\tau^L_j=\sigma^L_j+1$
as required.
In the case $i=t_{k+1}$ ($ {}=\sigma^L_j$), Lemma \ref{BDParamLem}
implies that $\delta^{p^{(i)},p^{(i)\prime}}_{a^{(i)},e^{(i)}}=1$.
The induction hypothesis then implies that $i+1=\tau^L_j$ as required.

An entirely analogous argument yields the next bulleted item.
 

We now tackle the remaining cases of the final two bulleted items together
(we only consider the cases with superscript \lq L\rq: those with
superscript \lq R\rq\ follow similarly).
The remaining cases arise when $d^L_0(i-1)<d^L_0(i)$ and
$\eta^L(i-1)>\eta^L(i)$.
In the former of these cases, necessarily $\sigma^L_1<i=t_n$ when
$d^L_0(i-1)=0$ and $d^L_0(i)=1$, and thus here we set $j=0$.
In the latter, necessarily $i=\tau^L_{\eta^L+1}$ when
$\eta^L(i-1)=\eta^L+1$ and $\eta^L(i)=\eta^L$, and thus here we set
$j=\eta^L$ (${}>0$).
In either case, the above analysis has shown that
$\mu^{(i-1)*}_j=a^{(i-1)}-\Delta^L_j=a^{(i)\prime}$ is interfacial.
It remains to show that
$\mu^{(i-1)}_j=a^{(i)\prime}+2\Delta^L_{j+1}$ is interfacial and that
$\rho^{(i-1)}(\mu_j^{(i-1)*})=\tilde\mu^{(i-1)*}_j$
and $\rho^{(i-1)}(\mu_j^{(i-1)})=\tilde\mu^{(i-1)}_j$.

In the $j=0$ case (when $\sigma^L_1<i=t_n$), the definitions
directly give $\mu_0^{(i-1)}=0$ if $\Delta^L_1=-1$ and
$\mu_0^{(i-1)}=p^{(i-1)\prime}$ if $\Delta^L_1=1$,
so that $\mu_0^{(i-1)}$ is certainly interfacial.
We then immediately have
$\rho^{(i-1)}(\mu_0^{(i-1)})=0=\tilde\mu^{(i-1)}_0$ if $\Delta^L_1=-1$
and $\rho^{(i-1)}(\mu_0^{(i-1)})=p^{(i-1)}=\tilde\mu^{(i-1)}_0$
if $\Delta^L_1=1$, as required.
It follows from Note \ref{DireNote} that
$\rho^{(i-1)}(\mu_0^{(i-1)*})=1=\tilde\mu^{(i-1)*}_0$ if $\Delta^L_1=-1$
and $\rho^{(i-1)}(\mu_0^{(i-1)*})=p^{(i-1)}-1=\tilde\mu^{(i-1)*}_0$
if $\Delta^L_1=1$, as required.

For the $j>0$ cases, if we can show that
$\rho^{(i-1)}(\mu_j^{(i-1)*})=\tilde\mu^{(i-1)*}_j$
and that $a^{(i)\prime}+2\Delta^L_{j+1}$ is interfacial then,
via Note \ref{DireNote}, it will follow that
$\rho^{(i-1)}(\mu_j^{(i-1)})=\rho^{(i-1)}(\mu_j^{(i-1)*})+\Delta^L_{j+1}
=\tilde\mu^{(i-1)*}_j+\Delta^L_{j+1}=\tilde\mu^{(i-1)}_j$,
the final equality following from the definition of $\tilde\mu^{(i-1)}_j$.
To show this, we consider various cases.

With $i=\tau^L_{j+1}$, that
$\{\tau^L_J,\sigma^L_J,\Delta^L_J\}_{J=1}^{d^L}$ is a reduced run
implies that $i=t_k$ unless $i=t_{n}-1\ne t_{n-1}$.

In the $i=t_{n}-1\ne t_{n-1}$ case, we have $j=1$, $t_n<\sigma^L_1<t$,
$\Delta^L_2\ne\Delta^L_1$ and
$p^{(i-1)}/p^{(i-1)\prime}=3+1/c_n$.
For $1\le r<c_n$, the $r$th odd band in the
$(p^{(i-1)},p^{(i-1)\prime})$-model lies between heights $3r$ and $3r+1$.
We readily obtain $\kappa_{r+2}^{(i-1)}=3r+1$ and
$\tkappa_{r+2}^{(i-1)}=r$.
Then, if $\Delta^L_1=-1$, we have
$a^{(i)\prime}=\kappa^{(i-1)}_{\sigma^L_1-i+1}=3(\sigma^L_1-t_n)+1$
and therefore
$\rho^{(i-1)}(\mu^{(i-1)*}_1)=\sigma^L_1-t_n=
\tkappa^{(i-1)}_{\sigma^L_1-t_n+2}=\tilde\mu^{(i-1)*}_1$.
In addition, $\mu^{(i-1)}_1=a^{(i)\prime}+2\Delta^L_2=a^{(i)\prime}+2$
is clearly interfacial.
If $\Delta^L_1=1$, we have
$a^{(i)\prime}=p^{(i-1)\prime}-\kappa^{(i-1)}_{\sigma^L_1-i+1}
=3(c_n-\sigma^L_1+t_n)$
and therefore
$\rho^{(i-1)}(\mu^{(i-1)*}_1)=c_n-\sigma^L_1+t_n=
p^{(i-1)}-\tkappa^{(i-1)}_{\sigma^L_1-t_n+2}=\tilde\mu^{(i-1)*}_1$.
In addition, $\mu^{(i-1)}_1=a^{(i)\prime}+2\Delta^L_2=a^{(i)\prime}-2$
is clearly interfacial.

We now tackle the cases for which $i=\tau^L_{j+1}=t_k$.
Since $i+1\ne\tau^L_j$ we obtain
$\delta^{p^{(i)},p^{(i)\prime}}_{a^{(i)},e^{(i)}}=0$.
Now set $\hat a=a^{(i)}+\Delta^L_{j+1}$, and set
$\hat e=e^{(i)}$ if
$\lfloor \hat a p^{(i)}/p^{(i)\prime}\rfloor=
 \lfloor a^{(i)}p^{(i)}/p^{(i)\prime}\rfloor$
and $\hat e=1-e^{(i)}$ otherwise.
We claim that $\delta^{p^{(i)},p^{(i)\prime}}_{\hat a,\hat e}=0$
(this will be established below).
Thereupon, invoking Lemma \ref{Dire2Lem} in the case
$\lfloor \hat a p^{(i)}/p^{(i)\prime}\rfloor=
 \lfloor a^{(i)}p^{(i)}/p^{(i)\prime}\rfloor$
or Lemma \ref{Dire1Lem} in the case
$\lfloor \hat a p^{(i)}/p^{(i)\prime}\rfloor\ne
 \lfloor a^{(i)}p^{(i)}/p^{(i)\prime}\rfloor$
shows that $a^{(i)\prime}+2\Delta^L_{j+1}$
(and $a^{(i)\prime}$) is
interfacial in the $(p^{(i-1)},p^{(i-1)\prime})$-model.
Using first the fact that $a^{(i)\prime}$ is interfacial in the
$(p^{(i-1)},p^{(i-1)\prime})$-model, then Lemma \ref{Aprime_AprimeLem}
noting that (\ref{DefEiEq}) implies $\Delta^L_{j+1}=2{e^{(i-1)}}-1$,
then the definitions of $\tilde a^{(i-1)}$ and $\tilde\mu^{(i-1)*}_j$,
yields:
\begin{equation*}
\begin{split}
\rho^{(i-1)}(a^{(i)\prime})&=
\lfloor a^{(i)\prime}p^{(i-1)}/p^{(i-1)\prime}\rfloor+1-e^{(i)}\\
&=a^{(i)}-\tilde a^{(i)}-e^{(i)}
 =\tilde a^{(i-1)}+e^{(i-1)}-\Delta^L_{j+1}\\
&=\tilde\mu^{(i-1)*}_j,
\end{split}
\end{equation*}
as required.

To establish our claim that
$\delta^{p^{(i)},p^{(i)\prime}}_{\hat a,\hat e}=0$
requires the consideration of a number of cases.
%
%
%
It is convenient to separately treat $c_k>1$ and $c_k=1$.
Note that because $p^{(i)\prime}/p^{(i)}$ has continued
fraction $[c_{k}+1,c_{k+1},\ldots,c_n]$, each neighbouring
pair of odd bands in the $(p^{(i)},p^{(i)\prime})$-model
is separated by either $c_k$ or $c_k+1$ even bands.

For $c_k>1$, first consider $\sigma^L_j>t_{k+1}$.
Here, since $\tau^L_{j+1}=t_k=i<t_{k+1}<\sigma^L_j$, we have
that $a^{(i)}$ is interfacial in the $(p^{(i)},p^{(i)\prime})$-model.
Furthermore, since $c_k>1$, each pair of odd bands is separated
by at least two even bands, whereupon 
$\delta^{p^{(i)},p^{(i)\prime}}_{\hat a,\hat e}=0$ immediately.

Now consider $c_k>1$ and $\sigma^L_j\le t_{k+1}$.
Since $k(i)=k^L(i)=k$, we obtain $\Delta^L_j=(-1)^{e^{(i)}}$
from (\ref{DefEiEq}).
We readily calculate
$\mu_{j-1}^{(i)}-a^{(i)}=\Delta^L_j\kappa^{(i)}_{\sigma^L_j-i}
=\Delta^L_j(\sigma^L_j-i+1)$ where we have relied on the above
continued fraction expansion to evaluate $\kappa^{(i)}_{\sigma^L_j-i}$.
Thus $2\le\vert\mu_{j-1}^{(i)}-a^{(i)}\vert\le c_k+1$.
Then $\delta^{p^{(i)},p^{(i)\prime}}_{\hat a,\hat e}=0$
follows because, on the one hand,
$\mu_{j-1}^{(i)}$ is interfacial by the induction hypothesis,
on the second hand,
neighbouring odd bands in the $(p^{(i)},p^{(i)\prime})$-model are
separated by either $c_k$ or $c_k+1$ even bands, and on the other hand,
if $\sigma^L_j=t_{k}+1$ then $\Delta^L_{j+1}\ne\Delta^L_j$.
Note that this reasoning applies even if $j=1$ when necessarily
$k\le n-1$, if we interpret both $0$ and $p^{(i)\prime}$ as
bordering odd bands.

For $c_k=1$, first consider $\sigma^L_j=t_{k+1}$.
Lemma \ref{BigHashLem} implies that $\Delta^L_{j+1}\ne\Delta^L_j$.
Since $k(i)=k^L(i)=k$, we have
$\Delta^L_{j+1}=-\Delta^L_j=-(-1)^{e^{(i)}}$ via (\ref{DefEiEq}).
If $a^{(i)}$ is interfacial then
$\delta^{p^{(i)},p^{(i)\prime}}_{a^{(i)},e^{(i)}}=0$
implies that
$\lfloor \hat a p^{(i)}/p^{(i)\prime}\rfloor\ne
 \lfloor a^{(i)}p^{(i)}/p^{(i)\prime}\rfloor$
and therefore
$\delta^{p^{(i)},p^{(i)\prime}}_{\hat a,\hat e}=0$.
If $a^{(i)}$ is not interfacial then
$\lfloor \hat a p^{(i)}/p^{(i)\prime}\rfloor=
 \lfloor a^{(i)}p^{(i)}/p^{(i)\prime}\rfloor$,
which since $\hat a=a^{(i)}+\Delta^L_{j+1}=a^{(i)}-(-1)^{\hat e}$
immediately implies that
$\delta^{p^{(i)},p^{(i)\prime}}_{\hat a,\hat e}=0$.

Now consider $c_k=1$ and $\sigma^L_j=t_{k+1}+1$.
Lemma \ref{BigHashLem} implies that $\Delta^L_{j+1}=\Delta^L_j$.
Since $k^L(i)=k+1$ and $k(i)=k$, we have
$\Delta^L_{j+1}=\Delta^L_j=-(-1)^{e^{(i)}}$.
Then $\delta^{p^{(i)},p^{(i)\prime}}_{\hat a,\hat e}=0$
as in the previous case.

Now consider $c_k=1$ and $\sigma^L_j>t_{k+1}+1$.
Lemma \ref{BigHashLem} implies that in fact, $\sigma^L_j>t_{k+2}$.
Here, since $\tau^L_{j+1}=t_k=i<t_{k+1}<\sigma^L_j$, the induction
hypothesis implies that $a^{(i)}$ is interfacial in the
$(p^{(i)},p^{(i)\prime})$-model.
Now, however, $\delta^{p^{(i)},p^{(i)\prime}}_{\hat a,\hat e}=0$
is immediate only in the
$\lfloor \hat a p^{(i)}/p^{(i)\prime}\rfloor\ne
 \lfloor a^{(i)}p^{(i)}/p^{(i)\prime}\rfloor$ case.
In the case
$\lfloor \hat a p^{(i)}/p^{(i)\prime}\rfloor=
 \lfloor a^{(i)}p^{(i)}/p^{(i)\prime}\rfloor$,
since $a^{(i)}$ is interfacial and
$\delta^{p^{(i)},p^{(i)\prime}}_{a^{(i)},e^{(i)}}=0$,
it follows that $\hat a=a^{(i)}+(-1)^{e^{(i)}}$.
Then $\Delta^L_{j+1}=(-1)^{e^{(i)}}=-(-1)^{e^{(i+1)}}$,
the final equality arising from (\ref{DefEiEq}) because
$i+1=t_{k+1}<\sigma^L_j$ implies that $k(i+1)=k+1$
and $k^L(i+1)=k^L(i)$.
Again via the induction hypothesis,
$\tau^L_{j+1}<t_{k+1}=i+1<t_{k+2}<\sigma^L_j$ implies that
$a^{(i+1)}$ is interfacial in the $(p^{(i)},p^{(i)\prime})$-model
and $i+2<\sigma^L_j$ implies that 
$\delta^{p^{(i+1)},p^{(i+1)\prime}}_{a^{(i+1)},e^{(i+1)}}=0$,
so that
$\lfloor (a^{(i+1)}+\Delta^L_{j+1}) p^{(i+1)}/p^{(i+1)\prime}\rfloor\ne
 \lfloor a^{(i+1)}p^{(i+1)}/p^{(i+1)\prime}\rfloor$.
Invoking Lemma \ref{Dire1Lem} with $p'=p^{(i+1)\prime}$,
$p=p^{(i)}$, and $\{a,\hat a\}=\{a^{(i+1)},a^{(i+1)}+\Delta^L_{j+1}\}$
now shows that
$\lfloor a^{(i)} p^{(i)}/p^{(i)\prime}\rfloor
=\lfloor(a^{(i)}+2\Delta^L_{j+1})p^{(i)}/p^{(i)\prime}\rfloor$,
whereupon, certainly
$\delta^{p^{(i)},p^{(i)\prime}}_{\hat a,\hat e}=0$.

This proves the ninth bulleted item for all $i$ with
$0\le i<\tau^L_{j+1}$ when $j>0$ and all $i$ with
$0\le i<t_n$ when $j=0$.
An entirely analogous argument gives the remaining bulleted item.

We now see that our proposition holds for $i$ replaced by $i-1$,
and thus the lemma is proved by induction.
\cqfd
\medskip


Before taking the final step to proving Theorem \ref{CoreThrm}, we need
the following result concerning the quantities defined in Section
\ref{ConSec}, for the vector $\boldu=\boldu^L+\boldu^R$.

\begin{lemma}\label{ParityLem}
For $0\le j\le t$,
\begin{equation*}
\begin{split}
\alpha_j''&\equiv Q_j\;(\mod2);\\
\beta_j'&\equiv Q_j-Q_{j+1}\;(\mod2).
\end{split}
\end{equation*}
\end{lemma}

\Proof Since $\alpha_t''=0$, $\beta_t'=0$ and $Q_t=Q_{t+1}=0$,
this result is manifest for $j=t$.

We now proceed by downward induction. Thus assume the result holds
for a particular $j>0$.
When $j\ne t_{k(j)}$, equations (\ref{Const1Eq}) and (\ref{Const3Eq})
imply that $\beta_{j-1}'\equiv\beta_j'+(\boldu^L)_j-(\boldu^R)_j$.
Equations (\ref{ParityDef}), (\ref{MNmatEq2})  and (\ref{MNEq2})
imply that $Q_{j-1}\equiv Q_{j+1}-(\boldu^L)_j-(\boldu^R)_j$.
Thus the induction hypothesis immediately gives
$\beta_{j-1}'\equiv Q_{j-1}-Q_{j}$ in this case.

When $j=t_{k(j)}$, equations (\ref{Const1Eq}) and (\ref{Const3Eq})
imply that $\beta_{j-1}'\equiv\alpha_j''-\beta_j'+(\boldu^L)_j-(\boldu^R)_j$.
Equations (\ref{ParityDef}), (\ref{MNmatEq2}) and (\ref{MNEq1})
imply that $Q_{j-1}\equiv Q_{j}+Q_{j+1}-(\boldu^L)_j-(\boldu^R)_j$.
Thus the induction hypothesis also gives
$\beta_{j-1}'\equiv Q_{j-1}-Q_{j}$ in this case.

In both cases, equations (\ref{Const2Eq}) and
(\ref{Const3Eq}) give $\alpha_{j-1}''=\alpha_j''+\beta_{j-1}'$,
whence the induction hypothesis together with the above result
immediately gives $\alpha_{j-1}''\equiv Q_{j-1}$ as required.
\cqfd
\medskip

In order to isolate the $i=0$ case of Lemma \ref{CoreIndLem},
we note that comparison of the definitions of
Sections \ref{CoreSec} and \ref{IndParamSec},
gives $a=a^{(0)}$, $b=b^{(0)}$, $e=e^{(0)}$, $f=f^{(0)}$,
$p=p^{(0)}$, $p'=p^{(0)\prime}$,
$\boldmu=\boldmu^{(0)}$, $\boldmu^*=\boldmu^{(0)*}$,
$\boldnu=\boldnu^{(0)}$, $\boldnu^*=\boldnu^{(0)*}$,
and since $k(0)=0$,
$\boldu^{L}_{(\flat,k(0))}=\boldu^{L}_{\flat}$
and
$\boldu^{R}_{(\sharp,k(0))}=\boldu^{R}_{\sharp}$.
We now readily obtain:

\begin{corollary}\label{CoreCor}
Let $p'>2p$,
$m_0\equiv\alpha^{p,p'}_{a,b}\;(\mod2)$ and
$m_1\equiv\alpha^{p,p'}_{a,b}+\beta^{p,p'}_{a,b,e,f}\;(\mod2)$.
Then:
\begin{equation*}
\begin{split}
&\mchi^{p,p'}_{a,b,e,f}(m_0,m_1+2\delta^{p,p^{\prime}}_{a,e})
\left\{
\begin{matrix} \boldmu^{\phantom{*}};\boldnu^{\phantom{*}}\\
                    \boldmu^*;\boldnu^* \end{matrix} \right\}\\
&\qquad=\sum
  q^{\frac{1}{4}\hat{\boldm}^T\boldC\hat{\boldm}-\frac{1}{4} m_0^2
   -\frac{1}{2}(\boldu^{L}_{\flat}+\boldu^{R}_{\sharp})\cdot\boldm
  +\frac{1}{4}\gamma_0}
  \prod_{j=1}^{t-1}
  \left[
  {m_j-\frac{1}{2}(\boldC^*\hat{\boldm}\!+\!\boldu^{L}\!+\!\boldu^{R})_j
       \atop m_j} \right]_q,
\end{split}
\end{equation*}
with the sum to be taken over all
$(m_2,m_3,\ldots,m_{t-1})\equiv(Q_2,Q_3,\ldots,Q_{t-1})$,
with $\boldm=(m_1,m_2,m_3,\ldots,m_{t-1})$
and $\hat{\boldm}=(m_0,m_1,m_2,m_3,\ldots,m_{t-1})$.
For the case $t=2$, the summation in the above expression is omitted.
For $t=1$, we have:
\begin{equation*}
\mchi^{p,p'}_{a,b,e,f}(m_0,m_1)
\left\{
\begin{matrix} \boldmu^{\phantom{*}};\boldnu^{\phantom{*}}\\
                    \boldmu^*;\boldnu^* \end{matrix} \right\}
 =q^{\frac{1}{4}m_0^2+\frac{1}{4}\gamma_0} \delta_{m_1,0}.
\end{equation*}
(In this $t=1$ case, necessarily $\delta^{p,p^{\prime}}_{a,e}=0$.)

In addition, $\boldmu,\boldmu^*,\boldnu,\boldnu^*$
are an interfacial mazy-four in the $(p,p')$-model sandwiching $(a,b)$,
with $\rho^{p,p'}(\mu_j)=\tilde\mu_j$ and
$\rho^{p,p'}(\mu_j^*)=\tilde\mu_j^*$ for $d^L_0\le j<d^L$, and
$\rho^{p,p'}(\nu_j)=\tilde\nu_j$ and
$\rho^{p,p'}(\nu_j^*)=\tilde\nu_j^*$ for $d^R_0\le j<d^R$.

Also, $\delta^{p,p^{\prime}}_{b,f}=1$ only if $\tau^R_{d^R}=1$
and $\nu_{d^R-1}=b-(-1)^f$.
If $\sigma^R_{d^R}>t_1$ then $b$ is interfacial in
the $(p,p')$-model.
If $\sigma^R_{d^R}=t_1$ and $b$ is not interfacial in the
$(p,p')$-model then $\tau^R_{d^R}=t_1+1$.

If $b$ is interfacial and $\sigma^R_{d^R}>0$ then
\begin{equation*}
\rho^{p,p'}(b)=\sum_{m=2}^{d^R}
        \Delta^R_m (\tkappa_{\tau^R_m}-\tkappa_{\sigma^R_m})
\:+\:\begin{cases}
        \tkappa_{\sigma^{R}_1} &\text{if $\Delta^R_1=-1$;}\\
        p-\tkappa_{\sigma^{R}_1} &\text{if $\Delta^R_1=+1$.}
  \end{cases}
\end{equation*}
If $b$ is not interfacial and $\sigma^R_{d^R}>0$ then
\begin{equation*}
\left\lfloor\frac{bp}{p'}\right\rfloor=\sum_{m=2}^{d^R}
        \Delta^R_m (\tkappa_{\tau^R_m}-\tkappa_{\sigma^R_m})
\:-\:\frac12(1-\Delta^R_{d^R})
\:+\:\begin{cases}
        \tkappa_{\sigma^{R}_1} &\text{if $\Delta^R_1=-1$;}\\
        p-\tkappa_{\sigma^{R}_1} &\text{if $\Delta^R_1=+1$.}
  \end{cases}
\end{equation*}
\end{corollary}

\Proof
Lemma \ref{CoreIndLem} gives
$\alpha^{p,p'}_{a,b}=\alpha_0''=\alpha_0$ and
$\beta^{p,p'}_{a,b,e,f}=\beta_0'=\beta_0$,
whereupon, the first statement follows from Lemma \ref{ParityLem}.
The statements in the second and third paragraphs follow from the $i=0$
case of Lemma \ref{CoreIndLem} after noting that $t_0<0=\tau^R_{d^R+1}<t_1$.
For the final paragraph, the $i=0$ case of Lemma \ref{CoreIndLem}
implies that if $\sigma^R_{d^R}>0$ then
\begin{equation*}
\left\lfloor\frac{bp}{p'}\right\rfloor
=\tilde b^{(0)}
=-f
+\sum_{m=2}^{d^R}
        \Delta^R_m (\tkappa_{\tau_m^{R}}-\tkappa_{\sigma_m^{R}})
\:+\:\begin{cases}
        \tkappa_{\sigma^{R}_1} &\text{if $\Delta^R_1=-1$;}\\
        p-\tkappa_{\sigma^{R}_1} &\text{if $\Delta^R_1=+1$.}
  \end{cases}
\end{equation*}
Note that $\sigma^R_{d^R}>0$ implies that $\tau^R_{d^R}>1$ which,
from above, implies that $\delta^{p,p'}_{b,f}=0$.
Then if $b$ is interfacial, $\rho^{p,p'}(b)=\lfloor bp/p'\rfloor+f$
giving the desired result in this case.
If $b$ is not interfacial then, from above, $\sigma^R_{d^R}\le t_1$
and therefore $k^R(0)=0$. Since $k(0)=0$, (\ref{DefFiEq}) yields
$\Delta^R_{d^R}=(-1)^f=1-2f$, from which the desired result follows.
\cqfd
\medskip

\medskip
\noindent {\it Proof of Theorem \ref{CoreThrm}: }
For reduced runs, Theorem \ref{CoreThrm}
now follows from Lemma \ref{MazyResPathGenLem} and Corollary
\ref{CoreCor} since $\alpha^{p,p'}_{a,b}=b-a$ and, via Lemma \ref{ParityLem},
$\alpha^{p,p'}_{a,b}+\beta^{p,p'}_{a,b,e,f}\equiv Q_1\;(\mod2)$.
That it holds in the case of arbitrary naive runs then follows from
Lemma \ref{ReducedLem}.
\cqfd
\medskip

For later convenience, we show here that the values
$\{\mu^*_j\}_{j=1}^{d^L-1}$ and $\{\nu^*_j\}_{j=1}^{d^R-1}$
that are obtained from {\em naive} runs are interfacial.
It is not possible using the results of this section to show that all
$\{\mu_j\}_{j=1}^{d^L-1}$ and $\{\nu_j\}_{j=1}^{d^R-1}$ 
are also interfacial.
However, that they actually are will be deduced in
Section \ref{AssimilateTrees}.

\begin{corollary}\label{StarInterCor}
Let $\{\tau^{L}_j,\sigma^{L}_j,\Delta^{L}_j\}_{j=1}^{d^{L}}$ and
$\{\tau^{R}_j,\sigma^{R}_j,\Delta^{R}_j\}_{j=1}^{d^{R}}$ be naive runs.

If $1\le j<d^L$ then $\mu^*_j$ is interfacial with
$\rho^{p,p'}(\mu^*_j)=\tilde\mu^*_j$.
In addition, if $1\le j<d^L$, and $\tau^{L}_{j+1}\ne\sigma^{L}_{j+1}$ then
$\mu_j$ is interfacial with $\rho^{p,p'}(\mu_j)=\tilde\mu_j$.

If $1\le j<d^R$ then $\nu^*_j$ is interfacial with
$\rho^{p,p'}(\nu^*_j)=\tilde\nu^*_j$.
In addition, if $1\le j<d^R$, and $\tau^{R}_{j+1}\ne\sigma^{R}_{j+1}$ then
$\nu_j$ is interfacial with $\rho^{p,p'}(\nu_j)=\tilde\nu_j$.
\end{corollary}

\Proof That this holds for reduced runs is immediate from Corollary
\ref{CoreCor}.
If $\tau^L_j=\sigma^L_j$ then $\mu^*_j=\mu^*_{j-1}$ and
$\tilde\mu^*_j=\tilde\mu^*_{j-1}$.
Since $\tau^L_j=\sigma^L_j$ cannot occur for $j=1$,
the required result follows in the $\mu^*_j$ case.
With $\tau^{L}_{j+1}\ne\sigma^{L}_{j+1}$, the result for the $\mu_j$
cases follows immediately from that for reduced runs.

The final paragraph follows in the same way.
\cqfd
\medskip

\subsection{Transferring to the original weighting}\label{TransSec}

In this and the following sections, we fix naive runs
$\run^L=\{\tau^{L}_j,\sigma^{L}_j,\Delta^{L}_j\}_{j=1}^{d^{L}}$ and
$\run^R=\{\tau^{R}_j,\sigma^{R}_j,\Delta^{R}_j\}_{j=1}^{d^{R}}$,
and make use of the definitions of Section \ref{CoreSec}.
The set
$\P^{p,p'}_{a,b,c}(L)\left\{
\begin{smallmatrix} \boldmu^{\phantom{*}}\!&;&\boldnu^{\phantom{*}}\\
                    \boldmu^*\!&;&\boldnu^* \end{smallmatrix} \right\}$
and its generating function 
$\ochi^{p,p'}_{a,b,c}(L)\left\{
\begin{smallmatrix} \boldmu^{\phantom{*}};\boldnu^{\phantom{*}}\\
                    \boldmu^*;\boldnu^* \end{smallmatrix} \right\}$
are defined in Section \ref{MazyOSec}.

The following lemma relates this generating function,
defined in terms of the original weight function (\ref{WtDef}),
to that considered in Theorem \ref{CoreThrm}
(the latter generating function was defined in (\ref{MazyPathGenDef})
in terms of the modified weight function (\ref{WtDef2})).

\begin{lemma}\label{TransLem}
Let $p'>2p$.
Set $c=b+\Delta^R_{d^R}$.

If $\sigma^R_{d^R}=0$ and either $\lfloor pb/p'\rfloor=\lfloor pc/p'\rfloor$
or $c\in\{0,p'\}$ then:
\begin{equation*}
\ochi^{p,p'}_{a,b,c}(L)\left\{
\begin{matrix} \boldmu^{\phantom{*}};\boldnu^{\phantom{*}}\\
                    \boldmu^*;\boldnu^* \end{matrix} \right\}
=
q^{-\frac12(L+\Delta^R_{d^R}(a-b))}
\mchi^{p,p'}_{a,b,e,f}(L)\left\{
\begin{matrix} \boldmu^{\phantom{*}};\boldnu^{\phantom{*}}\\
                    \boldmu^*;\boldnu^* \end{matrix} \right\}.
\end{equation*}
Otherwise:
\begin{equation*}
\ochi^{p,p'}_{a,b,c}(L)\left\{
\begin{matrix} \boldmu^{\phantom{*}};\boldnu^{\phantom{*}}\\
                    \boldmu^*;\boldnu^* \end{matrix} \right\}
=
\mchi^{p,p'}_{a,b,e,f}(L)\left\{
\begin{matrix} \boldmu^{\phantom{*}};\boldnu^{\phantom{*}}\\
                    \boldmu^*;\boldnu^* \end{matrix} \right\}.
\end{equation*}
\end{lemma}

\Proof
Let
$\tilde h\in\P^{p,p'}_{a,b,e,f}(L)\left\{
\begin{smallmatrix} \boldmu^{\phantom{*}};\boldnu^{\phantom{*}}\\
                    \boldmu^*;\boldnu^* \end{smallmatrix} \right\}$,
and let
$h\in\P^{p,p'}_{a,b,c}(L)\left\{
\begin{smallmatrix} \boldmu^{\phantom{*}};\boldnu^{\phantom{*}}\\
                    \boldmu^*;\boldnu^* \end{smallmatrix} \right\}$
have the same sequence of heights: $h_i=\tilde h_i$ for $0\le i\le L$.

First consider $\sigma^R_{d^R}=0$.
The $i=1$ case of Lemma \ref{Mu_ALem} implies that
$\nu_{d^R-1}=b+\Delta^R_{d^R}$ and hence
$\tilde h_{L-1}=b-\Delta^R_{d^R}$.
The definition of $f$ (in Section \ref{CoreSec})
gives $\Delta^R_{d^R}=-(-1)^f$.
Then, by (\ref{ModWtDef}), the $L$th vertex of $\tilde h$ is
scoring with weight $\frac12(L+\Delta^R_{d^R}(a-b))$.
On the other hand, with $c=b+\Delta^R_{d^R}$,
if $\lfloor pb/p'\rfloor=\lfloor pc/p'\rfloor$ and $0<c<p'$ then
Table \ref{WtsTable} shows that the $L$th vertex of $h$ is non-scoring.
This is also the case if $c\in\{0,p'\}$ because the $0$th and
$(p'-1)$th bands are defined to be even when $p'>2p$.
Therefore $\owt(h)=\mwt(\tilde h)-\frac12(L+\Delta^R_{d^R}(a-b))$.
In the case where $0<c<p'$ and $\lfloor pb/p'\rfloor\ne\lfloor pc/p'\rfloor$,
the $L$th vertex is scoring in both $\tilde h$ and $h$, and
therefore $\owt(h)=\mwt(\tilde h)$.
The required identities between generating functions then follow
whenever $\sigma^R_{d^R}=0$.

For $\sigma^R_{d^R}>0$, note that $\tau^R_{d^R}>\sigma^R_{d^R}$
implies that $\tau^R_{d^R}>1$ and hence, via Corollary \ref{CoreCor},
that $\delta^{p,p'}_{b,f}=0$.

In the case $0<\sigma^R_{d^R}\le t_1$, we have $k^R=0$
and therefore $\Delta^R_{d^R}=(-1)^f$.
{}From $\delta^{p,p'}_{b,f}=0$ then follows that $\owt(h)=\mwt(\tilde h)$,
and hence the required identity between the generating functions results.

If $\sigma^R_{d^R}>t_1$, then $b$ is interfacial by Corollary \ref{CoreCor}.
Again $\delta^{p,p'}_{b,f}=0$ implies that $\owt(h)=\mwt(\tilde h)$
whereupon the lemma follows.
\cqfd
\medskip

The following theorem specifies the subset of $\P^{p,p'}_{a,b,c}(L)$
(weighted as in (\ref{PathGenDef1})) for which the generating function is
$F(\boldu^{L},\boldu^{R},L)$.
It deals with both the cases $p'>2p$ and $p'<2p$.

\begin{theorem}\label{Core2Thrm}
If $b$ is interfacial, let $c\in\{b\pm1\}$.
Otherwise if $b$ is not interfacial, set $c=b+\Delta^R_{d^R}$.
If $L\equiv b-a\;(\mod2)$ then:
\begin{equation}\label{Core2Eq}
\ochi^{p,p'}_{a,b,c}(L)
\left\{
\begin{matrix} \boldmu^{\phantom{*}};\boldnu^{\phantom{*}}\\
                    \boldmu^*;\boldnu^* \end{matrix} \right\}
= F(\boldu^{L},\boldu^{R},L).
\end{equation}
\end{theorem}

\Proof First consider $p'>2p$ and $c=b+\Delta^R_{d^R}$.

If $\sigma^R_{d^R}=0$ and $\lfloor bp/p'\rfloor=\lfloor cp/p'\rfloor$
with $0<c<p'$
then $\gamma'=\gamma_0-2(L+\Delta^R_{d^R}(a-b))$ by the definition in
Section \ref{ConSec}.
This is also the case if $c=0$ (when $b=1$ and $\Delta^R_{d^R}=-1$)
or $c=p'$ (when $b=p'-1$ and $\Delta^R_{d^R}=1$).
Otherwise $\gamma'=\gamma_0$.
Expression (\ref{Core2Eq}) then follows from
Theorem \ref{CoreThrm} and Lemma \ref{TransLem}.

If $b$ is interfacial then Table \ref{WtsTable} shows that a path
$h\in\P^{p,p'}_{a,b,c}(L)$ has the same generating function
whether $c=b+1$ or $c=b-1$. The theorem then follows in the $p'>2p$ case.

Now consider $p'<2p$. If $h\in\P^{p,p'}_{a,b,c}(L)$, and
$\hat h\in\P^{p'-p,p'}_{a,b,c}(L)$ is defined by $\hat h_i=h_i$
for $0\le i\le L$, then $\owt(\hat h)=\frac14(L^2-\alpha^2)-\owt(h)$,
where $\alpha=b-a$, by a direct analogue of Lemma \ref{DresLem}
(here, as there, each scoring vertex of $\hat h$ corresponds to
a non-scoring vertex of $h$ and vice-versa).
Since $p'>2(p'-p)$, we thus obtain:
\begin{equation}\label{Core2bEq}
\ochi^{p,p'}_{a,b,c}(L)
\left\{
\begin{matrix} \boldmu^{\phantom{*}};\boldnu^{\phantom{*}}\\
                    \boldmu^*;\boldnu^* \end{matrix} \right\}
= q^{\frac14(L^2-\alpha^2)}
  F^{p'-p,p'}(\boldu^{L},\boldu^{R},L;q^{-1}),
\end{equation}
where $F^{p'-p,p'}(\boldu^{L},\boldu^{R},L;q^{-1})$
is $F(\boldu^{L},\boldu^{R},L;q^{-1})$ with all the quantities
employed in its definition by (\ref{FEq}) pertaining to the
continued fraction of $p'/(p'-p)$. Below, we make use of
$\boldC^{p'-p,p'}$ and $(\boldu^L_\flat+\boldu^R_\sharp)^{p'-p,p'}$
and $\gamma(\run^L,\run^R)^{p'-p,p'}$ defined similarly.

Now note that
$\hat{\boldm}^T\boldC^{p'-p,p'}\hat{\boldm}
 =\hat{\boldm}\boldC^{p,p'}\hat{\boldm}+L^2-2Lm_1$,
$$
(\boldu^L_\flat+\boldu^R_\sharp)^{p'-p,p'}\cdot\boldm
=(\boldu^L+\boldu^R)\cdot\boldm
-(\boldu^L_\flat+\boldu^R_\sharp)^{p,p'}\cdot\boldm,
$$
and
$\gamma(\run^L,\run^R)^{p,p'}
=-\gamma(\run^L,\run^R)^{p'-p,p'}-\alpha_0^2$ by (\ref{Const3Eq}).
{}From these expressions together with Lemma \ref{PartitionInvLem},
we obtain:
\begin{equation}\label{Core2cEq}
F^{p'-p,p'}(\boldu^{L},\boldu^{R},L;q^{-1})
= q^{-\frac14(L^2-\alpha_0^2)} F^{p,p'}(\boldu^{L},\boldu^{R},L;q).
\end{equation}
Together with (\ref{Core2bEq}), this proves (\ref{Core2Eq}) in the
$p'<2p$ case, having noted, via Lemma \ref{CoreIndLem}, that
$\alpha=\alpha^{p,p'}_{a,b}=\alpha_0$.
\cqfd
\medskip

\subsection{The other direction}\label{OtherSec}

Theorem \ref{Core2Thrm} deals with either $b$ being interfacial
or $c$ taking the specific value $b+\Delta^R_{d^R}$.
In this section, we deal with the other cases in which
$b$ is not interfacial and $c=b-\Delta^R_{d^R}$.
The generating functions then involve the fermion-like
expressions defined in Section \ref{FermLikeSec}.

As in the previous section, we fix naive runs
$\run^L=\{\tau^{L}_j,\sigma^{L}_j,\Delta^{L}_j\}_{j=1}^{d^{L}}$ and
$\run^R=\{\tau^{R}_j,\sigma^{R}_j,\Delta^{R}_j\}_{j=1}^{d^{R}}$,
and make use of the definitions of Section \ref{CoreSec}.
In addition, for convenience, we set $\sigma=\sigma^R_{d^R}$
and $\Delta=\Delta^R_{d^R}$. It will also be useful to set
\begin{equation}\label{tauDef2}
\tau=
\left\{ \begin{array}{ll}
             \tau^R_{d^R} & \quad\mbox{if } d^R>1;\\
             t_n & \quad\mbox{if } d^R=1 \mbox{ and } \sigma<t_n;\\
             t   & \quad\mbox{if } d^R=1 \mbox{ and } \sigma\ge t_n.
  \end{array} \right.
\end{equation}
Recall that $\sigma<t_n$ implies that $d^R_0=0$, and
$\sigma\ge t_n$ implies that $d^R_0=1$.

\begin{theorem}\label{FermLikeThrm}
Let all parameters be as in Section \ref{CoreSec}.
Let $b$ be non-interfacial and set $c=b-\Delta^R_{d^R}$.
If $L\equiv b-a\;(\mod2)$ then:
\begin{equation}\label{FermLikeEq}
\ochi^{p,p'}_{a,b,c}(L)
\left\{
\begin{matrix} \boldmu^{\phantom{*}};\boldnu^{\phantom{*}}\\
                    \boldmu^*;\boldnu^* \end{matrix} \right\}
= \widetilde F(\boldu^{L},\boldu^{R},L).
\end{equation}
\end{theorem}

\Proof We first consider $p'>2p$.
Since $b$ is not interfacial,
Corollary \ref{CoreCor} implies that $\sigma\le t_1$.

In the subcases for which $d^R>d^R_0$ note that $\boldnu\ne()\ne\boldnu^*$. 
Then, the definitions of Section \ref{CoreSec} yield
$\nu_{d^R-1}=b+\Delta\kappa_{\sigma}$
and $\nu^*_{d^R-1}=b-\Delta(\kappa_{\tau}-\kappa_{\sigma})$.
In the subcase for which $d^R=d^R_0$, note that $\boldnu=\boldnu^*=()$. 

For $\sigma=0$, when $d^R>d^R_0$ we have $\nu_{d^R-1}=b+\Delta$,
whereupon the first identity of Lemma \ref{RecurseLem}(2) gives:
\begin{equation}\label{FermLike1Eq}
\begin{split}
\ochi^{p,p'}_{a,b,c}(L)
\left\{
\begin{matrix} \boldmu^{\phantom{*}};\boldnu^{\phantom{*}}\\
                    \boldmu^*;\boldnu^* \end{matrix} \right\}
&=
q^{\frac12(L+\Delta(a-b))}
\ochi^{p,p'}_{a,b,b+\Delta}(L)
\left\{
\begin{matrix} \boldmu^{\phantom{*}};\boldnu^{\phantom{*}}\\
                    \boldmu^*;\boldnu^* \end{matrix} \right\}\\
&=
q^{\frac12(L+\Delta(a-b))} F(\boldu^{L},\boldu^{R},L),
\end{split}
\end{equation}
the second equality following from Theorem \ref{Core2Thrm}.
When $d^R=d^R_0$ so that $\boldnu=\boldnu^*=()$,
(\ref{FermLike1Eq}) also follows from Lemma \ref{RecurseLem}(2)
and Theorem \ref{Core2Thrm}.
The definition (\ref{Tilde1Def}) then gives the required result
in this $\sigma=0$ case.

Next, we consider $0\le\sigma<\tau-1$.  Set $b'=b-\Delta$.
When $d^R>d^R_0$, we find that $b'-\nu_{d^R-1}=-\Delta(\kappa_{\sigma}+1)$
and $b'-\nu^*_{d^R-1}=\Delta(\kappa_{\tau}-\kappa_{\sigma}-1)$.
In particular, $b'$ is strictly between $\nu_{d^R-1}$ and $\nu^*_{d^R-1}$.
The same is true of $b'-\Delta$ unless $\sigma+2=\tau\le t_1+1$,
in which case $b'-\Delta=\nu^*_{d^R-1}$.
Substituting $b\to b'$ and $L\to L+1$ into the second identity of
Lemma \ref{RecurseLem}(1) and rearranging, yields:
\begin{equation}\label{FermLike2Eq}
\begin{split}
&q^{\frac12(L+1-\Delta(a-b'))}
\ochi^{p,p'}_{a,b'+\Delta,b'}(L)
\left\{
\begin{matrix} \boldmu^{\phantom{*}};\boldnu^{\phantom{*}}\\
                    \boldmu^*;\boldnu^* \end{matrix} \right\} \\
&\hskip25mm=
\ochi^{p,p'}_{a,b',b'+\Delta}(L+1)
\left\{
\begin{matrix} \boldmu^{\phantom{*}};\boldnu^{\phantom{*}}\\
                    \boldmu^*;\boldnu^* \end{matrix} \right\}
-
\ochi^{p,p'}_{a,b'-\Delta,b'}(L)
\left\{
\begin{matrix} \boldmu^{\phantom{*}};\boldnu^{+\phantom{*}}\\
                    \boldmu^*;\boldnu^{+*} \end{matrix} \right\}.
\end{split}
\end{equation}
where if $b'-\Delta\ne\nu^*_{d^R-1}$ then $\boldnu^+=\boldnu$ and
$\boldnu^{+*}=\boldnu^*$;
and if $b'-\Delta=\nu^*_{d^R-1}$ then $d^+\ge d^R_0$ is the smallest value
such that $b'-\Delta=\nu^*_{d^+}$, and then
$\boldnu^+=(\nu_{d^R_0},\ldots,\nu_{d^+-1})$ and
$\boldnu^{+*}=(\nu^*_{d^R_0},\ldots,\nu^*_{d^+-1})$.
The term on the left of (\ref{FermLike2Eq}) is (apart from the prefactor)
that on the left of (\ref{FermLikeEq}).
The first term on the right is equal to
$F(\boldu^{L},\boldu^{R+},L+1)$ by Theorem \ref{Core2Thrm},
after noting that with $\run^{R+}$ in place of $\run^{R}$, the definitions
of Section \ref{CoreSec} lead to the same mazy-pair $\boldnu,\boldnu^*$.
If, on the one hand, $b'-\Delta\ne\nu^*_{d^R-1}$ the second term on the
right is likewise equal to $F(\boldu^{L},\boldu^{R++},L)$ by using
$\run^{R++}$ in place of $\run^{R}$.
Therefore:
\begin{equation}\label{FermLike3Eq}
q^{\frac12(L-\Delta(a-b))}
\ochi^{p,p'}_{a,b'+\Delta,b'}(L)
\left\{
\begin{matrix} \boldmu^{\phantom{*}};\boldnu^{\phantom{*}}\\
                    \boldmu^*;\boldnu^* \end{matrix} \right\}
=F(\boldu^{L},\boldu^{R+},L+1)-F(\boldu^{L},\boldu^{R++},L),
\hskip-2.5mm
\end{equation}
from which (\ref{FermLikeEq}) follows in this case via the
definition (\ref{Tilde1Def}).

On the other hand, if $b'-\Delta=\nu^*_{d^R-1}$, note
in the subcase for which $d^R>1$, that the mazy pair
$\boldnu^+,\boldnu^{+*}$ arises from the run
$\run'=\{\tau^{R}_j,\sigma^{R}_j,\Delta^{R}_j\}_{j=1}^{d^*}$
where $d^*=\max\{d^+,1\}$.
Also note that $\nu^*_{d^*}=\nu^*_{d^*+1}=\cdots=\nu^*_{d^R-1}=b-2\Delta$
implies that $\tau^R_j=\sigma^R_j$ for $d^*<j<d$ and
$\tau^R_{d^R}=\sigma^R_{d^R}+2$ whereupon $\boldu(\run')=\boldu^{R++}$
as defined in Section \ref{FermLikeSec}, and (\ref{DeltaEq}) gives
$\boldDelta(\run')=\boldDelta(\run^{++})$.
In the subcase for which $d^R=1$ and $d^R_0=0$,
the mazy pair $\boldnu^+=\boldnu^{+*}=()$ arises from
$\run'=\{\tau^{R}_1,t_n,\Delta^{R}_1\}$ whereupon again
$\boldu(\run')=\boldu^{R++}$ and $\boldDelta(\run')=\boldDelta(\run^{++})$.
Thus, even when $b'-\Delta=\nu^*_{d^R-1}$,
the second term on the right of (\ref{FermLike2Eq}) is equal to
$F(\boldu^{L},\boldu^{R++},L)$ by Theorem \ref{Core2Thrm},
and (\ref{FermLike3Eq}) and then (\ref{FermLikeEq}) follow.

For this $0\le \sigma<\tau-1$ case, it remains to consider
the subcase for which $d^R=d^R_0=1$.
Here $\boldnu=\boldnu^*=\boldnu^+=\boldnu^{+*}=()$.
The instance $\sigma=t_1$ may be excluded because then $n=1$
and $b$ may be shown to be interfacial.
Otherwise Lemma \ref{RecurseLem}(1) and Theorem \ref{Core2Thrm}
lead to (\ref{FermLike2Eq}), (\ref{FermLike3Eq}) and (\ref{FermLikeEq})
in a similar fashion to the $b'-\Delta\ne\nu^*_{d^R-1}$ case considered above
after noting that $\sigma<t-1$ ensures that $1\le b',b'-\Delta<p'$.

We now consider $0<\sigma=\tau-1$.
In the subcase where $d^R>d^R_0$ we have
             $b-\Delta=\nu^*_{d^R-1}$ and
             $b+\Delta\kappa_{\sigma}=\nu_{d^R-1}$,
so that certainly both $b$ and $b+\Delta$ are strictly between
$\nu^*_{d^R-1}$ and $\nu_{d^R-1}$.
Eliminating the terms that appear second on the right sides of the two
identities in Lemma \ref{RecurseLem}(1), yields:
\begin{equation}\label{FermLike4Eq}
\begin{split}
&q^{\frac12(L-\Delta(a-b))}
\ochi^{p,p'}_{a,b,b-\Delta}(L)
\left\{
\begin{matrix} \boldmu^{\phantom{*}};\boldnu^{\phantom{*}}\\
                    \boldmu^*;\boldnu^* \end{matrix} \right\} \\
&\hskip18mm=
\ochi^{p,p'}_{a,b,b+\Delta}(L)
\left\{
\begin{matrix} \boldmu^{\phantom{*}};\boldnu^{\phantom{*}}\\
                    \boldmu^*;\boldnu^* \end{matrix} \right\}
+(q^L-1)
\ochi^{p,p'}_{a,b-\Delta,b}(L-1)
\left\{
\begin{matrix} \boldmu^{\phantom{*}};\boldnu^{+\phantom{*}}\\
                    \boldmu^*;\boldnu^{+*} \end{matrix} \right\},
\hskip-2mm
\end{split}
\end{equation}
where if $d^+\ge d^R_0$ is the smallest value
such that $b-\Delta=\nu^*_{d^+}$, then
$\boldnu^+=(\nu_{d^R_0},\ldots,\nu_{d^+-1})$ and
$\boldnu^{+*}=(\nu^*_{d^R_0},\ldots,\nu^*_{d^+-1})$.

The first term on the right of (\ref{FermLike4Eq}) is equal to
$F(\boldu^{L},\boldu^{R},L)$ by Theorem \ref{Core2Thrm}.
Now note that in the subcase for which $d^R>1$, the mazy pair
$\boldnu^+,\boldnu^{+*}$ arises from the run
$\run'=\{\tau^{R}_j,\sigma^{R}_j,\Delta^{R}_j\}_{j=1}^{d^*}$
where $d^*=\max\{d^+,1\}$.
Also note that $\nu^*_{d^*}=\nu^*_{d^*+1}=\cdots=\nu^*_{d^R-1}=b-\Delta$
implies that $\tau^R_j=\sigma^R_j$ for $d^*<j<d^R$ and
$\tau^R_{d^R}=\sigma^R_{d^R}+1$ whereupon $\boldu(\run')=\boldu^{R+}$
as defined in Section \ref{FermLikeSec}, and (\ref{DeltaEq}) gives
$\boldDelta(\run')=\boldDelta(\run^{+})$.
In the subcase for which $d^R=1$ and $d^R_0=0$ the mazy pair
$\boldnu^+=\boldnu^{+*}=()$ arises from
$\run'=\{\tau^{R}_1,t_n,\Delta^{R}_1\}$ whereupon again
$\boldu(\run')=\boldu^{R+}$ and $\boldDelta(\run')=\boldDelta(\run^{+})$.
Thus whenever $d^R>d^R_0$, Theorem \ref{Core2Thrm} implies that the
second term on the right of (\ref{FermLike4Eq}) is equal to
$(q^L-1)F(\boldu^{L},\boldu^{R+},L-1)$.
Thus, (\ref{FermLike4Eq}) gives:
\begin{equation}\label{FermLike5Eq}
\begin{split}
q^{\frac12(L-\Delta(a-b))}
&\ochi^{p,p'}_{a,b,b-\Delta}(L)
\left\{
\begin{matrix} \boldmu^{\phantom{*}};\boldnu^{\phantom{*}}\\
                    \boldmu^*;\boldnu^* \end{matrix} \right\}\\
&=F(\boldu^{L},\boldu^{R},L) +(q^L-1) F(\boldu^{L},\boldu^{R+},L-1).
\end{split}
\end{equation}
The definition (\ref{Tilde1Def}) then yields the required (\ref{FermLikeEq}).
For the subcase for which $d^R=d^R_0=1$,
Lemma \ref{RecurseLem}(1) and Theorem \ref{Core2Thrm}
lead to (\ref{FermLike4Eq}), (\ref{FermLike5Eq}) and (\ref{FermLikeEq})
in a similar fashion after noting that
$\boldnu=\boldnu^*=\boldnu^+=\boldnu^{+*}=()$,
and $0<\sigma<t$ implies that $1\le b\pm1<p'$.

The proof for the $p'<2p$ cases differs from that of the $p'<2p$ cases
above only in that parts 3) and 4) of Lemma \ref{RecurseLem}
are used as appropriate, yielding different coefficients,
but yielding (\ref{FermLikeEq}) as required.
\cqfd

\begin{note}\label{SpuriousNote}
As indicated in the above proof, if $0=\sigma<\tau-1$ then the
polynomials defined by the first and second cases of (\ref{Tilde1Def})
are actually equal.
The same comment applies to the first and second cases of (\ref{Tilde2Def}).
\end{note}
\medskip

We prove the following result, which is not needed until
Section \ref{FermCharSec}, along lines similar to the proof of
Theorem \ref{FermLikeThrm} above.

\begin{lemma}\label{FermLikeBitsLem}
Let all parameters be as in Section \ref{CoreSec}.

1) Let $p'>2p$, let $0\le\sigma<\tau-1$ and
let $b-\Delta^R_{d^R}$ be non-interfacial.
If $b$ is interfacial then:
\begin{equation*}
F(\boldu^{L},\boldu^{R},L)
= \widetilde F(\boldu^{L},\boldu^{R+},L+1)
-q^{\frac12(L+2+\Delta(a-b))} F(\boldu^{L},\boldu^{R++},L).
\end{equation*}
If $b$ is non-interfacial then:
\begin{equation*}
\widetilde F(\boldu^{L},\boldu^{R},L)
= \widetilde F(\boldu^{L},\boldu^{R+},L+1)
-q^{\frac12(L+2+\Delta(a-b))} F(\boldu^{L},\boldu^{R++},L).
\end{equation*}

2) If $p'<2p$ and $\sigma_d=0$ then
\begin{equation*}
\widetilde F(\boldu^{L},\boldu^{R},L)
= F(\boldu^{L},\boldu^{R+},L-1).
\end{equation*}
\end{lemma}

\Proof If $d^R>d^R_0$, note that $\boldnu\ne()\ne\boldnu^*$. 
Then, the definitions of Section \ref{CoreSec} yield
$\nu_{d^R-1}=b+\Delta\kappa_{\sigma}$
and $\nu^*_{d^R-1}=b-\Delta(\kappa_{\tau}-\kappa_{\sigma})$.
If $d^R=d^R_0$, note that $\boldnu=\boldnu^*=()$. 

1) Set $b'=b-\Delta$.
When $d^R>d^R_0$ we have $b'-\nu_{d^R-1}=-\Delta(\kappa_{\sigma}+1)$ and
$b'-\nu^*_{d^R-1}=\Delta(\kappa_{\tau}-\kappa_{\sigma}-1)$.
In particular, $b'$ is strictly between $\nu_{d^R-1}$ and $\nu^*_{d^R-1}$.
The same is true of $b'-\Delta$ unless $\sigma+2=\tau\le t_1+1$,
in which case $b'-\Delta=\nu^*_{d^R-1}$.
Substituting $b\to b'$ and $L\to L+1$ into the first identity of
Lemma \ref{RecurseLem}(1) and rearranging, yields:
\begin{equation}\label{FermLikeBits1Eq}
\begin{split}
&\ochi^{p,p'}_{a,b'+\Delta,b'}(L)
\left\{
\begin{matrix} \boldmu^{\phantom{*}};\boldnu^{\phantom{*}}\\
                    \boldmu^*;\boldnu^* \end{matrix} \right\} \\
&\hskip3mm=
\ochi^{p,p'}_{a,b',b'-\Delta}(L+1)
\left\{
\begin{matrix} \boldmu^{\phantom{*}};\boldnu^{\phantom{*}}\\
                    \boldmu^*;\boldnu^* \end{matrix} \right\}
-
q^{\frac12(L+1+\Delta(a-b'))}
\ochi^{p,p'}_{a,b'-\Delta,b'}(L)
\left\{
\begin{matrix} \boldmu^{\phantom{*}};\boldnu^{+\phantom{*}}\\
                    \boldmu^*;\boldnu^{+*} \end{matrix} \right\}\!.
\hskip-1.5mm
\end{split}
\end{equation}
where if $b'-\Delta\ne\nu^*_{d^R-1}$ then $\boldnu^+=\boldnu$ and
$\boldnu^{+*}=\boldnu^*$;
and if $b'-\Delta=\nu^*_{d^R-1}$ then $d^+\ge d^R_0$ is the smallest value
such that $b'-\Delta=\nu^*_{d^+}$, and then
$\boldnu^+=(\nu_{d^R_0},\ldots,\nu_{d^+-1})$ and
$\boldnu^{+*}=(\nu^*_{d^R_0},\ldots,\nu^*_{d^+-1})$.
If $b$ is interfacial then the left side of (\ref{FermLikeBits1Eq})
is equal to $F(\boldu^{L},\boldu^{R},L)$ by Theorem \ref{Core2Thrm},
whereas if $b$ is non-interfacial, it is equal to
$\widetilde F(\boldu^{L},\boldu^{R},L)$ by Theorem \ref{FermLikeThrm}.
The first term on the right is equal to
$\widetilde F(\boldu^{L},\boldu^{R+},L+1)$ by Theorem \ref{FermLikeThrm},
after noting that with $\run^{R+}$ in place of $\run^{R}$, the definitions
of Section \ref{CoreSec} lead to the same mazy-pair $\boldnu,\boldnu^*$.
Regardless of whether $b'-\Delta\ne\nu^*_{d^R-1}$ or $b'-\Delta=\nu^*_{d^R-1}$
(see the corresponding argument which applies to the second term on
the right of (\ref{FermLike2Eq}) in the proof of Theorem \ref{FermLikeThrm}),
the second term on the right is equal to $F(\boldu^{L},\boldu^{R++},L)$
via Theorem \ref{Core2Thrm}.
The required result follows in the $d^R>d^R_0$ case.

In the subcase for which $d^R=d^R_0=1$, we have
$\boldnu=\boldnu^*=\boldnu^+=\boldnu^{+*}=()$.
Lemma \ref{RecurseLem}(1) and Theorems \ref{Core2Thrm} and \ref{FermLikeThrm}
then lead to (\ref{FermLikeBits1Eq}) and the required result
in a similar fashion to the $b'-\Delta\ne\nu^*_{d^R-1}$ case considered
above after noting that $\sigma<t-1$ ensures that $1\le b',b'-\Delta<p'$.

2) In the case for which $d^R>d^R_0$ we have $b+\Delta=\nu_{d^R-1}$ here,
whereupon the second identity of Lemma \ref{RecurseLem}(4) gives:
\begin{equation}\label{FermLikeBits2Eq}
\ochi^{p,p'}_{a,b,b-\Delta}(L)
\left\{
\begin{matrix} \boldmu^{\phantom{*}};\boldnu^{\phantom{*}}\\
                    \boldmu^*;\boldnu^* \end{matrix} \right\}
=\ochi^{p,p'}_{a,b-\Delta,b}(L-1)
\left\{
\begin{matrix} \boldmu^{\phantom{*}};\boldnu^{+\phantom{*}}\\
                    \boldmu^*;\boldnu^{+*} \end{matrix} \right\}.
\end{equation}
where if $b-\Delta\ne\nu^*_{d^R-1}$ then $\boldnu^+=\boldnu$ and
$\boldnu^{+*}=\boldnu^*$;
and if $b-\Delta=\nu^*_{d^R-1}$ then $d^+\ge d^R_0$ is the smallest value
such that $b-\Delta=\nu^*_{d^+}$, and then
$\boldnu^+=(\nu_{d^R_0},\ldots,\nu_{d^+-1})$ and
$\boldnu^{+*}=(\nu^*_{d^R_0},\ldots,\nu^*_{d^+-1})$.
The left side of (\ref{FermLikeBits2Eq}) is equal to
$\widetilde F(\boldu^{L},\boldu^{R},L)$ by Theorem \ref{FermLikeThrm}.
If $b-\Delta\ne\nu^*_{d^R-1}$, the right side of (\ref{FermLikeBits2Eq})
is equal to $F(\boldu^{L},\boldu^{R+},L-1)$ by Theorem \ref{Core2Thrm},
after noting that with $\run^{R+}$ in place of $\run^{R}$, the definitions
of Section \ref{CoreSec} lead to the same mazy-pair $\boldnu,\boldnu^*$.
If $b-\Delta=\nu^*_{d^R-1}$, the argument which applies to the second term on
the right of (\ref{FermLike4Eq}) in the proof of Theorem \ref{FermLikeThrm},
equally applies here to show that the term on the right of
(\ref{FermLikeBits2Eq}) is equal to $F(\boldu^{L},\boldu^{R+},L-1)$
via Theorem \ref{Core2Thrm}.
The required result follows in this $d^R>d^R_0$ case.

In the case for which $d^R=d^R_0$, we have $b+\Delta\in\{0,p'\}$
whereupon Lemma \ref{RecurseLem}(4) again implies (\ref{FermLikeBits2Eq}).
Again, Theorems \ref{Core2Thrm} and \ref{FermLikeThrm} yield the
required result.
\cqfd

\subsection{The $\boldm\boldn$-system}\label{MN2sysSec}

Equation (\ref{MNmatEq1}) defines the vector $\boldn$ in terms of
$\hat{\boldm}$ and $\boldu=\boldu^L+\boldu^R$. The summands of
(\ref{FEq}) have $m_i\equiv Q_i$ for $0\le i<t$ whereupon, via
(\ref{ParityDef}), it follows that $n_i\in\Z$ for $1\le i\le t$.
Moreover, on expressing (\ref{MNmatEq1}) in the form
$n_i=-\frac12(\boldC^*\hat{\boldm}-\boldu^L-\boldu^R)_i$,
we see that the non-zero terms in (\ref{FEq}) have each $n_i\ge0$.
(On examining the proof of Lemma \ref{CoreIndLem}, we find that
$n_i$ is the number of particles inserted at the $i$th stage of
the induction performed there.)

Equation (\ref{Particles1Eq}) follows immediately from the
following lemma.

\begin{lemma} For $0\le j<t$,
\begin{displaymath}
\sum_{i=1}^t l_iC_{ij}=-\delta_{j0}.
\end{displaymath}
\end{lemma}

\Proof In view of the definition (\ref{CijDef}), only $i=j-1$,
$i=j$ and $i=j+1$ contribute in the summation.
Since $C_{j+1,j}=-1$, the lemma holds for $j=0$.

Now consider $j\ge1$ and let $k$ be such that $t_k<j\le t_{k+1}$.
The definition of $\{l_i\}_{i=1}^t$ shows that if $j<t_{k+1}$ then
$l_{j+1}=l_j+y_k$ and if $j=t_{k+1}$ then $l_{j+1}=y_k$.
Since $C_{jj}=2$ in the former case and $C_{jj}=1$ in the latter case,
$l_jC_{jj}+l_{j+1}C_{j+1,j}=l_j-y_k$ in both cases.
For $j=1$, we readily find $l_1=y_k=1$ thus giving the required result.

For $j>1$, a similar argument shows that 
if $j-1>t_k$ then $l_j=l_{j-1}+y_k$ and
if $j-1=t_k$ then $l_j=y_{k-1}$.
In the former case, $C_{j-1,j}=-1$, whereupon
$\sum_{i=j-1}^{j+1} l_iC_{i,j}=0$.
In the latter case, $C_{j-1,j}=1$, whereupon
$\sum_{i=j-1}^{j+1} l_iC_{i,j}=l_{j-1}+y_{k-1}-y_k=0$,
with the final equality following because
$l_{t_k}=y_{k-2}+(t_k-t_{k-1}-1)y_{k-1}=y_k-y_{k-1}$.
\cqfd

\noindent
The summation over $\boldm\equiv\boldQ(\boldu^L+\boldu^R)$ in
(\ref{FEq}) may therefore be carried out by finding all solutions
to (\ref{Particles2Eq}) with $n_i\in\Z_{\ge0}$ for $1\le i\le t$
(the set of solutions is clearly finite), and then using (\ref{MNmatEq2})
to obtain the corresponding values $(m_1,m_2,\ldots,m_{t-1})$.

\newpage

\setcounter{section}{7}

\section{Collating the runs}\label{CollateSec}

In this section, we show that Takahashi trees for $a$ and $b$ enable
us to obtain a set of $F(\boldu^L,\boldu^R,L)$ whose sum is
$\chi^{p,p'}_{a,b,c}(L)$.

\subsection{Pulling the drawstrings}

This section provides a few preliminary results that are required
to make use of the Takahashi trees.

\begin{lemma}\label{GatherLem}
Let $\boldmu,\boldmu^*,\boldnu,\boldnu^*$ be a
mazy-four in the $(p,p')$-model sandwiching $(a,b)$.

If $d^L\ge2$ then
\begin{equation*}
\begin{split}
\ochi^{p,p'}_{a,b,c}&(L)\left\{
\begin{matrix}
  \mu_1,\ldots,\mu_{d^L-1};\boldnu^{\phantom{*}}\\
  \mu^*_1,\ldots,\mu^*_{d^L-1};\boldnu^*
\end{matrix} \right\}\\
&=
\ochi^{p,p'}_{a,b,c}(L)\left\{
\begin{matrix}
  \mu_1,\ldots,\mu_{d^L-1},\mu_{d^L};\boldnu^{\phantom{*}}\\
  \mu^*_1,\ldots,\mu^*_{d^L-1},\mu^*_{d^L};\boldnu^*
\end{matrix} \right\}
+
\ochi^{p,p'}_{a,b,c}(L)\left\{
\begin{matrix}
  \mu_1,\ldots,\mu_{d^L-1},\mu^*_{d^L};\boldnu^{\phantom{*}}\\
  \mu^*_1,\ldots,\mu^*_{d^L-1},\mu_{d^L};\boldnu^*
\end{matrix} \right\}.
\end{split}
\end{equation*}
If $d^R\ge2$ then
\begin{equation*}
\begin{split}
\ochi^{p,p'}_{a,b,c}&(L)\left\{
\begin{matrix}
  \boldmu^{\phantom{*}};\nu_1,\ldots,\nu_{d^R-1}\\
  \boldmu^*;\nu^*_1,\ldots,\nu^*_{d^R-1}
\end{matrix} \right\}\\
&=
\ochi^{p,p'}_{a,b,c}(L)\left\{
\begin{matrix}
  \boldmu^{\phantom{*}};\nu_1,\ldots,\nu_{d^R-1},\nu_{d^R}\\
  \boldmu^*;\nu^*_1,\ldots,\nu^*_{d^R-1},\nu^*_{d^R}
\end{matrix} \right\}
+
\ochi^{p,p'}_{a,b,c}(L)\left\{
\begin{matrix}
  \boldmu^{\phantom{*}};\nu_1,\ldots,\nu_{d^R-1},\nu^*_{d^R}\\
  \boldmu^*;\nu^*_1,\ldots,\nu^*_{d^R-1},\nu_{d^R}
\end{matrix} \right\}.
\end{split}
\end{equation*}
\end{lemma}

\Proof In the first case, since $a$ is between $\mu_{d^L}$ and $\mu^*_{d^L}$,
and $\mu_{d^L}$ and $\mu^*_{d^L}$ are both between
$\mu_{{d^L}-1}$ and $\mu^*_{{d^L}-1}$,
it follows that each path $h$ with $h_0=a$ that attains $\mu^*_{{d^L}-1}$,
necessarily attains $\mu_{d^L}$ or $\mu^*_{d^L}$. Thence:
\begin{equation*}
\begin{split}
\P^{p,p'}_{a,b,c}&(L)\left\{
\begin{matrix}
  \mu_1,\ldots,\mu_{d^L-1};\boldnu^{\phantom{*}}\\
  \mu^*_1,\ldots,\mu^*_{d^L-1};\boldnu^*
\end{matrix} \right\}\\
&=
\P^{p,p'}_{a,b,c}(L)\left\{
\begin{matrix}
  \mu_1,\ldots,\mu_{d^L-1},\mu_{d^L};\boldnu^{\phantom{*}}\\
  \mu^*_1,\ldots,\mu^*_{d^L-1},\mu^*_{d^L};\boldnu^*
\end{matrix} \right\}
\cup
\P^{p,p'}_{a,b,c}(L)\left\{
\begin{matrix}
  \mu_1,\ldots,\mu_{d^L-1},\mu^*_{d^L};\boldnu^{\phantom{*}}\\
  \mu^*_1,\ldots,\mu^*_{d^L-1},\mu_{d^L};\boldnu^*
\end{matrix} \right\}.
\end{split}
\end{equation*}
Then, since the latter two sets are disjoint, the first part of the
lemma follows. The second part follows in an analogous way.
\cqfd
\medskip

In the above lemma, it is essential that $d^L\ge2$ for the first expression,
and $d^R\ge2$ for the second expression, since the two expressions
don't necessarily hold if $d^L=1$ and $d^R=1$ respectively.
The following two lemmas deal with that extra case.

\begin{lemma}\label{GatherLast1Lem}
Let $0\le \mu^*<a<\mu\le p'$ and $0\le \nu^*<b<\nu\le p'$.
If $\mu\le\nu^*$ or $\nu\le\mu^*$ then:
\begin{equation*}
\begin{split}
\ochi^{p,p'}_{a,b,c}(L)
=
\ochi^{p,p'}_{a,b,c}(L)&\left\{
\begin{matrix}
  \mu^{\phantom*};\nu^{\phantom{*}}\\
  \mu^*;\nu^*
\end{matrix} \right\}
+
\ochi^{p,p'}_{a,b,c}(L)\left\{
\begin{matrix}
  \mu^{\phantom*};\nu^*\\
  \mu^*;\nu^{\phantom{*}}
\end{matrix} \right\} \\
&+
\ochi^{p,p'}_{a,b,c}(L)\left\{
\begin{matrix}
  \mu^*;\nu^{\phantom{*}}\\
  \mu^{\phantom*};\nu^*
\end{matrix} \right\}
+
\ochi^{p,p'}_{a,b,c}(L)\left\{
\begin{matrix}
  \mu^*;\nu^*\\
  \mu^{\phantom*};\nu^{\phantom{*}}
\end{matrix} \right\}.
\end{split}
\end{equation*}
\end{lemma}

\Proof
Let $h\in\P^{p,p'}_{a,b,c}(L)$ whence $h_0=a$ and $h_L=b$.
In the first case $a<\mu\le\nu^*<b$, so that $h_i=\mu$
for some $i$ and $h_{i'}=\nu^*$ for some $i'$.
Thus $h$ is certainly an element of one of
$\P^{p,p'}_{a,b,c}(L)\left\{
\begin{smallmatrix}
  \mu^{\phantom*};\nu^{\phantom{*}}\\
  \mu^*;\nu^*
\end{smallmatrix} \right\}$,
$\P^{p,p'}_{a,b,c}(L)\left\{
\begin{smallmatrix}
  \mu^{\phantom*};\nu^*\\
  \mu^*;\nu^{\phantom{*}}
\end{smallmatrix} \right\}$,
$\P^{p,p'}_{a,b,c}(L)\left\{
\begin{smallmatrix}
  \mu^*;\nu^{\phantom{*}}\\
  \mu^{\phantom*};\nu^*
\end{smallmatrix} \right\}$,
and
$\P^{p,p'}_{a,b,c}(L)\left\{
\begin{smallmatrix}
  \mu^*;\nu^*\\
  \mu^{\phantom*};\nu^{\phantom{*}}
\end{smallmatrix} \right\}$.
Since these four sets are disjoint, the lemma follows when
$a<\mu\le \nu^*<b$.

When $\nu\le \mu*$, the lemma follows in an analogous way.
\cqfd
\medskip

\begin{lemma}\label{GatherLast2Lem}
Let $0\le \mu^*<a,b<\mu\le p'$ and $\hat p'=\mu-\mu^*$.
If $\hat p$ is such that the $s$th band of the $(\hat p,\hat p')$-model
is of the same parity as the $(s+\mu^*)$th band of the $(p,p')$-model
for $1\le s\le\hat p'-2$,%
\footnote{In general, there is no guarantee that such a $\hat p$ exists.}
then:
\begin{equation*}
\begin{split}
\ochi^{p,p'}_{a,b,c}(L)
=
\ochi^{p,p'}_{a,b,c}(L)&\left\{
\begin{matrix}
  \mu^{\phantom*};\mu^{\phantom{*}}\\
  \mu^*;\mu^*
\end{matrix} \right\}
+
\ochi^{p,p'}_{a,b,c}(L)\left\{
\begin{matrix}
  \mu^{\phantom*};\mu^*\\
  \mu^*;\mu^{\phantom{*}}
\end{matrix} \right\} \\
&+
\ochi^{p,p'}_{a,b,c}(L)\left\{
\begin{matrix}
  \mu^*;\mu^{\phantom{*}}\\
  \mu^{\phantom*};\mu^*
\end{matrix} \right\}
+
\ochi^{p,p'}_{a,b,c}(L)\left\{
\begin{matrix}
  \mu^*;\mu^*\\
  \mu^{\phantom*};\mu^{\phantom{*}}
\end{matrix} \right\}
+
\ochi^{\hat p,\hat p'}_{\hat a,\hat b,\hat c}(L),
\end{split}
\end{equation*}
where $\hat a=a-u^*$, $\hat b=b-u^*$ and
\begin{equation*}
\hat c= \begin{cases}
       2&\text{if $c=u^*>0$ and 
          $\left\lfloor\frac{(b+1)p}{p'}\right\rfloor=
           \left\lfloor\frac{(b-1)p}{p'}\right\rfloor+1$;}\\[1.5mm]
       \hat p'-2&\text{if $c=u<p'$ and 
          $\left\lfloor\frac{(b+1)p}{p'}\right\rfloor=
           \left\lfloor\frac{(b-1)p}{p'}\right\rfloor+1$;}\\[1.5mm]
       c-u^* &\text{otherwise}.
         \end{cases} 
\end{equation*}
\end{lemma}

\Proof
Let $h\in\P^{p,p'}_{a,b,c}(L)$ whence $h_0=a$ and $h_L=b$.
If $h_i=\mu$ for some $i$ with $0\le i\le L$,
or $h_{i'}=\mu^*$ for some $i'$ with $0\le i'\le L$,
then $h$ is an element of one of the disjoint sets
$\P^{p,p'}_{a,b,c}(L)\left\{
\begin{smallmatrix}
  \mu^{\phantom*};\mu^{\phantom{*}}\\
  \mu^*;\mu^*
\end{smallmatrix} \right\}$,
$\P^{p,p'}_{a,b,c}(L)\left\{
\begin{smallmatrix}
  \mu^{\phantom*};\mu^*\\
  \mu^*;\mu^{\phantom{*}}
\end{smallmatrix} \right\}$,
$\P^{p,p'}_{a,b,c}(L)\left\{
\begin{smallmatrix}
  \mu^*;\mu^{\phantom{*}}\\
  \mu^{\phantom*};\mu^*
\end{smallmatrix} \right\}$,
and
$\P^{p,p'}_{a,b,c}(L)\left\{
\begin{smallmatrix}
  \mu^*;\mu^*\\
  \mu^{\phantom*};\mu^{\phantom{*}}
\end{smallmatrix} \right\}$.

Otherwise $\mu^*<h_i<\mu$ for $0\le i\le L$, whereupon setting
$h'_i=h_i-\mu^*$ for $0\le i\le L$, defines an element of
$\P^{\hat p,\hat p'}_{\hat a,\hat b,c-\mu^*}(L)$.
Clearly, every element of
$\P^{\hat p,\hat p'}_{\hat a,\hat b,c-\mu^*}(L)$
arises uniquely in this way.
The sequence of band parities between heights
$\mu^*+1$ and $\mu-1$ of the $(p,p')$-model being identical to that
between heights $1$ and $\hat p'-1$ of the $(\hat p,\hat p')$-model,
guarantees that the contribution of the $i$th vertex to the weight
of the path is identical in the two cases $h$ and $h'$,
except possibly when $i=L$ and either $c=u^*$ or $c=u$.
After noting the assignment of parities to the $0$th and
$(\hat p'-1)$th bands as specified in Note \ref{ExtraBandsNote},
we see that this discrepancy is rectified by switching the direction
of the final vertex whenever $b$ is interfacial and either
$c=u^*$ or $c=u$.
Then $\owt(h')=\owt(h)$. The lemma follows.
\cqfd
\medskip

\begin{note} Lemmas \ref{GatherLast1Lem} and \ref{GatherLast2Lem}
will be employed with $\mu^*=0$, $\mu=p'$, $\nu^*=0$ or $\nu=p'$.
Then, at least two of the terms on the right side of the
expressions in Lemmas \ref{GatherLast1Lem} and \ref{GatherLast2Lem}
are automatically zero. 
\end{note}

\subsection{Constraints of the Takahashi tree}\label{ConTreeSec}

Let $1\le a<p'$ and consider the Takahashi tree for $a$ that
is described in Section \ref{FormSec}.
Let $a_{i_1i_2\cdots i_d}$ be a particular leaf-node.
{}From this leaf-node, obtain
$\{\tau_j,\sigma_j,\Delta_j\}_{j=1}^d$ as in Section \ref{FormVSec}.
The next result shows that $\{\tau_j,\sigma_j,\Delta_j\}_{j=1}^d$
satisfies the definition of a naive run that was given
in Section \ref{CoreSec}.
We will then refer to $\{\tau_j,\sigma_j,\Delta_j\}_{j=1}^d$
as the naive run corresponding to (the leaf-node)
$a_{i_1i_2\cdots i_d}$.

\begin{lemma}\label{NaiveRunLem}
For $2\le j<d$, if $\Delta_j=\Delta_{j+1}$ then
$\sigma_j\le\tau_j$, and if $\Delta_j\ne\Delta_{j+1}$ then
$\sigma_j<\tau_j$. Also $\sigma_d<\tau_d$.

Let $d>1$. If $2\le j<d$, or both $j=1$ and $\sigma_1\le t_n$, then
\begin{displaymath}
\begin{array}{ll}
\tau_{j+1}=t_{\zeta(\sigma_j-1)}
 &\mbox{ if } \Delta_{j+1}=\Delta_j;\\
\tau_{j+1}=t_{\zeta(\sigma_j)}
 &\mbox{ if } \Delta_{j+1}\ne\Delta_j,
\end{array}
\end{displaymath}
and if $t_n<\sigma_1<t$ then
\begin{displaymath}
\begin{array}{ll}
\tau_{2}=t_{n-1}
 &\mbox{ if } \Delta_{2}=\Delta_1;\\
\tau_{2}=t_{n}-1
 &\mbox{ if } \Delta_{2}\ne\Delta_1\mbox{ and }c_n>1;\\
\tau_{2}=t_{n-2}
 &\mbox{ if } \Delta_{2}\ne\Delta_1\mbox{ and }c_n=1.
\end{array}
\end{displaymath}
\end{lemma}

\Proof 
For $2\le j<d$, the description of Section \ref{FormSec}
implies that there exists $x$ such that:
\begin{equation}\label{RunProofEq1}
\begin{split}
a_{i_1i_2\cdots i_{j-1}}-a_{i_1i_2\cdots i_{j-1}i_{j-1}}
&=(-1)^{i_{j-1}}\kappa_x;\\
a_{i_1i_2\cdots i_{j-1}}-a_{i_1i_2\cdots i_{j-1}\bari_{j-1}}
&=(-1)^{i_{j-1}}\kappa_{x+1}.
\end{split}
\end{equation}
The prescription of Section \ref{FormVSec} then implies that
$\sigma_j=x$ if $i_j\ne i_{j-1}$, and
$\sigma_j=x+1$ if $i_j= i_{j-1}$.
Since $\kappa_{\tau_j}=a_{i_1i_2\cdots i_{j-2}1}-a_{i_1i_2\cdots i_{j-2}0}$,
$a_{i_1i_2\cdots i_{j-2}0}<a<a_{i_1i_2\cdots i_{j-2}1}$,
and
$a_{i_1i_2\cdots i_{j-1}0}<a<a_{i_1i_2\cdots i_{j-1}1}$,
it follows that $\kappa_{x+1}\le\kappa_{\tau_j}$ and thus
that $x+1\le\tau_j$.
It now follows that if $\Delta_{j+1}\ne\Delta_{j}$ then $\sigma_j=x<\tau_j$,
and if $\Delta_{j+1}=\Delta_{j}$ then $\sigma_j=x+1\le\tau_j$.
That $\sigma_d\le\tau_d$ follows similarly: the case
$\sigma_d=\tau_d$ is then excluded since otherwise
$a_{i_1i_2\cdots i_{d-2}\bari_{d-1}}=a_{i_1i_2\cdots i_{d}}=a$,
and the former of these would label a leaf-node.

The prescription of Section \ref{FormVSec} also implies,
via (\ref{RunProofEq1}), that
$\kappa_{\tau_{j+1}}=\kappa_{x+1}-\kappa_{x}=y_{\zeta(x)}=
\kappa_{t_{\zeta(x)}}$, so that $\tau_{j+1}=t_{\zeta(x)}$.
This immediately yields:
\begin{equation*}
\tau_{j+1}=\left\{
  \begin{array}{ll}
    t_{\zeta(\sigma_j-1)}&\mbox{ if }\Delta_{j+1}=\Delta_j;\\
    t_{\zeta(\sigma_j)}&\mbox{ if }\Delta_{j+1}\ne\Delta_j,
  \end{array}
\right.
\end{equation*}
as required in this $2\le j<d$ case.

If $a<y_{n}+y_{n-1}=\kappa_{t_n+1}$ then
$a_0=\kappa_x$ and $a_1=\kappa_{x+1}$ for some $x$,
whereupon $\tau_2=t_{\zeta(x)}$.
So with $\Delta_2=-(-1)^{i_1}$ and $a_{\bari_1}=\kappa_{\sigma_1}$,
we obtain:
\begin{displaymath}
\tau_{2}=\left\{
  \begin{array}{ll}
    t_{\zeta(\sigma_1-1)}&\mbox{ if }\Delta_{2}=-1;\\
    t_{\zeta(\sigma_1)}&\mbox{ if }\Delta_{2}=+1.
  \end{array}
\right.
\end{displaymath}
Since $\Delta_1=-1$ in this case, this is as required. Note
that there is exactly one case here where $\sigma_1\not\le t_n$,
and that is when $i_1=0$ and $\sigma_1=t_n+1$, whereupon
$\Delta_2=-1$ and $\zeta(\sigma_1-1)=n-1$ as required.

If $a>(c_n-1)y_{n}=p'-\kappa_{t_n+1}$ then
$a_0=p'-\kappa_{x+1}$ and $a_1=p'-\kappa_{x}$ for some $x$,
whereupon $\tau_2=t_{\zeta(x)}$.
So with $\Delta_2=-(-1)^{i_1}$ and $a_{\bari_1}=p'-\kappa_{\sigma_1}$,
we obtain:
\begin{displaymath}
\tau_{2}=\left\{
  \begin{array}{ll}
    t_{\zeta(\sigma_1-1)}&\mbox{ if }\Delta_{2}=+1;\\
    t_{\zeta(\sigma_1)}&\mbox{ if }\Delta_{2}=-1.
  \end{array}
\right.
\end{displaymath}
Since $\Delta_1=+1$ in this case, this is as required. Note
that there is exactly one case here where $\sigma_1\not\le t_n$,
and that is when $i_1=1$ and $\sigma_0=t_n+1$, whereupon
$\Delta_2=+1$ and $\zeta(\sigma_1-1)=n-1$ as required.

If $ky_{n}+y_{n-1}<a<(k+1)y_{n}$ for $1\le k\le c_n-2$
then $a_0=\kappa_{t_n+k}$ and $a_1=p'-\kappa_{t_{n+1}-k-1}$
so that $a_0\in{\mathcal T}$ and $a_1\in{\mathcal T}^\prime$.
Then:
\begin{displaymath}
\kappa_{\tau_{2}}=y_{n}-y_{n-1}=\left\{
  \begin{array}{ll}
    \kappa_{t_{n}-1} &\mbox{ if } c_n>1;\\
    \kappa_{t_{n-2}} &\mbox{ if } c_n=1.
  \end{array}
\right.
\end{displaymath}
We now claim that $\Delta_1\ne\Delta_2$ in this case as required.
This follows because if
$i_1=0$ then $\Delta_1=+1$ and $\Delta_2=-(-1)^0=-1$ and if
$i_1=1$ then $\Delta_1=-1$ and $\Delta_2=-(-1)^1=+1$.

If $ky_{n}<a<ky_{n}+y_{n-1}$ for $2\le k\le c_n-2$
then $a_0=p'-\kappa_{t_{n+1}-k}$ and $a_1=\kappa_{t_{n}+k}$
so that $a_0\in{\mathcal T}^\prime$ and $a_1\in{\mathcal T}$.
Then $\kappa_{\tau_{2}}=y_{n-1}=\kappa_{t_{n-1}}$.
We now claim that $\Delta_1=\Delta_2$ in this case as required.
This follows because if
$i_1=0$ then $\Delta_1=-1$ and $\Delta_2=-(-1)^0=-1$ and if
$i_1=1$ then $\Delta_1=+1$ and $\Delta_2=-(-1)^1=+1$.
\cqfd
\medskip

\begin{lemma}\label{NodesMuLem}
Let $\{\tau_j,\sigma_j,\Delta_j\}_{j=1}^{d}$
be the naive run corresponding to the leaf-node
$a_{i_1i_2\cdots i_{d-1}0}$ of the Takahashi tree for $a$.
Then
\begin{align*}
a_{i_1i_2\cdots i_{k-1}\bari_k}&= \sum_{m=2}^{k}
        \Delta_m (\kappa_{\tau_m^{}}-\kappa_{\sigma_m^{}})\quad
+\quad\begin{cases}
        \kappa_{\sigma^{}_1} &\text{if $\Delta_1=-1$;}\\
        p'-\kappa_{\sigma^{}_1} &\text{if $\Delta_1=+1$,}
  \end{cases} \\
\intertext{for $1\le k\le d$, and}
a_{i_1i_2\cdots i_{k-1}i_{k}}&=
a_{i_1i_2\cdots i_{k-1}\bari_k}+\Delta_{k+1}\kappa_{\tau^{}_{k+1}},
\end{align*}
for $1\le k<d$.
\end{lemma}

\Proof 
The prescription of Section \ref{FormVSec} implies that
$\kappa_{\tau^{}_{k+1}}=
a_{i_1i_2\cdots i_{k-1}1}-a_{i_1i_2\cdots i_{k-1}0}$, whence
$\Delta_{k+1}=-(-1)^{i_k}$ implies that
$a_{i_1i_2\cdots i_{k}}=a_{i_1i_2\cdots \bari_{k}}
  +\Delta_{k+1}\kappa_{\tau^{}_{k+1}}$ for $1\le k<d$,
thus giving the second expression.

By definition, $a_{\bari_1}=\kappa_{\sigma^{}_1}$ if
$a_{\bari_1}\in{\mathcal T}$ and thus $\Delta_1=-1$;
and $a_{\bari_1}=p'-\kappa_{\sigma^{}_1}$ if
$a_{\bari_1}\in{\mathcal T}'$ and thus $\Delta_1=+1$.
This gives the $k=1$ case of the first expression.

The description of Section \ref{FormSec} implies that
$-(-1)^{i_{k-1}}(a_{i_1i_2\cdots i_{k-1}}-a_{i_1i_2\cdots i_{k-1}i})>0$
for $i\in\{0,1\}$ and $2\le k\le d$.
When $i=\bari_k$, Section \ref{FormVSec} specifies that this
expression equals $\kappa_{\sigma^{}_k}$.
Thereupon, we obtain
$a_{i_1i_2\cdots i_{k-1}\bari_k}
=a_{i_1i_2\cdots i_{k-1}}-\Delta_k\kappa_{\sigma^{}_k}$.
The lemma then follows by induction.
\cqfd
\medskip

\subsection{Gathering in the Takahashi tree}\label{GatherTree}

As described in Section \ref{FormVSec},
each leaf-node $a_{i_1i_2\cdots i_{d-1}0}$ of the Takahashi tree for
$a$ gives rise firstly to the naive run
$\{\tau_j,\sigma_j,\Delta_j\}_{j=1}^{d}$
and then to a vector $\boldu\in\mathcal U(a)$ by (\ref{uEq}).
For the particular leaf-node $a_{i_1i_2\cdots i_{d-1}0}$,
we will denote this corresponding vector by
$\boldu(a;i_1,i_2,\ldots,i_{d-1})$.
Note that
$\Delta(\boldu(a;i_1,i_2,\ldots,i_{d-1}))=\Delta_{d}$.

For vectors $\boldu(a;i_1,i_2,\ldots,i_{d^L-1})$ and
$\boldu(b;j_1,j_2,\ldots,j_{d^R-1})$ arising from leaf-nodes
of the Takahashi trees for $a$ and $b$ respectively, the following lemma
identifies a set of paths for which
$F(\boldu(a;i_1,i_2,\ldots,i_{d^L-1}),
        \boldu(b;j_1,j_2,\ldots,j_{d^R-1}),L)$ is the generating function.
The lemma makes use of the following definition:
\begin{displaymath}
\begin{array}{ll}
\displaystyle
\Xi_0(s)&=\quad
\displaystyle
\left\{
  \begin{array}{ll}
  \displaystyle
  0 \quad &\mbox{if }s<\kappa_{t_n};\\[1mm]
  \displaystyle
  p'\phantom{\displaystyle{}-\kappa_{t_n}}\quad &\mbox{if }s>p'-\kappa_{t_n};
  \end{array}
\right.\\[5mm]
\displaystyle
\Xi_1(s)&=\quad
\displaystyle
\left\{
  \begin{array}{ll}
  \displaystyle
  \kappa_{t_n} \quad &\mbox{if }s<\kappa_{t_n};\\[1mm]
  \displaystyle
  p'-\kappa_{t_n}\quad &\mbox{if }s>p'-\kappa_{t_n}.
  \end{array}
\right.
\end{array}
\end{displaymath}
(Outside of the ranges given here, $\Xi_0(s)$ and $\Xi_1(s)$ are
not defined and not needed.)

\begin{lemma}\label{OneNodeLem}
Let $1\le a,b<p'$,
let $a_{i_1i_2\cdots i_{d^L-1}0}$
be a leaf-node of the Takahashi tree of $a$, and
let $b_{j_1j_2\cdots j_{d^R-1}0}$
be a leaf-node of the Takahashi tree of $b$.
If $b$ is interfacial in the $(p,p')$-model, set $c\in\{b\pm1\}$.
If $b$ is not interfacial in the $(p,p')$-model,
set $c=b+\Delta(\boldu(b;j_1,j_2,\ldots,j_{d^R-1}))$.
\begin{enumerate}
\item If $\kappa_{t_n}\le a_{\bari_1}\le p'-\kappa_{t_n}$
      and $\kappa_{t_n}\le b_{\barj_1}\le p'-\kappa_{t_n}$
      then:
\begin{displaymath}
\begin{array}{l}
\displaystyle
\chi^{p,p'}_{a,b,c}(L)\left\{
{a_{i_1},\atop a_{\bari_1},}
{a_{i_1i_2},\atop a_{i_1\bari_2},}
{\ldots,\atop\ldots,}
{a_{i_1i_2\cdots i_{d^L-1}};\atop a_{i_1i_2\cdots\bari_{d^L-1}};}
{b_{j_1},\atop b_{\barj_1},}
{b_{j_1j_2},\atop b_{j_1\barj_2},}
{\ldots,\atop\ldots,}
{b_{j_1j_2\cdots j_{d^R-1}} \atop b_{j_1j_2\cdots\barj_{d^R-1}}}
\right\}\\[7mm]
\hskip25mm
\displaystyle
=F(\boldu(a;i_1,i_2,\ldots,i_{d^L-1}),
        \boldu(b;j_1,j_2,\ldots,j_{d^R-1}),L).
\end{array}
\end{displaymath}
\item If $a_{\bari_1}<\kappa_{t_n}$ or $a_{\bari_1}>p'-\kappa_{t_n}$,
      and $\kappa_{t_n}\le b_{\barj_1}\le p'-\kappa_{t_n}$
      then:
\begin{displaymath}
\begin{array}{l}
\displaystyle
\chi^{p,p'}_{a,b,c}(L)\left\{
{\Xi_0(a_{\bari_1}),\atop \Xi_1(a_{\bari_1}),}
{a_{i_1},\atop a_{\bari_1},}
{a_{i_1i_2},\atop a_{i_1\bari_2},}
{\ldots,\atop\ldots,}
{a_{i_1i_2\cdots i_{d^L-1}};\atop a_{i_1i_2\cdots\bari_{d^L-1}};}
{b_{j_1},\atop b_{\barj_1},}
{b_{j_1j_2},\atop b_{j_1\barj_2},}
{\ldots,\atop\ldots,}
{b_{j_1j_2\cdots j_{d^R-1}} \atop b_{j_1j_2\cdots\barj_{d^R-1}}}
\right\}\\[7mm]
\hskip25mm
\displaystyle
=F(\boldu(a;i_1,i_2,\ldots,i_{d^L-1}),
        \boldu(b;j_1,j_2,\ldots,j_{d^R-1}),L).
\end{array}
\end{displaymath}
\item If $\kappa_{t_n}\le a_{\bari_1}\le p'-\kappa_{t_n}$
      and, $b_{\barj_1}<\kappa_{t_n}$ or $b_{\barj_1}>p'-\kappa_{t_n}$
      then:
\begin{displaymath}
\begin{array}{l}
\displaystyle
\chi^{p,p'}_{a,b,c}(L)\left\{
{a_{i_1},\atop a_{\bari_1},}
{a_{i_1i_2},\atop a_{i_1\bari_2},}
{\ldots,\atop\ldots,}
{a_{i_1i_2\cdots i_{d^L-1}};\atop a_{i_1i_2\cdots\bari_{d^L-1}};}
{\Xi_0(b_{\barj_1}),\atop \Xi_1(b_{\barj_1}),}
{b_{j_1},\atop b_{\barj_1},}
{b_{j_1j_2},\atop b_{j_1\barj_2},}
{\ldots,\atop\ldots,}
{b_{j_1j_2\cdots j_{d^R-1}} \atop b_{j_1j_2\cdots\barj_{d^R-1}}}
\right\}\\[7mm]
\hskip25mm
\displaystyle
=F(\boldu(a;i_1,i_2,\ldots,i_{d^L-1}),
        \boldu(b;j_1,j_2,\ldots,j_{d^R-1}),L).
\end{array}
\end{displaymath}
\item If $a_{\bari_1}<\kappa_{t_n}$ or $a_{\bari_1}>p'-\kappa_{t_n}$,
      and, $b_{\barj_1}<\kappa_{t_n}$ or $b_{\barj_1}>p'-\kappa_{t_n}$
      then:
\begin{displaymath}
\begin{array}{l}
\displaystyle
\chi^{p,p'}_{a,b,c}(L)\!\left\{
\!{\Xi_0(a_{\bari_1}),\atop \Xi_1(a_{\bari_1}),}
{a_{i_1},\atop a_{\bari_1},}
{a_{i_1i_2},\atop a_{i_1\bari_2},}
{\ldots,\atop\ldots,}
{a_{i_1i_2\cdots i_{d^L-1}};\atop a_{i_1i_2\cdots\bari_{d^L-1}};}
{\Xi_0(b_{\barj_1}),\atop \Xi_1(b_{\barj_1}),}
{b_{j_1},\atop b_{\barj_1},}
{b_{j_1j_2},\atop b_{j_1\barj_2},}
{\ldots,\atop\ldots,}
{b_{j_1j_2\cdots j_{d^R-1}} \atop b_{j_1j_2\cdots\barj_{d^R-1}}}
\!\right\}\!\\[7mm]
\hskip25mm
\displaystyle
=F(\boldu(a;i_1,i_2,\ldots,i_{d^L-1}),
        \boldu(b;j_1,j_2,\ldots,j_{d^R-1}),L).
\end{array}
\end{displaymath}
\end{enumerate}
\end{lemma}

\Proof 
Let the leaf-node $a_{i_1i_2\cdots i_{d^L-1}0}$ of the Takahashi tree
for $a$ give rise to the run
$\{\tau^L_j,\sigma^L_j,\Delta^L_j\}_{j=1}^{d^L}$
as in Section \ref{FormVSec}.
Likewise, let the leaf-node $b_{j_1j_2\cdots j_{d^R-1}0}$
of the Takahashi tree for $b$ give rise to the run
$\{\tau^R_j,\sigma^R_j,\Delta^R_j\}_{j=1}^{d^R}$.
Lemma \ref{NaiveRunLem} states that
$\{\tau^L_j,\sigma^L_j,\Delta^L_j\}_{j=1}^{d^L}$ and
$\{\tau^R_j,\sigma^R_j,\Delta^R_j\}_{j=1}^{d^R}$ are both naive runs.

Using these naive runs, define
$\mu^*_k$ and $\mu_k$ for $0\le k<d^L$,
$\nu^*_k$ and $\nu_k$ for $0\le k<d^R$,
and $d^L_0$ and $d^R_0$ as in Section \ref{CoreSec}.
Then, Lemma \ref{NodesMuLem}
states that $a_{i_1i_2\cdots i_{k-1}\bari_k}=\mu^*_k$
and $a_{i_1i_2\cdots i_{k-1} i_k}=\mu_k$
for $1\le k<d^L$, and
$b_{j_1j_2\cdots j_{k-1}\barj_k}=\nu^*_k$
and $b_{j_1j_2\cdots j_{k-1} j_k}=\nu_k$ for $1\le k<d^R$.

Now consider $\kappa_{t_n}\le a_{\bari_1}\le p'-\kappa_{t_n}$.
If $a_{\bari_1}\in\mathcal T$ then $\kappa_{\sigma^L_1}\ge\kappa_{t_n}$
which implies that $\sigma^L_1\ge t_n$.
If $a_{\bari_1}\in\mathcal T'$ then
$p'-\kappa_{\sigma^L_1}\le p'-\kappa_{t_n}$ which also implies that
$\sigma^L_1\ge t_n$.
Similarly $\kappa_{t_n}\le b_{\barj_1}\le p'-\kappa_{t_n}$ implies that
$\sigma^R_1\ge t_n$.
Then $d^L_0=d^R_0=1$, whence the first part of the lemma follows
directly from Theorem \ref{Core2Thrm} after noting that
$\Delta^R_{d^R}=\Delta(\boldu(b;j_1,j_2,\ldots,j_{d^R-1}))$.

If $a_{\bari_1}<\kappa_{t_n}$ then $a_{\bari_1}\in\mathcal T$
whereupon $\sigma^L_1<t_n$, $\Delta^L_1=-1$ and $d^L_0=0$.
If $a_{\bari_1}>p'-\kappa_{t_n}$ then $a_{\bari_1}\in\mathcal T'$
whereupon $\sigma^L_1<t_n$, $\Delta^L_1=+1$ and $d^L_0=0$.
The second part of the lemma now follows directly from
Theorem \ref{Core2Thrm} after again noting that
$\Delta^R_{d^R}=\Delta(\boldu(b;j_1,j_2,\ldots,j_{d^R-1}))$.

The remaining two parts follow analogously.
\cqfd
\medskip

\begin{lemma}\label{DoubleNodeLem}
Let $1\le a<p'$, and $2\le b<p'-1$
with $b$ interfacial in the $(p,p')$-model,
$a\not\in\mathcal T\cup\mathcal T'$ and
$b\not\in\mathcal T\cup\mathcal T'$.
Set $c\in\{b\pm1\}$.
\begin{enumerate}
\item If $\kappa_{t_n}\le a_{\bari_1}\le p'-\kappa_{t_n}$
      and $\kappa_{t_n}\le b_{\barj_1}\le p'-\kappa_{t_n}$
      then:
\begin{displaymath}
\begin{array}{l}
\displaystyle
\chi^{p,p'}_{a,b,c}(L)\left\{
{a_{i_1};\atop a_{\bari_1};}
{b_{j_1}\atop b_{\barj_1}}
\right\}\\[7mm]
\hskip8mm
\displaystyle
= \sum_{{\scriptstyle \boldu(a;i'_1,i'_2,\ldots,i'_{d^L-1})
                      \in{\mathcal U}(a):i'_1=i_1\atop
         \scriptstyle \boldu(b;j'_1,j'_2,\ldots,j'_{d^R-1})
                      \in{\mathcal U}(b):j'_1=j_1}}
\hskip-12mm
F(\boldu(a;i'_1,i'_2,\ldots,i'_{d^L-1}),
        \boldu(b;j'_1,j'_2,\ldots,j'_{d^R-1}),L).
\end{array}
\end{displaymath}
\item If $a_{\bari_1}<\kappa_{t_n}$ or $a_{\bari_1}>p'-\kappa_{t_n}$,
      and $\kappa_{t_n}\le b_{\barj_1}\le p'-\kappa_{t_n}$
      then:
\begin{displaymath}
\begin{array}{l}
\displaystyle
\chi^{p,p'}_{a,b,c}(L)\left\{
{\Xi_0(a_{\bari_1}),\atop \Xi_1(a_{\bari_1}),}
{a_{i_1};\atop a_{\bari_1};}
{b_{j_1}\atop b_{\barj_1}}
\right\}\\[7mm]
\hskip8mm
\displaystyle
= \sum_{{\scriptstyle \boldu(a;i'_1,i'_2,\ldots,i'_{d^L-1})
                      \in{\mathcal U}(a):i'_1=i_1\atop
         \scriptstyle \boldu(b;j'_1,j'_2,\ldots,j'_{d^R-1})
                      \in{\mathcal U}(b):j'_1=j_1}}
\hskip-12mm
F(\boldu(a;i'_1,i'_2,\ldots,i'_{d^L-1}),
        \boldu(b;j'_1,j'_2,\ldots,j'_{d^R-1}),L).
\end{array}
\end{displaymath}
\item If $\kappa_{t_n}\le a_{\bari_1}\le p'-\kappa_{t_n}$
      and, $b_{\barj_1}<\kappa_{t_n}$ or $b_{\barj_1}>p'-\kappa_{t_n}$
      then:
\begin{displaymath}
\begin{array}{l}
\displaystyle
\chi^{p,p'}_{a,b,c}(L)\left\{
{a_{i_1};\atop a_{\bari_1};}
{\Xi_0(b_{\barj_1}),\atop \Xi_1(b_{\barj_1}),}
{b_{j_1}\atop b_{\barj_1}}
\right\}\\[7mm]
\hskip8mm
\displaystyle
= \sum_{{\scriptstyle \boldu(a;i'_1,i'_2,\ldots,i'_{d^L-1})
                      \in{\mathcal U}(a):i'_1=i_1\atop
         \scriptstyle \boldu(b;j'_1,j'_2,\ldots,j'_{d^R-1})
                      \in{\mathcal U}(b):j'_1=j_1}}
\hskip-12mm
F(\boldu(a;i'_1,i'_2,\ldots,i'_{d^L-1}),
        \boldu(b;j'_1,j'_2,\ldots,j'_{d^R-1}),L).
\end{array}
\end{displaymath}
\item If $a_{\bari_1}<\kappa_{t_n}$ or $a_{\bari_1}>p'-\kappa_{t_n}$,
      and, $b_{\barj_1}<\kappa_{t_n}$ or $b_{\barj_1}>p'-\kappa_{t_n}$
      then:
\begin{displaymath}
\begin{array}{l}
\displaystyle
\chi^{p,p'}_{a,b,c}(L)\left\{
{\Xi_0(a_{\bari_1}),\atop \Xi_1(a_{\bari_1}),}
{a_{i_1};\atop a_{\bari_1};}
{\Xi_0(b_{\barj_1}),\atop \Xi_1(b_{\barj_1}),}
{b_{j_1}\atop b_{\barj_1}}
\right\}\\[7mm]
\hskip8mm
\displaystyle
= \sum_{{\scriptstyle \boldu(a;i'_1,i'_2,\ldots,i'_{d^L-1})
                      \in{\mathcal U}(a):i'_1=i_1\atop
         \scriptstyle \boldu(b;j'_1,j'_2,\ldots,j'_{d^R-1})
                      \in{\mathcal U}(b):j'_1=j_1}}
\hskip-12mm
F(\boldu(a;i'_1,i'_2,\ldots,i'_{d^L-1}),
        \boldu(b;j'_1,j'_2,\ldots,j'_{d^R-1}),L).
\end{array}
\end{displaymath}
\end{enumerate}
(Note that in each of the above sums, the values of $d^L$ and $d^R$ vary
over the elements of ${\mathcal U}(a)$ and ${\mathcal U}(b)$.)
\end{lemma}

\Proof Since $a\not\in\mathcal T\cup\mathcal T'$ and
$b\not\in\mathcal T\cup\mathcal T'$, the Takahashi trees of $a$ and $b$
each have more than one leaf-node.
The results then follow from Lemma \ref{OneNodeLem}, and repeated use
of Lemma \ref{GatherLem}.
\cqfd
\medskip

\begin{note}\label{DoubleNodeNote}
If either $a\in\mathcal T\cup\mathcal T'$ or
$b\in\mathcal T\cup\mathcal T'$, then we obtain analogues of
the four results of Lemma \ref{DoubleNodeLem}
(obtaining sixteen cases in all), which we won't write out in full.
Each such case is obtained from the corresponding case of
Lemma \ref{DoubleNodeLem} by,
if $a\in\mathcal T\cup\mathcal T'$ (so that $a_0=a_1=a$),
omitting the column containing $a_{i_1}$ and $a_{\bari_1}$;
and if $b\in\mathcal T\cup\mathcal T'$ (so that $b_0=b_1=b$),
omitting the column containing $b_{j_1}$ and $b_{\barj_1}$.
For example, if $a\in\mathcal T\cup\mathcal T'$ and
$b\not\in\mathcal T\cup\mathcal T'$, then the appropriate
analogue of Lemma \ref{DoubleNodeLem}(2) states that:

\begin{quote}
      If $a<\kappa_{t_n}$ or $a>p'-\kappa_{t_n}$,
      and $\kappa_{t_n}\le b_{\barj_1}\le p'-\kappa_{t_n}$
      then:
\begin{displaymath}
\begin{array}{l}
\displaystyle
\chi^{p,p'}_{a,b,c}(L)\left\{
{\Xi_0(a);\atop \Xi_1(a);}
{b_{j_1}\atop b_{\barj_1}}
\right\}\\[7mm]
\hskip8mm
\displaystyle
= \sum_{{\scriptstyle \boldu(a;i'_1,i'_2,\ldots,i'_{d^L-1})
                      \in{\mathcal U}(a):i'_1=i_1\atop
         \scriptstyle \boldu(b;j'_1,j'_2,\ldots,j'_{d^R-1})
                      \in{\mathcal U}(b):j'_1=j_1}}
\hskip-12mm
F(\boldu(a;i'_1,i'_2,\ldots,i'_{d^L-1}),
        \boldu(b;j'_1,j'_2,\ldots,j'_{d^R-1}),L).
\end{array}
\end{displaymath}
\end{quote}

Note further that if $a\in\mathcal T\cup\mathcal T'$
then $\mathcal U(a)$ contains just one element, and likewise,
if $b\in\mathcal T\cup\mathcal T'$ then $\mathcal U(b)$ contains
just one element.
\end{note}

The following two theorems make use of the following definition
of $\hat c$ in terms of parameters that will be specified in those
theorems:
\begin{equation}\label{ExtraCEq1}
\hat c= \begin{cases}
       2&\text{if $c=\xi_{\eta}>0$ and 
          $\left\lfloor\frac{(b+1)p}{p'}\right\rfloor=
           \left\lfloor\frac{(b-1)p}{p'}\right\rfloor+1$;}\\[1.5mm]
       \hat p'-2&\text{if $c=\xi_{\eta+1}<p'$ and 
          $\left\lfloor\frac{(b+1)p}{p'}\right\rfloor=
           \left\lfloor\frac{(b-1)p}{p'}\right\rfloor+1$;}\\[1.5mm]
       c-\xi_{\eta} &\text{otherwise}.
         \end{cases} 
\end{equation}
Note that if $b$ is interfacial, then certainly
$\left\lfloor\frac{(b+1)p}{p'}\right\rfloor=
           \left\lfloor\frac{(b-1)p}{p'}\right\rfloor+1$.

\begin{theorem}\label{AllNode1Thrm}
Let $1\le a<p'$, and $2\le b<p'-1$
with $b$ interfacial in the $(p,p')$-model.
Let $\eta$ be such that $\xi_{\eta}\le a<\xi_{\eta+1}$ and
$\eta^\prime$ be such that $\xi_{\eta^\prime}\le b<\xi_{\eta^\prime+1}$.
Set $c\in\{b\pm1\}$.

\noindent
1. If $\eta=\eta^\prime$ and
      $\xi_{\eta}<a$ and $\xi_{\eta}<b$ then:
\begin{equation*}
\chi^{p,p'}_{a,b,c}(L)
= \sum_{{\scriptstyle \boldu^L\in{\mathcal U}(a)\atop
         \scriptstyle \boldu^R\in{\mathcal U}(b)}}
F(\boldu^L,\boldu^R,L)
\quad+\quad\ochi^{\hat p,\hat p'}_{\hat a,\hat b,\hat c}(L),
\end{equation*}
where $\hat p'=\xi_{\eta+1}-\xi_{\eta}$,
$\hat p=\tilde\xi_{\eta+1}-\tilde\xi_{\eta}$,
$\hat a=a-\xi_{\eta}$, $\hat b=b-\xi_{\eta}$ and
$\hat c$ is given by (\ref{ExtraCEq1}).

\noindent
2. Otherwise:
\begin{equation*}
\chi^{p,p'}_{a,b,c}(L)
= \sum_{{\scriptstyle \boldu^L\in{\mathcal U}(a)\atop
         \scriptstyle \boldu^R\in{\mathcal U}(b)}}
F(\boldu^L,\boldu^R,L).
\end{equation*}
\end{theorem}

\Proof There are many cases to consider here.

For the moment, assume that $a\not\in\mathcal T\cup\mathcal T'$
and $b\not\in\mathcal T\cup\mathcal T'$ so that the Takahashi trees
for $a$ and $b$ each have at least two leaf-nodes.
In particular, we have $a_0<a<a_1$ and $b_0<b<b_1$, and
certainly $\xi_{\eta}<a$ and $\xi_{\eta'}<b$.
If $\kappa_{t_n}<a<p'-\kappa_{t_n}$ and $\kappa_{t_n}<b<p'-\kappa_{t_n}$
then
$\kappa_{t_n}\le a_0<a_1\le p'-\kappa_{t_n}$
and
$\kappa_{t_n}\le b_0<b_1\le p'-\kappa_{t_n}$.
Thereupon, if $\eta\ne\eta'$,
combining Lemma \ref{GatherLast1Lem} and Lemma \ref{DoubleNodeLem}(1)
gives the required result.
If $\eta=\eta'$, combining Lemma \ref{GatherLast2Lem} and
Lemma \ref{DoubleNodeLem}(1) gives the required result after
noting that Lemma \ref{SegmentLem} implies that the band structure of the
$(p,p')$-model between heights $\xi_{\eta}+1$ and $\xi_{\eta+1}-1$
is identical to the band structure of the $(\hat p,\hat p')$-model.

If $a_0<a_1<\kappa_{t_n}$
then $\Xi_0(a_0)=\Xi_0(a_1)=0$ and $\Xi_1(a_0)=\Xi_1(a_1)=\kappa_{t_n}$.
When $\kappa_{t_n}\le b_0<b_1\le p'-\kappa_{t_n}$,
use of Lemma \ref{DoubleNodeLem}(2) and Lemma \ref{GatherLem} yields:
\begin{displaymath}
\begin{array}{l}
\displaystyle
\chi^{p,p'}_{a,b,c}(L)\left\{
{0;\atop \kappa_{t_n};}
{b_{j_1}\atop b_{\barj_1}}
\right\}\\[7mm]
\hskip8mm
\displaystyle
= \sum_{{\scriptstyle \boldu(a;i'_1,i'_2,\ldots,i'_{d^L-1})
                      \in{\mathcal U}(a)\atop
         \scriptstyle \boldu(b;j'_1,j'_2,\ldots,j'_{d^R-1})
                      \in{\mathcal U}(b):j'_1=j_1}}
\hskip-12mm
F(\boldu(a;i'_1,i'_2,\ldots,i'_{d^L-1}),
        \boldu(b;j'_1,j'_2,\ldots,j'_{d^R-1}),L).
\end{array}
\end{displaymath}
Then, since $0<a<\kappa_{t_n}<b$, Lemma \ref{GatherLast1Lem}
gives the required result.

If $a_0<a_1=\kappa_{t_n}$ then we can use the same argument as in
the above case, after noting that:
\begin{displaymath}
\begin{array}{l}
\displaystyle
\chi^{p,p'}_{a,b,c}(L)\left\{
{0,\atop \kappa_{t_n},}
{a_0;\atop a_1;}
{b_{j_1}\atop b_{\barj_1}}
\right\}
=
\chi^{p,p'}_{a,b,c}(L)\left\{
{a_0;\atop a_1;}
{b_{j_1}\atop b_{\barj_1}}
\right\}\\[7mm]
\hskip8mm
\displaystyle
= \sum_{{\scriptstyle \boldu(a;i'_1,i'_2,\ldots,i'_{d^L-1})
                      \in{\mathcal U}(a):i'_1=0\atop
         \scriptstyle \boldu(b;j'_1,j'_2,\ldots,j'_{d^R-1})
                      \in{\mathcal U}(b):j'_1=j_1}}
\hskip-12mm
F(\boldu(a;i'_1,i'_2,\ldots,i'_{d^L-1}),
        \boldu(b;j'_1,j'_2,\ldots,j'_{d^R-1}),L),
\end{array}
\end{displaymath}
where the first equality is a consequence of (the analogue of)
Lemma \ref{CutParam2Lem}, and the second follows from
Lemma \ref{DoubleNodeLem}(1).

The case $p'-\kappa_{t_n}\le a_0<a_1$ is dealt with in a similar fashion,
again using Lemmas \ref{DoubleNodeLem}(2), \ref{GatherLem},
and \ref{GatherLast1Lem}, and additionally if $a_0=p'-\kappa_{t_n}$,
Lemmas \ref{CutParam2Lem} and \ref{DoubleNodeLem}(1).

Each of the remaining cases for which
$a\not\in\mathcal T\cup\mathcal T'$
and $b\not\in\mathcal T\cup\mathcal T'$ is dealt with in a similar
fashion, using the appropriate part(s) of Lemma \ref{DoubleNodeLem}.

Now if either $a\in\mathcal T\cup\mathcal T'$ or
$b\in\mathcal T\cup\mathcal T'$ (so that $d^L=1$ or $d^R=1$ respectively),
then the result also follows in a similar way, making use of the
appropriate analogue of Lemma \ref{DoubleNodeLem} as described
in Note \ref{DoubleNodeNote}.
\cqfd
\medskip

\begin{note}\label{FunnyNote}
1. An analogue of Lemma \ref{OneNodeLem} that deals with the case of
$b$ being non-interfacial in the $(p,p')$-model and
$c=b-\Delta(\boldu(b;j_1,j_2,\ldots,j_{d^R-1}))$ is obtained
upon replacing (in each case)
$F(\boldu(a;i_1,i_2,\ldots,i_{d^L-1}),
        \boldu(b;j_1,j_2,\ldots,j_{d^R-1}),L)$
by
$\tilde F(\boldu(a;i_1,i_2,\ldots,i_{d^L-1}),$ $
        \boldu(b;j_1,j_2,\ldots,j_{d^R-1}),L)$.
The proof of this result is identical to that of Lemma \ref{OneNodeLem}
except that Theorem \ref{FermLikeThrm} is used instead of
Theorem \ref{Core2Thrm}.

2. An analogue of Lemma \ref{DoubleNodeLem} that deals with the case of
$b$ being non-interfacial in the $(p,p')$-model or $b\in\{1,p'-1\}$,
is obtained upon, whenever the second summation index is such that
$\Delta(\boldu(b;j'_1,j'_2,\ldots,j'_{d^R-1}),L)=b-c$,
replacing the summand
$F(\,\boldu(a\,;i'_1,i'_2,\ldots,i'_{d^L-1}\,),
        \,\boldu(b\,;j'_1,j'_2,\ldots,j'_{d^R-1}\,),\,L\,)$
by
$\tilde F(\,\boldu(a\,;i'_1,i'_2,\ldots,i'_{d^L-1}\,),$\linebreak $
        \boldu(b\,;j'_1,j'_2,\ldots,j'_{d^R-1}\,),\,L\,)$.
This follows from the use of Lemma \ref{OneNodeLem}, the modified form
of Lemma \ref{OneNodeLem} described above, and repeated use of
Lemma \ref{GatherLem} as in the proof of Lemma \ref{DoubleNodeLem}.
The analogues of those results described in Note \ref{DoubleNodeNote}
also follow.
\end{note}

\begin{theorem}\label{AllNode2Thrm}
Let $1\le a,b<p'$ with $b$ non-interfacial in the $(p,p')$-model.
Let $\eta$ be such that $\xi_{\eta}\le a<\xi_{\eta+1}$ and
$\eta^\prime$ be such that $\xi_{\eta^\prime}\le b<\xi_{\eta^\prime+1}$.
Set $c\in\{b\pm1\}$.

\noindent
1. If $\eta=\eta^\prime$ and
      $\xi_{\eta}<a$ and $\xi_{\eta}<b$ then:
\begin{equation*}
\chi^{p,p'}_{a,b,c}(L)=
  \sum_{{\scriptstyle \boldu^L\in{\mathcal U}(a)\atop
         \scriptstyle \boldu^R\in{\mathcal U}^{c-b}(b)}}
\hskip-2mm
F(\boldu^L,\boldu^R,L)
\hskip2mm+
  \sum_{{\scriptstyle \boldu^L\in{\mathcal U}(a)\atop
         \scriptstyle \boldu^R\in{\mathcal U}^{b-c}(b)}}
\hskip-2mm
\widetilde F(\boldu^L,\boldu^R,L)
\hskip2mm+
\hskip2mm
\ochi^{\hat p,\hat p'}_{\hat a,\hat b,\hat c}(L),
\end{equation*}
where $\hat p'=\xi_{\eta+1}-\xi_{\eta}$,
$\hat p=\tilde\xi_{\eta+1}-\tilde\xi_{\eta}$,
$\hat a=a-\xi_{\eta}$, $\hat b=b-\xi_{\eta}$ and
$\hat c$ is given by (\ref{ExtraCEq1}).

\noindent
2. Otherwise:
\begin{equation*}
\chi^{p,p'}_{a,b,c}(L)=
  \sum_{{\scriptstyle \boldu^L\in{\mathcal U}(a)\atop
         \scriptstyle \boldu^R\in{\mathcal U}^{c-b}(b)}}
\hskip-2mm
F(\boldu^L,\boldu^R,L)
\quad+
  \sum_{{\scriptstyle \boldu^L\in{\mathcal U}(a)\atop
         \scriptstyle \boldu^R\in{\mathcal U}^{b-c}(b)}}
\hskip-2mm
\widetilde F(\boldu^L,\boldu^R,L).
\end{equation*}
\end{theorem}

\Proof With Note \ref{FunnyNote} in mind, this theorem is proved
in the same way as Theorem \ref{AllNode1Thrm}.
\cqfd
\medskip

\noindent
Theorems \ref{AllNode1Thrm} and \ref{AllNode2Thrm} prove the expressions
given in (\ref{FermEq}) and (\ref{Ferm2Eq}) respectively.

\newpage

\setcounter{section}{8}

\section{Fermionic character expressions}\label{FermCharSec}

In this section we obtain fermionic expressions for the
Virasoro characters $\chi^{p,p'}_{r,s}$.
In the first instance, in Section \ref{FermLimSec},
these character expressions are in terms of the Takahashi trees
of $s$ and $b$, where $b$ is related to $r$ by (\ref{groundstatelabel})
for some $c\in\{b\pm1\}$.
{}From these expressions, in Section \ref{AssimilateTrees}, simpler
fermionic expressions for $\chi^{p,p'}_{r,s}$ in terms of the Takahashi
tree for $s$ and the truncated Takahashi tree for $r$ are derived.

\subsection{Taking the limit $L\to\infinity$.}\label{FermLimSec}

In order to obtain fermionic expressions for the characters
$\chi^{p,p'}_{r,s}$, we take the $L\to\infinity$ limit
in (\ref{FermEq}) and (\ref{Ferm2Eq}).
However, this cannot be done immediately:
it is first necessary to change variables in (\ref{FEq}).
To this end, given $\run=\{\tau_j,\sigma_j,\Delta_j\}_{j=1}^d$,
we define $\boldN(\run)=(N_1,N_2,\ldots,N_{t_1})$, where
for $1\le i\le t_1$:
\begin{equation*}
N_i=\begin{cases}
      0 & \text{if }1\le i\le\sigma_d;\\
      i-\sigma_d & \text{if }\sigma_d\le i\le\tau_d;\\
      \tau_d-\sigma_d & \text{if }\tau_d\le i\le t_1.
    \end{cases}
\end{equation*}
Then define:
\begin{equation*}\label{ChangeEq1}
\begin{split}
F^*(\boldu^L,\boldu^R,L)=
  &\sum
\hskip2mm
  q^{\tilde{\boldn}^T\boldB\tilde{\boldn}
   -\boldN(\run^R)\cdot\tilde{\boldn}
   +\frac{1}{4}\hat{\boldm}^T\boldC\hat{\boldm}
   -\frac{1}{2}(\boldu^L_\flat+\boldu^R_\sharp)\cdot{\boldm}
   +\frac{1}{4}\gamma(\run^L,\run^R)}\\
&\hskip9mm\times\hskip2mm
  \prod_{\makebox[10mm]{$\scriptstyle j=1$}}^{t_1}
  \left[
  {n_j+m_{t_1}+2\sum_{i=j+1}^{t_1}(i-j)\tilde n_i\atop n_j}
  \right]_q \\
&\hskip9mm\times\hskip2mm
  \prod_{\makebox[10mm]{$\scriptstyle j=t_1+1$}}^{t-1}
  \left[
  {m_j-\frac{1}{2}\left(\sum_{i=t_1}^{t-1}
      C_{ji}m_i\!-\!(\boldu^L)_j\!-\!(\boldu^R)_j\right)\atop m_j}
  \right]_q,
\end{split}
\end{equation*}
where the sum is over all
$(m_{t_1},m_{t_1+1},\ldots,m_{t-1})\equiv(Q_{t_1},Q_{t_1+1},\ldots,Q_{t-1})$,
and over all $\boldn=(n_1,n_2,\ldots,n_{t_1})$ such that
\begin{equation}\label{ConstraintEq}
\sum_{i=1}^{t_1} i(2n_i-u_i)=L-m_{t_1}-t_1m_{t_1+1},
\end{equation}
where $(u_1,u_2,\ldots,u_{t})=\boldu^L+\boldu^R$.
Here, $\tilde{\boldn}$ is defined by (\ref{ntildeEq}), and the
$t$-dimensional vector $\hat{\boldm}$ and the $(t-1)$-dimensional
vector ${\boldm}$ are defined by:
\begin{equation}\label{mextendEq}
\begin{split}
\hat{\boldm}=(0,0,0,\ldots,0,&m_{t_1+1},m_{t_1+2},\ldots,m_{t-1});\\
{\boldm}=(0,0,\ldots,0,&m_{t_1+1},m_{t_1+2},\ldots,m_{t-1}),
\end{split}
\end{equation}
where the appropriate number of zeros have been included.
In addition, $\run^L=\run(a,\boldu^L)$ and $\run^R=\run(b,\boldu^R)$,
as defined in Section \ref{FormVSec}.

For $\boldu\in\U(b)$ and $\boldu=\boldu(\run)$ with
$\run=\{\tau_j,\sigma_j,\Delta_j\}_{j=1}^d$, we also
define $\sigma(\boldu)=\sigma_d$ and $\Delta(\boldu)=\Delta_d$.

The results of this section involve the subset
$\overline{\mathcal U}(b)$ of ${\mathcal U}(b)$ defined by:
\begin{equation}\label{ExcludeEq}
\begin{array}{rl}
\overline{\mathcal U}(b)=
\{\boldu\in{\mathcal U}(b):&\sigma(\boldu)=0,\\
&0<b+\Delta(\boldu)<p'\implies
\lfloor bp/p'\rfloor\ne\lfloor(b+\Delta(\boldu))p/p'\rfloor,\\
&b+\Delta(\boldu)\in\{0,p'\}\implies p'<2p \}.\\
\end{array}
\end{equation}
Correspondingly, we define
$\overline{\mathcal U}{}^\Delta(b)
      ={\mathcal U}^\Delta(b)\cap\overline{\mathcal U}(b)$.

The results that follow are valid in both cases $p'>2p$ and $p'<2p$.

\begin{lemma}\label{ChangeLem}
Let $1\le a,b<p'$, and let
$\boldu^L\in{\mathcal U}(a)$ and
$\boldu^R\in{\mathcal U}(b)$.
Then set $(Q_1,Q_2,\ldots,Q_{t-1})=\boldQ(\boldu^L+\boldu^R)$.

1. If $\boldu^R\in\overline{\mathcal U}(b)$ then:
\begin{equation*}
F(\boldu^L,\boldu^R,L)=q^{\frac{1}{2}L} F^*(\boldu^L,\boldu^R,L).
\end{equation*}

2. If $\boldu^R\notin\overline{\mathcal U}(b)$ then:
\begin{equation*}
F(\boldu^L,\boldu^R,L)=F^*(\boldu^L,\boldu^R,L).
\end{equation*}
\end{lemma}

\Proof
In the expression (\ref{FEq}) for $F(\boldu^L,\boldu^R,L)$,
the summation is over variables $\{m_i\}_{i=1}^{t-1}$ with $m_0=L$.
Here we exchange $\{m_i\}_{i=0}^{t_1-1}$ for
$\{n_i\}_{i=1}^{t_1}$ using the $\boldm\boldn$-system.
Since $m_0=L$ is fixed, a constraint amongst the $\{n_i\}_{i=1}^{t_1}$
is required.
In accordance with (\ref{MNEq1}) and (\ref{MNEq2}),
set $n_i=\frac12(m_{i-1}-2m_i+m_{i+1}+u_i)$
for $1\le i< t_1$,
and set $n_{t_1}=\frac12(m_{{t_1}-1}-m_{t_1}-m_{{t_1}+1}+u_{t_1})$.
For $1\le j\le t_1$, we immediately obtain:
\begin{equation*}\label{ChangeProof1}
m_{j-1}-m_j=m_{t_1+1}+\sum_{i=j}^{t_1}(2n_i-u_i)
=2\sum_{i=j}^{t_1}\tilde n_i.
\end{equation*}
Thereupon, for $0\le j\le t_1$,
\begin{equation}\label{ChangeProof2}
m_j=m_{t_1}+2\sum_{i=j+1}^{t_1} (i-j)\tilde n_i.
\end{equation}
The $j=0$ case here yields the constraint (\ref{ConstraintEq}).
Additionally, if $\hat{\boldm}_*=(m_0,m_1,m_2,$ $\ldots,m_{t-1})$,
we obtain:
\begin{equation*}
\hat{\boldm}_*^{\smash{T}}\boldC\hat{\boldm}_*^{\smash{\phantom{T}}}-L^2
=\sum_{k=1}^{t_1}(m_{k-1}-m_k)^2+\hat{\boldm}^T\boldC\hat{\boldm}
=4\tilde{\boldn}^T\boldB\tilde{\boldn}+\hat{\boldm}^T\boldC\hat{\boldm}
\end{equation*}

Now restrict to the case $p'>2p$ so that $t_1>0$ (all the above is
valid in the $t_1=0$ case, although it's mainly degenerate).
Let $\run^R=\{\tau_j,\sigma_j,\Delta_j\}_{j=1}^d$
and $\boldN=\boldN(\run^R)$.
If $\sigma_d\le t_1\le\tau_d$ then
$(\boldu^L_\flat+\boldu^R_\sharp)_i
 =\delta_{\sigma_d,i}-\delta_{t_1,i}$ for $1\le i\le t_1$.
If $\sigma_d<\tau_d\le t_1$ then
$(\boldu^L_\flat+\boldu^R_\sharp)_i
 =\delta_{\sigma_d,i}-\delta_{\tau_d,i}$ for $1\le i\le t_1$.
Thereupon, for $0<\sigma_d\le t_1$,
\begin{align*}
(\boldu^L_\flat+\boldu^R_\sharp)\cdot(m_1,m_2,\ldots,m_{t-1})
&=(\boldu^L_\flat+\boldu^R_\sharp)\cdot{\boldm}+m_{\sigma_d}
 -\begin{cases}
     m_{\tau_d} &\text{if }\tau_d\le t_1\\
     m_{t_1}    &\text{if }\tau_d\ge t_1
  \end{cases}
\\
&=(\boldu^L_\flat+\boldu^R_\sharp)\cdot{\boldm}
    +2\boldN\cdot\tilde{\boldn},\\
\intertext{and for $\sigma_d=0$,}
(\boldu^L_\flat+\boldu^R_\sharp)\cdot(m_1,m_2,\ldots,m_{t-1})
&=(\boldu^L_\flat+\boldu^R_\sharp)\cdot{\boldm}
 -\begin{cases}
     m_{\tau_d} &\text{if }\tau_d\le t_1\\
     m_{t_1}    &\text{if }\tau_d\ge t_1
  \end{cases}
\\
&=(\boldu^L_\flat+\boldu^R_\sharp)\cdot{\boldm}
    +2\boldN\cdot\tilde{\boldn}-L,
\end{align*}
having used (\ref{ChangeProof2}) in each case.
In the case $p'>2p$ and $\sigma_d\le t_1$,
the lemma then follows from the definition of $F(\boldu^L,\boldu^R,L)$,
after noting that the definition of Section \ref{ConSec} implies that
$\gamma'(\run^L,\run^R)=\gamma(\run^L,\run^R)$ unless both $\sigma_d=0$
and $\boldu^R\notin{\mathcal U}(b)$, in which case
$\gamma'(\run^L,\run^R)=\gamma(\run^L,\run^R)-2L$.

In the case $p'>2p$ and $\sigma_d>t_1$, the lemma follows after noting that
\begin{equation}\label{Change2Eq}
(\boldu^L_\flat+\boldu^R_\sharp)\cdot(m_1,m_2,\ldots,m_{t-1})
=(\boldu^L_\flat+\boldu^R_\sharp)\cdot{\boldm},
\end{equation}
$\boldN=0$, and $\gamma'(\run^L,\run^R)=\gamma(\run^L,\run^R)$
in this case.

In the case $p'<2p$ so that $t_1=0$, eq.\ (\ref{Change2Eq}) also holds.
The lemma then follows in this case after noting from
Section \ref{ConSec} that
$\gamma'(\run^L,\run^R)=\gamma(\run^L,\run^R)+2L$
if $\sigma_d=0$ and either
$\lfloor (b+\Delta_d)p/p'\rfloor\ne\lfloor bp/p'\rfloor$
or $b+\Delta_d\in\{0,p'\}$, whereas otherwise,
$\gamma'(\run^L,\run^R)=\gamma(\run^L,\run^R)$.
\cqfd
\medskip

Note that although the set of vectors ${\mathcal U}(b)$ is
identical in the cases of the $(p,p')$-model and the $(p'-p,p')$-model,
those ${\boldu}^R$ which fall into the first category in
Lemma \ref{ChangeLem} differ in the two cases.

We can now obtain $\lim_{L\to\infinity} F(\boldu^L,\boldu^R,L)$.
The result involves $F^*(\boldu^L,\boldu^R)$ defined by:
\begin{equation}\label{ChangeEq2}
\begin{split}
F^*(\boldu^L,\boldu^R)=
  &\sum
\hskip2mm
  q^{\tilde{\boldn}^T\boldB\tilde{\boldn}
   -\boldN(\run^R)\cdot\tilde{\boldn}
   +\frac{1}{4}\hat{\boldm}^T\boldC\hat{\boldm}
   -\frac{1}{2}(\boldu^L_\flat+\boldu^R_\sharp)\cdot{\boldm}
   +\frac{1}{4}\gamma(\run^L,\run^R)}\\
&\hskip3mm\times\hskip2mm
  \frac{1}{(q)_{m_{t_1+1}}}
  \prod_{j=1}^{t_1} \frac{1}{(q)_{n_j}}
  \prod_{j=t_1+2}^{t-1}
  \left[
 {m_j-\frac{1}{2}
     ({\boldCC}^*{\oboldm}\!-\!\oboldu^L\!-\!\oboldu^R)_j\atop m_j}
  \right]_q,
\end{split}
\end{equation}
where the sum is over all vectors
$\oboldm=(m_{t_1+1},m_{t_1+2},\ldots,m_{t-1})$ such that
$\oboldm\equiv(Q_{t_1+1},Q_{t_1+2},\ldots,Q_{t-1})$,
and over all $\boldn=(n_1,n_2,\ldots,n_{t_1})$, and as above,
we use (\ref{ntildeEq}) and (\ref{mextendEq}).
We set $\run^L=\run(a,\boldu^L)$.
When $\boldu^R\in{\mathcal U}(b)$, we set $\run^R=\run(b,\boldu^R)$,
but if $\boldu^R=\boldu^+$ with $\boldu^+\in{\mathcal U}(b)$ then
we set $\run^R=\run^+$ where $\run=\run(b,\boldu)$.
We also set $m_t=0$.

\begin{lemma}\label{LimitLem}
Let $1\le a,b<p'$,
$\boldu^L\in{\mathcal U}(a)$ and
$\boldu^R\in{\mathcal U}(b)$.

1. If $\boldu^R\in\overline{\mathcal U}(b)$ then:
\begin{equation*}
\lim_{L\to\infinity} F(\boldu^L,\boldu^R,L)=0.
\end{equation*}

2. If $\boldu^R\notin\overline{\mathcal U}(b)$ then:
\begin{equation*}
\lim_{L\to\infinity} F(\boldu^L,\boldu^R,L)=F^*(\boldu^L,\boldu^R).
\end{equation*}
\end{lemma}

\Proof We use the notation in the proof of Lemma \ref{ChangeLem}.
If $t_1+1\le j<t$ note that $\sum_{i=t_1}^{t-1} C_{ji}m_i$
depends on $m_{t_1}$ only if $j=t_1+1$, and then $C_{jt_1}=-1$.
Otherwise, for $t_1+1<j<t$, we have
$\sum_{i=t_1}^{t-1} C_{ji}m_i=(\boldCC^*\overline{\boldm})_j$.

In view of (\ref{ChangeProof2}) and since $m_0=L$,
taking the limit $L\to\infty$ effects $m_{t_1}\to\infinity$.
The current theorem then follows from Lemma \ref{ChangeLem}
after noting that
\begin{equation*}
\lim_{m\to\infinity}
\left[{m+n}\atop n\right]_q=
\lim_{m\to\infinity}
\left[{m+n}\atop m\right]_q=
\frac{1}{(q)_n}.
\end{equation*}
\cqfd
\medskip

\begin{theorem}\label{Limit1Thrm}
Let $1\le a<p'$, and $2\le b<p'-1$
with $b$ interfacial in the $(p,p')$-model.
Set $r=\lfloor(b+1)p/p'\rfloor$.
Let $\eta$ be such that $\xi_{\eta}\le a<\xi_{\eta+1}$ and $\eta^\prime$
be such that $\tilde\xi_{\eta^\prime}\le r<\tilde\xi_{\eta^\prime+1}$.

\noindent
1. If $\eta=\eta^\prime$ and
      $\xi_{\eta}<a$ and $\tilde\xi_{\eta}<r$ then:
\begin{equation}\label{Limit1Eq1}
\chi^{p,p'}_{r,a}
= \sum_{\begin{subarray}{c}
           \boldu^L\in{\mathcal U}(a)\\
           \boldu^R\in{\mathcal U}(b)\backslash\overline{\mathcal U}(b)
        \end{subarray}}
F^*(\boldu^L,\boldu^R)
\quad+\quad\ochi^{\hat p,\hat p'}_{\hat r,\hat a},
\end{equation}
where $\hat p'=\xi_{\eta+1}-\xi_{\eta}$,
$\hat p=\tilde\xi_{\eta+1}-\tilde\xi_{\eta}$,
$\hat a=a-\xi_{\eta}$, $\hat r=r-\tilde\xi_{\eta}$.

\noindent
2. Otherwise:
\begin{equation}\label{Limit1Eq2}
\chi^{p,p'}_{r,a}
= \sum_{\begin{subarray}{c}
           \boldu^L\in{\mathcal U}(a)\\
           \boldu^R\in{\mathcal U}(b)\backslash\overline{\mathcal U}(b)
        \end{subarray}}
F^*(\boldu^L,\boldu^R).
\end{equation}
\end{theorem}

\Proof To prove these results, we take the $L\to\infinity$ limit
of the results given in Theorem \ref{AllNode1Thrm}.
Since $b$ is interfacial, we may choose either $c\in\{b\pm1\}$.
After checking that $0<r<p$, (\ref{ChiLimEq}) gives
$\lim_{L\to\infinity}\chi^{p,p'}_{a,b,c}(L)=\chi^{p,p'}_{r,a}$.

Now note that if $\boldu^R\in\overline{\mathcal U}(b)$,
Lemma \ref{LimitLem}(1) implies that
$\lim_{L\to\infinity} F(\boldu^L,\boldu^R,L)=0$.
For the remaining terms where
$\boldu^R\in{\mathcal U}(b)\backslash\overline{\mathcal U}(b)$,
Lemma \ref{LimitLem}(2) implies that
$\lim_{L\to\infinity} F(\boldu^L,\boldu^R,L)=F^*(\boldu^L,\boldu^R)$.
Then in the case $a=\xi_{\eta}$,
taking the $L\to\infinity$ limit of Theorem \ref{AllNode1Thrm}(2)
immediately yields the required (\ref{Limit1Eq2}).

Hereafter, we assume that $a>\xi_{\eta}$.
Define $\eta''$ such that $\xi_{\eta''}\le b<\xi_{\eta''+1}$
and set $\hat b=b-\xi_{\eta''}$.
Since, by Lemma \ref{SegmentHalfLem}, $\xi_i$ is interfacial
and neighbours the $\tilde\xi_i$th odd band for $1\le i\le2c_n-2$,
it follows that
$\tilde\xi_{\eta''}\le r\le\tilde\xi_{\eta''+1}$.

Consider the case $\eta''\ne\eta$.
Now if $\eta'=\eta$ then $r=\tilde\xi_{\eta''+1}$ and thus
$\eta'=\eta''+1$ so that $r=\tilde\xi_{\eta}$.
So whether $\eta'\ne\eta$ or $\eta'=\eta$, we are required to
prove (\ref{Limit1Eq2}).
This follows immediately from taking the $L\to\infinity$ limit
of Theorem \ref{AllNode1Thrm}(2) which holds because $\eta''\ne\eta$.

The case $\eta''=\eta$ spawns three subcases:
$r=\tilde\xi_{\eta}$;
$r=\tilde\xi_{\eta+1}$;
$\tilde\xi_{\eta}<r<\tilde\xi_{\eta+1}$.
We consider each in turn.

If $r=\tilde\xi_{\eta}$ then it is required to prove (\ref{Limit1Eq2}).
Since $b$ is interfacial and neighbours the $r$th odd band,
then $\eta>0$ and either $b=\xi_{\eta}$ or $b=\xi_{\eta}+1$.
In the former case, we take the $L\to\infinity$ limit
of Theorem \ref{AllNode1Thrm}(2) to obtain the desired (\ref{Limit1Eq2}).
In the latter case, $\hat b=1$ and from (\ref{ExtraCEq1}), $\hat c=2$.
It then follows from (\ref{ChiLim0Eq}) that
$\lim_{L\to\infinity}\chi^{\hat p,\hat p'}_{\hat a,\hat b,\hat c}(L)=0$,
and thus Theorem \ref{AllNode1Thrm}(1) yields the desired (\ref{Limit1Eq2}).

If $r=\tilde\xi_{\eta+1}$ then $\eta'=\eta+1$, so that again
it is required to prove (\ref{Limit1Eq2}).
Here, we necessarily have $b=\xi_{\eta+1}-1$, and so
$\hat b=\xi_{\eta+1}-1-\xi_{\eta}=\hat p'-1$.
Then $\hat c=\hat p'-2$ from (\ref{ExtraCEq1}).
Then (\ref{ChiLim0Eq}) again implies that
$\lim_{L\to\infinity}\chi^{\hat p,\hat p'}_{\hat a,\hat b,\hat c}(L)=0$,
and thus Theorem \ref{AllNode1Thrm}(1) yields the desired (\ref{Limit1Eq2}).

If $\tilde\xi_{\eta}<r<\tilde\xi_{\eta+1}$ then $\eta'=\eta$,
and $\xi_{\eta}<b<\xi_{\eta+1}$.
Thus Theorem \ref{AllNode1Thrm}(1) applies in this case,
and we are required to prove (\ref{Limit1Eq1}).
It suffices to show that
$\lim_{L\to\infinity}\chi^{\hat p,\hat p'}_{\hat a,\hat b,\hat c}(L)
=\chi^{\hat p,\hat p'}_{\hat r,\hat a}$.
That this is so follows from (\ref{ChiLimEq}) because there are
$\tilde\xi_{\eta}$ odd bands in the $(p,p')$-model below height
$\xi_{\eta}+1$, and therefore $\hat b$ neighbours the
$(r-\tilde\xi_{\eta})$th odd band in the $(\hat p,\hat p')$-model.
\cqfd
\medskip

\noindent
Theorem \ref{Limit1Thrm} certainly provides a fermionic expression
for $\chi^{p,p'}_{r,s}$ when $b$ with the required property can be
found. However a more succinct expression for $\chi^{p,p'}_{r,s}$,
and which works for all $r$, is obtained in Section \ref{AssimilateTrees}
in terms of the truncated Takahashi tree of $r$.

\begin{lemma}\label{Tilde1LimLem}
Let $p'>2p$ and $1\le a,b<p'$ with
$b$ non-interfacial in the $(p,p')$-model.
If $\boldu^L\in{\mathcal U}(a)$ and $\boldu^R\in{\mathcal U}(b)$ then:
\begin{displaymath}
\lim_{L\to\infinity} \widetilde F(\boldu^L,\boldu^R,L) = 0.
\end{displaymath}
\end{lemma}

\Proof
Let $\{\tau_j,\sigma_j,\Delta_j\}_{j=1}^d=\run(b,\boldu^R)$ so
that $\sigma_d=\sigma(\boldu^R)$ and $\Delta_d=\Delta(\boldu^R)$.
In the case $\sigma_d=0$, the required result follows
from the definition (\ref{Tilde1Def}) after noting via
Theorem \ref{Core2Thrm} that $F(\boldu^L,\boldu^R,L)$ is the
generating function for certain paths.

Now consider $\sigma_d>0$.
Let $\run^-=\{\tau_j,\sigma^-_j,\Delta_j\}_{j=1}^d$,
with $\sigma_j^-=\sigma_j$ for $1\le j<d$ and $\sigma_d^-=\sigma_d-1$.
Then let $\boldu^{R-}=\boldu(\run^-)$, using (\ref{uEq}).
Note that $\sigma(\boldu^{R-})=\sigma(\boldu^R)-1$.
We then have $0\le\sigma(\boldu^{R-})<\tau-1$,
where $\tau$ is specified by (\ref{tauDef}).
Substituting $\boldu^R\to\boldu^{R-}$
(so that $\boldu^{R+}\to\boldu^R$ and $b\to b+\Delta$)
into Lemma \ref{FermLikeBitsLem}(1) and rearranging, yields:
\begin{equation}\label{Tilde1LimEq}
\lim_{L\to\infty}
\widetilde F(\boldu^{L},\boldu^{R},L)
=
\begin{cases}
\lim_{L\to\infty} \widetilde F(\boldu^{L},\boldu^{R-},L)
&\text{if $b+\Delta$ is non-interfacial;}\\
\lim_{L\to\infty} F(\boldu^{L},\boldu^{R-},L)
&\text{if $b+\Delta$ is interfacial.}
\end{cases}
\end{equation}

Since $b$ is not interfacial,
Corollary \ref{CoreCor} implies that $\sigma_d\le t_1$.
Since $p'/p$ has continued fraction $[c_0,c_1,\ldots,c_n]$,
and $t_1=c_0-1$, neighbouring odd bands in the $(p,p')$-model
are separated by $t_1$ or $t_1+1$ even bands.
If $b_{j_1j_2\cdots j_{d-1}0}$ is the leaf-node of the
Takahashi tree for $b$ that corresponds to the run
$\{\tau_j,\sigma_j,\Delta_j\}_{j=1}^d$ then,
from Lemma \ref{NodesMuLem},
$b=b_{j_1j_2\cdots j_{d-1}}-\Delta\kappa_{\sigma_d}
        =b_{j_1j_2\cdots j_{d-1}}-\Delta(\sigma_d+1)$,
and $b_{j_1j_2\cdots j_{d-1}}$ is interfacial.
It follows that if $b$ is non-interfacial,
then $b+\Delta$ is interfacial only if $\sigma_d=1$.
In this case that $b+\Delta$ is interfacial then
$\sigma(\boldu^{R-})=0$ whereupon it follows that
$\boldu^{R-}\in\overline{\mathcal U}(b+\Delta)$
and Lemma \ref{LimitLem}(1) gives
$\lim_{L\to\infty} F(\boldu^{L},\boldu^{R-},L)=0$.

In the case that $\sigma_d=1$ and $b$ is non-interfacial, the use
of (\ref{Tilde1Def}) again yields
$\lim_{L\to\infty} \widetilde F(\boldu^{L},\boldu^{R-},L)=0$.
Thus the required result holds whenever $\sigma(\boldu^R)=1$.
For $\sigma(\boldu^R)>1$, the required result then follows by
repeated use of (\ref{Tilde1LimEq}).
\cqfd
\medskip

\begin{lemma}\label{Tilde2LimLem}
Let $p'<2p$ and $1\le a,b<p'$ with
$b$ non-interfacial in the $(p,p')$-model.
Let $\boldu^L\in{\mathcal U}(a)$ and $\boldu^R\in{\mathcal U}(b)$.
Then:
\begin{equation*}
\lim_{L\to\infinity} \widetilde F(\boldu^L,\boldu^R,L)=
F^*(\boldu^L,\boldu^{R+}).
\end{equation*}
\end{lemma}

\Proof In the case $\sigma(\boldu^R)>0$, this follows directly
from the definition (\ref{Tilde2Def}) and Lemma \ref{LimitLem}(2).
In the case $\sigma(\boldu^R)=0$, it follows from
Lemma \ref{FermLikeBitsLem}(2) combined with Lemma \ref{LimitLem}(2),
on noting that $\sigma(\boldu^{R+})>0$
\cqfd
\medskip

\begin{theorem}\label{Limit2Thrm}
Let $p'<2p$ and $1\le a,b<p'$ with $b$ non-interfacial in the $(p,p')$-model.
Let $c\in\{b\pm1\}$, let $r$ be given by (\ref{groundstatelabel2})
and assume that $0<r<p$.
Let $\eta$ be such that $\xi_{\eta}\le a<\xi_{\eta+1}$ and $\eta^\prime$
be such that $\tilde\xi_{\eta^\prime}\le r<\tilde\xi_{\eta^\prime+1}$.

\noindent
1. If $\eta=\eta^\prime$ and
      $\xi_{\eta}<a$ and $\tilde\xi_{\eta}<r$ then:
\begin{equation}\label{Limit2Eq1}
\chi^{p,p'}_{r,a}
= \sum_{\begin{subarray}{c}
      \boldu^L\in{\mathcal U}(a)\\
      \boldu^R\in{\mathcal U}^{c-b}(b)\backslash
                                      \overline{\mathcal U}{}^{c-b}(b)
   \end{subarray}}
\hskip-2mm
F^*(\boldu^L,\boldu^R)
\hskip2mm
+ \sum_{\begin{subarray}{c}
      \boldu^L\in{\mathcal U}(a)\\
      \boldu^R\in{\mathcal U}^{b-c}(b)
   \end{subarray}}
\hskip-2mm
F^*(\boldu^L,\boldu^{R+})
\hskip2mm+
\hskip2mm
\ochi^{\hat p,\hat p'}_{\hat r,\hat a}.
\end{equation}
where $\hat p'=\xi_{\eta+1}-\xi_{\eta}$,
$\hat p=\tilde\xi_{\eta+1}-\tilde\xi_{\eta}$,
$\hat a=a-\xi_{\eta}$, $\hat r=r-\tilde\xi_{\eta}$.

\noindent
2. Otherwise:
\begin{equation}\label{Limit2Eq2}
\chi^{p,p'}_{r,a}
= \sum_{\begin{subarray}{c}
      \boldu^L\in{\mathcal U}(a)\\
      \boldu^R\in{\mathcal U}^{c-b}(b)\backslash
                                      \overline{\mathcal U}{}^{c-b}(b)
   \end{subarray}}
\hskip-2mm
F^*(\boldu^L,\boldu^R)
\hskip2mm
+ \sum_{\begin{subarray}{c}
      \boldu^L\in{\mathcal U}(a)\\
      \boldu^R\in{\mathcal U}^{b-c}(b)
   \end{subarray}}
\hskip-2mm
F^*(\boldu^L,\boldu^{R+}).
\end{equation}
\end{theorem}

\Proof To prove these results, we take the $L\to\infinity$ limit
of the results given in Theorem \ref{AllNode2Thrm}.
With $r$ as specified, (\ref{ChiLimEq}) gives
$\lim_{L\to\infinity}\chi^{p,p'}_{a,b,c}(L)=\chi^{p,p'}_{r,a}$.

If $\boldu^R\in\overline{\mathcal U}{}^{c-b}(b)$
(so that $c=b+\Delta(\boldu^R)$) then Lemma \ref{LimitLem}(1) implies
that $\lim_{L\to\infinity} F(\boldu^L,\boldu^R,L)=0$.
If $\boldu^R\in{\mathcal U}^{c-b}(b)\backslash
    \overline{\mathcal U}{}^{c-b}(b)$,
then Lemma \ref{LimitLem}(2) implies that
$\lim_{L\to\infinity} F(\boldu^L,\boldu^R,L)=F^*(\boldu^L,\boldu^R)$.
For the remaining cases, where $\boldu^R\in{\mathcal U}^{c-b}(b)$,
Lemma \ref{Tilde2LimLem} implies that
$\lim_{L\to\infinity} \widetilde F(\boldu^L,\boldu^R,L)=
F^*(\boldu^L,\boldu^{R+})$.
Then in the case $a=\xi_{\eta}$,
taking the $L\to\infinity$ limit of Theorem \ref{AllNode2Thrm}(2)
immediately yields the required (\ref{Limit2Eq2}).

Hereafter, we assume that $a>\xi_{\eta}$. The proof is now very similar
to that of Theorem \ref{Limit1Thrm}.
Note that since $p'<2p$ and $b$ is non-interfacial,
$b$ lies between two odd bands.
Define $\eta''$ such that $\xi_{\eta''}\le b<\xi_{\eta''+1}$
and set $\hat b=b-\xi_{\eta''}$.
Since, by Lemma \ref{SegmentHalfLem}, $\xi_i$ is interfacial
and neighbours the $\tilde\xi_i$th odd band for $1\le i\le2c_n-2$,
it follows that
$\tilde\xi_{\eta''}\le r\le\tilde\xi_{\eta''+1}$
and $\xi_{\eta''}<b<\xi_{\eta''+1}$.

Consider the case $\eta''\ne\eta$.
Now if $\eta'=\eta$ then $r=\tilde\xi_{\eta''+1}$ and thus
$\eta'=\eta''+1$ so that $r=\tilde\xi_{\eta}$.
So whether $\eta'\ne\eta$ or $\eta'=\eta$, we are required to
prove (\ref{Limit2Eq2}).
This follows immediately from taking the $L\to\infinity$ limit
of Theorem \ref{AllNode2Thrm}(2) which holds because $\eta''\ne\eta$.

For the case $\eta''=\eta$, Theorem \ref{AllNode2Thrm}(1) holds.
We consider the three subcases:
$r=\tilde\xi_{\eta}$;
$r=\tilde\xi_{\eta+1}$;
$\tilde\xi_{\eta}<r<\tilde\xi_{\eta+1}$, in turn.

If $r=\tilde\xi_{\eta}$ then it is required to prove (\ref{Limit2Eq2}).
Since $b$ neighbours the $r$th odd band, then $b=\xi_{\eta}+1$ and
$c=\xi_{\eta}$ so that $\hat b=1$ and from (\ref{ExtraCEq1}), $\hat c=0$.
It then follows from (\ref{groundstatelabel2}) and (\ref{ChiLim0Eq}) that
$\lim_{L\to\infinity}\chi^{\hat p,\hat p'}_{\hat a,\hat b,\hat c}(L)=0$,
and thus Theorem \ref{AllNode2Thrm}(1) yields the desired (\ref{Limit2Eq2}).

If $r=\tilde\xi_{\eta+1}$ then $\eta'=\eta+1$, so that again
it is required to prove (\ref{Limit2Eq2}).
Here, we necessarily have $b=\xi_{\eta+1}-1$ and $c=\xi_{\eta+1}$, and so
$\hat b=\xi_{\eta+1}-1-\xi_{\eta}=\hat p'-1$ and
$\hat c=\hat p'$ from (\ref{ExtraCEq1}).
Again,
$\lim_{L\to\infinity}\chi^{\hat p,\hat p'}_{\hat a,\hat b,\hat c}(L)=0$
follows from (\ref{groundstatelabel2}) and (\ref{ChiLim0Eq}),
and thus Theorem \ref{AllNode2Thrm}(1) yields the desired (\ref{Limit2Eq2}).

If $\tilde\xi_{\eta}<r<\tilde\xi_{\eta+1}$ then $\eta'=\eta$.
Thus we are required to prove (\ref{Limit2Eq1}).
Note that (\ref{ExtraCEq1}) implies that $\hat c=c-\xi_{\eta}$.
To prove (\ref{Limit2Eq1}), it suffices to show that
$\lim_{L\to\infinity}\chi^{\hat p,\hat p'}_{\hat a,\hat b,\hat c}(L)
=\chi^{\hat p,\hat p'}_{\hat r,\hat a}$.
That this is so follows from (\ref{ChiLimEq}) after noting that because
$b$ and $c$ straddle the $r$th odd band in the $(p,p')$-model and there
are $\tilde\xi_{\eta}$ odd bands in the $(p,p')$-model below height
$\xi_{\eta}+1$, then $\hat b$ and $\hat c$ straddle the
$(r-\tilde\xi_{\eta})$th odd band in the $(\hat p,\hat p')$-model.
\cqfd
\medskip

\noindent
Theorem \ref{Limit2Thrm} provides fermionic expressions for those
$\chi^{p,p'}_{r,s}$ that are not covered by Theorem \ref{Limit1Thrm}.
In Section \ref{AssimilateTrees}, these two Theorems are subsumed into
a more succinct expression for $\chi^{p,p'}_{r,s}$ which makes use
of the truncated Takahashi tree of $r$.

\begin{lemma}\label{ShiftSigmaLem}
Let $1\le a,b<p'$.
If $\boldu^L\in\mathcal U(a)$ and
$\boldu\in{\mathcal U}(b)\backslash\overline{\mathcal U}(b)$ with
$\sigma(\boldu)\le t_1$ then:
\begin{equation*}
F^*(\boldu^L,\boldu)=F^*(\boldu^L,\boldu^+).
\end{equation*}
\end{lemma}

\Proof Let $\run^L=\run(a,\boldu^L)$, $\boldDelta^L=\boldDelta(\run^L)$,
$\run=\run(b,\boldu)$, $\sigma=\sigma(\boldu)$ and $\Delta=\Delta(\boldu)$,
and let $\run^+$ and $\boldu^+$ be defined as in Section \ref{FermLikeSec}.
Using (\ref{DeltaEq}), set $\boldDelta=\boldDelta(\run)$ and 
$\boldDelta^+=\boldDelta(\run^+)$.
On exchanging $\boldu^R=\boldu$ and $\run^R=\run$ for
$\boldu^R=\boldu^+$ and $\run^R=\run^+$ in the definition
(\ref{ChangeEq2}), the only possible changes are in the
terms $\tilde{\boldn}^T\boldB\tilde{\boldn}$,
$\boldN(\boldu^R)\cdot\tilde{\boldn}$ and $\gamma(\run^L,\run^R)$
(by virtue of $\boldDelta^R$ changing from
$\boldDelta(\run)$ to $\boldDelta(\run^+)$).

In the case $\sigma=t_1$, the first two of these are clearly unchanged.
{}From (\ref{Const1Eq}), (\ref{Const2Eq}) and (\ref{Const3Eq}),
it is easily shown that
$\gamma_{t_1}=-\gamma_{t_1+1}-\alpha_{t_1+1}^2
              -2\alpha_{t_1+1}(\beta_{t_1+1}+(\boldDelta^L)_{t_1+1})$.
Then since $(\boldDelta^R)_j=0$ for $0\le j\le t_1$ and $\beta'_j$
is unchanged for $0\le j\le t_1$, it follows that $\gamma_0$ is unchanged.
The lemma follows in the $\sigma=t_1>0$ case.
When $\sigma=t_1=0$ (here $p'<2p$), the lemma follows after noting that
$\boldu\in{\mathcal U}(b)\backslash\overline{\mathcal U}(b)$ implies that
$\gamma(\run^L,\run)=\gamma_0$ and thus
$\gamma(\run^L,\run)=\gamma(\run^L,\run^+)$.

For $0\le\sigma<t_1$, the value of $\tilde{\boldn}^T\boldB\tilde{\boldn}$
is seen to change by $\frac14-(\tilde n_{\sigma+1}+\cdots+\tilde n_{t_1})$,
and the value of $\boldN(\boldu^R)\cdot\tilde{\boldn}$ is seen to change by
$-(\tilde n_{\sigma+1}+\cdots+\tilde n_{t_1})$.
Thus to prove the lemma it remains to show that
$\gamma(\boldu^L,\boldu^R)$ changes by $-1$.

For $0\le j<t_1$, we have
$\gamma_j=\gamma_{j+1}+2\alpha_{j+1}(\boldDelta^R)_{j+1}-\beta_j^2$.
Combining the $j=\sigma$ and $j=\sigma-1$ instances of this and
noting that $\alpha_j=\alpha_{j+1}+\beta_j$, gives:
$$
\gamma_{\sigma-1}=\gamma_{\sigma+1}
+2\alpha_{\sigma+1}((\boldDelta^R)_{\sigma+1}+(\boldDelta^R)_{\sigma})
+2\beta_{\sigma}(\boldDelta^R)_{\sigma}
-\beta_{\sigma}^2-\beta_{\sigma-1}^2
$$
for $0<\sigma<t_1$. 
In passing from $\boldu$ to $\boldu^+$, we have that
$(\boldDelta^R)_{\sigma+1}+(\boldDelta^R)_{\sigma}$ is unchanged,
$(\boldDelta^R)_{\sigma}$ changes from $\Delta$ to $0$
and $\beta_{\sigma}$ decreases by $\Delta$, with
$\beta_j$ unchanged for $0\le j<\sigma$.
Thereupon $\gamma(\run^L,\run^R)=\gamma_0$ changes by $-1$
for $0<\sigma<t_1$.

By the definition of Section \ref{ConSec}, for $\sigma=0$ and $p'>2p$,
we have on the one hand
$\gamma(\run^L,\run)=\gamma_0+2\Delta\alpha_0$
whereas
$\gamma(\run^L,\run^+)=\gamma_0$ on the other.
As above,
$\gamma_0=\gamma_1+2\alpha_1(\boldDelta^R)_1-\beta_0^2$.
We then see that $\gamma_0$ changes by
$2\alpha_1\Delta+2\beta_0\Delta-1=2\alpha_0\Delta-1$, whereupon
the required change of $-1$ in $\gamma(\run^L,\run^R)$ follows.
\cqfd
\medskip

\subsection{Assimilating the trees}\label{AssimilateTrees}

In Theorems \ref{Limit1Thrm} and \ref{Limit2Thrm} of the previous section,
we obtained fermionic expressions for the characters $\chi^{p,p'}_{r,s}$
in terms of the Takahashi trees for $s$ and $b$, where $b$ is related to
$r$ by (\ref{groundstatelabel}) for some $c\in\{b\pm1\}$.
In this section, we obviate the need to find such $b$ and $c$,
by obtaining more succinct expressions for $\chi^{p,p'}_{r,s}$ in terms
of the truncated Takahashi tree for $r$ and the Takahashi tree for $s$.
We do this via a detailed examination of the relationship between the
truncated Takahashi tree for $r$ and the Takahashi tree for $b$,
where $r$ and $b$ are related as above.

\begin{figure}[ht]
\includegraphics[scale=1.00]{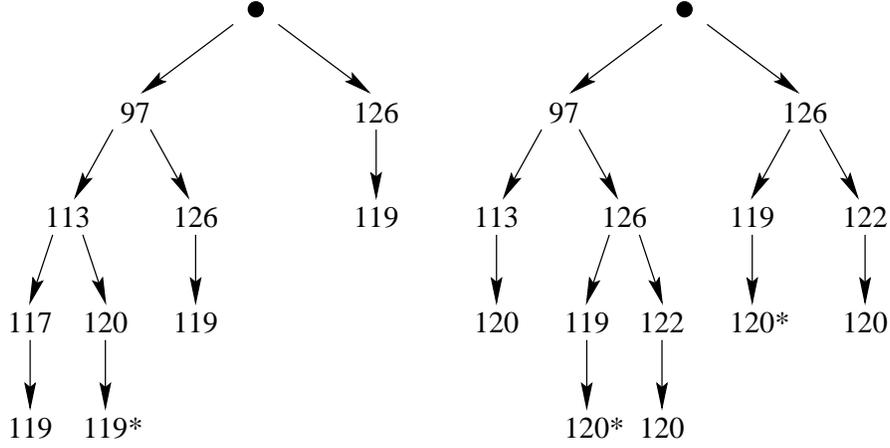}
\caption{\hbox{Takahashi trees for $b=119$ and $b=120$ when
         $(p,p')=(69,223)$.}}
\label{TreesbFig}
\medskip
\end{figure}

As an example, consider Fig.~\ref{TreesbFig}, where we give the Takahashi
trees for both $b=119$ and $b=120$ in the case where $(p,p')=(69,223)$.
Each of these values of $b$ is interfacial and, in fact,
$\rho^{69,223}(b)=37$ in both cases.
As will be shown in this section, for any $s$, the fermionic expressions
for $\chi^{69,223}_{37,s}$ that result from these two trees
are identical (up to trivial differences), and furthermore are identical
to the fermionic expression which results from the use of the truncated
Takahashi tree for $r=37$ that was given in Fig.~\ref{TruncTreeFig}.
(The leaf-nodes in Fig.~\ref{TreesbFig} that are affixed with
an asterisk are those for which the corresponding vector
$\boldu^R\in\overline{\mathcal U}(b)$. By virtue of Theorem \ref{Limit1Thrm},
these nodes do not contribute terms to the fermionic expression
(\ref{Limit1Eq1}) or (\ref{Limit1Eq2}) for $\chi^{69,223}_{37,s}$.)

In this section, we will sometimes refer to the Takahashi tree for $b$
as the $b$-tree and the truncated Takahashi tree for $r$ as the $r$-tree.
Since $p$ and $p'$ will be fixed throughout this section,
we abbreviate $\rho^{p,p'}(\cdot)$ to $\rho(\cdot)$.

First, we need to extend the results on the $b$-tree contained in
Lemma \ref{NodesMuLem}.

\begin{lemma}\label{NodesRhoLem}
Let $\{\tau_j,\sigma_j,\Delta_j\}_{j=1}^{d}$
be the naive run corresponding to the leaf-node
$b_{i_1i_2\cdots i_{d-1}0}$ of the Takahashi tree for $b$.
Then, for $1\le k<d$, both
$b_{i_1i_2\cdots i_{k-1}\bari_k}$ and
$b_{i_1i_2\cdots i_{k-1}i_{k}}$ are interfacial.
Moreover,
\begin{align*}
\rho(b_{i_1i_2\cdots i_{k-1}\bari_k})&= \sum_{m=2}^{k}
        \Delta_m (\tkappa_{\tau_m^{}}-\tkappa_{\sigma_m^{}})
\:+\:\begin{cases}
        \tkappa_{\sigma^{}_1} &\text{if $\Delta_1=-1$;}\\
        p-\tkappa_{\sigma^{}_1} &\text{if $\Delta_1=+1$,}
  \end{cases} \\
\intertext{and}
\rho(b_{i_1i_2\cdots i_{k-1}i_{k}})&=
\rho(b_{i_1i_2\cdots i_{k-1}\bari_k})
+\Delta_{k+1}\tkappa_{\tau^{}_{k+1}}.
\intertext{If $b$ is interfacial and $\sigma_d>0$ then:}
\rho(b)&= \sum_{m=2}^{d}
        \Delta_m (\tkappa_{\tau_m^{}}-\tkappa_{\sigma_m^{}})
\:+\:\begin{cases}
        \tkappa_{\sigma^{}_1} &\text{if $\Delta_1=-1$;}\\
        p-\tkappa_{\sigma^{}_1} &\text{if $\Delta_1=+1$.}
  \end{cases}
\end{align*}
If $b$ is not interfacial and $\sigma_d>0$ then:
\begin{equation*}
\left\lfloor\frac{bp}{p'}\right\rfloor=\sum_{m=2}^{d}
        \Delta_m (\tkappa_{\tau_m^{}}-\tkappa_{\sigma_m^{}})
\:-\:\frac12(1\pm\Delta_d)
\:+\:\begin{cases}
        \tkappa_{\sigma^{}_1} &\text{if $\Delta_1=-1$;}\\
        p-\tkappa_{\sigma^{}_1} &\text{if $\Delta_1=+1$,}
  \end{cases}
\end{equation*}
where the \lq${}-{}$\rq\ sign applies when $p'>2p$ and the
\lq${}+{}$\rq\ sign applies when $p'<2p$.
\end{lemma}

\Proof First consider the case $p'>2p$.

On comparing the definitions of Section \ref{CoreSec} with
the expressions stated in Lemma \ref{NodesMuLem}, we see that
$b_{i_1i_2\cdots i_{k-1}\bari_k}=\nu^*_k$ and
$b_{i_1i_2\cdots i_{k-1}i_k}=\nu_k$ for $1\le k<d$.
Use of Corollary \ref{StarInterCor} then shows that
$b_{i_1i_2\cdots i_{k-1}\bari_k}$ is interfacial and
$\rho(b_{i_1i_2\cdots i_{k-1}\bari_k})=\tilde\nu^*_k$,
the definition of which in Section \ref{CoreSec} proves the
required expression for $\rho(b_{i_1i_2\cdots i_{k-1}\bari_k})$.

In the case of $b_{i_1i_2\cdots i_{k-1}i_{k}}$, first consider
$\sigma_{k+1}\ne\tau_{k+1}$.
Then Corollary \ref{StarInterCor} shows that
$b_{i_1i_2\cdots i_{k-1}\bari_k}$ is interfacial and
$\rho(b_{i_1i_2\cdots i_{k-1}i_{k}})
 -\rho(b_{i_1i_2\cdots i_{k-1}\bari_k})
 =\tilde\nu_k-\tilde\nu^*_k
 =\Delta_{k+1}\tkappa_{\tau^{}_{k+1}}$
as required.

If $\tau_{k+1}=\sigma_{k+1}$, consider any leaf-node of the form
$b_{i_1\cdots i_k\bari_{k+1}i'_{k+2}\cdots i'_{d'}}$,
where necessarily $d'>k+1$.
Using Lemma \ref{NodesMuLem}, it is easily seen that the
corresponding naive run 
$\{\tau'_j,\sigma'_j,\Delta'_j\}_{j=1}^{d'}$ is such that $\tau'_j=\tau_j$
for $1\le j\le k+1$ and $\sigma'_j=\sigma_j$ for $1\le j\le k$.
Moreover, $\sigma'_{k+1}\ne\sigma_{k+1}$
and so $\sigma'_{k+1}\ne\tau'_{k+1}$.
Then the above argument applied to this case shows that
$b_{i_1i_2\cdots i_{k-1}i_{k}}$ is interfacial and yields the
required expression for $\rho(b_{i_1i_2\cdots i_{k-1}i_{k}})$.

When $\sigma_d>0$, Corollary \ref{StarInterCor}
immediately gives the required expressions for $\rho(b)$
and $\lfloor bp/p'\rfloor$.

In the case where $p'<2p$, note that $p'>2(p'-p)$.
We are then able to use the result established above for
$(p,p')\to(p'-p,p')$.
If $\{\kappa_j'\}_{j=1}^t$ and $\{\tkappa_j'\}_{j=1}^t$ are respectively
the Takahashi lengths and the truncated Takahashi lengths for the
$(p'-p,p')$-model then Lemma \ref{DmodelLem} yields
$\kappa_j'=\kappa_j$ and $\tkappa_j'=\kappa_j-\tkappa_j$
for $1\le j\le t$.
Lemma \ref{InterLem}(2) implies that
$b_{i_1i_2\cdots i_{k-1}\bari_k}$ is interfacial in the $(p,p')$-model,
and that
\begin{align*}
\rho^{p,p'}(b_{i_1i_2\cdots i_{k-1}\bari_k})&=
b_{i_1i_2\cdots i_{k-1}\bari_k}-
\rho^{p'-p,p'}(b_{i_1i_2\cdots i_{k-1}\bari_k})\\
&\hskip-8mm
=b_{i_1i_2\cdots i_{k-1}\bari_k}
-\sum_{m=2}^{k}
        \Delta_m (\tkappa'_{\tau_m^{}}-\tkappa'_{\sigma_m^{}})
-\begin{cases}
        \tkappa'_{\sigma^{}_1} &\text{if $\Delta_1=-1$;}\\
        p'-p-\tkappa'_{\sigma^{}_1} &\text{if $\Delta_1=+1$.}
  \end{cases}
\end{align*}
The desired expression for
$\rho^{p,p'}(b_{i_1i_2\cdots i_{k-1}\bari_k})$ then results
from the expression for $b_{i_1i_2\cdots i_{k-1}\bari_k}$
given in Lemma \ref{NodesMuLem} and
$\kappa_j'=\kappa_j$ and $\tkappa_j'=\kappa_j-\tkappa_j$
for $1\le j\le t$ (the $j=0$ case is not required because $k<d$).

The other expressions follow in a similar way, noting in the
final case that $\lfloor bp/p'\rfloor=b-1-\lfloor b(p'-p)/p'\rfloor$.
\cqfd
\medskip

In particular, this lemma implies that every node of the $b$-tree that is
not a leaf-node is interfacial.

Now consider a leaf-node $r_{i_1i_2\cdots i_d}$ of the truncated
Takahashi tree for $r$, and let 
$\{\tilde\tau_j,\tilde\sigma_j,\tilde\Delta_j\}_{j=1}^d$ be the
corresponding run. That it is actually a naive run follows
from the appropriate analogue of Lemma \ref{NaiveRunLem}.
We also have an analogue of Lemma \ref{NodesMuLem}, which
for ease of reference, we state here in full:

\begin{lemma}\label{NodesTildeLem}
Let $\{\tilde\tau_j,\tilde\sigma_j,\tilde\Delta_j\}_{j=1}^{d}$
be the naive run corresponding to the leaf-node
$r_{i_1i_2\cdots i_{d-1}0}$ of the truncated Takahashi tree for $r$.
Then
\begin{align*}
r_{i_1i_2\cdots i_{k-1}\bari_k}&= \sum_{m=2}^{k}
        \tilde\Delta_m (\tkappa_{\tilde\tau_m^{}}
                           -\tkappa_{\tilde\sigma_m^{}})\quad
+\quad\begin{cases}
        \tkappa_{\tilde\sigma^{}_1} &\text{if $\tilde\Delta_1=-1$;}\\
        p-\tkappa_{\tilde\sigma^{}_1} &\text{if $\tilde\Delta_1=+1$,}
  \end{cases} \\
\intertext{for $1\le k\le d$, and}
r_{i_1i_2\cdots i_{k-1}i_{k}}&=
r_{i_1i_2\cdots i_{k-1}\bari_k}
      +\tilde\Delta_{k+1}\tkappa_{\tilde\tau^{}_{k+1}},
\end{align*}
for $1\le k<d$.
\end{lemma}

\Proof This is proved in the same way as Lemma \ref{NodesMuLem}.
\cqfd
\medskip

For Lemmas \ref{BBLem} through \ref{GGLem}, and Theorem \ref{CombinedThrm}
we let $1\le r<p$ and select $b$ and $c$ such that
$b,c\in\{\lfloor rp'/p\rfloor,\lfloor rp'/p\rfloor+1\}$ with $b\ne c$.
This ensures that (\ref{groundstatelabel}) is satisfied,
and that $b$ and $c$ straddle the $r$th odd band.
If $b$ is interfacial then $\rho(b)=r$, and in this case
we could actually choose either $c=b\pm1$.
If $b$ is non-interfacial then necessarily $p'<2p$.
It will be useful to note that if
$\{\tau_j,\sigma_j,\Delta_j\}_{j=1}^d$ is the naive run
corresponding to a leaf-node of the $b$-tree then $\sigma_d\le t_2$.
This follows from Corollary \ref{CoreCor} (replacing $t_1$ with
$t_2$ because we are dealing with the $p'<2p$ case).
It will be also be useful to note that $\kappa_j=j+1$
and $\tkappa_j=j$ for $1\le j\le t_2+1$ in this $p'<2p$ case.

The next few lemmas will enable us to compare the $b$-tree with
the $r$-tree.

\begin{lemma}\label{BBLem}
For $1\le r<p$, let $b,c\in\{\lfloor rp'/p\rfloor,\lfloor rp'/p\rfloor+1\}$
with $b\ne c$.
For $k\ge1$, let $b_{i_1i_2\cdots i_k}$ be a non-leaf-node
of the $b$-tree with $\rho(b_{i_1i_2\cdots i_k})=r$.
Then $b_{i_1i_2\cdots i_k}=b-(-1)^{i_k}$ and
$b_{i_1i_2\cdots i_k}$ is a through-node of the $b$-tree.
If $\boldu\in{\mathcal U}(b)$ is the vector corresponding to the leaf-node
$b_{i_1i_2\cdots i_k0}$, then $\boldu\in\overline{\mathcal U}(b)$.
Moreover, if $b$ is not interfacial then
$b-c=(-1)^{i_k}$ and $\boldu\in\overline{\mathcal U}{}^{c-b}(b)$
\end{lemma}

\Proof $\rho(b_{i_1i_2\cdots i_k})=r$ implies that
$b_{i_1i_2\cdots i_k}$ borders the $r$th odd band.
So does $b$ and thus $\vert b_{i_1i_2\cdots i_k}-b\vert=1$.
The definition of the $b$-tree then stipulates that
$b_{i_1i_2\cdots i_k}$ is a through-node,
$b_{i_1i_2\cdots i_k0}$ is a leaf-node
and $b_{i_1i_2\cdots i_k}-b=-(-1)^{i_k}\kappa_{\sigma_{k+1}}$.
Thereupon, $\kappa_{\sigma_{k+1}}=1$ so that $\sigma_{k+1}=0$.

For the vector $\boldu$, we then have $\sigma(\boldu)=0$
and $\Delta(\boldu)=-(-1)^{i_k}$.
Since $b$ and $b_{i_1i_2\cdots i_k}$ are either side of the
$r$th odd band, we have that
$\lfloor bp/p'\rfloor\ne\lfloor(b+\Delta(\boldu))p/p'\rfloor$.
Then $\boldu\in\overline{\mathcal U}(b)$ by definition.
In the case that $b$ is not interfacial, then necessarily
$c=b_{i_1i_2\cdots i_k}$ so that $c-b=-(-1)^{i_k}=\Delta(\boldu)$,
giving $\boldu\in\overline{\mathcal U}{}^{c-b}(b)$ as required.
\cqfd
\medskip

\begin{lemma}\label{DDLem}
For $1\le r<p$, let $b,c\in\{\lfloor rp'/p\rfloor,\lfloor rp'/p\rfloor+1\}$
with $b\ne c$.

1) For $d\ge2$, let $b_{i_1i_2\cdots i_{d-1}}$ be a through-node
of the $b$-tree with $\rho(b_{i_1i_2\cdots i_{d-1}})\ne r$,
and let $\boldu\in{\mathcal U}(b)$ be the vector corresponding to
the leaf-node $b_{i_1i_2\cdots i_{d-1}0}$.
If $b$ is interfacial then $\boldu\not\in\overline{\mathcal U}(b)$,
and if $b$ is not interfacial then
$\boldu\not\in\overline{\mathcal U}{}^{c-b}(b)$.
In addition, if $r_{i_1i_2\cdots i_{d-1}}$ is a node of the $r$-tree
with $r_{i_1i_2\cdots i_{d-1}}=\rho(b_{i_1i_2\cdots i_{d-1}})$ then
$r_{i_1i_2\cdots i_{d-1}}$ is a through-node.

2) Let $b_0$ be a leaf-node of the $b$-tree,
and let $\boldu\in{\mathcal U}(b)$ be the corresponding vector.
If $b$ is interfacial then $\boldu\not\in\overline{\mathcal U}(b)$,
and if $b$ is not interfacial then
$\boldu\not\in\overline{\mathcal U}{}^{c-b}(b)$.
In addition, $r_0$ is a leaf-node of the $r$-tree.
\end{lemma}

\Proof
1) Lemma \ref{NodesMuLem} implies that
$b_{i_1i_2\cdots i_{d-1}}-b=\Delta_d\kappa_{\sigma_d}$,
for $\Delta_{d}=-(-1)^{i_k}$ and some $\sigma_d$.
Note that $\sigma(\boldu)=\sigma_{d}$ and $\Delta(\boldu)=\Delta_{d}$.
If $\sigma_{d}>0$ then $\boldu\not\in\overline{\mathcal U}(b)$ by definition.
If $\sigma_{d}=0$ and $b$ is interfacial, we claim that
$\lfloor bp/p'\rfloor=\lfloor(b+\Delta_d)p/p'\rfloor$.
Otherwise, $\sigma_{d}=0$ and
$\lfloor bp/p'\rfloor\ne\lfloor(b+\Delta_d)p/p'\rfloor$
imply that $b_{i_1i_2\cdots i_{d-1}}$ neighbours the
$r$th odd band, as does $b$, whereupon noting that
$b_{i_1i_2\cdots i_{d-1}}$ is interfacial (by Lemma \ref{NodesRhoLem}),
we would have $\rho(b_{i_1i_2\cdots i_{d-1}})=r$, which is not the case.
So $\boldu\not\in\overline{\mathcal U}(b)$.
If $\sigma_{d}=0$ and $b$ is not interfacial, we claim that
$c=b-\Delta$. Otherwise $c=b+\Delta=b_{i_1i_2\cdots i_{d-1}}$,
again implying that $\rho(b_{i_1i_2\cdots i_{d-1}})=r$, which is not the case.
Therefore $\boldu\in{\mathcal U}^{b-c}(b)$ so that certainly
$\boldu\not\in\overline{\mathcal U}{}^{c-b}(b)$ as required.

For $\sigma_d>0$ and $b$ interfacial, Lemma \ref{NodesRhoLem} implies
that $r_{i_1i_2\cdots i_{d-1}}-r=\Delta_d\tkappa_{\sigma_d}$.
For $\sigma_d>0$ and $b$ non-interfacial (when necessarily $p'<2p$),
Lemma \ref{NodesRhoLem} implies
that $r_{i_1i_2\cdots i_{d-1}}-r=\Delta_d\tkappa_{\sigma_d}$ or
$r_{i_1i_2\cdots i_{d-1}}-r
 =\Delta_d(\tkappa_{\sigma_d}+1)=\Delta_d(\tkappa_{\sigma_d+1})$.
In either case, the definition of the $r$-tree implies that
$r_{i_1i_2\cdots i_{d-1}}$ is a through-node of the $r$-tree.
For $\sigma_d=0$,
$b_{i_1i_2\cdots i_{d-1}}-b=\Delta_d$ implies that
$r_{i_1i_2\cdots i_{d-1}}-r=\Delta_d$ whence
$r_{i_1i_2\cdots i_{d-1}}$ is a through-node of the $r$-tree.

2) This case arises only if $b\in{\mathcal T}\cup{\mathcal T}'$ so
that either $b=\kappa_{\sigma_1}$ or $b=p'-\kappa_{\sigma_1}$ for
some $\sigma_1$.
If $\sigma_1>0$ then immediately $\boldu\not\in\overline{\mathcal U}(b)$.
If $\sigma_1>0$ and $b$ is interfacial then
Lemma \ref{NodesRhoLem} also implies that either
$r=\rho(b)=\tkappa_\sigma$ or $r=\rho(b)=p-\tkappa_\sigma$.
It follows that $r_0$ is a leaf-node of the $r$-tree.
If $\sigma_1>0$ and $b$ is not interfacial then
Lemma \ref{NodesRhoLem} implies that either
$\lfloor bp/p'\rfloor=\tkappa_{\sigma_1}$
or $\lfloor bp/p'\rfloor=p-1-\tkappa_{\sigma_1}$.
Since $r=\lfloor bp/p'\rfloor$ or $r=\lfloor bp/p'\rfloor+1$,
and noting that $p'<2p$ and $\sigma_1\le t_2$ here,
it immediately follows that $r_0$ is a leaf-node of the $r$-tree.

If $\sigma_1=0$ then either $b=1$ and $\Delta_1=-1$,
or $b=p'-1$ and $\Delta_1=1$
(this case only arises when $p'<2p$).
That $1\le r<p$ forces $c=2$ or $c=p'-2$ respectively:
i.e.~$c=b-\Delta_d$ so that $\boldu\in{\mathcal U}^{b-c}(b)$ and
thus $\boldu\not\in\overline{\mathcal U}{}^{c-b}(b)$ as required.
We also have $r=1$ or $r=p-1$ respectively so that
$r_0$ is a leaf-node of the $r$-tree.
\cqfd
\medskip

\begin{lemma}\label{EELem}
For $1\le r<p$, let $b,c\in\{\lfloor rp'/p\rfloor,\lfloor rp'/p\rfloor+1\}$
with $b\ne c$.

(1) For $k\ge1$, let $b_{i_1i_2\cdots i_k}$ be a branch-node of
the $b$-tree and $r_{i_1i_2\cdots i_k}$ a node of the $r$-tree
such that $\rho(b_{i_1i_2\cdots i_k})=r_{i_1i_2\cdots i_k}\ne r$.
If $\rho(b_{i_1i_2\cdots i_kh})=r$ for some $h\in\{0,1\}$,
then $r_{i_1i_2\cdots i_k}$ is a through-node of the $r$-tree.

(2) Let the $b$-tree nodes $b_0$ and $b_1$ be such that $b_0\ne b_1$
with $\rho(b_{h})=r$ for either $h=0$ or $h=1$.
Then $r\in\tilde{\mathcal T}\cup\tilde{\mathcal T}'$ and
$r_0$ is a leaf-node. Set $k=0$ for what follows.

(3) Let $h,k,i_1,i_2,\ldots,i_k$ be as in (1) or (2) above.
There exists a unique $d\ge k+2$ such that if we set
$i_{k+1}=\cdots=i_{d-1}=\barh$ then
$b_{i_1i_2\cdots i_{d-1}0}$ is a leaf-node.
The vector $\boldu\in\mathcal U(b)$ corresponding to this leaf-node is
such that $\boldu\not\in\overline{\mathcal U}(b)$ when $b$ is interfacial,
and $\boldu\in{\mathcal U}{}^{b-c}(b)$ when $b$ is not interfacial.
In addition, if $k+2\le j\le d$ then $b_{i_1i_2\cdots i_{j-2}h}=b-(-1)^h$
and $b_{i_1i_2\cdots i_{j-2}h0}$ is a leaf-node.
In each of these instances, the corresponding vector
$\boldu\in\mathcal U(b)$ is such that $\boldu\in\overline{\mathcal U}(b)$,
and also such that $\boldu\in\overline{\mathcal U}{}^{c-b}(b)$
when $b$ is not interfacial.

If $\{\tau_j,\sigma_j,\Delta_j\}_{j=1}^{d}$ is the naive run
corresponding to the leaf-node $b_{i_1i_2\cdots i_{d-1}0}$
then $\tau_j=\sigma_j$ for $k+2\le j<d$ and $\Delta_j=(-1)^h$
for $k+2\le j\le d$.
In addition, $\sigma_d+1=\tau_d\le t_1+1$ in the $p'>2p$ case,
and $\sigma_d+1=\tau_d\le t_2+1$ in the $p'<2p$ case.
If $b$ is interfacial then $t_1\le\sigma_d+1$.
\end{lemma}

\Proof (1) Consider a leaf-node of the form
$b_{i_1\cdots i_k\barh\cdots{}}$.
Lemma \ref{NodesMuLem} shows that
$b_{i_1i_2\cdots i_kh}
=b_{i_1i_2\cdots i_k}+(-1)^{i_k}\kappa_{\sigma}$
for some $\sigma$.
Lemma \ref{NodesRhoLem} then shows that
$b_{i_1i_2\cdots i_kh}$ and $b_{i_1i_2\cdots i_k}$
are both interfacial with
$\rho(b_{i_1i_2\cdots i_kh})
=\rho(b_{i_1i_2\cdots i_k})+(-1)^{i_k}\tkappa_{\sigma}$.
Since $r_{i_1i_2\cdots i_k}=\rho(b_{i_1i_2\cdots i_k})$ and
$\rho(b_{i_1i_2\cdots i_kh})=r$, the definition of
the $r$-tree stipulates that $r_{i_1i_2\cdots i_k}$ is a through-node.

(2) Consider a leaf-node of the form $b_{\barh\cdots{}}$.
Lemma \ref{NodesMuLem} shows that either
$b_{i_1}=\kappa_\sigma$ or $b_{i_1}=p'-\kappa_\sigma$ for some $\sigma$.
Lemma \ref{NodesRhoLem} shows that $b_{i_1}$ is interfacial and
$\rho(b_{i_1})=\tkappa_\sigma$ or $\rho(b_{i_1})=p-\tkappa_\sigma$
respectively.
Then $r=\rho(b_{i_1})\in\tilde{\mathcal T}\cup\tilde{\mathcal T}'$
and the definition of the $r$-tree stipulates that
$r_0$ is a leaf-node.

(3) $b_{i_1i_2\cdots i_kh}$ is a non-leaf-node of the $b$-tree.
Use of Lemma \ref{BBLem} shows that it is a through-node
and that the vector $\boldu\in\mathcal U(b)$ corresponding to the
leaf-node $b_{i_1i_2\cdots i_kh0}$ is such that
$\boldu\in\overline{\mathcal U}(b)$, and also such that
$\boldu\in\overline{\mathcal U}{}^{c-b}(b)$ when $b$ is not interfacial.
With $\Delta=(-1)^{h}$, Lemma \ref{BBLem} also shows that
$b_{i_1i_2\cdots i_kh}=b-\Delta$,
and if $b$ is non-interfacial then $\Delta=b-c$.

$b_{i_1i_2\cdots i_k\barh}$ is a non-leaf-node of the $b$-tree and
is thus interfacial.
$\rho(b_{i_1i_2\cdots i_k\barh})\ne r$ since otherwise one of
$b_{i_1i_2\cdots i_kh}$ and $b_{i_1i_2\cdots i_k\barh}$ would
necessarily be equal to $b$. So either $b_{i_1i_2\cdots i_k\barh}$
is a through-node or a branch-node.
In the former case, let $\boldu\in\mathcal U(b)$ be the vector
corresponding to the leaf-node $b_{i_1i_2\cdots i_k\barh0}$.
If $b$ is interfacial, Lemma \ref{DDLem} implies that
$\boldu\not\in\overline{\mathcal U}(b)$, and
if $b$ is non-interfacial then $\Delta(\boldu)=-(-1)^{\barh}=\Delta=b-c$
so that $\boldu\in{\mathcal U}{}^{b-c}(b)$, as required.

In the case that $b_{i_1i_2\cdots i_k\barh}$ is a branch-node,
there exists no $\sigma$ such that
$b_{i_1i_2\cdots i_k\barh}-\Delta\kappa_{\sigma}=b$.
However, if $\tau$ is such that
$\kappa_\tau=b_{i_1i_2\cdots i_k1}-b_{i_1i_2\cdots i_k0}$,
then $b_{i_1i_2\cdots i_k\barh}-\Delta\kappa_\tau=
      b_{i_1i_2\cdots i_kh}=b-\Delta$,
whereupon the definition of the $b$-tree stipulates that
$b_{i_1i_2\cdots i_k\barh h}=b_{i_1i_2\cdots i_kh}$
and $b_{i_1i_2\cdots i_k\barh\,\barh}=
     b_{i_1i_2\cdots i_k\barh}-\Delta\kappa_{\tau-1}$.

Since $\rho(b_{i_1i_2\cdots i_k\barh h})=r$,
the node $b_{i_1i_2\cdots i_k\barh}$ now satisfies the same criteria
as did $b_{i_1i_2\cdots i_k}$.
We then iterate the above argument until we encounter a through-node
of the form $b_{i_1i_2\cdots i_k\barh\,\barh\cdots\barh}$.
This must eventually occur because the $b$-tree is finite.
Then as above, the vector $\boldu\in\mathcal U(b)$
corresponding to the leaf-node
$b_{i_1i_2\cdots i_k\barh\,\barh\cdots\barh0}$ is such that
$\boldu\not\in\overline{\mathcal U}(b)$ when $b$ is interfacial, and
$\boldu\in{\mathcal U}{}^{b-c}(b)$ when $b$ is not interfacial.
Set $d$ to be the number of levels that this leaf-node occurs
below the root node.
Then also as above, for $k+2\le j\le d$
we have $b_{i_1i_2\cdots i_ki_{k+1}\cdots i_{j-2}h}=b-\Delta$
and that the vector $\boldu\in\mathcal U(b)$ corresponding to
the leaf-node $b_{i_1i_2\cdots i_ki_{k+1}\cdots i_{j-2}h0}$
is such that $\boldu\in\overline{\mathcal U}(b)$,
and also satisfies $\boldu\in\overline{\mathcal U}{}^{c-b}(b)$
when $b$ is not interfacial.
This accounts for all descendents of $b_{i_1i_2\cdots i_k}$.

Since $b_{i_1\cdots i_{j-2}\bari_{j-1}}=b_{i_1\cdots i_{j-1}\bari_{j}}$
for $k+2\le j<d$, it follows from Lemma \ref{NodesMuLem} that
$\tau_j=\sigma_j$.
For $k+2\le j\le d$, that $\Delta_j=(-1)^h$ follows immediately
from $i_{j-1}=\barh$

Finally, we have
$b_{i_1\cdots i_{d-2}\bari_{d-1}}=b_{i_1\cdots i_{k}h}=b-\Delta$,
and from Lemma \ref{NodesMuLem}, we have
$b_{i_1\cdots i_{d-1}}-b_{i_1\cdots i_{d-2}\bari_{d-1}}
 =\Delta\kappa_{\tau_d}$ and
$b=b_{i_1\cdots i_{d-1}}-\Delta\kappa_{\sigma_d}$.
Together, these imply that $\kappa_{\sigma_d}=\kappa_{\tau_d}-1$.
In the $p'>2p$ case this can only occur if $\sigma_d+1=\tau_d\le t_1+1$,
whereas in the $p'<2p$ case this can only occur if
$\sigma_d+1=\tau_d\le t_2+1$.
In addition, since $b_{i_1\cdots i_{d-1}}$ is interfacial
with $\rho(b_{i_1\cdots i_{d-1}})\ne r$, and odd bands in the
$(p,p')$-model are separated by at least $t_1$ even bands,
$b=b_{i_1\cdots i_{d-1}}-\Delta\kappa_{\sigma_d}$ implies that
if $b$ is interfacial then
$t_1\le\kappa_{\sigma_d}=\sigma_d+1$ which completes the proof.
\cqfd
\medskip

\begin{lemma}\label{AALem}
For $1\le r<p$, let $b,c\in\{\lfloor rp'/p\rfloor,\lfloor rp'/p\rfloor+1\}$
with $b\ne c$.

(1) For $k\ge1$, let $b_{i_1i_2\cdots i_k}$ be a branch-node of
the $b$-tree and $r_{i_1i_2\cdots i_k}$ a node of the $r$-tree
such that $\rho(b_{i_1i_2\cdots i_k})=r_{i_1i_2\cdots i_k}\ne r$.
If
$$\rho(b_{i_1i_2\cdots i_k0})<r<\rho(b_{i_1i_2\cdots i_k1})$$
then $r_{i_1i_2\cdots i_k}$ is a branch-node of the $r$-tree
and moreover
$\rho(b_{i_1i_2\cdots i_k0})=r_{i_1i_2\cdots i_k0}$ and
$\rho(b_{i_1i_2\cdots i_k1})=r_{i_1i_2\cdots i_k1}$.

(2) Let $b_0$ and $b_1$ be nodes of the $b$-tree such that
$\rho(b_0)<r<\rho(b_1)$.
Then $r_0$ and $r_1$ are nodes of the $r$-tree with
$\rho(b_0)=r_0$ and $\rho(b_1)=r_1$.
\end{lemma}

\Proof (1) The definition of the $b$-tree implies that there exists
$x$ such that:
\begin{align*}
b_{i_1i_2\cdots i_k0}&=b_{i_1i_2\cdots i_k}+(-1)^{i_k}\kappa_{x+i_k},\\
b_{i_1i_2\cdots i_k1}&=b_{i_1i_2\cdots i_k}+(-1)^{i_k}\kappa_{x+1-i_k}.
\end{align*}
By choosing an appropriate leaf-node, Lemmas \ref{NodesMuLem}
and \ref{NodesRhoLem} then imply that
\begin{align*}
\rho(b_{i_1i_2\cdots i_k0})=\rho(b_{i_1i_2\cdots i_k})
 +(-1)^{i_k}\tkappa_{x+i_k},\\
\rho(b_{i_1i_2\cdots i_k1})=\rho(b_{i_1i_2\cdots i_k})
 +(-1)^{i_k}\tkappa_{x+1-i_k}.
\end{align*}
By hypothesis, $r$ is strictly between these values.
Since $\rho(b_{i_1i_2\cdots i_k})=r_{i_1i_2\cdots i_k}\ne r$,
the definition of the $r$-tree stipulates that
$r_{i_1i_2\cdots i_k}$ is a branch-node and that
$r_{i_1i_2\cdots i_k0}=\rho(b_{i_1i_2\cdots i_k0})$ and
$r_{i_1i_2\cdots i_k1}=\rho(b_{i_1i_2\cdots i_k1})$.

(2) There are four cases to consider depending on which of the sets
${\mathcal T}$ and ${\mathcal T}'$ contains $b_0$ and $b_1$.
We consider only the case where $b_0,b_1\in{\mathcal T}$:
the other cases are similar.

The definition of the $b$-tree shows that there exists $\sigma_1$
such that $b_0=\kappa_{\sigma_1}$ and $b_1=\kappa_{\sigma_1+1}$
Lemmas \ref{NodesMuLem} and \ref{NodesRhoLem} imply that
$\rho(b_0)=\tkappa_{\sigma_1}$ and $\rho(b_1)=\tkappa_{\sigma_1+1}$.
We claim that $r_0=\tkappa_{\sigma_1}$ and $r_1=\tkappa_{\sigma_1+1}$.
This is so because $\tkappa_{\sigma_1}<r<\tkappa_{\sigma_1+1}$ and if
there were any elements of $\tilde{\mathcal T}\cup\tilde{\mathcal T}'$
strictly between $\tkappa_{\sigma_1}$ and $\tkappa_{\sigma_1+1}$ then,
via Lemma \ref{OddPosLem}(1), there would be elements of
${\mathcal T}\cup{\mathcal T}'$ strictly between $\kappa_{\sigma_1}$ and
$\kappa_{\sigma_1+1}$ whereupon the description of the $b$-tree would
not then yield $b_0=\kappa_{\sigma_1}$ and $b_1=\kappa_{\sigma_1+1}$.
\cqfd
\medskip

\begin{lemma}\label{GGLem}
For $1\le r<p$, let $b,c\in\{\lfloor rp'/p\rfloor,\lfloor rp'/p\rfloor+1\}$
with $b\ne c$.
For $d\ge1$, let $b_{i_1i_2\cdots i_{d-1}i_d}$ be a leaf-node of
the $b$-tree and let $\{\tau_j,\sigma_j,\Delta_j\}_{j=1}^d$ be
the corresponding naive run.
Assume that $r_{i_1i_2\cdots i_{\tilde d-1}i_{\tilde d}}$
is a leaf-node of the $r$-tree for some $\tilde d$ with $1\le\tilde d\le d$,
and let $\{\tilde\tau_j,\tilde\sigma_j,\tilde\Delta_j\}_{j=1}^{\tilde d}$
be the corresponding naive run.
Then $\tilde\tau_j=\tau_j$ and $\tilde\Delta_j=\Delta_j$ for
$1\le j\le\tilde d$, and $\tilde\sigma_j=\sigma_j$ for $1\le j<\tilde d$.

If either both $\tilde d<d$ and
$\rho(b_{i_1i_2\cdots i_{\tilde d-1}\bari_{\tilde d}})=r$,
or both $\tilde d=d$ and $b$ is interfacial, then
$\tilde\sigma_{\tilde d}=\max\{\sigma_{\tilde d},t_1+1\}$.
If $\tilde d=d$ and $b$ is not interfacial then
\begin{equation*}
\tilde\sigma_d=
 \begin{cases}
   \sigma_d&\text{if $c-b=\Delta_d$;}\\
   \sigma_d+1&\text{if $c-b\ne\Delta_d$.}
 \end{cases}
\end{equation*}
\end{lemma}

\Proof
We first claim that $\rho(b_{i_1\cdots i_k0})\ne r$ and
$\rho(b_{i_1\cdots i_k1})\ne r$ for $0\le k\le \tilde d-2$.
Otherwise, let $k$ be the smallest value for which the claim doesn't
hold, so that $\rho(b_{i_1\cdots i_kh})=r$ for some $h\in\{0,1\}$.
If $k=0$ (in which case $\tilde d\ge2$), then Lemma \ref{EELem}(2)
shows that $r_0$ is a leaf-node, whereas if $k\ge1$ then
Lemma \ref{EELem}(1) shows that $r_{i_1\cdots i_k0}$ is a leaf-node.
This contradicts the hypothesis that $r_{i_1i_2\cdots i_{\tilde d-1}0}$
is a leaf-node, thus establishing the claim.

Lemma \ref{AALem} now shows that
$\rho(b_{i_1\cdots i_k0})=r_{i_1\cdots i_k0}$ and
$\rho(b_{i_1\cdots i_k1})=r_{i_1\cdots i_k1}$ for
$0\le k\le \tilde d-2$.
Comparison of Lemmas \ref{NodesRhoLem} and \ref{NodesTildeLem}
then shows that $\tilde\tau_j=\tau_j$ for $1\le j\le\tilde d$, and
$\tilde\sigma_j=\sigma_j$ for $1\le j<\tilde d$.
That $\tilde\Delta_j=\Delta_j$ for $1\le j\le\tilde d$ is immediate
from their definitions. It remains to prove the expressions for
$\tilde\sigma_{\tilde d}$.

If $d>\tilde d>1$ then
Lemma \ref{NodesRhoLem} implies that
$\vert\rho(b_{i_1\cdots i_{\tilde d-1}})
 -\rho(b_{i_1\cdots i_{\tilde d-1}\bari_{\tilde d}})\vert
 =\tkappa_{\sigma_{\tilde d}}$.
Thereupon, $\vert r_{i_1\cdots i_{\tilde d-1}}-r\vert
       =\tkappa_{\sigma_{\tilde d}}$
so that $\tkappa_{\tilde\sigma_{\tilde d}}=\tkappa_{\sigma_{\tilde d}}$.
On noting that $\tkappa_i=1$ if and only if $0\le i\le t_1+1$,
the required result follows here.

If $d=\tilde d>1$ with $b$ interfacial (and thus $\rho(b)=r$),
a similar argument applies when $\sigma_d>0$.
If $\sigma_{d}=0$ then $\vert b_{i_1\cdots i_{d-1}}-b\vert=1$
so that, noting that $r_{i_1\cdots i_{\tilde d-1}}\ne r$,
then necessarily $\vert r_{i_1\cdots i_{\tilde d-1}}-r\vert=1$
and thus $\tilde\sigma_{\tilde d}=t_1+1$.

If $d>\tilde d=1$ then
Lemma \ref{NodesRhoLem} implies that
$\rho(b_{\bari_1})=\tkappa_{\sigma_1}$ or
$\rho(b_{\bari_1})=p-\tkappa_{\sigma_1}$.
Since $r=\rho(b_{\bari_1})$, we then immediately obtain
$\tkappa_{\tilde\sigma_{\tilde d}}=\tkappa_{\sigma_{\tilde d}}$,
from which again the desired result follows.

If $d=\tilde d=1$ with $b$ interfacial (and thus $\rho(b)=r$),
a similar argument applies when $\sigma_1>0$.
The case $\sigma_1=0$ is excluded here since then either
$b_{0}=1$ or $b_{0}=p'-1$, neither of which is interfacial.

Now consider $\tilde d=d$ and $b$ being non-interfacial
(so that $p'<2p$ and $\sigma_d\le t_2$).
Here $r=\lfloor bp/p'\rfloor+\frac12(c-b+1)$.
If $\sigma_d>0$, Lemma \ref{NodesRhoLem} then yields:
\begin{equation*}
r=\rho(b_{i_1\cdots i_{d-1}})-\Delta_d\tkappa_{\sigma_d}
                    +\frac12(c-b-\Delta_d).
\end{equation*}
Then, since
$r=r_{i_1\cdots i_{d-1}}-\Delta_d\tkappa_{\tilde\sigma_d}
  =\rho(b_{i_1\cdots i_{d-1}})-\Delta_d\tkappa_{\tilde\sigma_d}$,
we obtain
$\Delta_d\tkappa_{\tilde\sigma_d}=\Delta_d\tkappa_{\sigma_d}-
 \frac12(c-b-\Delta_d)$, from which the required expressions
for $\tilde\sigma_d$ follow after noting that $\sigma_d\le t_2$.

If $\sigma_d=0$ then $b_{i_1\cdots i_{d-1}}=b+\Delta_d$,
and we claim that $c=b-\Delta_d$.
Otherwise $b_{i_1\cdots i_{d-1}}=b-\Delta_d=c$ which is interfacial
and thus $\rho(b_{i_1\cdots i_{d-1}})=r$, contradicting the fact that
$r_{i_1\cdots i_{d-1}}=\rho(b_{i_1\cdots i_{d-1}})$ is not a leaf-node.
It now follows that $|r_{i_1\cdots i_{d-1}}-r|=1$
and so $\tilde\sigma_{d}=t_1+1=1$ in this case, as required.
\cqfd
\medskip

\begin{theorem}\label{CombinedThrm}
For $1\le r<p$, let $b,c\in\{\lfloor rp'/p\rfloor,\lfloor rp'/p\rfloor+1\}$
with $b\ne c$.

1) If $b$ is interfacial then there is a
bijection between ${\mathcal U}(b)\backslash\overline{\mathcal U}(b)$
and $\widetilde{\mathcal U}(r)$ such that if $\boldu$ maps to
$\tilde{\boldu}$ under this bijection, then:
\begin{equation}\label{BijectiveEq1}
F^*(\boldu^L,\tilde{\boldu})=F^*(\boldu^L,\boldu),
\end{equation}
for all $\boldu^L\in{\mathcal U}(a)$ with $1\le a<p'$.

2) If $b$ is not interfacial then there is a
bijection between ${\mathcal U}(b)\backslash\overline{\mathcal U}{}^{c-b}(b)$
and $\widetilde{\mathcal U}(r)$ such that if $\boldu$ maps
to $\tilde{\boldu}$ under this bijection, then:
\begin{equation}\label{BijectiveEq2}
F^*(\boldu^L,\tilde{\boldu})
=\begin{cases}
  F^*(\boldu^L,\boldu)&    \text{if $\boldu\in{\mathcal U}^{c-b}(b)$;}\\
  F^*(\boldu^L,\boldu^+)&  \text{if $\boldu\in{\mathcal U}^{b-c}(b)$,}
\end{cases}
\end{equation}
for all $\boldu^L\in{\mathcal U}(a)$ with $1\le a<p'$.
\end{theorem}

\Proof
We will describe a traversal of the $b$-tree level by level (breadth first)
starting at level one. During the traversal, the $r$-tree is constructed
level by level, and certain nodes of the $b$-tree are eliminated
from further consideration.

At the $k$th level ($k\ge1$)
this process will ensure that we only examine non-leaf $b$-tree
nodes $b_{i_1\cdots i_k}$ for which $r_{i_1\cdots i_k}$ is a node of
the $r$-tree and for which $\rho(b_{i_1\cdots i_k})\ne r$.
Consider such a node $b_{i_1\cdots i_k}$. This node is either a
through-node, a branch-node for which $\rho(b_{i_1\cdots i_kh})=r$
for some $h\in\{0,1\}$, or a branch-node for which
$\rho(b_{i_1\cdots i_k0})<r<\rho(b_{i_1\cdots i_k1})$.
We consider these three cases separately below.
However, we will also deal simultaneously with the analogous
three cases that arise when $k=0$: the root node being a through-node,
$\rho(b_h)=r$ for some $h\in\{0,1\}$, or $\rho(b_0)<r<\rho(b_1)$.

In the case that $k=0$ and the root node is a through-node, or
$k>0$ and $b_{i_1\cdots i_k}$ is a through-node, set $d=k+1$ so
that $b_{i_1\cdots i_{d-1}0}$ is a leaf-node of the $b$-tree and,
by Lemma \ref{DDLem}, $r_{i_1\cdots i_{d-1}0}$ is a leaf-node
of the $r$-tree.
Let $\{\tau_j,\sigma_j,\Delta_j\}_{j=1}^d$ and
$\{\tilde\tau_j,\tilde\sigma_j,\tilde\Delta_j\}_{j=1}^d$
be the respective corresponding naive runs, and let
$\boldu\in{\mathcal U}(b)$ and $\tilde{\boldu}\in\widetilde{\mathcal U}(r)$
be the respective corresponding vectors.
Lemma \ref{DDLem} shows that $\boldu\not\in\overline{\mathcal U}(b)$
if $b$ is interfacial, and
$\boldu\not\in\overline{\mathcal U}{}^{c-b}(b)$
if $b$ is not interfacial.
For the required bijection, we map $\boldu\mapsto\tilde{\boldu}$
with (\ref{BijectiveEq1}) and (\ref{BijectiveEq2}) yet to be demonstrated.

Lemma \ref{GGLem} shows that $\tau_j=\tilde\tau_j$ for $1\le j\le d$
and $\sigma_j=\tilde\sigma_j$ for $1\le j<d$.
In the case that $b$ is interfacial, Lemma \ref{GGLem} gives
$\tilde\sigma_d=\max\{\sigma_d,t_1+1\}$.
If $\sigma_d\ge t_1+1$ then $\tilde{\boldu}=\boldu$ and
(\ref{BijectiveEq1}) holds trivially.
Otherwise, (\ref{BijectiveEq1}) follows by using
Lemma \ref{ShiftSigmaLem} (possibly more than once).
In the case that $b$ is not interfacial,
note first that $\boldu\in{\mathcal U}^{\Delta_d}(b)$.
For $\Delta_d=c-b$, Lemma \ref{GGLem} shows that
$\tilde\sigma_d=\sigma_d$ whence $\tilde{\boldu}=\boldu$ and
(\ref{BijectiveEq2}) holds trivially.
For $\Delta_d=b-c$, Lemma \ref{GGLem} shows that
$\tilde\sigma_d=\sigma_d+1$ whence $\tilde{\boldu}=\boldu^+$ and
(\ref{BijectiveEq2}) holds trivially.
We exclude examining the leaf-node $b_{i_1\cdots i_{d-1}0}$
in our subsequent traversing of the $b$-tree.

We now deal with the cases in which either
both $k=0$ and the root node is a branch-node,
or both $k>0$ and $b_{i_1\cdots i_k}$ is a branch-node.
If, in this case,
$\rho(b_{i_1\cdots i_kh})=r$ for some $h\in\{0,1\}$,
Lemma \ref{EELem}(1,2) shows that $r_{i_1\cdots i_{\tilde d-1}0}$
is a leaf-node of the $r$-tree where $\tilde d=k+1$.
Let $\{\tilde\tau_j,\tilde\sigma_j,\tilde\Delta_j\}_{j=1}^{\tilde d}$
be the corresponding naive run, and let
$\tilde{\boldu}\in\widetilde{\mathcal U}(r)$ be the corresponding vector.
Lemma \ref{EELem}(3) shows that there is a unique $d\ge k+2$ for which
$b_{i_1\cdots i_{d-1}0}$ is a leaf-node with
$i_{k+1}=\cdots=i_{d-1}=\barh$.
Let $\{\tau_j,\sigma_j,\Delta_j\}_{j=1}^d$ be the corresponding naive run,
and let $\boldu\in{\mathcal U}(b)$ be the corresponding vector.
Lemma \ref{EELem}(3) shows that $\boldu\not\in\overline{\mathcal U}(b)$
if $b$ is interfacial, and
$\boldu\in{\mathcal U}{}^{b-c}(b)$ if $b$ is not interfacial.
For the required bijection, we map $\boldu\mapsto\tilde{\boldu}$
with (\ref{BijectiveEq1}) and (\ref{BijectiveEq2}) yet to be demonstrated.

Lemma \ref{GGLem} shows that $\tau_j=\tilde\tau_j$ for $1\le j\le\tilde d$
and $\sigma_j=\tilde\sigma_j$ for $1\le j<\tilde d$ with
$\tilde\sigma_{\tilde d}=\max\{\sigma_{\tilde d},t_1+1\}$,
on noting that $d>k+1=\tilde d$.
Then $\tilde\sigma_{\tilde d}=\sigma_{\tilde d}$ because otherwise
$d>\tilde d$ implies that $\tau_d<\sigma_{\tilde d}\le t_1$ which
cannot occur in a naive run.
Lemma \ref{EELem}(3) also shows that
$\tau_j=\sigma_j$ for $\tilde d+1\le j<d$, and
$\sigma_d+1=\tau_d\le t_1+1$ for $p'>2p$ and
$\sigma_d+1=\tau_d\le t_2+1$ for $p'<2p$.
Therefore, $\boldu=\tilde{\boldu}+\boldu_{\sigma_d,\sigma_d+1}$
so that $\boldu^+=\tilde{\boldu}$.
In the case of $b$ being non-interfacial,
(\ref{BijectiveEq2}) follows immediately.
In the case of $b$ being interfacial and $\sigma_d\le t_1$,
the use of Lemma \ref{ShiftSigmaLem} proves (\ref{BijectiveEq1}).
We now claim that the case of $b$ being interfacial and $\sigma_d>t_1$
does not arise here. First note that necessarily $p'<2p$ so that $t_1=0$.
Lemma \ref{EELem}(3) implies that $b_{i_1\cdots i_{d-2}h}=b-\Delta_d$
so that both $b$ and $b-\Delta_d$ are interfacial, and separated
by an odd band. This implies that $t_2=1$
(by considering the model obtained by toggling the parity of each band).
So $\sigma_d=t_2=1$ and $\tau_d=2$.
Now $b_{i_1\cdots i_{d-2}\barh}=b+\Delta_d\kappa_{\sigma_d}=
b+2\Delta_d$ is interfacial.
The definition of a naive run implies that either $\tau_d=t_k$ for some $k$,
or $\tau_d=t_n-1$.
Thus either $t_3=2$ or both $t_3=3$ and $n=3$.
Lemma \ref{AnnoyingLem} now confirms the above claim.

Lemma \ref{EELem}(3) also shows that all leaf-nodes other than
$b_{i_1\cdots i_{d-1}0}$ that are descendents of $b_{i_1\cdots i_k}$
yield vectors $\boldu\in\overline{\mathcal U}(b)$ when $b$ is interfacial,
and vectors $\boldu\in\overline{\mathcal U}{}^{c-b}(b)$ when $b$ is
non-interfacial.
Thus we exclude examining all nodes that are descendents of
$b_{i_1\cdots i_k}$ in our subsequent traversing of the $b$-tree.

In the case in which either
both $k=0$ and the root node is a branch-node,
or both $k>0$ and $b_{i_1\cdots i_k}$ is a branch-node,
we consider
$\rho(b_{i_1\cdots i_k0})<r<\rho(b_{i_1\cdots i_k1})$,
Lemma \ref{AALem} shows that both
$r_{i_1\cdots i_k0}$ and $r_{i_1\cdots i_k1}$ are non-leaf-nodes of the
$r$-tree with the value of neither equal to $r$.
Each will be examined in traversing the next level of the $b$-tree.

Once this recursive procedure has finished, the $r$-tree will have been
completely constructed from the $b$-tree.
In addition, in the case of $b$ being interfacial, an explicit bijection
will have been established
between the leaf-nodes of the $b$-tree which yield
vectors $\boldu\in{\mathcal U}(b)\backslash\overline{\mathcal U}(b)$
and the leaf-nodes of the $r$-tree.
Since the corresponding bijection between
$\boldu\in{\mathcal U}(b)\backslash\overline{\mathcal U}(b)$ and
$\tilde{\boldu}\in\widetilde{\mathcal U}(r)$ satisfies (\ref{BijectiveEq1})
in each case, the theorem is proved for $b$ interfacial.
Similarly, in the case of $b$ being non-interfacial, an explicit bijection
will have been established
between the leaf-nodes of the $b$-tree which yield
vectors $\boldu\in{\mathcal U}(b)\backslash\overline{\mathcal U}{}^{c-b}(b)$
and the leaf-nodes of the $r$-tree.
Since the corresponding bijection between
$\boldu\in{\mathcal U}(b)\backslash\overline{\mathcal U}{}^{c-b}(b)$ and
$\tilde{\boldu}\in\widetilde{\mathcal U}(r)$ satisfies (\ref{BijectiveEq2})
in each case, the theorem is proved for $b$ non-interfacial.
\cqfd
\medskip

We are now in a position to prove the fermionic expressions
for all Virasoro characters $\chi^{p,p'}_{r,s}$ that were stated
in (\ref{FermCEq}). For convenience, the cases of the extra term
$\ochi^{\hat p,\hat p'}_{\hat r,\hat s}$ being present or not
are stated separately.

\begin{theorem}\label{FinalThrm}
Let $1\le s<p'$ and $1\le r<p$.
Let $\eta$ be such that $\xi_{\eta}\le s<\xi_{\eta+1}$ and $\eta^\prime$
be such that $\tilde\xi_{\eta^\prime}\le r<\tilde\xi_{\eta^\prime+1}$.

\noindent
1. If $\eta=\eta^\prime$ and
      $\xi_{\eta}<s$ and $\tilde\xi_{\eta}<r$ then:
\begin{equation}\label{Final1Eq}
\chi^{p,p'}_{r,s}
= \sum_{\begin{subarray}{c}
           \boldu^L\in{\mathcal U}(s)\\
           \boldu^R\in\widetilde{\mathcal U}(r)
        \end{subarray}}
F(\boldu^L,\boldu^R)
\quad+\quad\ochi^{\hat p,\hat p'}_{\hat r,\hat s},
\end{equation}
where $\hat p'=\xi_{\eta+1}-\xi_{\eta}$,
$\hat p=\tilde\xi_{\eta+1}-\tilde\xi_{\eta}$,
$\hat s=s-\xi_{\eta}$, $\hat r=r-\tilde\xi_{\eta}$.

\noindent
2. Otherwise:
\begin{equation}\label{Final2Eq}
\chi^{p,p'}_{r,s}
= \sum_{\begin{subarray}{c}
           \boldu^L\in{\mathcal U}(s)\\
           \boldu^R\in\widetilde{\mathcal U}(r)
        \end{subarray}}
F(\boldu^L,\boldu^R).
\end{equation}
\end{theorem}

\Proof
First note that if $\boldu^R\in\widetilde{\mathcal U}(r)$
then $\sigma(\boldu^R)\ge t_1+1$.
In the notation used for (\ref{ChangeEq2}), this guarantees that
$\boldN(\boldu^R)=0$ and 
$(\boldu^L_\flat+\boldu^R_\sharp)\cdot{\boldm}
=(\oboldu^L_\flat+\oboldu^R_\sharp)\cdot{\oboldm}$.
Thereupon $F^*(\boldu^L,\boldu^R)=F(\boldu^L,\boldu^R)$.

Let $b$ be equal to one of $\lfloor rp'/p\rfloor$ and
$\lfloor rp'/p\rfloor+1$, and let $c$ be equal to the other.
If $b$ is interfacial, the application of Theorem \ref{CombinedThrm}(1)
to the two cases of Theorem \ref{Limit1Thrm} now yields
(\ref{Final1Eq}) and (\ref{Final2Eq}).
Similarly, if $b$ is non-interfacial, the application of
Theorem \ref{CombinedThrm}(2) to the two cases of Theorem \ref{Limit2Thrm}
also yields (\ref{Final1Eq}) and (\ref{Final2Eq}).
\cqfd

\newpage

\setcounter{section}{9}

\section{Discussion}\label{DissSec}

In this paper, we have shown how to write down constant-sign fermionic
expressions for all Virasoro characters $\chi^{p,p'}_{r,s}$,
where $1<p<p'$ with $p$ and $p'$ coprime, $1\le r<p$ and $1\le s<p'$.
The first step in this process constructs a tree for $s$ from
the set $\T\cup\T'$ where $\T$ is the set of Takahashi lengths associated
with the continued fraction of $p'/p$, and $\T'$ is a set of values
complementary to $\T$.
A tree for $r$ is constructed in a similar way using the set
$\tT\cup\tT'$, where $\tT$ is the set of truncated
Takahashi lengths associated with the continued fraction of $p'/p$,
and $\tT'$ is a set of values complementary to $\tT$.
A set $\U(s)$ of vectors is obtained from the leaf-nodes of the first
of these trees, with each leaf-node giving rise to precisely one vector.
Similarly, a set $\tU(r)$ of vectors is obtained from the leaf-nodes
of the second of these trees.
Each pair of vectors $\boldu^L\in\U(s)$ and $\boldu^R\in\tU(r)$
gives rise to a fundamental fermionic form $F(\boldu^L,\boldu^R)$.
In most cases, taking the sum of such terms over all pairs
$\boldu^L,\boldu^R$ yields the required fermionic expression for
$\chi^{p,p'}_{r,s}$. This expression thus comprises
$|\U(s)|\cdot|\tU(r)|$ fundamental fermionic forms.
In the remaining cases, in addition to these terms, we require
a further character $\chi^{\hat p,\hat p'}_{\hat r,\hat s}$.
The process described above is now applied to this character.
Thus in general, writing down the
fermionic expression for $\chi^{p,p'}_{r,s}$ is a recursive procedure. 
It is interesting to note that the number of fundamental
fermionic forms comprising the resulting fermionic expressions for
$\chi^{p,p'}_{r,s}$, varies erratically as $r$ and $s$ run over
their permitted ranges. A similar observation was made in \cite{bms}
for the expressions given there.

Expressions of a similar nature are provided for the finitized
characters $\chi^{p,p'}_{a,b,c}(L)$, where $1\le p<p'$ with $p$ and $p'$
coprime, $1\le a,b<p'$ and $c=b\pm1$. These finitized characters
are generating functions for length $L$ Forrester-Baxter paths.
The fermionic expressions here involve a sum over all pairs
$\boldu^L,\boldu^R$ with $\boldu^L\in\U(a)$ and $\boldu^R\in\U(b)$:
the sets $\U(a)$ and $\U(b)$ having both been produced using the
set $\T\cup\T'$.
If $b$ satisfies (\ref{interface}), each term in the sum is the
fundamental fermionic form $F(\boldu^L,\boldu^R,L)$.
For other values of $b$, each term in the sum is either
$F(\boldu^L,\boldu^R,L)$ or $\widetilde F(\boldu^L,\boldu^R,L)$,
where the latter is itself a linear combination of terms
$F(\boldu^L,\boldu,L')$ for various $\boldu$ and $L'$.
The fermionic expression for $\chi^{p,p'}_{a,b,c}(L)$ also sometimes
involves a term $\chi^{\hat p,\hat p'}_{\hat a,\hat b,\hat c}(L)$
whence, like the fermionic expression for $\chi^{p,p'}_{r,s}$,
it is recursive in nature.
For $b$ satisfying (\ref{interface}), the expressions obtained for
$\chi^{p,p'}_{a,b,c}(L)$ are genuinely fermionic in that they are
positive sums over manifestly positive definite $F(\boldu^L,\boldu^R,L)$.
For $b$ not satisfying (\ref{interface}), this is not the case.
Nonetheless, whatever the value of $b$, after some work, the fermionic
expressions for $\chi^{p,p'}_{r,s}$ described above emerge from these
expressions on taking the $L\to\infty$ limit.

By using the Takahashi trees that give rise to the sets $\U(a)$ and $\U(b)$,
it is possible to characterise the set of paths for which
$F(\boldu^L,\boldu^R,L)$ or $\widetilde F(\boldu^L,\boldu^R,L)$ is
the generating function. In the former case, this set of paths is
specified in Lemma \ref{OneNodeLem} as those which attain certain
heights in a certain order. It is a similar set of paths
in the latter case, as described in Note \ref{FunnyNote}.
The corresponding term $F(\boldu^L,\boldu^R)$ that forms a component of
the fermionic expression for $\chi^{p,p'}_{r,s}$ may then be seen to be
the generating function of a certain set of paths of infinite length.
However, it is somewhat unfortunate that characterising these latter
paths requires the construction of the Takahashi tree for a suitable $b$,
whereas in forming the fermionic expression for $\chi^{p,p'}_{r,s}$,
the construction of such a tree has been superseded by the construction
of the truncated Takahashi tree for $r$.

In fact, this and many aspects of the proof of the fermionic
expressions for $\chi^{p,p'}_{r,s}$ hint that perhaps
$\chi^{p,p'}_{a,b,c}(L)$ is not the most natural finitization
of $\chi^{p,p'}_{r,s}$. These aspects include the need to temporarily 
redefine the weight function to that of (\ref{WtDef2}), the frequent need
to treat the cases $p'>2p$ and $p'<2p$ separately, the appearance of
$\delta^{p,p'}_{a,e}$ and $\delta^{p'-p,p'}_{a,e}$
in Corollaries \ref{MazyBijCor} and \ref{MazyDijCor},
the appearance of terms $F(\boldu^L,\boldu^R,L)$ in the fermionic
expression for $\chi^{p,p'}_{a,b,c}(L)$ that vanish as $L\to\infinity$,
and the difference between the definitions of $\gamma(\run^L,\run^R)$
and $\gamma'(\run^L,\run^R)$ in Section \ref{ConSec}.
A more natural finitization may lead to a better combinatorial description
of the characters $\chi^{p,p'}_{r,s}$. This, in turn, would be of
great benefit both in the analysis of the expressions given in this
paper as well as in the construction of the corresponding irreducible
representations of the Virasoro algebra.

It is also intriguing that the constant term $\frac14\gamma(\run^L,\run^R)$
varies relatively little on running through the leaf-nodes of the
Takahashi tree for $s$ and the leaf-nodes of the truncated Takahashi tree
for $r$.
Moreover, these values are also close to the negative of the modular anomaly
$\Delta^{p,p'}_{r,s}-\frac1{24}c^{p,p'}$.
This hints at the fundamental fermionic forms $F(\boldu^L,\boldu^R)$
themselves having interesting properties under the action of the
modular group.

Equating the fermionic expression for each $\chi^{p,p'}_{r,s}$
with the corresponding bosonic expression given by (\ref{RochaEq})
yields bosonic-fermionic $q$-series identities.
Furthermore, (\ref{ProdEq1}) or (\ref{ProdEq3}) provides a product
expression for the character $\chi^{p,p'}_{r,s}$ in certain cases.
Equating each of these with the corresponding fermionic
expression given here yields identities of the Rogers-Ramanujan type.
In addition, equating the fermionic expression for each
finitized character $\chi^{p,p'}_{a,b,c}(L)$ with the
corresponding bosonic expression (\ref{FinRochaEq}) yields
bosonic-fermionic polynomial identities.

Finally, we remark that further fermionic expressions and identities
for $\chi^{p,p'}_{r,s}$ may be produced by
instead of using the set $\T\cup\T'$ to obtain the Takahashi tree
for $s$ and the set $\tT\cup\tT'$ to obtain the truncated Takahashi tree
for $r$, using either just the set $\T\cup\{\kappa_t\}$ or the
set $\T'\cup\{p'-\kappa_t\}$ in the former case, or using just
the set $\tT\cup\{\tkappa_t\}$ or the set $\tT'\cup\{p-\tkappa_t\}$
in the latter case.
However, such expressions for $\chi^{p,p'}_{r,s}$
are not as efficient as those given in Section \ref{PrologueSec},
in that they involve as least as many terms
$F(\boldu^L,\boldu^R)$ and sometimes more.
In addition, an extra term $\chi^{\hat p,\hat p'}_{\hat r,\hat s}$
might be present when the expression of Section \ref{PrologueSec}
requires no such term.
Similar comments pertain to the finitized characters
$\chi^{p,p'}_{a,b,c}(L)$.
The definitions of Section \ref{ExtraSec} must be suitably modified to
accommodate these additional expressions.

Looking ahead, the question now arises as to whether the techniques
developed in this paper can be extended to deal with the characters of
the $W_n$ algebras, the $n=2$ case of which is the Virasoro algebra.
The RSOS models that pertain to these algebras are described
in \cite{jmo,nakanishi,kns}. As yet, fermionic expressions for the
$W_n$ algebra characters are available only in very special cases
\cite{fow,schilling,dasmahapatra,warnaar-pearce,
schilling-shimozono,asw,warnaar3}.
In particular, the expressions of \cite{asw,warnaar3} have led to
a number of elegant new identities of Rogers-Ramanujan type.

\addtocontents{toc}{\SkipTocEntry}
\subsection*{Acknowledgement} I wish to thank Omar Foda very much for his
support and encouragement during this work.
I also wish to thank Omar, Alexander Berkovich, Barry McCoy, Will Orrick,
Anne Schilling and Ole Warnaar for useful discussions.

\newpage

\begin{appendix}

\section{Examples}\label{ExampleSec}

Here, we provide three examples of the constructions of the fermionic
expressions that are described in Section \ref{PrologueSec}.

\addtocontents{toc}{\SkipTocEntry}
\subsection{Example 1}

\begin{figure}[ht]
\includegraphics[scale=0.75]{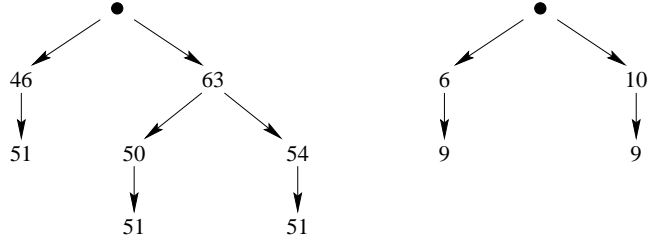}
\caption{Takahashi tree for $s=51$ and truncated Takahashi tree for $r=9$
         when $p'=109$ and $p=26$.}
\label{FormationEx1aFig}
\medskip
\end{figure}

Let $p'=109$, $p=26$, $r=9$ and $s=51$.
Then $p'/p$ has continued fraction $[4,5,5]$, so that $n=2$,
$\{t_k\}_{k=0}^3=\{-1,3,8,13\}$, $t=12$,
$\{y_k\}_{k=-1}^3=\{0,1,4,21,109\}$ and
$\{z_k\}_{k=-1}^3=\{1,0,1,5,26\}$.
The sets $\T$, $\T'$, $\tT$ and $\tT'$ of Takahashi lengths,
complementary Takahashi lengths, truncated Takahashi lengths
and complementary Takahashi lengths are:
\begin{align*}
\T=\{\kappa_i\}_{i=0}^{11}&=\{1,2,3,4,5,9,13,17,21,25,46,67\};\\
\T'=\{p'-\kappa_i\}_{i=0}^{11}
         &=\{108,107,106,105,104,100,96,92,88,84,63,42\};\\
\tT=\{\tkappa_i\}_{i=4}^{11}&=\{1,2,3,4,5,6,11,16\};\\
\tT'=\{p-\kappa_i\}_{i=4}^{11}&=\{25,24,23,22,21,20,15,10\}.
\end{align*}
Using these, we obtain the Takahashi tree for $s=51$ and the truncated
Takahashi tree for $r=9$ shown in Fig.~\ref{FormationEx1aFig}.

Since these trees have three and two leaf-nodes respectively,
$|\U(s)|=3$ and $|\tU(r)|=2$.
Thus for the current example, the summation in (\ref{FermCEq})
will be over six terms $F(\boldu^L,\boldu^R)$.
To determine whether the additional term 
$\chi^{\hat p,\hat p'}_{\hat r,\hat s}$ is present,
we obtain following the prescription of Section \ref{ExtraSec},
\begin{align*}
\{\xi_\ell\}_{\ell=0}^9&=\{0,21,25,42,46,63,67,84,88,109\};\\
\{\tilde\xi_\ell\}_{\ell=0}^9&=\{0,5,6,10,11,15,16,20,21,26\}.
\end{align*}
Then, since $\eta(s)=4\ne\tilde\eta(r)=2$,
that additional term is not present.

As described in Section \ref{FormVSec} each leaf-node of the
Takahashi tree for $s$ gives rise to a run $\run$.
For the three leaf-nodes present when $s=51$, we denote
the corresponding runs by $\run^L_1,\run^L_2,\run^L_3$.
Similarly, as described in Section \ref{FormVSec}, each leaf
node of the truncated Takahashi tree for $r$ gives rise to
a run $\run$: for the two leaf-nodes present when $r=9$,
we denote the corresponding runs by $\run^R_1,\run^R_2$.
These sets are given by:
\begin{align*}
\run^L_1&= \{ \{13,7\},\{10,4\},\{1,-1\} \}
&\run^R_1&= \{ \{13,7\},\{11,6\},\{1,-1\} \},\\
\run^L_2&= \{ \{13,7,3\},\{10,5,0\},\{-1,1,-1\} \}
&\run^R_2&= \{ \{13,7\},\{9,4\},\{-1,1\} \},\\
\run^L_3&= \{ \{13,7,3\},\{10,6,2\},\{-1,1,1\} \}.
\end{align*}

{}From these, using (\ref{uEq}), we obtain
$\boldu^L_i=\boldu(\run^L_i)$ for $i=1,2,3$,
and $\boldu^R_j=\boldu(\run^R_j)$ for $j=1,2$:
\begin{align*}
\boldu^L_1&= (0, 0, 0, 1, 0, 0, -1, 0, 0, 1, 0, 1),
&\boldu^R_1&= (0, 0, 0, 0, 0, 1, -1, 0, 0, 0, 1, 1),\\
\boldu^L_2&= (0, 0, -1, 0, 1, 0, -1, 0, 0, 1, 0, 0),
&\boldu^R_2&= (0, 0, 0, 1, 0, 0, -1, 0, 1, 0, 0, 0),\\
\boldu^L_3&= (0, 1, -1, 0, 0, 1, -1, 0, 0, 1, 0, 0).
\end{align*}
So $\U(s)=\{\boldu^L_1,\boldu^L_2,\boldu^L_3\}$
and $\tU(r)=\{\boldu^R_1,\boldu^R_2\}$.

Setting $\boldDelta^L_i=\boldDelta(\run^L_i)$ and
$\boldDelta^R_j=\boldDelta(\run^R_j)$, we also obtain:
\begin{align*}
\boldDelta^L_1&= (0, 0, 0, -1, 0, 0, 1, 0, 0, 1, 0, -1),
&\boldDelta^R_1&= (0, 0, 0, 0, 0, -1, 1, 0, 0, 0, 1, -1),\\
\boldDelta^L_2&= (0, 0, 1, 0, 1, 0, -1, 0, 0, -1, 0, 0),
&\boldDelta^R_2&= (0, 0, 0, 1, 0, 0, -1, 0, -1, 0, 0, 0),\\
\boldDelta^L_3&= (0, 1, -1, 0, 0, 1, -1, 0, 0, -1, 0, 0).
\end{align*}
For $i\in\{1,2,3\}$ and $j\in\{1,2\}$, let
$\boldm^{(1)}_{ij}=(\oboldu^L_{i\flat}+\oboldu^R_{j\sharp})\cdot\boldm$
for $\boldm=(m_4,m_5,\ldots,m_{11})$,
and
$\gamma_{ij}=\gamma(\run^L_i,\run^R_j)$.
These are respectively minus twice the linear term and four times the
constant term that appears in the exponent in (\ref{FEq}).
Via the prescriptions of Section \ref{LinSec} and \ref{ConSec},
we obtain:
\begin{align*}
&\boldm^{(1)}_{11}=m_4-m_7+m_{11},
&&\gamma_{11}=-42;
\\
&\boldm^{(1)}_{21}=m_5-m_7+m_{11},
&&\gamma_{21}=-42;
\\
&\boldm^{(1)}_{31}=m_6-m_7+m_{11},
&&\gamma_{31}=-43;
\\
&\boldm^{(1)}_{12}=m_4-m_7+m_9,
&&\gamma_{12}=-42;
\\
&\boldm^{(1)}_{22}=m_5-m_7+m_9,
&&\gamma_{22}=-42;
\\
&\boldm^{(1)}_{32}=m_6-m_7+m_9,
&&\gamma_{32}=-43.
\end{align*}

Let $\boldn^{(2)}=\tilde{\boldn}^T\boldB\tilde{\boldn}
+\frac14\boldm^T\boldC\boldm$
with $\tilde{\boldn}=(\tilde n_1,\tilde n_2,\tilde n_3)$
and $\boldm=(m_4,m_5,\ldots,m_{11})$.
Then, via the descriptions of $\boldCC$ and $\boldB$ given in
Section \ref{QuadSec}, we obtain:
\begin{equation*}
\begin{split}
\boldn^{(2)}
&=(\tilde n_1+\tilde n_2+\tilde n_3)^2
 +(\tilde n_2+\tilde n_3)^2+\tilde n_3^2 \\
&\hskip5mm +\genfrac{}{}{}114\left(
m_4^2 +(m_4-m_5)^2 +(m_5-m_6)^2 +(m_6-m_7)^2 +(m_7-m_8)^2\right.\\
&\hskip15mm
\left. +m_9^2 +(m_9-m_{10})^2 +(m_{10}-m_{11})^2 +m_{11}^2\right).
\end{split}
\end{equation*}

Now let $\boldu_{ij}=\boldu^L_i+\boldu^R_j$ for $i\in\{1,2,3\}$
and $j\in\{1,2\}$.
These together with the values of $\boldQQ_{ij}=\boldQQ(\boldu_{ij})$
obtained using $\boldC^*$ or $\boldCC^*$ 
(or even via (\ref{MNEq1}) and (\ref{MNEq2}))
as described in Section \ref{MNsysSec} are:
\begin{align*}
&\oboldu_{11}=(0, 1, -2, 0, 0, 1, 1, 2),
&&\boldQQ_{11}=(1, 1, 1, 0, 1, 1, 1, 0);
\\
&\oboldu_{21}=(1, 1, -2, 0, 0, 1, 1, 1),
&&\boldQQ_{21}=(0, 0, 1, 1, 1, 0, 1, 1);
\\
&\oboldu_{31}=(0, 2, -2, 0, 0, 1, 1, 1),
&&\boldQQ_{31}=(1, 1, 1, 1, 1, 0, 1, 1);
\\
&\oboldu_{12}=(0, 0, -2, 0, 1, 1, 0, 1),
&&\boldQQ_{12}=(1, 1, 1, 1, 1, 0, 0, 1);
\\
&\oboldu_{22}=(1, 0, -2, 0, 1, 1, 0, 0),
&&\boldQQ_{22}=(0, 0, 1, 0, 1, 1, 0, 0);
\\
&\oboldu_{32}=(0, 1, -2, 0, 1, 1, 0, 0),
&&\boldQQ_{32}=(1, 1, 1, 0, 1, 1, 0, 0).
\end{align*}
In this case, (\ref{FermCEq}) yields:
\begin{equation*}
\chi^{26,109}_{9,51}=
\sum_{i=1}^3\sum_{j=1}^2 F_{ij},
\end{equation*}
where
\begin{equation*}
\begin{split}
F_{ij}
=\sum_{\begin{subarray}{c}
           n_1,n_2,n_3\\ \boldm\equiv\boldQQ_{ij}
        \end{subarray}}
&\frac{q^{\boldn^{(2)}-\frac12\boldm^{(1)}_{ij}+\frac14\gamma_{ij}}}
{(q)_{n_1}(q)_{n_2}(q)_{n_3}(q)_{m_4}}
\qbinom{\frac12(m_4+m_6+(\boldu_{ij})_5)}{m_5}\\
&\times\qbinom{\frac12(m_5+m_7+(\boldu_{ij})_6)}{m_6}
\qbinom{\frac12(m_6+m_8+(\boldu_{ij})_7)}{m_7}\\
&\times\qbinom{\frac12(m_7+m_8-m_9+(\boldu_{ij})_8)}{m_8}
\qbinom{\frac12(m_8+m_{10}+(\boldu_{ij})_9)}{m_9}\\
&\times\qbinom{\frac12(m_9+m_{11}+(\boldu_{ij})_{10})}{m_{10}}
\qbinom{\frac12(m_{10}+(\boldu_{ij})_{11})}{m_{11}},
\end{split}
\end{equation*}
where $\boldm=(m_4,m_5,\ldots,m_{11})$ and
$\tilde n_1=n_1-\frac12(\boldu_{ij})_1$,
$\tilde n_2=n_2-\frac12(\boldu_{ij})_2$,
$\tilde n_3=n_3-\frac12(\boldu_{ij})_3+\frac12m_4$.
The summation here extends over an infinite number of terms because
the values that $n_1,n_2,n_3,m_4$ can take are unbounded.

We now seek a fermionic expression for a finitization
$\chi^{26,109}_{a,b,c}(L)$ of $\chi^{26,109}_{9,51}$.
Certainly, we require $a=51$, but a number of values of
$b$ and $c$ satisfy (\ref{groundstatelabel}).
The simplest fermionic expressions (\ref{FermEq}) are obtained when
$b$ is interfacial.
In this case, either $b=37$ or $b=38$ is interfacial,
whereupon either $c=b\pm1$ satisfies (\ref{groundstatelabel}). 
Here we select $b=38$ and $c=39$.
Thus we use (\ref{FermEq}) to obtain a fermionic expression for
$\chi^{26,109}_{51,38,39}(L)$.

\begin{figure}[ht]
\includegraphics[scale=0.75]{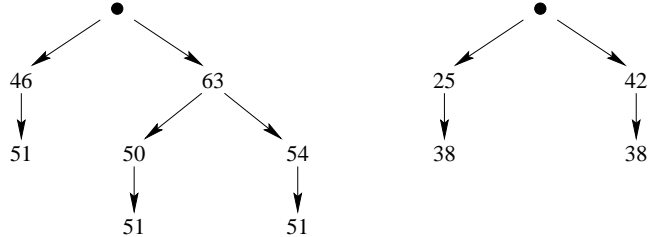}
\caption{Takahashi trees for $a=51$ and $b=38$ when $p'=109$ and $p=26$.}
\label{FormationEx1bFig}
\medskip
\end{figure}

The Takahashi trees for $a$ and $b$ are given
in Fig.\ \ref{FormationEx1bFig}.

The three leaf-nodes of the Takahashi tree for $a=51$ give rise
to the runs $\run^L_1,\run^L_2,\run^L_3$ which are, of course,
exactly as above.
The two leaf-nodes of the Takahashi tree for $b=38$ give rise
to runs $\run^R_1,\run^R_2$, which are {\em not} exactly as above.
\begin{align*}
\run^L_1&= \{ \{13,7\},\{10,4\},\{1,-1\} \}
&\run^R_1&= \{ \{13,7\},\{11,6\},\{1,-1\} \},\\
\run^L_2&= \{ \{13,7,3\},\{10,5,0\},\{-1,1,-1\} \}
&\run^R_2&= \{ \{13,7\},\{9,3\},\{-1,1\} \},\\
\run^L_3&= \{ \{13,7,3\},\{10,6,2\},\{-1,1,1\} \}.
\end{align*}

{}From these, using (\ref{uEq}), we obtain
$\boldu^L_i=\boldu(\run^L_i)$ for $i=1,2,3$,
and $\boldu^R_j=\boldu(\run^R_j)$ for $j=1,2$:
\begin{align*}
\boldu^L_1&= (0, 0, 0, 1, 0, 0, -1, 0, 0, 1, 0, 1),
&\boldu^R_1&= (0, 0, 0, 0, 0, 1, -1, 0, 0, 0, 1, 1),\\
\boldu^L_2&= (0, 0, -1, 0, 1, 0, -1, 0, 0, 1, 0, 0),
&\boldu^R_2&= (0, 0, 0, 0, 0, 0, -1, 0, 1, 0, 0, 0),\\
\boldu^L_3&= (0, 1, -1, 0, 0, 1, -1, 0, 0, 1, 0, 0).
\end{align*}
So $\U(a)=\{\boldu^L_1,\boldu^L_2,\boldu^L_3\}$
and $\U(b)=\{\boldu^R_1,\boldu^R_2\}$.

For $i\in\{1,2,3\}$ and $j\in\{1,2\}$, defining
$\boldm^{(1)}_{ij}=(\boldu^L_{i\flat}+\boldu^R_{j\sharp})\cdot\boldm$
for $\boldm=(L,m_1,m_2,\ldots,m_{11})$,
and
$\gamma'_{ij}=\gamma(\run^L_i,\run^R_j)$,
we obtain exactly the values of $\boldm^{(1)}_{ij}$ obtained above,
and $\gamma'_{ij}=\gamma_{ij}$ given above.

Let $\boldm^{(2)}=\hat{\boldm}^T\boldC\hat{\boldm}-L^2$.
Then, having calculated $\boldC$ as in Section \ref{QuadSec}, we obtain:
\begin{equation*}
\begin{split}
\boldm^{(2)}
&=(L-m_1)^2 +(m_1-m_2)^2 +(m_2-m_3)^2 +m_4^2 \\
&\hskip5mm +(m_4-m_5)^2 +(m_5-m_6)^2 +(m_6-m_7)^2 +(m_7-m_8)^2 +m_9^2\\
&\hskip5mm +(m_9-m_{10})^2 +(m_{10}-m_{11})^2 +m_{11}^2.
\end{split}
\end{equation*}
Then take $\boldu_{ij}=\boldu^L_i+\boldu^R_j$ for $i\in\{1,2,3\}$
and $j\in\{1,2\}$.
These together with the values of $\boldQ_{ij}=\boldQ(\boldu_{ij})$
obtained using $\boldC^*$ from Section \ref{QuadSec}
(or alternatively, obtained via (\ref{MNEq1}) and (\ref{MNEq2}))
are:
\begin{align*}
&\boldu_{11}=(0, 0, 0, 1, 0, 1, -2, 0, 0, 1, 1, 2),
&&\boldQ_{11}=(0, 1, 0, 1, 1, 1, 0, 1, 1, 1, 0);
\\
&\boldu_{21}=(0, 0, -1, 0, 1, 1, -2, 0, 0, 1, 1, 1),
&&\boldQ_{21}=(0, 1, 0, 0, 0, 1, 1, 1, 0, 1, 1);
\\
&\boldu_{31}=(0, 1, -1, 0, 0, 2, -2, 0, 0, 1, 1, 1),
&&\boldQ_{31}=(0, 1, 1, 1, 1, 1, 1, 1, 0, 1, 1);
\\
&\boldu_{12}=(0, 0, 0, 1, 0, 0, -2, 0, 1, 1, 0, 1),
&&\boldQ_{12}=(0, 1, 0, 1, 1, 1, 1, 1, 0, 0, 1);
\\
&\boldu_{22}=(0, 0, -1, 0, 1, 0, -2, 0, 1, 1, 0, 0),
&&\boldQ_{22}=(0, 1, 0, 0, 0, 1, 0, 1, 1, 0, 0);
\\
&\boldu_{32}=(0, 1, -1, 0, 0, 1, -2, 0, 1, 1, 0, 0),
&&\boldQ_{32}=(0, 1, 1, 1, 1, 1, 0, 1, 1, 0, 0).
\end{align*}
Putting all of this into (\ref{FermEq}) and (\ref{FEq}) produces:
\begin{equation*}
\chi^{26,109}_{51,38,37}(L)=
\chi^{26,109}_{51,38,39}(L)=
\sum_{i=1}^3\sum_{j=1}^2
F_{ij}(L),
\end{equation*}
where
\begin{equation*}
\begin{split}
F_{ij}(L)
&=\sum_{\boldm\equiv\boldQ_{ij}}
q^{\frac14\boldm^{(2)}-\frac12\boldm^{(1)}_{ij}+\frac14\gamma'_{ij}}
\qbinom{\frac12(L\!+\!m_2\!+\!(\boldu_{ij})_1)}{m_1}
\qbinom{\frac12(m_1\!+\!m_3\!+\!(\boldu_{ij})_2)}{m_2}\\
&\times\qbinom{\frac12(m_2\!+\!m_3\!-\!m_4\!+\!(\boldu_{ij})_3)}{m_3}
\qbinom{\frac12(m_3\!+\!m_5\!+\!(\boldu_{ij})_4)}{m_4}
\qbinom{\frac12(m_4\!+\!m_6\!+\!(\boldu_{ij})_5)}{m_5}\!\\
&\times\qbinom{\frac12(m_5\!+\!m_7\!+\!(\boldu_{ij})_6)}{m_6}
\qbinom{\frac12(m_6\!+\!m_8\!+\!(\boldu_{ij})_7)}{m_7}
\qbinom{\frac12(m_7\!+\!m_8\!-\!m_9\!+\!(\boldu_{ij})_8)}{m_8}\!\\
&\times\qbinom{\frac12(m_8\!+\!m_{10}\!+\!(\boldu_{ij})_9)}{m_9}
\qbinom{\frac12(m_9\!+\!m_{11}\!+\!(\boldu_{ij})_{10})}{m_{10}}
\qbinom{\frac12(m_{10}\!+\!(\boldu_{ij})_{11})}{m_{11}}.
\end{split}
\end{equation*}
Note that the summation here comprises a finite number of non-zero terms.
It might be better carried out by finding all solutions $\{n_i\}_{i=1}^{12}$
of (\ref{Particles2Eq}), where in this case,
$\{l_i\}_{i=1}^{12}=\{1, 2, 3, 1, 5, 9, 13, 17, 4, 25, 46, 67\}$,
and then obtaining  $\{m_i\}_{i=1}^{12}$ via (\ref{MNmatEq2}),
or via (\ref{MNEq1}) and (\ref{MNEq2}).
Note that the Gaussian polynomial terms are then best expressed in
the form $\qbinom{m_j+n_j}{n_j}$.

The $L\to\infty$ limit of the above expression yields
the Virasoro character $\chi^{26,109}_{9,51}$.
In fact, we have $\lim_{L\to\infty} F_{ij}(L)=F_{ij}$ for
$i=1,2,3$ and $j=1,2$.

\addtocontents{toc}{\SkipTocEntry}
\subsection{Example 2}

Let $p'=75$, $p=53$, $r=17$ and $s=72$.
Then $p'/p$ has continued fraction $[1,2,2,2,4]$, so that $n=4$,
$\{t_k\}_{k=0}^5=\{-1,0,2,4,6,10\}$, $t=9$,
$\{y_k\}_{k=-1}^5=\{0,1,1,3,7,17,75\}$ and
$\{z_k\}_{k=-1}^5=\{1,0,1,2,5,12,53\}$.
The set $\T$ of Takahashi lengths is then
$\T=\{\kappa_i\}_{i=0}^{8}=\{1,2,3,4,7,10,17,24,41\}$
and the set $\tT$ of truncated Takahashi lengths is
$\tT=\{\kappa_i\}_{i=1}^{8}=\{1,2,3,5,7,12,17,29\}$.
Since $p'-s=3\in\T$, the Takahashi tree for $s=72$ is trivial,
comprising just one node in addition to the root node.
This leaf-node gives rise to the run $\run^L=\{ \{10\},\{2\},\{1\} \}$.
Similarly, since $17\in\tT$, the truncated Takahashi tree for $r=17$
is trivial.
Its single leaf-node gives rise to the run $\run^R=\{ \{10\},\{7\},\{-1\} \}$.
{}From these, we obtain
$\boldu^L=\boldu(\run^L)$ and $\boldu^R=\boldu(\run^R)$ given by:
\begin{align*}
&\boldu^L= (0, 0, 0,-1, 0,-1, 0, 0, 1),
&&\boldu^R= (0, 0, 0, 0, 0, 0, 1, 0, 0);\\
\intertext{and $\boldDelta^L=\boldDelta(\run^L)$ and
$\boldDelta^R=\boldDelta(\run^R)$ given by:}
&\boldDelta^L= (0, 0, 0,-1, 0,-1, 0, 0, 1),
&&\boldDelta^R= (0, 0, 0, 0, 0, 0, -1, 0, 0).
\end{align*}
The procedure of Section \ref{ConSec} then yields
$\gamma(\run^L,\run^R)=-1624$.
The linear term is specified by
$\boldm^{(1)}=(\boldu^L_{\flat}+\boldu^R_{\sharp})\cdot\boldm$
for $\boldm=(m_1,m_2,\ldots,m_8)$, whereupon $\boldm^{(1)}=-m_6+m_7$.

Let $\boldn^{(2)}=\frac14\boldm^T\boldC\boldm$.
Then, via the descriptions of $\boldCC$ and $\boldB$ given in
Section \ref{QuadSec}, we obtain:
\begin{equation*}
\begin{split}
\boldn^{(2)}
&=\genfrac{}{}{}114\left( m_1^2 + (m_1-m_2)^2 + m_3^2 + (m_3-m_4)^2\right.\\
&\hskip20mm\left. + m_5^2 + (m_5-m_6)^2 + m_7^2 +(m_7-m_8)^2 + m_8^2\right).
\end{split}
\end{equation*}

Setting $\boldu=\boldu^L+\boldu^R$ and $\boldQQ=\boldQQ(\boldu)$
obtained using $\boldC^*$ or $\boldCC^*$ 
(or even via (\ref{MNEq1}) and (\ref{MNEq2})), gives:
\begin{equation*}
\oboldu=(0,0,-1,0,-1,1,0,1),\qquad
\boldQQ=(0,0,0,0,1,0,0,1).
\end{equation*}

We can now substitute all these values into (\ref{FermCEq}).
We note that the sum comprises just one term $F(\boldu^L,\boldu^R)$.
To determine whether the extra term is present, we calculate:
\begin{equation*}
\{\xi_\ell\}_{\ell=0}^9=\{0,17,24,34,41,51,58,75\};\qquad
\{\tilde\xi_\ell\}_{\ell=0}^9=\{0,12,17,24,29,36,41,53\}.
\end{equation*}
Then, since $\eta(s)=6\ne\tilde\eta(r)=2$,
that additional term is not present.
Thereupon:
\begin{equation*}
\begin{split}
\chi^{53,75}_{17,72}
=\sum_{\begin{subarray}{c}
           n_1,n_2\\ \boldm\equiv\boldQQ
        \end{subarray}}
&\frac{q^{\boldn^{(2)}-\frac12\boldm^{(1)}-\frac14\cdot1624}}
{(q)_{m_1}}
\qbinom{\frac12(m_1+m_2-m_3)}{m_2}
\qbinom{\frac12(m_2+m_4)}{m_3}\\
&\times\qbinom{\frac12(m_3+m_4-m_5-1)}{m_4}
\qbinom{\frac12(m_4+m_6)}{m_5}\\
&\times\qbinom{\frac12(m_5+m_6-m_7-1)}{m_6}
\qbinom{\frac12(m_6+m_8+1)}{m_7}
\qbinom{\frac12m_7}{m_8},
\end{split}
\end{equation*}
where $\boldm=(m_1,m_2,\ldots,m_8)$.
The summation here extends over an infinite number of terms because
the value of $m_1$ is unbounded.

\begin{figure}[ht]
\includegraphics[scale=0.75]{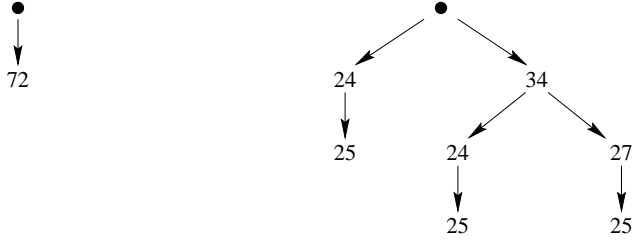}
\caption{Formations for $a=72$ and $b=25$ when $p'=75$ and $p=53$.}
\label{FormationEx2bFig}
\medskip
\end{figure}

To obtain a finitization $\chi^{53,75}_{a,b,c}(L)$ of $\chi^{53,75}_{17,72}$,
we set $a=72$, and can set $(b,c)=(24,23),(24,25),(25,24)$.
The most interesting case is the latter, so we choose that.
The Takahashi trees for $a=72$ and $b=25$ are given
in Fig.\ \ref{FormationEx2bFig}.

The single leaf-node of the Takahashi tree for $a=72$ gives rise to
the run $\run^L=\{ \{10\},\{2\},\{1\} \}$.
The three leaf-nodes of the Takahashi tree for $b=25$ gives rise
to the runs:
\begin{align*}
\run^R_1&=\{ \{10,5\},\{8,0\},\{1,-1\} \},\\
\run^R_2&=\{ \{10,5,2\},\{7,4,0\},\{-1,1,-1\} \},\\
\run^R_3&=\{ \{10,5,2\},\{7,5,1\},\{-1,1,1\} \}.
\end{align*}
Setting $\boldu^L=\boldu(\run^L)$,
$\boldDelta^L=\boldDelta(\run^L)$,
and $\boldu^R_j=\boldu(\run^R_j)$ and
$\boldDelta^R_j=\boldDelta(\run^R_j)$ for $j=1,2,3$,
we have:
\begin{align*}
\boldu^L_1&= (0, 0, 0,-1, 0,-1, 0, 0, 1),
&\boldDelta^L_1&= (0, 0, 0,-1, 0,-1, 0, 0, 1),\\
\boldu^R_1&= (0, -1, 0, -1, -1, 0, 0, 1, 1),
&\boldDelta^R_1&= (0, 1, 0, 1, 1, 0, 0, 1, -1),\\
\boldu^R_2&= (0, -1, 0, 0, -1, 0, 1, 0, 0),
&\boldDelta^R_2&= (0, 1, 0, 0, -1, 0, -1, 0, 0),\\
\boldu^R_3&= (1, -1, 0, 0, 0, 0, 1, 0, 0),
&\boldDelta^R_3&= (1, -1, 0, 0, 0, 0, -1, 0, 0).
\end{align*}
Now note that $\Delta(\boldu^R_1)=\Delta(\boldu^R_1)=-1=c-b$,
and that $\Delta(\boldu^R_1)=1=b-c$.
Therefore, since $b$ does not satisfy (\ref{interface}),
expression (\ref{Ferm2Eq}) applies in this case to give:
\begin{align*}
\chi^{53,75}_{72,25,24}(L)
&=F(\boldu^L,\boldu^R_1,L)+F(\boldu^L,\boldu^R_2,L)
                        +\tilde F(\boldu^L,\boldu^R_3,L)\\
&=F(\boldu^L,\boldu^R_1,L)+F(\boldu^L,\boldu^R_2,L)\\
&\qquad+q^{\frac12(L-47)}F(\boldu^L,\boldu^R_3,L)
          +(1-q^L)F(\boldu^L,\boldu^{R+}_3,L-1),
\end{align*}
the second equality being via the third case of (\ref{Tilde2Def}).

Using the definitions of Section \ref{FermLikeSec},
$\run^{R+}_3$ is given by $\{ \{10\},\{7\},\{-1\} \}=\run^R$,
as defined above. So $\boldu^{R+}_3=\boldu^R$.
For convenience, we set $\run^R_4=\run^R$ and $\boldu^R_4=\boldu^R$,
so that:
\begin{equation*}
\boldu^R_4= (0, 0, 0, 0, 0, 0, 1, 0, 0);\qquad
\boldDelta^R_4= (0, 0, 0, 0, 0, 0, -1, 0, 0).
\end{equation*}

For $j=1,2,3,4$, define
$\boldm^{(1)}_{1j}=(\boldu^L_{\flat}+\boldu^R_{j\sharp})\cdot\boldm$
for $\boldm=(m_1,m_2,\ldots,m_8)$,
and $\gamma'_{1j}=\gamma'(\run^L,\run^R_j)$.
Then:
\begin{align*}
&\boldm^{(1)}_{11}=-m_4-m_6+m_8,
&&\gamma'_{11}=-1624+2L;
\\
&\boldm^{(1)}_{12}=-m_6+m_7,
&&\gamma'_{12}=-1624+2L;
\\
&\boldm^{(1)}_{13}=-m_6+m_7,
&&\gamma'_{13}=-1530;
\\
&\boldm^{(1)}_{14}=-m_6+m_7,
&&\gamma'_{14}=-1624.
\end{align*}

The quadratic term $\boldm^{(2)}=\hat{\boldm}^T\boldC\hat{\boldm}-L^2$,
is given by $\boldm^{(2)}=4\boldn^{(2)}$ as specified above.

With $\boldu_{1j}=\boldu^L+\boldu^R_j$, and
$\boldQ_{1j}=\boldQ(\boldu_{1j})$, we have:
\begin{align*}
\boldu_{11}&= (0,-1, 0,-2,-1,-1, 0, 1, 2),
&\boldQ_{11}&= (1, 1, 1, 1, 0, 0, 1, 0),\\
\boldu_{12}&= (0,-1, 0,-1,-1,-1, 1, 0, 1),
&\boldQ_{12}&= (1, 1, 1, 1, 1, 0, 0, 1),\\
\boldu_{13}&= (1,-1, 0,-1, 0,-1, 1, 0, 1),
&\boldQ_{13}&= (1, 0, 0, 0, 1, 0, 0, 1),\\
\boldu_{14}&= (0, 0, 0,-1, 0,-1, 1, 0, 1),
&\boldQ_{14}&= (0,0,0,0,1,0,0,1).
\end{align*}

Then, using (\ref{FEq}), the above expression for
$\chi^{53,75}_{72,25,24}(L)$ produces:
\begin{equation*}
\chi^{53,75}_{72,25,24}(L)=
F_1(L) + F_2(L) + q^{\frac12(L-47)}F_3(L) + (1-q^L)F_4(L-1),
\end{equation*}
where for $j=1,2,3,4$, we set:
\begin{equation*}
\begin{split}
F_j(L)
=\sum_{\boldm\equiv\boldQ_{1j}}
&q^{\frac14\boldm^{(2)}-\frac12\boldm^{(1)}_{1j}+\frac14\gamma'_{1j}}
\qbinom{\frac12(L+m_2+(\boldu_{1j})_1)}{m_1}\\
&\times\qbinom{\frac12(m_1+m_2-m_3+(\boldu_{1j})_2)}{m_2}
\qbinom{\frac12(m_2+m_4+(\boldu_{1j})_3)}{m_3}\\
&\times\qbinom{\frac12(m_3+m_4-m_5+(\boldu_{1j})_4)}{m_4}
\qbinom{\frac12(m_4+m_6+(\boldu_{1j})_5)}{m_5}\\
&\times\qbinom{\frac12(m_5+m_6-m_7+(\boldu_{1j})_6)}{m_6}
\qbinom{\frac12(m_6+m_8+(\boldu_{1j})_7)}{m_7}\\
&\times
\qbinom{\frac12(m_7+(\boldu_{1j})_8)}{m_8},
\end{split}
\end{equation*}
where $\boldm=(m_1,m_2,\ldots,m_8)$.

Note that for $j=1,2$, the presence of $2L$ in $\gamma'_{1j}$
implies that $\lim_{L\to\infty} F_j(L)=0$.
Thus, as $L\to\infty$, only the term $F_4(L-1)$ remains.
This yields precisely the fermionic expression given above
for the character $\chi^{53,75}_{17,72}$.

\addtocontents{toc}{\SkipTocEntry}
\subsection{Example 3}

In this example, we provide a case when the term
on the right of (\ref{FermCEq}) appears.
We will not give every detail because the
expression involves many terms (thirty-one!), but will
concentrate on how the extra term(s) are dealt with.

\begin{figure}[ht]
\includegraphics[scale=0.75]{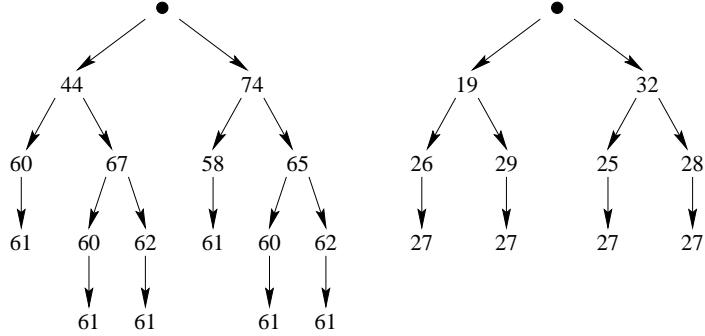}
\caption{Takahashi tree for $s=61$ and truncated Takahashi tree for
         $r=27$ when $p'=118$ and $p=51$.}
\label{FormationEx3aFig}
\medskip
\end{figure}

Let $p'=118$, $p=51$, $r=27$ and $s=61$.
Then $p'/p$ has continued fraction $[2,3,5,3]$, so that $n=3$,
$\{t_k\}_{k=0}^3=\{-1,1,4,9,12\}$, $t=11$,
$\{y_k\}_{k=-1}^3=\{0,1,2,7,37,118\}$ and
$\{z_k\}_{k=-1}^3=\{1,0,1,3,16,51\}$.
Then:
\begin{align*}
\T&=\{1,2,3,5,7,9,16,23,30,37,44\};\\
\T'&=\{117,116,115,113,111,109,102,95,88,81,74\};\\
\tT&=\{1,2,3,4,7,10,13,16,19\};\\
\tT'&=\{50,49,48,47,44,41,38,35,32\}.
\end{align*}
Using these, we obtain the Takahashi tree for $s=61$ and the truncated
Takahashi tree for $r=27$ shown in Fig.~\ref{FormationEx3aFig}.

We see that the Takahashi tree for $s=61$ and the truncated
Takahashi tree for $r=27$ have six and four leaf-nodes respectively.
Thus $|\U(s)|=6$ and $|\tU(r)|=4$.
Therefore, for $\chi^{51,118}_{27,61}$, the sum in expression
(\ref{FermCEq}) runs over 24 terms $F(\boldu^L,\boldu^R)$.
We will not compute these terms here, but just denote this
sum by $\sum^{(11)}$, where the superscript indicates that each element
of $\U(s)$ and $\tU(r)$ is an $11$-dimensional vector.

{}From the prescription of Section \ref{ExtraSec}, we obtain
$\{\xi_\ell\}_{\ell=0}^5=\{0,37,44,74,81,118\}$
and
$\{\tilde\xi_\ell\}_{\ell=0}^5=\{0,16,19,32,35,51\}$.
Then $\eta(61)=\tilde\eta(27)=2$ and $\hat s=61-44=17$, $\hat r=27-19=8$.
The additional term in (\ref{FermCEq}) is therefore
present in this case.
We obtain $\hat p'=74-44=30$ and $\hat p=32-19=13$.
Thereupon:
\begin{equation*}
\chi^{51,118}_{27,61}=\chi^{13,30}_{8,17}+\sum{}^{(11)}.
\end{equation*}

\begin{figure}[ht]
\includegraphics[scale=0.75]{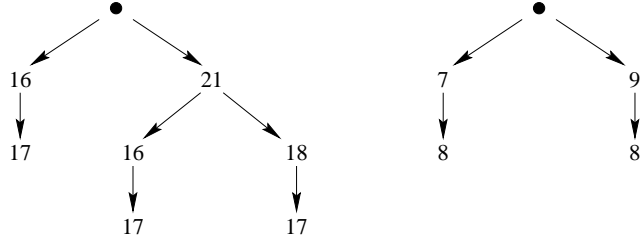}
\caption{Takahashi tree for $s=17$ and truncated Takahashi tree for
         $r=8$ when $p'=30$ and $p=13$.}
\label{FormationEx3bFig}
\medskip
\end{figure}

We now use (\ref{FermCEq}) to determine a fermionic expression for
$\chi^{13,30}_{8,17}$. So now set $p'=30$ and $p=13$.
The continued fraction of $p'/p$ is $[2,3,4]$.
Note that this may be obtained directly from $(\ref{ExtraCfEq})$.
{}From this, we obtain $n=2$,
$\{t_k\}_{k=0}^3=\{-1,1,4,8\}$, $t=7$,
$\{y_k\}_{k=-1}^3=\{0,1,2,7,30\}$,
$\{z_k\}_{k=-1}^3=\{1,0,1,3,13\}$,
$\T=\{1,2,3,5,7,9,16\}$,
$\T'=\{29,28,27,25,23,21,14\}$,
$\tT=\{1,2,3,4,7\}$ and
$\tT'=\{12,11,10,9,6\}$.
Resetting $s=\hat s=17$ and $r=\hat r=8$, we use these values
to obtain the Takahashi tree for $s=17$ and the truncated
Takahashi tree for $r=8$ that are shown in Fig.~\ref{FormationEx3bFig}.

Since these trees have three and two nodes respectively,
we have $|\U(s)|=3$ and $|\tU(r)|=2$, whereupon the
sum in expression (\ref{FermCEq}) for $\chi^{13,30}_{8,17}$
runs over six terms $F(\boldu^L,\boldu^R)$.
Again, we will not compute these terms here, but just denote this
sum by $\sum^{(7)}$, where the superscript is used as above.

In this $p'=30$, $p=13$ case, the prescription of Section \ref{ExtraSec}
yields $\{\xi_\ell\}_{\ell=0}^7=\{0,7,9,14,16,21,23,30\}$
and $\{\tilde\xi_\ell\}_{\ell=0}^7=\{0,3,4,6,7,9,10,13\}$.
Then here, $\eta(s)=\eta(r)=4$ and $\hat s=s-16=1$ and $\hat r=r-6=1$.
Therefore, we again have to include the extra term in (\ref{FermCEq}).
After calculating $\hat p'=21-16=5$ and $\hat p=9-7=2$, we therefore
have:
\begin{equation*}
\chi^{51,118}_{27,61}=\chi^{2,5}_{1,1}+\sum{}^{(7)}+\sum{}^{(11)}.
\end{equation*}
A further use of (\ref{FermCEq}) expresses $\chi^{2,5}_{1,1}$
in the fermionic form given in (\ref{RRinfferm}).
We have thus expressed the Virasoro character $\chi^{51,118}_{27,61}$
in fermionic form, as a sum of thirty-one fundamental fermionic forms.

A finitization of this character is provided by
$\chi^{51,118}_{a,b,c}$ with $a=61$ and either $(b,c)=(62,63)$,
$(62,61)$, $(63,64)$ or $(63,62)$.
We choose the latter and give a very brief outline
of how the prescription of Section \ref{PrologueSec} produces a fermionic
expression for $\chi^{51,118}_{61,63,62}(L)$.

With $p=51$ and $p'=118$, the values pertaining to the continued
fraction of $p'/p$ calculated above yield Takahashi trees for
$a=61$ and $b=63$ that contain six and five leaf-nodes respectively.
Thus $|{\mathcal U}(a)|=6$ and $|{\mathcal U}(b)|=5$.
Since $b$ satisfies (\ref{interface}), expression (\ref{FermEq})
is the appropriate expression to use for $\chi^{51,118}_{61,63,62}(L)$.
The sum in this expression thus runs over thirty terms
$F(\boldu^L,\boldu^R,L)$.
We denote this sum by $\sum^{(11)}(L)$.

We also find that $\eta(a)=\eta(b)=2$ and $\hat a=a-44=17$,
$\hat b=b-44=19$ and $\hat c=c-44=18$.
The additional term in (\ref{FermEq}) is therefore present in this case.
Again we have $\hat p'=74-44=30$ and $\hat p=32-19=13$, whereupon:
\begin{equation*}
\chi^{51,118}_{61,63,62}(L)=\chi^{13,30}_{17,19,18}(L)+\sum{}^{(11)}(L).
\end{equation*}

With $p=13$ and $p'=30$, $b=19$ satisfies (\ref{interface}) and therefore
we use (\ref{FermEq}) once more to express $\chi^{13,30}_{17,19,18}(L)$
in fermionic form.
The values pertaining to the continued
fraction of $30/13$ calculated above yield Takahashi trees for
$17$ and $19$ that contain three and two leaf-nodes respectively.
Therefore, for $\chi^{13,30}_{17,19,18}(L)$, the sum in expression
(\ref{FermEq}) runs over six terms $F(\boldu^L,\boldu^R,L)$.
We denote their sum by $\sum^{(7)}(L)$.

In this $p'=30$, $p=13$ case, we find that $\eta(17)=\eta(19)=4$
and $\hat a=17-16=1$, $\hat b=19-16=3$ and $\hat c=18-16=2$.
Therefore, we again require the extra term from (\ref{FermEq}).
With $\hat p'=21-16=5$ and $\hat p=9-7=2$, we thus have:
\begin{equation*}
\chi^{51,118}_{61,63,62}(L)=\chi^{2,5}_{1,3,2}(L)
                   +\sum{}^{(7)}(L)+\sum{}^{(11)}(L).
\end{equation*}

We may now use (\ref{FermEq}) yet again to express $\chi^{2,5}_{1,3,2}(L)$
in fermionic form.
The result is equivalent to that given in (\ref{RRfinferm}).
In this way, $\chi^{51,118}_{61,63,62}(L)$ is expressed in fermionic form,
as a sum of thirty-seven fundamental fermionic forms.

In the $L\to\infty$ limit, six of the fundamental fermionic forms
that comprise $\sum{}^{(11)}(L)$ tend to zero.
The remaining thirty-one terms reproduce the fermionic expression
for $\chi^{51,118}_{27,61}$ calculated above.

\newpage

\section{Obtaining the bosonic generating function}\label{BosonicSec}

In Section \ref{WeightSec}, the path generating
function $\chi^{p,p'}_{a,b,c}(L)$ was defined by:
$$
\chi^{p,p'}_{a,b,c}(L)=\sum_{p\in\P^{p,p'}_{a,b,c}(L)} q^{wt(h)}
$$
where $wt(h)$ was defined in (\ref{WtDef}).
Here we prove the bosonic expression for $\chi^{p,p'}_{a,b,c}(L)$
stated in (\ref{FinRochaEq}).
Our proof is similar to that given in \cite{forrester-baxter}
for what, via the bijection given in \cite{flpw}, turns out
to be an equivalent result.

Let $1\le p<p'$,  $1\le a\le p'-1$ and $L>0$.
If $2\le b\le p'-2$, then
the definition (\ref{WtDef}) implies the following recurrence 
relations for $\chi^{p,p'}_{a,b,c}(L)$:
if $\lfloor bp/p'\rfloor=\lfloor (b+1)p/p'\rfloor$ then
\begin{equation}\label{RecEq1}
\chi^{p,p'}_{a,b,b+1}(L)
 =q^{\frac12(L-a+b)}\chi^{p,p'}_{a,b+1,b}(L-1)
 +\chi^{p,p'}_{a,b-1,b}(L-1);
\end{equation}
if $\lfloor bp/p'\rfloor\ne\lfloor (b+1)p/p'\rfloor$ then
\begin{equation}\label{RecEq2}
\chi^{p,p'}_{a,b,b+1}(L)
 =\chi^{p,p'}_{a,b+1,b}(L-1)
 +q^{\frac12(L-a+b)}\chi^{p,p'}_{a,b-1,b}(L-1);
\end{equation}
if $\lfloor bp/p'\rfloor=\lfloor (b-1)p/p'\rfloor$ then
\begin{equation}\label{RecEq3}
\chi^{p,p'}_{a,b,b-1}(L)
 =q^{\frac12(L+a-b)}\chi^{p,p'}_{a,b-1,b}(L-1)
 +\chi^{p,p'}_{a,b+1,b}(L-1);
\end{equation}
if $\lfloor bp/p'\rfloor\ne\lfloor (b-1)p/p'\rfloor$ then
\begin{equation}\label{RecEq4}
\chi^{p,p'}_{a,b,b-1}(L)
 =\chi^{p,p'}_{a,b-1,b}(L-1)
 +q^{\frac12(L+a-b)}\chi^{p,p'}_{a,b+1,b}(L-1).
\end{equation}
In addition, (\ref{WtDef}) implies that:
if $\lfloor p/p'\rfloor=\lfloor 2p/p'\rfloor$ then
\begin{equation}\label{BoundEq1}
\chi^{p,p'}_{a,1,2}(L)
 =q^{\frac12(L-a+1)}\chi^{p,p'}_{a,2,1}(L-1);
\end{equation}
if $\lfloor p/p'\rfloor\ne\lfloor 2p/p'\rfloor$ then
\begin{equation}\label{BoundEq2}
\chi^{p,p'}_{a,1,2}(L)
 =\chi^{p,p'}_{a,2,1}(L-1);
\end{equation}
if $\lfloor (p'-1)p/p'\rfloor=\lfloor (p'-2)p/p'\rfloor$ then
\begin{equation}\label{BoundEq3}
\chi^{p,p'}_{a,p'-1,p'-2}(L)
 =q^{\frac12(L+a-p'+1)}\chi^{p,p'}_{a,p'-2,p'-1}(L-1);
\end{equation}
if $\lfloor (p'-1)p/p'\rfloor\ne\lfloor (p'-2)p/p'\rfloor$ then
\begin{equation}\label{BoundEq4}
\chi^{p,p'}_{a,p'-1,p'-2}(L)
 =\chi^{p,p'}_{a,p'-2,p'-1}(L-1).
\end{equation}
Finally, we obtain:
\begin{equation}\label{StartEq}
\chi^{p,p'}_{a,b,b\pm1}(0)=\delta_{a,b}.
\end{equation}

Since expressions (\ref{RecEq1})-(\ref{BoundEq4}) enable
$\chi^{p,p'}_{a,b,c}(L)$, for $1\le a,b,c<p'$ with $b-c=\pm1$ and $L>0$,
to be expressed in terms of various $\chi^{p,p'}_{a,b',c'}(L-1)$,
with $1\le b',c'<p'$ and $b'-c'=\pm1$,
the above expressions together with (\ref{StartEq})
determine $\chi^{p,p'}_{a,b,c}(L)$ uniquely.

For $1\le p<p'$ and $1\le a,c<p'$ and $0\le b\le p'$,
with $L\equiv a-b\,(\mod2)$ and $L\ge0$, define:
\begin{align*}
\pchi^{p,p'}_{a,b,c}(L)&=
\sum_{\lambda=-{\infinity}}^{\infinity}
q^{\lambda^{2} p p'+ \lambda (p'r(b,c)-pa)}
\left[ {L \atop {\frac{L+a-b}{2}}-p'\lambda} \right]_q\\[0.5mm]
&\qquad\qquad\qquad
-\sum_{\lambda=-\infinity}^\infinity
q^{(\lambda p+r(b,c))(\lambda p'+a)}
\left[ {L \atop {\frac{L-a-b}{2}}-p'\lambda} \right]_q,
\end{align*}
where
\begin{equation}\label{Appgroundstatelabel1}
r(b,c)=\lfloor pc/p'\rfloor+(b-c+1)/2.
\end{equation}

Once we establish that expressions (\ref{RecEq1})-(\ref{StartEq})
are satisfied with $\pchi^{p,p'}_{a,b,c}(L)$ in place
of $\chi^{p,p'}_{a,b,c}(L)$, we can then conclude that
$\chi^{p,p'}_{a,b,c}(L)=\pchi^{p,p'}_{a,b,c}(L)$.

First consider the case $c=b+1$ with
$\lfloor bp/p'\rfloor=\lfloor cp/p'\rfloor$.
Set $r=r(b,b+1)=\lfloor bp/p'\rfloor$, whence
$r(b+1,b)=r+1$ and $r(b-1,b)=r$.
Then, on using the recurrence
$\left[{A\atop B}\right]_q
=q^{A-B} \left[{A-1\atop B-1}\right]_q+\left[{A-1\atop B}\right]_q$
for Gaussian polynomials, we obtain:
\begin{align*}
\pchi^{p,p'}_{a,b,b+1}(L)&=
q^{\frac12(L-a+b)}
\left(
\sum_{\lambda=-{\infinity}}^{\infinity}
q^{\lambda^{2} p p'+ \lambda (p'(r+1)-pa)}
\left[ {L-1 \atop {\frac{(L-1)+a-(b+1)}{2}}-p'\lambda} \right]_q\right.
\\[0.5mm]
&\qquad\qquad\qquad\qquad
\left.
-\sum_{\lambda=-\infinity}^\infinity
q^{(\lambda p+r+1)(\lambda p'+a)}
\left[ {L-1 \atop {\frac{(L-1)-a-(b+1)}{2}}-p'\lambda} \right]_q\right)
\\[0.5mm]
&\qquad
+\left(
\sum_{\lambda=-{\infinity}}^{\infinity}
q^{\lambda^{2} p p'+ \lambda (p'r-pa)}
\left[ {L-1 \atop {\frac{(L-1)+a-(b-1)}{2}}-p'\lambda} \right]_q\right.
\\[0.5mm]
&\qquad\qquad\qquad\qquad
\left.
-\sum_{\lambda=-\infinity}^\infinity
q^{(\lambda p+r)(\lambda p'+a)}
\left[ {L-1 \atop {\frac{(L-1)-a-(b-1)}{2}}-p'\lambda} \right]_q\right)
\\[0.5mm]
&=
q^{\frac12(L-a+b)}
\pchi^{p,p'}_{a,b+1,b}(L-1)+\pchi^{p,p'}_{a,b-1,b}(L-1).
\end{align*}
Thus, (\ref{RecEq1}) holds with $\pchi^{p,p'}_{a,b,c}(L)$
in place of $\chi^{p,p'}_{a,b,c}(L)$.

Using $r(0,1)=0$, we calculate:
\begin{align*}
\pchi^{p,p'}_{a,0,1}(L)&=
\sum_{\lambda=-{\infinity}}^{\infinity}
q^{\lambda^{2} p p'- \lambda pa}
\left[ {L \atop {\frac{L+a}{2}}-p'\lambda} \right]_q
-\sum_{\lambda=-\infinity}^\infinity
q^{\lambda p(\lambda p'+a)}
\left[ {L \atop {\frac{L-a}{2}}-p'\lambda} \right]_q \\[0.5mm]
&=0,
\end{align*}
after changing the sign of the second summation parameter and using
$\left[{A\atop B}\right]_q=
\left[{A\atop A-B}\right]_q$.
On substituting this into the previous expression when $b=1$, we obtain
$\pchi^{p,p'}_{a,1,2}(L)=
q^{\frac12(L-a+1)} \pchi^{p,p'}_{a,2,1}(L-1)$,
and thus (\ref{BoundEq1}) holds with $\pchi^{p,p'}_{a,b,c}(L)$
in place of $\chi^{p,p'}_{a,b,c}(L)$.

In a similar way, perhaps making use of the recurrence
$\left[{A\atop B}\right]_q
= \left[{A-1\atop B-1}\right]_q+q^B\left[{A-1\atop B}\right]_q$
for Gaussian polynomials,
we obtain each of (\ref{RecEq2})-(\ref{RecEq4})
and (\ref{BoundEq2})-(\ref{BoundEq4}) with
$\pchi^{p,p'}_{a,b,c}(L)$ in place of $\chi^{p,p'}_{a,b,c}(L)$.

Finally, with $1\le a,b,c<p$, we obtain:
\begin{align*}
\pchi^{p,p'}_{a,b,c}(0)&=
\sum_{\lambda=-{\infinity}}^{\infinity}
q^{\lambda^{2} p p'+ \lambda (p'r(b,c)-pa)}
\left[ {0 \atop {\frac{a-b}{2}}-p'\lambda} \right]_q\\[0.5mm]
&\qquad\qquad\qquad
-\sum_{\lambda=-\infinity}^\infinity
q^{(\lambda p+r(b,c))(\lambda p'+a)}
\left[ {0 \atop {\frac{a+b}{2}}+p'\lambda} \right]_q\\[0.5mm]
&=\delta_{a,b},
\end{align*}
since $\left[{0\atop B}\right]_q=\delta_{B,0}$ and
$\vert\frac12(a\pm b)\vert<p'$.

We thus conclude that if $1\le a,b,c<p'$ with $c=b\pm1$, then
$\chi^{p,p'}_{a,b,c}(L)=\pchi^{p,p'}_{a,b,c}(L)$, thus
verifying (\ref{FinRochaEq}).

In Section \ref{WeightSec}, we require an extension of this result to
where $c=0$ or $c=p'$. In accordance with (\ref{groundstatelabel2}),
extend the definition (\ref{Appgroundstatelabel1}) to:
\begin{equation}\label{Appgroundstatelabel2}
r(b,c)= \begin{cases}
       \lfloor pc/p'\rfloor+(b-c+1)/2 &\text{if $1\le c<p'$;}\\
             1 &\text{if $c=0$ and $p'>2p$;}\\
             0 &\text{if $c=0$ and $p'<2p$;}\\
             p-1 &\text{if $c=p'$ and $p'>2p$;}\\
             p &\text{if $c=p'$ and $p'<2p$.}
   \end{cases}
\end{equation}
In the case $p'>2p$, that the $0$th and
$(p'-1)$th bands are even implies that:
\begin{equation*}
\chi^{p,p'}_{a,1,0}(L)=q^{-\frac12(L+1-a)}\chi^{p,p'}_{a,1,2}(L);
\quad
\chi^{p,p'}_{a,p'-1,p'}(L)=q^{-\frac12(L-p'+1+a)}\chi^{p,p'}_{a,p'-1,p'-2}(L).
\end{equation*}
Since $r(1,0)=1$, $r(1,2)=0$, $r(p'-1,p')=p-1$ and $r(p'-1,p'-2)=p$,
it follows in this $p'>2p$ case that
$\chi^{p,p'}_{a,b,c}(L)=\pchi^{p,p'}_{a,b,c}(L)$ for
$1\le a,b<p'$ and $0\le c\le p'$ with $b-c=\pm1$ and $L\ge0$.

In the case $p'<2p$, that the $0$th and $(p'-1)$th bands are odd implies that:
\begin{equation*}
\chi^{p,p'}_{a,1,0}(L)=q^{\frac12(L+1-a)}\chi^{p,p'}_{a,1,2}(L);
\quad
\chi^{p,p'}_{a,p'-1,p'}(L)=q^{\frac12(L-p'+1+a)}\chi^{p,p'}_{a,p'-1,p'-2}(L).
\end{equation*}
Since $r(1,0)=0$, $r(1,2)=1$, $r(p'-1,p')=p$ and $r(p'-1,p'-2)=p-1$,
we conclude that
$\chi^{p,p'}_{a,b,c}(L)=\pchi^{p,p'}_{a,b,c}(L)$ for
$1\le a,b<p'$ and $0\le c\le p'$ with $b-c=\pm1$ and $L\ge0$, as required.

\newpage

\section{Bands and the floor function}\label{AppBSec}

In this section, we derive some basic results concerning the floor function
$\lfloor\cdot\rfloor$ and how it relates to the $(p,p')$-model and other
closely related models.
$\rho^{p,p'}(a)$ is defined in Section \ref{BandSec}.
$\delta^{p,p'}_{a,e}$ is defined in Section \ref{PparamSec}.
$\omega^{p,p'}(a)$ is defined at the beginning of Section \ref{MazySec}.

\begin{lemma}\label{StartPtLem}
Let $1\le p<p'$. If $1\le a<p'$,
$e\in\{0,1\}$ and $a'=a+\lfloor ap/p'\rfloor+e$
then $\lfloor a'p/(p'+p)\rfloor=\lfloor ap/p'\rfloor$ and
$\delta^{p,p'+p}_{a',e}=0$.
On the other hand, if $1\le a'<p'+p$, $e\in\{0,1\}$,
$\delta^{p,p'+p}_{a',e}=0$ and
$a=a'-\lfloor a'p/(p'+p)\rfloor-e$ then
$\lfloor ap/p'\rfloor=\lfloor a'p/(p'+p)\rfloor$.

In addition, if $a$ is interfacial in the $(p,p')$-model and
$\delta^{p,p'}_{a,e}=0$, or if $a$ is multifacial in the $(p,p')$-model,
then $a'$ is interfacial in the $(p,p'+p)$-model.
\end{lemma}

\Proof Let $r=\lfloor ap/p'\rfloor$ whence $p'r\le pa<p'(r+1)$.
Then, for $x\in\{0,1\}$, we have
$(p'+p)r\le p(a+r+x)<(p'+p)r+p'+xp$, so that
$\lfloor (a+r+x)p/(p'+p)\rfloor=r$.
In particular, $\lfloor a'p/(p'+p)\rfloor=r$,
and $\lfloor (a+r+e+(-1)^e)p/(p'+p)\rfloor=r$.
Thus $r=\lfloor a'p/(p'+p)\rfloor=\lfloor (a'+(-1)^e)p/(p'+p)\rfloor$
which gives the first results.

For the second statement, let $r=\lfloor a'p/(p'+p)\rfloor$.
With $r=\lfloor (a'+(-1)^e)p/(p'+p)\rfloor$, this implies that
$(p'+p)r\le p(a+r+x)<(p'+p)(r+1)$ for both $x\in\{0,1\}$.
That $\lfloor ap/p'\rfloor=\lfloor a'p/(p'+p)\rfloor$ then follows.

Now, if $a$ is interfacial in the $(p,p')$-model and $\delta^{p,p'}_{a,e}=0$
then $\lfloor (a+(-1)^e)p/p'\rfloor=\lfloor ap/p'\rfloor
\ne\lfloor (a-(-1)^e)p/p'\rfloor$.
Thus $r\ne\lfloor (a-(-1)^e)p/p'\rfloor$.
Clearly, this inequality also holds if $a$ is multifacial.
When $e=0$, we have $(a-1)p<rp'$ and thus $(a+r-1)p<r(p'+p)$ so that
$\lfloor (a'-1)p/(p'+p)\rfloor<r$, which when compared to the above
result implies that $a'$ is interfacial in the $(p,p'+p)$-model.
Similarly $e=1$ gives $(a+1)p\ge(r+1)p'$ whence
$\lfloor (a'+1)p/(p'+p)\rfloor\ge r+1$, which when compared to the above
result also implies that $a'$ is interfacial in the $(p,p'+p)$-model.
\cqfd
\medskip

\begin{lemma}\label{DParamLem}
Let $1\le p<p'$ with $p$ coprime to $p'$ and $1\le a<p'$.
Then $\lfloor a(p'-p)/p'\rfloor=a-1-\lfloor ap/p'\rfloor$.

If, in addition, $a$ is interfacial in the $(p,p')$-model
and $\delta^{p,p'}_{a,e}=0$ then
$a$ is interfacial in the $(p'-p,p')$-model
and $\delta^{p'-p,p'}_{a,1-e}=0$.
\end{lemma}

\Proof Since $p$ and $p'$ are coprime, $\lfloor ap/p'\rfloor<ap/p'$.
Hence $\lfloor ap/p'\rfloor+\lfloor a(p'-p)/p'\rfloor=a-1$.

Since the $(p,p')$-model differs from the $(p'-p,p')$-model only in
that corresponding bands are of the opposite parity, $a$ being interfacial
in one model implies that it also is in the other.
The final part then follows immediately.
\cqfd
\medskip

\begin{lemma}\label{BDParamLem}
Let $1\le p<p'<2p$ with $p$ coprime to $p'$,
$1\le a<p'$ and $e\in\{0,1\}$ and set
$a'=2a-e-\lfloor a(p'-p)/p'\rfloor$.
Then $\delta^{p,p'+p}_{a',1-e}=0$ and
$\lfloor a'p/(p'+p)\rfloor=a-1-\lfloor a(p'-p)/p'\rfloor$.

Furthermore, if $\delta^{p'-p,p'}_{a,e}=0$ then $a'$ is interfacial
in the $(p,p'+p)$-model.
\end{lemma}

\Proof Lemma \ref{DParamLem} implies that
$a'=a+1-e+\lfloor a'p/p'\rfloor$.
Lemma \ref{StartPtLem} immediately implies that $\delta^{p,p'+p}_{a',1-e}=0$.
Using first Lemma \ref{StartPtLem} and then Lemma \ref{DParamLem} yields
$\lfloor a'p/(p'+p)\rfloor=\lfloor ap/p'\rfloor
=a-1-\lfloor a(p'-p)/p'\rfloor$.

If $\delta^{p'-p,p'}_{a,e}=0$ then $a$ either is interfacial
in the $(p'-p,p')$-model or lies between two even bands.
In the former case, Lemma \ref{DParamLem} implies that $a$ is
interfacial in the $(p,p')$-model and $\delta^{p,p'}_{a,1-e}=0$.
In the latter case, $a$ is multifacial in the $(p,p')$-model.
In either case, Lemma \ref{StartPtLem} then shows that $a'$ is
interfacial in the $(p,p'+p)$-model.
\cqfd
\medskip

\begin{lemma}\label{InterLem}
1. Let $0\le a\le p'$ with $a$ interfacial in the $(p,p')$-model.
Let $r=\rho^{p,p'}(a)$.
Then $a+r$ is interfacial in the $(p,p'+p)$-model and
$\rho^{p,p'+p}(a+r)=r$.
Moreover, $\omega^{p,p'+p}(a+r)=\omega^{p,p'}(a)$.

2. Let $p'<2p$ and $0\le a\le p'$ with $a$ interfacial in the
$(p'-p,p')$-model.
Let $r=\rho^{p'-p,p'}(a)$.
Then $a$ is interfacial in the $(p,p')$-model and
$\rho^{p,p'}(a)=a-r$.
Moreover, if $R=\lfloor rp/(p'-p)\rfloor$, then
$\omega^{p'-p,p'}(a)=r^-\implies \omega^{p,p'}(a)=R^+$,
and
$\omega^{p'-p,p'}(a)=r^+\implies \omega^{p,p'}(a)=(R+1)^-$.
\end{lemma}

\Proof 1.
Since $a$ is interfacial in the $(p,p')$-model, we have
$r-1\le (a-1)p/p'< r\le (a+1)p/p'<r+1$ whence
$(r-1)p'\le (a-1)p< rp'\le (a+1)p<(r+1)p'$ and
$r(p'+p)-p'\le (a+r-1)p< r(p'+p)\le (a+r+1)p<r(p'+p)+p'$ so that
$r-p'/(p'+p)\le (a+r-1)p/(p'+p)< r\le (a+r+1)p/(p'+p)<r+p'/(p'+p)$ and
therefore $a+r$ is interfacial in the $(p,p'+p)$-model
and $\rho^{p,p'+p}(a+r)=r$.
Now, $\lfloor r(p'+p)/p\rfloor=\lfloor rp'/p\rfloor+r$,
whence by definition, $\omega^{p,p'+p}(a+r)=\omega^{p,p'}(a)$.

2. Since corresponding bands in the $(p,p')$- and the
$(p'-p,p')$-models are of opposite parity, $a$ being interfacial
in the $(p'-p,p')$-model implies that $a$ is interfacial in the
$(p,p')$-model.
Then $\omega^{p'-p,p'}(a)=r^-$ implies that
$r=\lfloor (p'-p)a/p'\rfloor+1=a-1-\lfloor pa/p'\rfloor+1$.
Thus $\lfloor pa/p'\rfloor=a-r=\lfloor rp'/(p'-p)\rfloor-r=R$,
and since $a$ has an odd band immediately below it in the $(p,p')$-model,
$\rho^{p,p'}(a)=a-r$ and $\omega^{p,p'}(a)=R^+$.
In a similar way, $\omega^{p'-p,p'}(a)=r^+$ implies that
$\rho^{p,p'}(a)=a-r=R+1$ and $\omega^{p,p'}(a)=(R+1)^-$.
\cqfd
\medskip

\begin{lemma}\label{BadInterLem}
Let $0<a'<p'+p-1$ with both $a'$ and $a'+1$ interfacial in the
$(p,p'+p)$-model and set $a=a'-\rho^{p,p'+p}(a')$.
If $\rho^{p,p'+p}(a'+1)=\rho^{p,p'+p}(a')+1$ then
$a$ is multifacial in the $(p,p')$-model.
\end{lemma}

\Proof Let $r=\rho^{p,p'+p}(a')$.
Then the $r$th band in the $(p,p'+p)$-model lies between heights
$a'-1$ and $a'$ so that $a'-1<r(p'+p)/p<a'$ whereupon
$a'-r-1<rp'/p<a'-r$ so that the $r$th odd band in the $(p,p')$-model
lies between heights $a-1$ and $a$.
In a similar way, we see that the $(r+1)$th odd band in the $(p,p')$-model
lies between heights $a$ and $a+1$. The lemma follows.
\cqfd
\medskip

\newpage

\section{Bands on the move}\label{DireApp}

In the main proof that is presented in Section \ref{ProofIndSec},
whenever Lemma \ref{ExtGen1Lem} (or similarly Lemma \ref{ExtGen2Lem})
is invoked to obtain a generating function by extending paths,
it is sometimes necessary to prove that
(in the notation of Lemma \ref{ExtGen1Lem})
$a$ and $a+2\Delta$ are both interfacial in the $(p,p')$-model.
In the case of $a+2\Delta$, this is somewhat tricky.
In this section, we present a number of quite technical auxiliary
results which facilitate this proof.
Each result examines how the relative displacement of two particular
startpoints (or similarly endpoints) in a certain model
changes under certain combinations of the $\B$- and $\D$-transforms.

The following result examines the differing effects of a
$\B\D$-transform on two adjacent points which are separated
by an odd band in the $(p'-p,p')$-model and for which the
associated pre-segments are in even bands.

\begin{figure}[ht]
  \psfrag{a}{\footnotesize$a$}
  \psfrag{c}{\footnotesize$\hat a$}
  \psfrag{x}{\footnotesize$a'$}
  \psfrag{z}{\footnotesize$\hat a'$}
\includegraphics[scale=1.00]{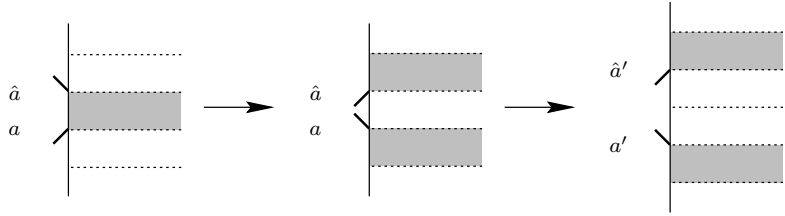}
\caption{Situation in Lemma \ref{Dire1Lem}.}
\label{Dire1Fig}
\medskip
\end{figure}
 
\begin{lemma}\label{Dire1Lem}
Let $p'<2p$ and $1\le a\le p'-2$, with
$\lfloor a(p'-p)/p'\rfloor\ne \lfloor (a+1)(p'-p)/p'\rfloor$,
so that $a$ and $\hat a=a+1$ are both interfacial in the
$(p'-p,p')$-model. 
Let $e,\hat e\in\{0,1\}$ be such that
$\delta^{p'-p,p'}_{a,e}=\delta^{p'-p,p'}_{\hat a,\hat e}=0$.
Let $a'=a+1-e+\lfloor ap/p'\rfloor$ and
$\hat a'=\hat a+1-\hat e+\lfloor\hat ap/p'\rfloor$.
Then $a'$ and $\hat a'$ are both interfacial in the $(p,p'+p)$-model
and $\hat a'=a'+2$.
In addition, $\lfloor a'p/(p'+p)\rfloor=\lfloor\hat a'p/(p'+p)\rfloor$.
\end{lemma}

\Proof Since $a$ and $\hat a$ are separated by an even band in
the $(p,p')$-model,
$\lfloor ap/p'\rfloor=\lfloor\hat ap/p'\rfloor$, whence
$\hat a'=a'+2$ immediately.
In addition, $\delta^{p'-p,p'}_{a,e}=\delta^{p'-p,p'}_{\hat a,\hat e}=0$
implies that $e=1$ and $\hat e=0$.
Then $\delta^{p,p'}_{a,0}=0$ which, by Lemma \ref{StartPtLem},
implies that $a'$ is interfacial in the $(p,p'+p)$-model.
Similarly, $\delta^{p,p'}_{\hat a,1}=0$ implies that
$\hat a'$ is interfacial in the $(p,p'+p)$-model.
Lemma \ref{StartPtLem} also implies that
$\delta^{p,p'+p}_{a',0}=\delta^{p,p'}_{\hat a',1}=0$,
so that $\lfloor a'p/(p'+p)\rfloor=
\lfloor (a'+1)p/(p'+p)\rfloor=
\lfloor\hat a'p/(p'+p)\rfloor$.
\cqfd
\medskip

\noindent
The conditions in the premise of Lemma \ref{Dire1Lem} of course
imply that $e=1$ and $\hat e=0$. However, it will be more
readily applied in the stated format.

The following result examines the differing effects of a
$\B\D$-transform on two adjacent points which are separated
by an even band in the $(p'-p,p')$-model and for which the
associated pre-segments are in even bands and in the same direction.

\begin{figure}[ht]
  \psfrag{a}{\footnotesize$a$}
  \psfrag{c}{\footnotesize$\hat a$}
  \psfrag{x}{\footnotesize$a'$}
  \psfrag{z}{\footnotesize$\hat a'$}
\includegraphics[scale=1.00]{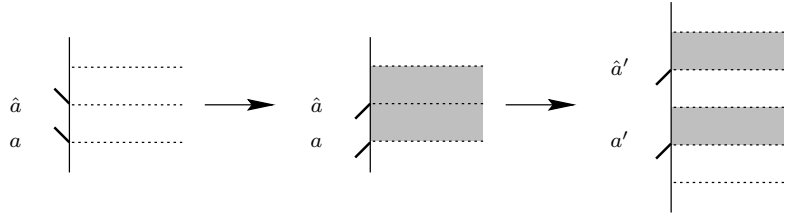}
\caption{Situation in Lemma \ref{Dire2Lem}.}
\label{Dire2Fig}
\medskip
\end{figure}

\begin{lemma}\label{Dire2Lem}
Let $p'<2p$ and $1\le a\le p'-2$, with
$\lfloor a(p'-p)/p'\rfloor=\lfloor \hat a(p'-p)/p'\rfloor$
and $\hat a=a+1$, and let $e\in\{0,1\}$ be such that
$\delta^{p'-p,p'}_{a,e}=\delta^{p'-p,p'}_{\hat a,e}=0$.
Let $a'=2a-e-\lfloor a(p'-p)/p'\rfloor$ and
$\hat a'=2\hat a-e-\lfloor\hat a(p'-p)/p'\rfloor$.
Then $a'$ and $\hat a'$ are both interfacial in the $(p,p'+p)$-model,
and $\hat a'=a'+2$.
\end{lemma}

\Proof Lemma \ref{BDParamLem} implies that both $a'$ and $\hat a'$
are interfacial in the $(p,p'+p)$-model.
That $\hat a'=a'+2$ is immediate.
\cqfd
\medskip

The following result examines the differing effects of a $\B\D$-transform
followed by a $\B$-transform on a single non-interfacial point in the
$(p'-p,p')$-model with the two possible directions of the pre-segment.

\begin{figure}[ht]
  \psfrag{a}{\footnotesize$a$}
  \psfrag{c}{\footnotesize$\hat a$}
  \psfrag{x}{\footnotesize$a'$}
  \psfrag{z}{\footnotesize$\hat a'$}
  \psfrag{r}{\footnotesize$a''$}
  \psfrag{s}{\footnotesize$\hat a''$}
\includegraphics[scale=1.00]{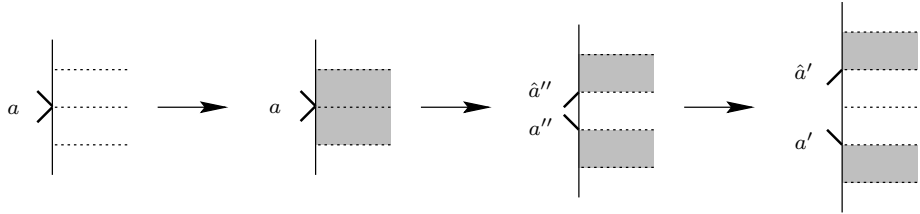}
\caption{Situation in Lemma \ref{Dire3Lem}.}
\label{Dire3Fig}
\medskip
\end{figure}

\begin{lemma}\label{Dire3Lem}
Let $p'<2p$ and $2\le a\le p'-2$, with
$\lfloor (a-1)(p'-p)/p'\rfloor=\lfloor (a+1)(p'-p)/p'\rfloor$,
so that $a$ is not interfacial in the $(p'-p,p')$-model. 
Let $a'=a+2\lfloor ap/p'\rfloor$ and
$\hat a'=a+2+2\lfloor ap/p'\rfloor$.
Then $a'$ and $\hat a'$ are both interfacial in the $(p,p'+2p)$-model.
\end{lemma}

\Proof $a$ is multifacial in the $(p,p')$-model, whence by Lemma
\ref{StartPtLem}, both $a''=a+\lfloor ap/p'\rfloor$ and
$\hat a''=a+1+\lfloor ap/p'\rfloor$ are interfacial
in the $(p,p'+p)$-model,
$\delta^{p,p'+p}_{a'',0}=\delta^{p,p'+p}_{\hat a'',1}=0$,
and
$\lfloor a''p/(p'+p)\rfloor=\lfloor\hat a''p/(p'+p)\rfloor=
\lfloor ap/p'\rfloor$.
A further application of Lemma \ref{StartPtLem} then proves the current lemma.
\cqfd
\medskip

\begin{note}\label{DireNote}
Once it is determined that $a$ and $a+2$ are both interfacial in the
$(p,p')$-model, it is readily seen that $\rho^{p,p'}(a+2)=\rho^{p,p'}(a)+1$.
This holds even if $a=0$ or $a+2=p'$.
\end{note}

\medskip

\newpage
\section{Combinatorics of the Takahashi lengths}\label{TakApp}

In this appendix, we obtain some basic but important
results concerning the continued fraction of $p'/p$ described in
Section \ref{ContFSec}, and the Takahashi lengths $\{\kappa_i\}_{i=0}^t$
and truncated Takahashi lengths $\{\tkappa_i\}_{i=0}^t$
that are derived from it as in Section \ref{TakSec}.

For fixed $p$ and $p'$, let $\{\xi_\ell\}_{\ell=0}^{2c_n-1}$
be as defined in Section \ref{ExtraSec}.
In Section \ref{SegmentSec}, we show that these values partition the
$(p,p')$-model into pieces that resemble specific other models.

\addtocontents{toc}{\SkipTocEntry}
\subsection{Model comparisons}\label{ModComSec}

Here, we relate the parameters associated with the $(p,p')$-model
for which the continued fraction is $[c_0,c_1,\ldots,c_n]$ to those
associated with certain \lq simpler\rq\ models.
In particular, if $c_0>1$, we compare them with those associated
with the $(p,p'-p)$-model and, if $c_0=1$, we compare them with
those associated with the $(p'-p,p')$-model.

In the following two lemmas, the parameters associated with those
simpler models will be primed to distinguish them from those
associated with the $(p,p')$-model.
In particular if $c_0>1$, $(p'-p)/p$ has
continued fraction $[c_0-1,c_1,\ldots,c_n]$, so that in this case,
$t'=t-1$, $n'=n$ and $t_k'=t_k-1$ for $1\le k\le n$.
If $c_0=1$, $p'/(p'-p)$ has continued fraction
$[c_1+1,c_2,\ldots,c_n]$, so that in this case,
$t'=t$, $n'=n-1$ and $t_k'=t_{k+1}$ for $1\le k\le n'$.

\begin{lemma}\label{BmodelLem}
Let $c_0>1$. For $1\le k\le n$ and $0\le j\le t$, let $y_k$, $z_k$,
$\kappa_j$ and $\tkappa_j$ be the parameters associated with
the $(p,p')$-model as defined in Section \ref{TakSec}.
For $1\le k\le n$ and $0\le j\le t'$, let $y_k'$, $z_k'$,
$\kappa_j'$ and $\tkappa_j'$ be the corresponding parameters
for the $(p,p'-p)$-model.
Then:
\begin{itemize}
\item $y_k=y_k'+z_k'\quad (0\le k\le n)$;
\item $z_k=z_k'\quad (0\le k\le n)$;
\item $\kappa_j=\kappa_{j-1}'+\tkappa_{j-1}' \quad (1\le j\le t)$;
\item $\tkappa_j=\tkappa_{j-1}' \quad (1\le j\le t)$.
\end{itemize}
\end{lemma}

\Proof Straightforward.
\cqfd
\medskip

\begin{lemma}\label{DmodelLem}
Let $c_0=1$. For $1\le k\le n$ and $0\le j\le t$, let $y_k$, $z_k$,
$\kappa_j$ and $\tkappa_j$ be the parameters associated with
the $(p,p')$-model as defined in Section \ref{TakSec}.
For $1\le k\le n'$ and $0\le j\le t$, let $y_k'$, $z_k'$,
$\kappa_j'$ and $\tkappa_j'$ be the corresponding parameters
for the $(p'-p,p')$-model.
Then:
\begin{itemize}
\item $y_k=y_{k-1}'\quad (1\le k\le n)$;
\item $z_k=y_{k-1}'-z_{k-1}'\quad (1\le k\le n)$;
\item $\kappa_j=\kappa_j' \quad (1\le j\le t)$;
\item $\tkappa_j=\kappa_j'-\tkappa_j' \quad (1\le j\le t)$.
\end{itemize}
\end{lemma}

\Proof Straightforward.
\cqfd
\medskip

\addtocontents{toc}{\SkipTocEntry}
\subsection{Segmenting the model}\label{SegmentSec}

Recall that the zone $\zeta(j)=k$ of an index $j$ satisfying $0\le j\le t$
is such that $t_k<j\le t_{k+1}$. We now express an arbitrary integer
as a sum of Takahashi lengths.

\begin{lemma}\label{TakExpLem}
Let $1\le s<p'$. Then, there is an expression:
\begin{equation}\label{sEq}
s=\sum_{i=1}^g \kappa_{\mu_i},
\end{equation}
with each $\mu_i<t_{n+1}$,
such that $0\le\zeta(\mu_1)<\zeta(\mu_2)<\cdots<\zeta(\mu_g)\le n$,
and, in addition, for each $i<g$ such that $\mu_i=t_{\zeta(\mu_i)}$,
we have $\mu_{i+1}>t_{\zeta(\mu_i)+1}$. 
Moreover, the expression is unique.
\end{lemma}
 
\Proof Use the \lq greedy algorithm\rq: let $\mu\le t$ be the largest
value such that $\kappa_\mu\le s$.
Then repeat with $s-\kappa_\mu$, if not zero.
\cqfd
\medskip

Informally, this expresses $s$ as a sum of Takahashi lengths with at
most one from each zone,
and if one of these Takahashi lengths is the last element of a zone, the
Takahashi lengths from the following zone are also excluded.
Often, having formed the expression (\ref{sEq}), we set $x=0$ if
$\zeta(\mu_1)>0$ and $x=\kappa_{\mu_1}$ if $\zeta(\mu_1)=0$ to obtain
an expression of the form:
\begin{equation}\label{sxEq}
s=x+\sum_{i=1}^g \kappa_{\mu_i},
\end{equation}
with $0\le x\le c_0$ and
$1\le\zeta(\mu_1)<\zeta(\mu_2)<\cdots<\zeta(\mu_g)\le n$,
and in addition, for each $i<g$ such that $\mu_i=t_{\zeta(\mu_i)}$,
we have $\zeta(\mu_{i+1})\ge\zeta(\mu_i)+2$, and if $x=c_0$ then
$\zeta(\mu_1)\ge2$. 

By examining the definitions in Sections \ref{ContFSec} and \ref{TakSec},
we see that $\{\tkappa_j: t_1<j\le t\}$ is the set of Takahashi
lengths for $p/(p'-c_0p)$. Therefore, we have:

\begin{corollary}\label{TakExpCor}
Let $1\le r<p$. Then, there is a unique expression:
\begin{equation}\label{rEq}
r=\sum_{i=1}^g \tkappa_{\mu_i},
\end{equation}
with $t_1<\mu_i<t_{n+1}$ for $1\le i\le g$,
such that $1\le\zeta(\mu_1)<\zeta(\mu_2)<\cdots<\zeta(\mu_g)\le n$,
and, in addition, for each $i<g$ such that $\mu_i=t_{\zeta(\mu_i)}$,
we have $\mu_{i+1}>t_{\zeta(\mu_i)+1}$. 
\end{corollary}

As we saw earlier, if $1\le r<p$, the $r$th odd band in the
$(p,p')$-model lies between heights
$\lfloor rp'/p\rfloor$ and $\lfloor rp'/p\rfloor+1$.
In the following result, we use Lemma \ref{TakExpLem} and
Corollary \ref{TakExpCor} to provide a different prescription of this.
This result was stated in \cite[p333]{bms}.

\begin{lemma}\label{OddPosLem}
1. Let $1\le r<p$ and use Corollary \ref{TakExpCor} to express $r$ in
the form (\ref{rEq}). Then
$$
\sum_{i=1}^g \kappa_{\mu_i} =
\left\lfloor \frac{rp'}{p}\right\rfloor+\delta^{(2)}_{\zeta(\mu_1),1}.
$$

\noindent
2. Let $1\le s<p'$ and use Lemma \ref{TakExpLem} to express $s$ in
the form (\ref{sxEq}). Then
$$
\sum_{i=1}^g \tkappa_{\mu_i} =
\left\lfloor \frac{sp}{p'}\right\rfloor
  +\delta_{x,0}\delta^{(2)}_{\zeta(\mu_1),0}.
$$
\end{lemma}

\Proof We first prove that the two statements are equivalent.
First assume statement 1, with $r=\sum_{i=1}^g \tkappa_{\mu_i}$.
Now $\lfloor ps/p'\rfloor=r$ if and only if
$\lfloor rp'/p\rfloor<s\le \lfloor (r+1)p'/p\rfloor$.
Then, on using Corollary \ref{TakExpCor} to write
$r+1=\sum_{i=1}^{g'} \tkappa_{\nu_i}$, we have
$\lfloor ps/p'\rfloor=r$
if and only if $\sum_{i=1}^g \kappa_{\mu_i}-\delta^{(2)}_{\zeta(\mu_1),1}<
s\le \sum_{i=1}^{g'} \kappa_{\nu_i}-\delta^{(2)}_{\zeta(\nu_1),1}$.
If $\delta^{(2)}_{\zeta(\nu_1),1}=0$ then $s=\sum_{i=1}^{g'} \kappa_{\nu_i}$
is in this range. In this case, we obtain
$\lfloor ps/p'\rfloor=r=(r+1)-1=\sum_{i=1}^{g'} \kappa_{\nu_i}-1$
via statement 1, thus giving statement 2 here.
Otherwise, $\sum_{i=1}^g \kappa_{\mu_i}-\delta^{(2)}_{\zeta(\mu_1),1}<
s<\sum_{i=1}^{g'} \kappa_{\nu_i}$.
We then claim that $s=x+\sum_{i=1}^g \kappa_{\mu_i}$ for $0\le x<c_1$,
since if this is not the case then there exists $x'+\sum_{i=1}^{g''}
\kappa_{\lambda_i}$ in this range, implying that
$s'=\sum_{i=1}^{g''}\kappa_{\lambda_i}$ also is in the range.
Thence if $r''=\sum_{i=1}^{g''}\tkappa_{\lambda_i}$ then,
using statement 1, $\lfloor p'r''/p\rfloor=s'-\delta^{(2)}_{\lambda_1,1}$.
This implies that $r<r''<r+1$, which is absurd.
So $a=x+\sum_{i=1}^g \kappa_{\mu_i}$ for $0\le x<c_1$, and
$x=0$ only if $\delta^{(2)}_{\mu_1,1}=1$.
Statement 2 then follows.
By reversing the reasoning, it may be shown that statement 2 implies
statement 1.

We now prove the two statements by using induction on the height plus
the rank of $p'/p$.
Let $p'/p$ have continued fraction $[c_0,c_1,\ldots,c_n]$ so that the
height is $n$, and the rank is $t=c_0+\cdots+c_n-2$.
Assume that the two results hold in the case where rank plus
height is $n+t-1$.

First consider the case where $c_0>1$. For $1\le r<p$, write
$r=\sum_{i=1}^g \tkappa_{\mu_i}$. Then, in terms of the parameters
of the $(p,p'-p)$-model, $r=\sum_{i=1}^g \tkappa_{\mu_i-1}'$.
Thereupon, by the induction assumption,
$\sum_{i=1}^g \kappa_{\mu_i-1}'
  +\delta^{(2)}_{\zeta'(\mu_1-1),1}=
\lfloor r(p'-p)/p\rfloor=\lfloor rp'/p\rfloor-r$.
Thus, since $\zeta'(\mu_1-1)=\zeta(\mu_1)$, we obtain
$\lfloor rp'/p\rfloor=\sum_{i=1}^g \kappa_{\mu_i-1}'
+\delta^{(2)}_{\zeta(\mu_1),1}+
\sum_{i=1}^g \tkappa_{\mu_i-1}'=
\sum_{i=1}^g \kappa_{\mu_i}+\delta^{(2)}_{\zeta(\mu_1),1}$, as required to
prove statement 1 for $c_0>1$. Statement 2 then follows in this case.

Now consider the case where $c_0=1$ and let $1\le s<p'$.
In terms of the parameters of the $(p'-p,p')$-model, we obtain
$s=x'+\sum_{i=1}^g \kappa_{\mu_i}'$.
Then, by the induction hypothesis,
$\lfloor (p'-p)s/p'\rfloor+\delta_{x',0}\delta^{(2)}_{\zeta'(\mu_1),0}
=\sum_{i=1}^g \tkappa_{\mu_i}'$.
Thereupon,
$s-1-\lfloor ps/p'\rfloor+\delta_{x',0}\delta^{(2)}_{\zeta(\mu_1),1}
=\sum_{i=1}^g \tkappa_{\mu_i}'$,
where we have made use of $\zeta(\mu_1)=\zeta'(\mu_1)+1$.
Thence,
$\lfloor ps/p'\rfloor
=\sum_{i=1}^g(\kappa_{\mu_i}'-\tkappa_{\mu_i}')
+x'-1+\delta_{x',0}\delta^{(2)}_{\zeta(\mu_1),1}
=\sum_{i=1}^g \tkappa_{\mu_i}+x'-1
+\delta_{x',0}\delta^{(2)}_{\zeta(\mu_1),1}$.
If $x'=0$ (so that $s=\sum_{i=1}^g \kappa_{\mu_i}$) then
$\lfloor ps/p'\rfloor=\sum_{i=1}^g \tkappa_{\mu_i}
-\delta^{(2)}_{\zeta(\mu_1),0}$ as required.
If $x'=1$ (so that $s=1+\sum_{i=1}^g \kappa_{\mu_i}$) then
$\lfloor ps/p'\rfloor=\sum_{i=1}^g \tkappa_{\mu_i}$,
as required.
If $x'>1$ (so that $s=\kappa_{x-1}+\sum_{i=1}^g \kappa_{\mu_i}$) then
$\lfloor ps/p'\rfloor=\tkappa_{x-1}
+\sum_{i=1}^g \tkappa_{\mu_i}$.
We thus see that statement 2 holds for the $c_0=1$ case.
Statement 1 then follows in this case, and the lemma
is thus proved.
\cqfd
\medskip

Since $\lfloor rp'/p\rfloor$ is the height of the lowermost edge of the
$r$th odd band, this theorem states that with $s=\sum_{i=1}^g \kappa_{\mu_i}$,
then if $\zeta(\mu_1)\equiv0\,(\mod2)$, the $r$th odd band
lies between heights $s$ and $s+1$;
and if $\zeta(\mu_1)\equiv1\,(\mod2)$, the $r$th odd band
lies between heights $s-1$ and $s$.

\begin{lemma}\label{OneShift1Lem}
Let $1<k\le n$ and $1\le r<z_{k}$.
If $t_k<j\le t$ then
\begin{equation*}
\left\lfloor\frac{p'(r+\tkappa_j)}{p}\right\rfloor
=\left\lfloor\frac{p'r}{p}\right\rfloor+\kappa_j.
\end{equation*}
\end{lemma}

\Proof
Write $r$ in the form (\ref{rEq}).
{}From $1\le r<z_{k}$, it follows that $g\ge1$, $\mu_g<t_k$
and $r+\tkappa_j<p$.
Set $\mu_{g+1}=j$, whereupon $\sum_{i=1}^{g+1}\tkappa_{\mu_i}$
is an expression of the form (\ref{rEq}) for $r+\tkappa_j$.
Then:
\begin{equation*}
\left\lfloor\frac{p'(r+\tkappa_j)}{p}\right\rfloor
=\sum_{i=1}^{g+1}\kappa_{\mu_i}
     -\delta^{(2)}_{\zeta(\mu_1),1}
=\sum_{i=1}^{g}\kappa_{\mu_i}
     -\delta^{(2)}_{\zeta(\mu_1),1}+\kappa_j
=\left\lfloor\frac{p'r}{p}\right\rfloor+\kappa_j,
\end{equation*}
where the first and third equalities follow from Lemma \ref{OddPosLem}(1).
\cqfd
\medskip

\begin{lemma}\label{OneShift2Lem}
Let $1\le k\le n$ and $1\le r<z_{k+1}$.
If $r+z_k<p$ then
\begin{equation*}
\left\lfloor\frac{p'(r+z_k)}{p}\right\rfloor
=\left\lfloor\frac{p'r}{p}\right\rfloor+y_k.
\end{equation*}
\end{lemma}

\Proof
Let $r$ be expressed in the form (\ref{rEq}).
Note that $g\ge1$, $\mu_1>t_1$ and $\mu_g<t_{k+1}$.
Depending on the value of $\mu_g$, we reexpress
$\kappa_{\mu_g}+y_k$ and $\tkappa_{\mu_g}+z_k$ according
to one of the following five cases:
i)   $t_k<\mu_g<t_{k+1}$.
     Here $\kappa_{\mu_g}+y_k=y_{k-1}+(\mu_g-t_k)y_k+y_k=\kappa_{\mu_g+1}$.
     Similarly, $\tkappa_{\mu_g}+z_k=\tkappa_{\mu_g+1}$.
ii)  $t_{k-1}+1<\mu_g\le t_{k}$.
     Here $\kappa_{\mu_g}+y_k=y_{k-2}+(\mu_g-t_{k-1})y_{k-1}+y_k
         =y_{k-2}+(\mu_g-t_{k-1}-1)y_{k-1}+(y_{k-1}+y_k)
         =\kappa_{\mu_g-1}+\kappa_{t_k+1}$.
     Similarly, $\tkappa_{\mu_g}+z_k=\tkappa_{\mu_g-1}+\tkappa_{t_k+1}$.
iii) $\mu_g=t_{k-1}+1$.
     Here $\kappa_{\mu_g}+y_k=y_{k-2}+(\mu_g-t_{k-1})y_{k-1}+y_k
         =y_{k-2}+y_{k-1}+y_{k}=y_{k-2}+\kappa_{t_k+1}$.
     Similarly, $\tkappa_{\mu_g}+z_k=z_{k-2}+\tkappa_{t_k+1}$.
iv)  $\mu_g=t_{k-1}$.
     Here $\kappa_{\mu_g}+y_k=y_{k-1}+y_k=\kappa_{t_k+1}$.
     Similarly, $\tkappa_{\mu_g}+z_k=\tkappa_{t_k+1}$.
v)   $\mu_g<t_{k-1}$. Here, no reexpressing is required.

We now use these to examine $\lfloor p'(r+z_k)/p\rfloor$,
tackling each of the five cases in turn (but leaving case iii) until last).
For convenience, we set $r'=\sum_{i=1}^{g-1}\tkappa_{\mu_i}$
so that $r+z_k=r'+\tkappa_{\mu_g}+z_k$ in which the last two
terms will be reexpressed as above.

i) $t_k<\mu_g<t_{k+1}$. Here,
\begin{equation*}
\begin{split}
\left\lfloor\frac{p'(r+z_k)}{p}\right\rfloor
&=\left\lfloor\frac{p'(r'+\tkappa_{\mu_g+1})}{p}\right\rfloor
=\sum_{i=1}^{g-1}\kappa_{\mu_i}+\kappa_{\mu_g+1}
                 -\delta^{(2)}_{\zeta(\mu_1),1}\\
&=\sum_{i=1}^{g}\kappa_{\mu_i}+y_k
                 -\delta^{(2)}_{\zeta(\mu_1),1}
=\left\lfloor\frac{p'r}{p}\right\rfloor+y_k,
\end{split}
\end{equation*}
where the first and third equalities follow from the case i) reexpressings
above, and the second and fourth follow from Lemma \ref{OddPosLem}(1),
noting when $g=1$ that $\zeta(\mu_1+1)=k=\zeta(\mu_1)$.

ii) $t_{k-1}+1<\mu_g\le t_{k}$. Here,
\begin{equation*}
\begin{split}
\left\lfloor\frac{p'(r+z_k)}{p}\right\rfloor
&=\left\lfloor\frac{p'(r'+\tkappa_{\mu_g-1}+\tkappa_{t_k+1})}{p}\right\rfloor
=\sum_{i=1}^{g-1}\kappa_{\mu_i}+\kappa_{\mu_g-1}+\kappa_{t_k+1}
                 -\delta^{(2)}_{\zeta(\mu_1),1}\\
&=\sum_{i=1}^{g}\kappa_{\mu_i}+y_k
                 -\delta^{(2)}_{\zeta(\mu_1),1}
=\left\lfloor\frac{p'r}{p}\right\rfloor+y_k,
\end{split}
\end{equation*}
where the first and third equalities follow from the case ii) reexpressings
above, and the second and fourth follow from Lemma \ref{OddPosLem}(1),
noting when $g=1$ that $\zeta(\mu_1-1)=k-1=\zeta(\mu_1)$.

iv) $\mu_g=t_{k-1}$. If $g>1$ then
\begin{equation*}
\begin{split}
\left\lfloor\frac{p'(r+z_k)}{p}\right\rfloor
&=\left\lfloor\frac{p'(r'+\tkappa_{t_k+1})}{p}\right\rfloor
=\sum_{i=1}^{g-1}\kappa_{\mu_i}+\kappa_{t_k+1}
                 -\delta^{(2)}_{\zeta(\mu_1),1}\\
&=\sum_{i=1}^{g}\kappa_{\mu_i}+y_k
                 -\delta^{(2)}_{\zeta(\mu_1),1}
=\left\lfloor\frac{p'r}{p}\right\rfloor+y_k,
\end{split}
\end{equation*}
where the first and third equalities follow from the case iv) reexpressings
above, and the second and fourth follow from Lemma \ref{OddPosLem}(1).
If $g=1$ then the argument is similar (with $r'=0$), noting that
$\delta^{(2)}_{\zeta(t_k+1),1}=\delta^{(2)}_{\zeta(t_{k-1}),1}$
because $\zeta(t_k+1)=k$ and $\zeta(t_{k-1})=k-2$.

v) $\mu_g<t_{k-1}$. Since $z_k=\tkappa_{t_k}$ and $\zeta(t_k)=k-1$,
we may apply Lemma \ref{OddPosLem}(1) directly. Then
\begin{equation*}
\begin{split}
\left\lfloor\frac{p'(r+z_k)}{p}\right\rfloor
&=\left\lfloor\frac{p'(r+\tkappa_{t_k})}{p}\right\rfloor
=\sum_{i=1}^{g}\kappa_{\mu_i}+\kappa_{t_k}
                 -\delta^{(2)}_{\zeta(\mu_1),1}\\
&=\sum_{i=1}^{g}\kappa_{\mu_i}+y_k
                 -\delta^{(2)}_{\zeta(\mu_1),1}
=\left\lfloor\frac{p'r}{p}\right\rfloor+y_k,
\end{split}
\end{equation*}
as required.

Case iii) requires the whole proof to be entombed in an induction argument.
We first prove the lemma for $k=1$ and $k=2$.
If $k=1$ then $g=1$ and necessarily $t_1<\mu_1<t_2$.
This $k=1$ case then follows from i) as above.
If $k=2$ then necessarily $t_{1}<\mu_g<t_{3}$.
Except for $\mu_g=t_{1}+1$ this is dealt with by i) and ii) above.
For $\mu_g=t_{1}+1$, we necessarily have $g=1$ so that $r=1$.
By definition, $z_2=c_1$, $y_1=c_0$ and $y_2=c_0c_1+1$.
Since $p'/p=c_0+1/(c_1+\epsilon)$ where $0<\epsilon<1$, it follows
that $\lfloor p'(r+z_2)/p\rfloor=c_0c_1+c_0+1=\lfloor p'r/p \rfloor+y_2$,
as required.

Having established the lemma for $k\le2$, we now consider $k>2$.
For the purposes of induction, assume that the lemma holds for $k$
replaced by $k-2$. For all cases except where $\mu_g=t_{k-1}+1$,
the lemma is immediately proved by i), ii), iv) or v) above.
In the remaining case $\mu_g=t_{k-1}+1$, we immediately
obtain that $\tkappa_{\mu_g+1}=\tkappa_{\mu_g}+z_{k-1}$.
Then, because $r<\tkappa_{\mu_g+1}$, we obtain $r'<z_{k-1}$
and furthermore $r'+z_{k-2}<z_{k-2}+z_{k-1}\le z_k$.
Then, if $r'>0$:
\begin{equation*}
\begin{split}
\left\lfloor\frac{p'(r+z_k)}{p}\right\rfloor
&=\left\lfloor\frac{p'(r'+z_{k-2}+\tkappa_{t_k+1})}{p}\right\rfloor
=\left\lfloor\frac{p'(r'+z_{k-2})}{p}\right\rfloor+\kappa_{t_k+1}\\
&=\left\lfloor\frac{p'r'}{p}\right\rfloor+y_{k-2}+\kappa_{t_k+1}
=\left\lfloor\frac{p'r'}{p}\right\rfloor+\kappa_{\mu_g}+y_k\\
&=\left\lfloor\frac{p'(r'+\tkappa_{\mu_g})}{p}\right\rfloor+y_k
=\left\lfloor\frac{p'r}{p}\right\rfloor+y_k,
\end{split}
\end{equation*}
where the first and fourth equalities follow from the case iii) reexpressings
above, the second and fifth equalities follow from Lemma \ref{OneShift1Lem},
and the third follows from the induction hypothesis.
The same string of equalities may be used in the case $r'=0$ provided
that we subtract $\delta^{(2)}_{k-3,1}$ ($=\delta^{(2)}_{\mu_g,1}$)
from the fourth and fifth expressions: the third and fifth
equalities then follow from Lemma \ref{OddPosLem}(1).
This completes the induction step and the lemma follows.
\cqfd
\medskip

In some of the following results, we state that certain quantities
have continued fractions of a certain form $[d_0,d_1,\ldots,d_{m-1},d_{m}]$.
The case $d_m=1$ may arise for $m>0$. This continued fraction would then
not be of the standard form stated in Section \ref{ContFSec}.
In such a case, the continued fraction should be
reinterpreted as $[d_0,d_1,\ldots,d_{m-1}+1]$.

\begin{lemma}\label{CFLem}
Let $p'/p$ have continued fraction $[c_0,c_1,\ldots,c_n]$.
For $0\le j\le t$, set $k$ such that $t_k<j\le t_{k+1}$.
Then:
\par\noindent 1. $y_k\tkappa_j-z_k\kappa_j=(-1)^k$;
\par\noindent 2. $\kappa_j$ and $\tkappa_j$ are coprime;
\par\noindent 3. $\kappa_j/\tkappa_j$ has continued fraction
$[c_0,c_1,\ldots,c_{k-1},j-t_k]$.
\end{lemma}

\Proof If $k=0$ then $y_0=1$, $z_0=0$, whereupon $\kappa_j=j+1$
and $\tkappa_j=1$. Each part is then seen to hold when $k=0$.
For the purposes of induction, let $k>0$ and assume that each part
holds for $k$ replaced by $k-1$.
Note that $\kappa_{t_k}=y_k$ and $\tkappa_{t_k}=z_k$.
The induction hypothesis then implies that
$y_{k-1}z_k-z_{k-1}y_k=(-1)^{k-1}$.
Statement 1.\ is now verified by using
$\kappa_j=y_{k-1}+(j-t_k)y_k$ and $\tkappa_j=z_{k-1}+(j-t_k)z_k$.
Statement 1.\ then follows for general $k$ by induction.
Statement 2.\ follows immediately from statement 1.
In the case of statement 3., the induction hypothesis states that:
\begin{equation*}
\frac{c_{k-1}y_{k-1}+y_{k-2}}{c_{k-1}z_{k-1}+z_{k-2}}=
\frac{y_k}{z_k}=
{c_0+\frac{\displaystyle\strut 1}{\phantom{1+}
\genfrac{}{}{}0{\lower-5pt\hbox{$\vdots$}}
{\strut c_{k-2} +
\frac{\displaystyle\strut 1}{\displaystyle\strut c_{k-1}}}}}.
\end{equation*}
Both sides here are rational functions of $c_{k-1}$ which are
positive when $c_{k-1}$ is positive.
Therefore, we may replace $c_{k-1}$ with $c_{k-1}+1/(j-t_k)$.
The right side yields the continued fraction
$[c_0,c_1,\ldots,c_{k-1},j-t_k]$, while the left side yields:
\begin{equation*}
\frac{(c_{k-1}+1/(j-t_k))y_{k-1}+y_{k-2}}
{(c_{k-1}+1/(j-t_k))z_{k-1}+z_{k-2}}
=\frac{y_k+y_{k-1}/(j-t_k)}{z_k+z_{k-1}/(j-t_k)}
=\frac{(j-t_k)y_k+y_{k-1}}{(j-t_k)z_k+z_{k-1}}
=\frac{\kappa_j}{\tkappa_j}.
\end{equation*}
This completes the induction step, whereupon the lemma follows.
\cqfd
\medskip

\begin{lemma}\label{ExCFLem}
Let $p'/p$ have continued fraction $[c_0,c_1,\ldots,c_n]$.
Let $0<\ell<2c_n-1$, and set $\hat p'=\xi_{\ell+1}-\xi_{\ell}$
and $\hat p=\tilde\xi_{\ell+1}-\tilde\xi_{\ell}$.

\noindent
1. $\hat p'\ge2$ if and only if one of the following cases applies:
\par i) $\ell$ is odd, $n\ge2$ and if $n=2$ then $c_0>1$.
        In this case $\hat p'/\hat p$ has continued fraction
        $[c_0,c_1,\ldots,c_{n-2}]$;
\par ii) $\ell\in\{0,2c_n-2\}$, $n\ge1$ and if $n=1$ then $c_0>1$.
        In this case $\hat p'/\hat p$ has continued fraction
        $[c_0,c_1,\ldots,c_{n-1}]$;
\par iii) $\ell$ is even, $0<\ell<2c_n-2$, $c_{n-1}>1$, $n\ge1$ and
        if $n=1$ then $c_0>2$.
        In this case $\hat p'/\hat p$ has continued fraction
        $[c_0,c_1,\ldots,c_{n-1}-1]$;
\par iv) $\ell$ is even, $0<\ell<2c_n-2$, $c_{n-1}=1$, $n\ge3$ and
        if $n=3$ then $c_0>1$.
        In this case $\hat p'/\hat p$ has continued fraction
        $[c_0,c_1,\ldots,c_{n-3}]$.

\noindent
2. Assume that $\hat p'\ge2$, and let $\{\kappa_i\}_{i=0}^t$ and
$\{\tkappa_i\}_{i=0}^t$ be respectively
the Takahashi lengths and the truncated Takahashi lengths for $p'/p$.
Let $\{\kappa_i'\}_{i=0}^{t'}$ and $\{\tkappa'_i\}_{i=0}^{t'}$
be the corresponding values for $\hat p'/\hat p$.
Then $\kappa_i'=\kappa_i$ and $\tkappa_i'=\tkappa_i$ for $0\le i\le t'$.
\end{lemma}

\Proof First note, using the definitions of Section \ref{TakSec},
that (when defined) $y_{-1}=0$, $y_0=1$, $y_1=c_0$, $y_2=c_1c_0+1$,
$z_{-1}=1$, $z_0=0$, $z_1=1$, and $z_2=c_1$.
We now split into the following four cases to calculate
$\hat p'$ and $\hat p$ using the definitions of Section \ref{ExtraSec}.

i).\ $\ell$ is odd: here $\hat p'=y_{n-1}$ and $\hat p=z_{n-1}$.
So we certainly require $n\ge2$. Since $y_{n-1}=\kappa_{t_{n-1}}$ and
$z_{n-1}=\tkappa_{t_{n-1}}$, Lemma \ref{CFLem}(3) implies that
$\hat p'/\hat p$ has continued fraction $[c_0,c_1,\ldots,c_{n-2}]$.
Then to ensure that $\hat p'\ge2$, we must exclude the case
$n=2$ and $c_0=1$.

ii).\ $\ell\in\{0,2c_n-2\}$: here $\hat p'=y_{n}$ and $\hat p=z_{n}$
whereupon this case follows similarly to the first.

iii).\ $\ell$ is even, $0<\ell<2c_n-2$ and $c_{n-1}>1$:
here $\hat p'=y_n-y_{n-1}$ and $\hat p=z_n-z_{n-1}$.
So here we certainly require $n\ge1$.
Then $\hat p'=y_{n-2}+(c_{n-1}-1)y_{n-1}=\kappa_{t_n-1}$
and similarly $\hat p=\tkappa_{t_n-1}$.
Lemma \ref{CFLem}(3) now implies that
$\hat p'/\hat p$ has continued fraction $[c_0,c_1,\ldots,c_{n-1}-1]$,
whereupon we see that to ensure $\hat p'\ge2$, we must exclude the case 
$n=1$ and $c_0=2$.

iv).\ $\ell$ is even, $0<\ell<2c_n-2$ and $c_{n-1}=1$:
here, $\hat p'=y_n-y_{n-1}=y_{n-2}+(c_{n-1}-1)y_{n-1}=y_{n-2}$ and
$\hat p=z_{n-2}$, whereupon this case follows similarly to the first.

Statement 2 now follows immediately via the definitions of
Section \ref{TakSec}.
\cqfd
\medskip

\begin{theorem}\label{SegmentHalfLem}
Let $p'/p$ have continued fraction $[c_0,c_1,\ldots,c_n]$
with $n\ge1$, and if $n=1$ then $c_0>1$.
If $0<\ell<2c_n-1$, then $\xi_{\ell}$ is interfacial in the
$(p,p')$-model and neighbours the $\tilde\xi_{\ell}$th odd band.
\end{theorem}

\Proof
For $n=1$, the $k$th odd band lies between heights $kc_0$ and $kc_0+1$
for $1\le k<c_0$.
In this case, $\xi_{2k-1}=kc_0$, $\xi_{2k}=kc_0+1$,
and $\tilde\xi_{2k-1}=\tilde\xi_{2k}=k$ for $1\le k<c_0$,
whereupon this case follows immediately.

Hereafter, assume that $n\ge2$.
Note that it is sufficient to prove the $p'<2p$ case: the
$p'<2p$ case then follows because on passing from the
$(p,p')$-model to the $(p'-p,p')$ model, each band changes parity,
and via Lemma \ref{DmodelLem},
the values $\xi_{\ell}$ are unchanged and the values of
$\tilde\xi_{\ell}$ change to $\xi_{\ell}-\tilde\xi_{\ell}$.

If $\ell$ is odd,
$\ell=2k-1$ for $1\le k\le c_n-1$, and
$\tilde\xi_{\ell}=kz_n$ and $\xi_{\ell}=ky_n$.
Then:
\begin{equation*}
\left\lfloor\frac{p'\tilde\xi_\ell}{p}\right\rfloor
=\left\lfloor\frac{p'(kz_n)}{p}\right\rfloor
=\left\lfloor\frac{p'z_{n}}{p}\right\rfloor+(k-1)y_n
=y_{n}-\delta^{(2)}_{n,0}+(k-1)y_n
=\xi_\ell-\delta^{(2)}_{n,0},
\end{equation*}
where the second equality follows by repeated use of Lemma \ref{OneShift2Lem}
on noting that $kz_n<z_{n-1}+c_nz_n=z_{n+1}=p$,
and the third follows from Lemma \ref{OddPosLem}(1).
This implies that either the upper or lower edge of the
$\tilde\xi_{\ell}$th odd band is at height $\xi_{\ell}$.
For $p'>2p$ this immediately implies that $\xi_{\ell}$ is interfacial.

If $\ell$ is even,
$\ell=2k$ for $1\le k\le c_n-1$, and
$\tilde\xi_{\ell}=kz_n+z_{n-1}$ and $\xi_{\ell}=ky_n+y_{n-1}$.
Then if $n>2$, we obtain in a similar way to the odd $\ell$ case above,
\begin{equation*}
\left\lfloor\frac{p'\tilde\xi_\ell}{p}\right\rfloor
=\left\lfloor\frac{p'(z_{n-1}+kz_n)}{p}\right\rfloor
=\left\lfloor\frac{p'z_{n-1}}{p}\right\rfloor+ky_n
=y_{n-1}-\delta^{(2)}_{n,1}+ky_n
=\xi_\ell-\delta^{(2)}_{n,1}.
\end{equation*}
This calculation is also valid for $n=2$, when the validity of the
third equality follows from $z_1=1$ and $\lfloor p'/p\rfloor=c_0=y_1$.
As above, for $p'>2p$ this immediately implies that $\xi_{\ell}$ is
interfacial.
The lemma then follows.
\cqfd
\medskip

\begin{theorem}\label{SegmentLem}
Let $p'/p$ have continued fraction $[c_0,c_1,\ldots,c_n]$ with $n\ge1$.
Let $0\le\ell\le 2c_n-2$,
let $\hat p'=\xi_{\ell+1}-\xi_{\ell}$, and let
$\hat p=\tilde\xi_{\ell+1}-\tilde\xi_{\ell}$.
Then for $1\le s\le \hat p'-2$,
the parity of the $s$th band of the $(\hat p,\hat p')$-model
is equal to the parity of the $(\xi_\ell+s)$th band of the $(p,p')$-model.
\end{theorem}

\Proof
%
Let $1\le r<\hat p$.
We claim that
\begin{equation*}
\left\lfloor\frac{p'(r+\tilde\xi_\ell)}{p}\right\rfloor
=\left\lfloor\frac{p'r}{p}\right\rfloor+\xi_\ell,
\end{equation*}
In the $\ell=0$ case this holds trivially because
$\xi_\ell=\tilde\xi_\ell=0$.
There are two other subcases: $\ell$ odd and $\ell$ even.

For $\ell=2k-1$ with $1\le k\le c_n-1$, we have
$\xi_\ell=ky_n$, $\tilde\xi_\ell=kz_n$ and
$\hat p=\tilde\xi_{\ell+1}-\tilde\xi_\ell=z_{n-1}$.
Then:
\begin{equation*}
\left\lfloor\frac{p'(r+\tilde\xi_\ell)}{p}\right\rfloor
=\left\lfloor\frac{p'(r+kz_n)}{p}\right\rfloor
=\left\lfloor\frac{p'r}{p}\right\rfloor+ky_n,
=\left\lfloor\frac{p'r}{p}\right\rfloor+\xi_\ell,
\end{equation*}
where the second equality is established by repeated use of
Lemma \ref{OneShift2Lem}, the validity of which is guaranteed because
$1\le r+kz_n<z_{n-1}+c_nz_n=z_{n+1}=p$.

For $\ell=2k$ with $1\le k\le c_n-1$, we have
$\xi_\ell=y_{n-1}+ky_n$ and $\tilde\xi_\ell=z_{n-1}+kz_n$.
For $k<c_n-1$, we have
$\hat p=\tilde\xi_{\ell+1}-\tilde\xi_\ell=z_{n}-z_{n-1}$,
whereas for $k=c_n-1$, we have
$\hat p=\tilde\xi_{\ell+1}-\tilde\xi_\ell=z_{n}$.
Since $z_1=1$, we may restrict consideration here to $n>1$.
Then:
\begin{equation*}
\begin{split}
\left\lfloor\frac{p'(r+\tilde\xi_\ell)}{p}\right\rfloor
&=\left\lfloor\frac{p'(r+z_{n-1}+kz_n)}{p}\right\rfloor\\
&=\left\lfloor\frac{p'(r+z_{n-1})}{p}\right\rfloor+ky_n
=\left\lfloor\frac{p'r}{p}\right\rfloor+y_{n-1}+ky_n
=\left\lfloor\frac{p'r}{p}\right\rfloor+\xi_\ell,
\end{split}
\end{equation*}
where the second equality is established by repeated use of
Lemma \ref{OneShift2Lem}, the validity of which is guaranteed because
$1\le r+z_{n-1}+kz_n<z_{n-1}+c_nz_n=z_{n+1}=p$;
and the third equality follows from a single use of
Lemma \ref{OneShift2Lem}, and is valid because $1\le r<z_n$.

Lemma \ref{ExCFLem}(2) states that the sequence of truncated Takahashi
lengths for $\hat p'/\hat p$ begins the sequence of truncated Takahashi
lengths for $p'/p$. This implies that the expression (\ref{rEq})
for $r$ in terms of the truncated Takahashi lengths of $p'/p$ is
identical to that in terms of the truncated Takahashi lengths of
$\hat p'/\hat p$.
It then follows from Lemma \ref{OddPosLem}(1) that
$\lfloor\hat p'r/\hat p\rfloor=\lfloor p'r/p\rfloor$
and therefore
\begin{equation*}
\left\lfloor\frac{p'(r+\tilde\xi_\ell)}{p}\right\rfloor
=\left\lfloor\frac{\hat p'r}{\hat p}\right\rfloor+\xi_\ell.
\end{equation*}
Thus, for $1\le r<\hat p$, the height of the $r$th odd band in the
$(\hat p,\hat p')$-model is precisely $\xi_\ell$ less that the height of
$(r+\tilde\xi_\ell)$th odd band in the $(p,p')$-model.
The theorem is then proved if it can be shown that there are no other
odd bands between heights $\xi_\ell+1$ and $\xi_{\ell+1}-1$ in
the $(p,p')$-model.
This is so because,
by Lemma \ref{SegmentHalfLem}, the lowermost edge of the $\tilde\xi_\ell$th 
odd band in the $(p,p')$-model (if there is one)
is at height $\xi_\ell$ or $\xi_\ell-1$,
and the lowermost edge of the $\tilde\xi_{\ell+1}$th odd band
in the $(p,p')$-model (if there is one) is at height
$\xi_{\ell+1}$ or $\xi_{\ell+1}-1$.
%
\cqfd
\medskip

The above result shows that for $0\le\ell\le 2c_n-2$, the
band structure between heights $\xi_{\ell}+1$ and $\xi_{\ell+1}-1$
in the $(p,p')$-model is identical to the band structure of a
certain smaller model (which is specified in Lemma \ref{ExCFLem}).
This result proves useful in that a certain subset of the set
of paths in the $(p,p')$-model will be shown to be the set of paths that
lie between heights $\xi_{\ell}+1$ and $\xi_{\ell+1}-1$ for some $\ell$.
Their generating function will thus be given by that of the smaller model.

\begin{figure}[ht]
\includegraphics[scale=1.00]{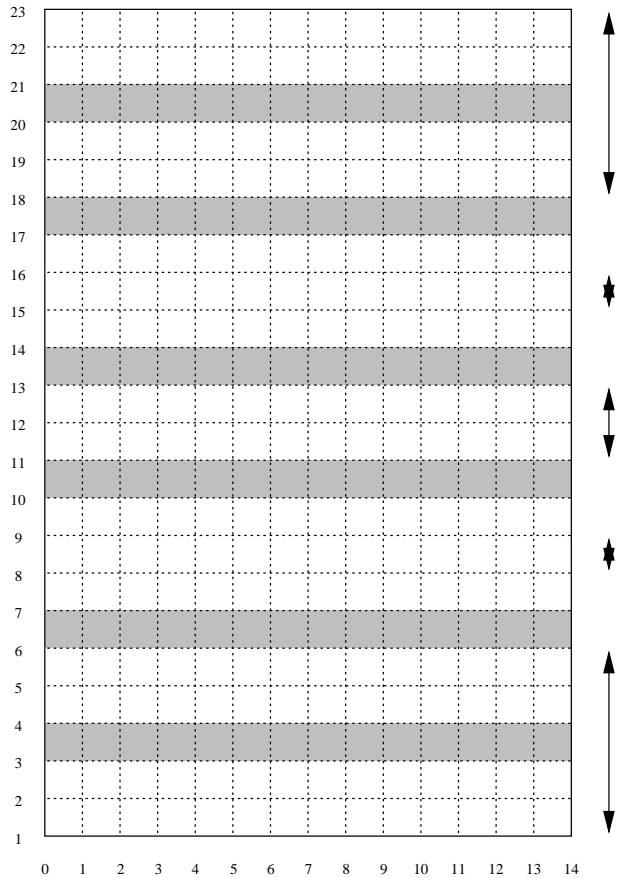}
\caption{$(7,24)$-model.}
\label{Model7_24}
\medskip
\end{figure}

We demonstrate Lemma \ref{SegmentLem} using the $(7,24)$-model
which is depicted in Fig.~\ref{Model7_24}.

\noindent Since $24/7$ has continued fraction $[3,2,3]$, we
obtain $y_{n-1}=3$ and $y_n=7$.
Thereupon $(\xi_0,\xi_1,\xi_2,\xi_3,\xi_4,\xi_5)=(0,7,10,14,17,24)$.
Lemma \ref{SegmentLem} states that for $\ell=0,1,2,3,4,$ the band structure
of the $(7,24)$-model between heights $\xi_{\ell}+1$ and $\xi_{\ell+1}-1$
is identical to the band structure of a specific smaller model.
The arrowed lines to the right of the grid in Fig.~8 show the
extent of these smaller models. For $\ell=0,1,2,3,4,$, they are
seen (by inspection, or by using Lemma \ref{SegmentLem}) to be the
$(2,7)$-, $(1,3)$-, $(1,4)$-, $(1,3)$-, and $(2,7)$-models
respectively.

Finally, we make use of Lemma \ref{OneShift2Lem} to provide a result
required in the proof of Theorem \ref{CombinedThrm}.

\begin{lemma}\label{AnnoyingLem}
For $1\le p<p'$, let $p'/p$ have continued fraction $[c_0,c_1,\ldots,c_n]$
with either $n\ge3$ and $c_0=c_1=c_2=1$, or $n=3$ and $c_0=c_1=1$
and $c_2=2$.
Then, for $\Delta=\pm1$, there exists no value $a$ with
$1\le a,a+5\Delta<p'$ such that each element of
$\{a+\Delta,a+2\Delta,a+4\Delta\}$ is interfacial and
$\rho^{p,p'}(a+\Delta)=\rho^{p,p'}(a+2\Delta)$.
\end{lemma}

\Proof First note that since odd bands are separated by either $c_0$
or $c_0-1$ even bands, then necessarily $a+3\Delta$ is interfacial.
Then, without loss of generality, we may assume that $\Delta=+1$, 
whereupon the specified configuration arises only if
$\lfloor (a+5)p/p'\rfloor=\lfloor ap/p'\rfloor+2$.

In the case where $c_0=c_1=c_2=1$, we have $p'/p=1+1/(1+1/(1+\epsilon))$
with $0<\epsilon<1$. Then $5p/p'=\frac52(1+1/(3+2\epsilon))>3$ implying
that $\lfloor (a+5)p/p'\rfloor\ge\lfloor ap/p'\rfloor+3$ so that
this case is excluded.

Now consider the case where $p'/p$ has continued fraction $[1,1,2,c_3]$.
Here, $y_3=5$, $p'=y_4=5c_3+2$, $z_3=3$, and $p=z_4=3c_3+1$.
Use of Lemma \ref{OneShift2Lem} for $k=3$ yields
$\lfloor (r+3)p'/p\rfloor=\lfloor p'r/p\rfloor+5$ for $1\le r<p-3$.
This shows that the $(r+1)$th, $(r+2)$th and $(r+3)$th
odd bands lie within five consecutive bands.
This also holds for $r=0$ because
$\lfloor 3p'/p\rfloor =\lfloor (15c_3+6)/(3c_3+1)\rfloor=5$.
Thus, for every three consecutive odd bands, there are necessarily
two that are adjacent.
The required result then follows.
\cqfd
\medskip

\end{appendix}

\end{document}